\documentclass{amsart}
\usepackage{amsmath,amssymb,verbatim}
\usepackage{verbatim,a4wide}
\usepackage{enumerate}
\newtheorem{theorem}{Theorem}[section]
\newtheorem{lemma}[theorem]{Lemma}
\newtheorem{proposition}[theorem]{Proposition}
\newtheorem{corollary}[theorem]{Corollary}
\newtheorem{claim}[theorem]{Claim}
\theoremstyle{definition}

\theoremstyle{remark}
\newtheorem{remark}{Remark}[section]
\newtheorem{example}{Example}[section]
\numberwithin{equation}{section}
\DeclareMathOperator{\sgn}{sgn}
\DeclareMathOperator{\sech}{sech}

\DeclareMathOperator{\spann}{span}
\DeclareMathOperator{\diag}{diag}
\DeclareMathOperator*{\res}{Res}
\newcommand{\R}{\mathbb{R}}
\newcommand{\C}{\mathbb{C}}
\newcommand{\N}{\mathbb{N}}
\newcommand{\Z}{\mathbb{Z}}
\newcommand{\la}{\langle}
\newcommand{\ra}{\rangle}
\newcommand{\pd}{\partial}
\newcommand{\eps}{\epsilon}
\newcommand{\bc}{\mathbf{c}}
\newcommand{\bx}{\mathbf{x}}
\newcommand{\bxd}{\mathbf{x'}}

\newcommand{\tbx}{\mathbf{\tilde{x}}}
\newcommand{\tbxd}{\mathbf{\tilde{x}'}}
\newcommand{\tbxdd}{\mathbf{\tilde{x}''}}
\newcommand{\bz}{\mathbf{z}}
\newcommand{\tg}{\tilde{g}}
\newcommand{\tk}{\tilde{k}}
\newcommand{\tpsi}{\tilde{\psi}}
\newcommand{\tu}{\tilde{u}}
\newcommand{\tv}{\tilde{v}}

\newcommand{\mL}{\mathcal{L}}
\newcommand{\mM}{\mathcal{M}}

\newcommand{\wmL}{\widetilde{\mathcal{L}}}
\newcommand{\wmM}{\widetilde{\mathcal{M}}}
\newcommand{\wT}{\widetilde{T}}
\newcommand{\calX}{\mathcal{X}}
\newcommand{\tfg}{\tilde{\mathfrak{g}}}
\renewcommand{\a}{\alpha} 
\renewcommand{\k}{\kappa}
\begin{document}
\title{Linear stability of elastic $2$-line solitons for the KP-II equation}
\author{Tetsu Mizumachi}
\address{Division of Mathematical and Information Sciences,
Hiroshima University, Kagamiyama 1-7-1, 739 - 8521 Japan}
\email{tetsum@hiroshima-u.ac.jp}
\dedicatory{Dedicated to Professor Robert L.~Pego's 65th birthday.}
\begin{abstract}
  The KP-II equation was derived by Kadomtsev and Petviashvili to
  explain stability of line solitary waves of shallow water.  Using
  the Darboux transformations, we study linear stability of $2$-line
  solitons whose line solitons interact elastically each other.  Time
  evolution of resonant continuous eigenfunctions is described by a
  damped wave equation in the transverse variable which is supposed to
  be a linear approximation of the local phase shifts of modulating
  line solitons.
\end{abstract}
\maketitle
\section{Introduction}
The KP-II equation
\begin{equation}\label{eq:KPII}
\partial_x(4\pd_tu+\pd_x^3u+3\pd_x(u^2))+3\partial_y^2u=0
\quad\text{for $t>0$ and $(x,y)\in \R^2$,}
\end{equation}
has been derived by Kadomtsev and Petviashvili (\cite{KP}) to explain transverse stability
of $1$-soliton solutions of the KdV equation in the media of negative dispersion.
\par
For the $1$-line solitons for the the KP-II equation posed on $\R^2$,
its spectral stability was proved by \cite{APS} and nonlinear
stability was proved by \cite{Miz15,Miz18,Miz19}. See also \cite{MT}
which proves stability in the $\R_x\times\mathbb{T}_y$ case and
\cite{HLP} which proves transverse linear stability of periodic waves.
Cauchy problems of \eqref{eq:KPII} have been studied by
\cite{Bourgain,Hadac,HHK,IM,Tak,TT,Tz,Ukai} and by \cite{MST} on the background of
line solitons.
\par
The KP-II equation has exact multi-line soliton solutions.  Most of
them have web patterns in the middle where line solitary waves
interact inelastically and the web pattern enlarges as time goes on.
In the present paper, we will study linear stability of elastic
$2$-line solitons of $P$-type (see Theorems~\ref{thm:AS} and \ref{thm:AS2})
and of $O$-type (see Theorems~\ref{thm:AS3} and \ref{thm:AS4}).
They tell us line solitons are
stable in exponentially weighted space and in the sense of linear
approximation, the modulations of line solitons are described by
$1$-dimensional damped wave equations in the transverse variable.  We
use exponentially weighted spaces not only they are useful to detect
directions of wave motion as they are for the KdV equation
(\cite{PW}) but also continuous eigenfunctions that describe wave
propagation along the crest of line solitons grow exponentially in
the rear of line solitons (see Section~\ref{sec:res2}).
\par
Our approach is to compare the spectrum of linearized operator around a $2$-line
soliton solution and that around a $1$-line soliton solution by using the
Darboux transformations and to prove that there is no spectrum in the unstable half
plane or in the vicinity of the imaginary axis except for a curve of resonant
continuous eigenvalues. 
\par
Darboux transformations or B\"acklund transformations become
powerful tools to investigate linear and nonlinear stability of
solitons and breathers of integrable systems as a PDE method
after the work by Merle and Vega \cite{MV}.  See
\cite{Alejo-Munoz13,Alejo-Munoz15,Contreras-Pelinovsky,MP,MPel,Munoz-Placios}.
Moreover, the B\"acklund transformation has been recently used to
prove non-degeneracy of lump solutions of the KP-I equation
(\cite{Liu-Wei}).
\par
The difficulty to use the Darboux transformation in our problem comes
from the fact that the endpoints of the curve of continuous eigenfunctions
that have to do with modulations of line solitons lie on the essential
spectrum of the linearized operator around the null solution. See
Lemma~\ref{lem:lambda2-est} and Remark~\ref{rem:lambda1-est}. Because
of this reason, it is difficult to construct an isomorphism between
solution spaces of linearized equations around a $1$-line soliton and
a $2$-line soliton.  However, as far as we consider elastic line
solitons, unstable modes of linearized equations should be exponentially
localized in a longitudinal direction if they exist.
We can prove linear stability of
line solitons by using the Darboux transformation and an extra exponential
decay in a longitudinal direction that unstable modes should have.
\par
Kadomtsev and Petviashvili also conclude instability of the KdV $1$-soliton
in the media of positive dispersion by using the KP-I equation.
Its rigorous instability results were proved by \cite{Brides},
\cite{Rousset-Tzvetkov1,Rousset-Tzvetkov2} and \cite{Z}.
\par
Using the inverse scattering method, Villarroel and Ablowitz \cite{VA} obtained
a solution formula for solutions around $1$-line soliton.
The inverse scattering method for solutions around multi-line solitons
has been studied by by Boiti, Pempinelli and Pogrebkov. See e.g. \cite{BPP} and the
reference therein.
Recently, Derchyi Wu has justified their method and proves that their solution formula
gives a classical solution of \eqref{eq:KPII} and its deviation from
the exact line soliton solution remains small in the $L^\infty$-norm for all the time.
\par
Our paper extends the result by \cite{Burtsev} to an elastic $2$-line
solitons.  Although our paper is limited to a special class of
$2$-line solitons, our linear stability results might be used
to investigate persistence of multi-line solitary waves of non
integrable $2$D long wave models studied by \cite{AbCu13} as the
linear stability of $1$-line solitons for the KP-II equation was
crucial to prove stability of line solitary waves for the Benney Luke
equation (\cite{Miz-Shimabukuro17,Miz-Shimabukuro20}).

\bigskip

\section{Preliminaries}
\subsection{Line soliton solutions}
It is well known that multi-line solutions of the KP-II
equation can be expressed by the Wronskian formulation.
Let $f_i$ $(1\le i\le N)$ be linearly independent solutions of
the system
\begin{equation}
  \label{eq:f1-N}
 \frac{\pd f}{\pd y}=\frac{\pd^2 f}{\pd x^2}\,,
\quad \frac{\pd f}{\pd t}=-\frac{\pd^3 f}{\pd x^3}\,,
\end{equation}
and 
\begin{equation}
  \label{eq:tau}
\tau(t,x,y)=\operatorname{Wr}(f_1,\cdots,f_N)=
\begin{vmatrix}
f_1 & f_2 & \ldots & f_N\\
\pd_xf_1 & \pd_xf_2 & \ldots & \pd_xf_N \\
\vdots & \vdots & \ddots & \vdots \\
\pd_x^{N-1}f_1 & \pd_x^{N-1}f_2 & \ldots & \pd_x^{N-1}f_N
\end{vmatrix}\,.  
\end{equation}
Then
\begin{equation}
  \label{eq:exact}
u(t,x,y)=2\frac{\pd_x^2}{\pd x^2}\log \tau(t,x,y)  
\end{equation}
gives an exact solution of \eqref{eq:KPII}.
See e.g. \cite{Hirota,Fr-Nim}.

The $\tau$-function of line soliton solutions of \eqref{eq:KPII} can be
expressed by the Wronskian determinant of
$$f_n=\sum_{m=1}^Ma_{mn}e^{\theta_m}\,,
\quad  \theta_m(t,x,y)=\k_mx+\k_m^2y-\k_m^3t\,,
\quad 1\le n\le N\,,$$
where $a_{mn}$ are real constants,
and $\{\k_m\}_{m=1}^M$ are constants satisfying
\begin{equation}
  \label{eq:orderkm}
 \k_1<\k_2<\cdots<\k_M\,. 
\end{equation}
Note that $f_n$ $(1\le n\le N)$ satisfy \eqref{eq:f1-N}.
Let $\Theta=\diag(e^{\theta_1},\cdots,e^{\theta_M})$ and $A=(a_{ij})$.
Then
\begin{gather}
\label{eq:tau-linesol}
\tau(t,x,y)=\operatorname{Wr}(f_1,\cdots,f_N)=K\Theta A\,,\quad
K=
\begin{pmatrix}
1 & 1 & \ldots & 1 \\
\k_1 & \k_2 & \ldots & \k_M \\
\vdots & \vdots & \ddots & \vdots \\
\k_1^{N-1} & \k_2^{N-1} & \ldots & \k_M^{N-1}
\end{pmatrix}\,.
\end{gather}
By the Binet-Cauchy formula,
$$
\tau(t,x,y)=\sum_{1\le m_1<\cdots<m_N}
A(m_1,\cdots,m_N)\exp(\theta_{m_1}+\cdots+\theta_{m_N})
\prod_{1\le s<r\le N}(\k_{m_r}-\k_{m_s})\,,$$
where $A(m_1,\cdots,m_N)$ is the $N\times N$ minor of $A$ obtained from
the $m_j$-th rows $(1\le j\le N)$.
To have regular line soliton solutions, all $N\times N$ minors of $A$
must be nonnegative and $\mathrm{rank}(A)=N$. 
\begin{remark}
\begin{enumerate}
\item If $N>M$, then $\tau(t,x,y)=0$.
\item Since \eqref{eq:exact} with \eqref{eq:tau-linesol} is invariant
under the transformation $A\mapsto AC$ with $C\in GL(N,\R)$,
we may assume that $A^{\mathrm{T}}$ is in the reduced row-echelon form without loss
of generality.
If $M=N$, we have $\tau=\exp(\sum_{m=1}^N\theta_m)$ and $u(t,x,y)\equiv0$.
\item
Suppose that $f_n=e^{\theta_m}$ for some $m$ and $n$.
That is, the $n$-th column of $A$ has only one nonzero entry
$a_{mn}=1$. Then $\tau(t,x,y)=e^{\theta_m}\tau_0(t,x,y)$, where $\tau_0$
is a $\tau$-function consists of $e^{\theta_j}$ $(j\ne m)$.
It is clear from \eqref{eq:exact} that $\tau(t,x,y)$ and $\tau_0(t,x,y)$
generates the same solution of \eqref{eq:KPII}.
  \end{enumerate}
\end{remark}
To have regular line soliton solutions and to avoid redundant expressions,
\begin{enumerate}
\item $\mathrm{rank}(A)=N<M$ and all minors of $A$ must be nonnegative,
\item each row of $A$ contains at least one nonzero element and each
row of $A^{\mathrm{T}}$ contains at least one nonzero element in addition to the
pivot.
\end{enumerate}
Then \eqref{eq:exact} with \eqref{eq:tau-linesol} is a multi-line soliton
solution that has $M-N$ solitons for $y\ll -1$ and $N$ solitons for $y\gg 1$.
For example, if $M=2$, $N=1$ and $\tau(t,x,y)=e^{\theta_1}+e^{\theta_2}$,
then we have a $1$-line soliton solution
$$ u(t,x,y)=\frac12(\k_2-\k_1)^2\sech^2\frac{\theta_2-\theta_1}{2}\,.$$
Since
$\theta_2-\theta_1=(\k_2-\k_1)\{x+(\k_1+\k_2)y-(\k_1^2+\k_1\k_2+\k_2^2)t\}$,
the line soliton moves toward $(1,\k_1+\k_2)$.
We refer the reder to \cite{CK} and a book \cite{K-book} and
the references therein.
\bigskip

\subsection{Elastic $2$-line solitons}
Generically, $2$-line soliton solutions
have a web like pattern in the middle thanks to
resonant interactions among line solitons.
Here we consider special cases that do not have resonant 
interactions. 

First, we consider $2$-line soliton solutions of P-type.
Let
\begin{equation}
\label{eq:def-tau2line}
  \Theta=\diag(e^{\theta_1},e^{\theta_2},e^{\theta_3},e^{\theta_4})\,,\quad
 K=\begin{pmatrix}1 & 1 & 1 & 1
   \\ \k_1 & \k_2 & \k_3 & \k_4 \end{pmatrix}\,,
\end{equation}
and 
\begin{gather}
\label{eq:KA-P}
 A^{\mathrm{T}}=\begin{pmatrix}
  1 & 0 & 0 & -1 \\ 0 & 1 & 1 &0
\end{pmatrix}\,,
\\
\label{eq:P-type}
\tau=\det(K\Theta A)=
\begin{vmatrix}
e^{\theta_2}+e^{\theta_3} &   e^{\theta_4}-e^{\theta_1}\\
\k_2e^{\theta_2}+\k_3e^{\theta_3} &   \k_4e^{\theta_4}-\k_1e^{\theta_1}
\end{vmatrix}\,.  
\end{gather}
Let
\begin{equation}
  \label{eqdef:a,c,omega}
a_{ij}=\frac{\k_i+\k_j}{2}\,,\quad c_{ij}=\frac{(\k_i-\k_j)^2}{4}\,,\quad
\omega_{ij}=\k_i^2+\k_i\k_j+\k_j^2\,.  
\end{equation}
Then for $\pm (a_{14}-a_{23})y\gg 1$,
the asymptotics of $u=2\pd_x^2\log\tau$ is given by
\begin{equation}
  \label{eq:u2-asymptotics}\left\{
  \begin{gathered}
u(t,x,y)\simeq 2c_{23}\sech^2Z_{1,\pm}+ 2c_{14}\sech^2 Z_{2,\mp}\,,
    \\
Z_{1,\pm}=\frac{\theta_4-\theta_1}{2}+\mu_{1,\pm}\,,\quad
Z_{2,\pm}=\frac{\theta_3-\theta_2}{2}+\mu_{2,\pm}\,,
\end{gathered}\right.
\end{equation}
where
\begin{gather*}
\mu_{1,+}=\frac{1}{2}\log\frac{\k_4-\k_3}{\k_4-\k_2}\,,\enskip
\mu_{1,-}=\frac12\log\frac{\k_3-\k_1}{\k_2-\k_1}\,,\\
\mu_{2,+}=\frac12\log\frac{\k_4-\k_3}{\k_3-\k_1}\,,\enskip
\mu_{2,-}=\frac12\log\frac{\k_4-\k_2}{\k_2-\k_1}\,.  
\end{gather*}
Since $\theta_i-\theta_j=(\k_i-\k_j)(x+2a_{ij}y-\omega_{ij}t)$,
we see from \eqref{eq:u2-asymptotics} that $u$ consists of line solitons
moving in the directions of $(1,2a_{23})$ and $(1,2a_{14})$.
We call each of them as $[2,3]$-soliton and $[1,4]$-soliton, respectively.
Both $[2,3]$-soliton and $[1,4]$-soliton of the $2$-soliton solution
of P-type can be placed in nearly parallel to the $y$-axis and the
each height of solitons cannot be the same.
Let
\begin{equation}
  \label{eq:b1-b2P}
X=x-b_1t\,,\quad Y=y-b_2t\quad\text{with}\quad b_1+2a_{23}b_2=\omega_{23}\,,\quad b_1+2a_{14}b_2=\omega_{14}\,.
\end{equation}
That is
\begin{equation}
  \label{eq:b1-b2P'}
b_1=-\left(4a_{23}a_{14}+\frac{a_{23}\k_1\k_4-a_{14}\k_2\k_3}{a_{23}-a_{14}}\right)\,,
\quad  
b_2=2(a_{23}+a_{14})+\frac{\k_1\k_4-\k_2\k_3}{2(a_{23}-a_{14})}\,.
\end{equation}
Then
\begin{equation}
  \label{eq:t3-t2}
\theta_3-\theta_2=(\k_3-\k_2)(X+2a_{23}Y)\,,
\quad\theta_4-\theta_1=(\k_4-\k_1)(X+2a_{14}Y)\,,  
\end{equation}
and the $2$-soliton solution $u=2\pd_x^2\log\tau$ with
\eqref{eq:P-type} is a stationary solution in the moving coordinate
$(X,Y)$.

Next, we consider $2$-line soliton solution of O-type.
Let $\Theta$ and $K$ be as \eqref{eq:def-tau2line} and
\begin{align*}
 A^{\mathrm{T}}=\begin{pmatrix}
  1 & 1 & 0 & 0 \\ 0 & 0 & 1 &1
\end{pmatrix}\,.
\end{align*}
Then
\begin{equation}
  \label{eq:O-type}
\tau=\det(K\Theta A)=
\begin{vmatrix}
e^{\theta_1}+e^{\theta_2} &   e^{\theta_3}+e^{\theta_4}\\
\k_1e^{\theta_1}+\k_2e^{\theta_2} &   \k_3e^{\theta_3}+\k_4e^{\theta_4}
\end{vmatrix}\,,
\end{equation}
and for $\pm y\gg 1$,
the asymptotics of $u=2\pd_x^2\log\tau$ is given by
\begin{equation}
  \label{eq:u2-asymptoticsO}
\left\{  \begin{gathered}
u(t,x,y) \simeq 
2c_{12}\sech^2Z_{1,\pm}+ 2c_{34}\sech^2Z_{2,\mp}
\\
Z_{1,\pm}=\frac{\theta_2-\theta_1}{2}+\mu_{1,\pm}\,,\quad
Z_{2,\pm}=\frac{\theta_4-\theta_3}{2}+\mu_{2,\pm}\,,
\end{gathered}\right.
\end{equation}
where
\begin{gather*}
\mu_{1,+}=\frac{1}{2}\log\frac{\k_4-\k_2}{\k_4-\k_1}\,,\enskip
\mu_{1,-}=\frac12\log\frac{\k_3-\k_2}{\k_3-\k_1}\,,\\
\mu_{2,+}=\frac12\log\frac{\k_4-\k_2}{\k_3-\k_2}\,,\enskip
\mu_{2,-}=\frac12\log\frac{\k_4-\k_1}{\k_3-\k_1}\,.  
\end{gather*}
The asymptotic line solitons of $u=2\pd_x^2\log\tau$ for $\pm y\gg 1$
move toward $(1,2a_{12})$ and $(1,2a_{34})$.  We
call them $[1,2]$-soliton and $[3,4]$-soliton, respectively.  The
$2$-soliton solution of O-type becomes singular if both
$[1,2]$-soliton and $[3,4]$-soliton are nearly parallel to the
$y$-axis.
Let
\begin{equation}
  \label{eq:b1-b2O}
X=x-b_1t\,,\quad Y=y-b_2t\quad\text{with}\quad
b_1+2a_{12}b_2=\omega_{12}\,,\quad b_1+2a_{34}b_2=\omega_{34}\,.
\end{equation}
Then $\theta_2-\theta_1=(\k_2-\k_1)(X+2a_{12}Y)$, $\theta_4-\theta_3=(\k_4-\k_3)(X+2a_{34}Y)$
and $u=2\pd_x^2\log\tau$ with \eqref{eq:O-type} is a stationary solution
in the moving coordinate $(X,Y)$.
\bigskip

\subsection{B\"acklund transformations}
Let $u=\pd_xw$. Then \eqref{eq:KPII} can be transformed into
\begin{equation}
  \label{eq:pkp}
4\pd_tw+\pd_x^3w+3(\pd_xw)^2+3\pd_x^{-1}\pd_y^2w=0\,.  
\end{equation}
The auto B\"acklund transformation of \eqref{eq:pkp} can be written as
\begin{equation}
  \label{eq:BT1}
\pd_x^{-1}\pd_y(w_2-w_1)=\pd_x(w_2+w_1)+\frac{1}{2}(w_2-w_1)^2\,,
\end{equation}
\begin{equation}
  \label{eq:BT2}
  \begin{split}
    & (4\pd_t+\pd_x^3)(w_2-w_1)+3\pd_x\pd_y(w_2+w_1)
\\ \quad & + \pd_x\{3(w_2-w_1)\pd_x(w_2+w_1)+(w_2-w_1)^3\}=0\,.
  \end{split}
\end{equation}
If $w_1$ and $w_2$ satisfy \eqref{eq:BT1} and $w_1$ is a solution of
\eqref{eq:pkp}, then $w_2$ is a solution of \eqref{eq:pkp}
See e.g. \cite{Chen75, Hirota}.
In the bilinear form, the B\"acklund transformation \eqref{eq:BT1} can be written as
\begin{equation}
\label{eq:BT-f}
\tau_1\pd_x^2\tau_2-2\pd_x\tau_2\pd_x\tau_1+\tau_2\pd_x^2\tau_1
+\tau_2\pd_y\tau_1-\tau_1\pd_y\tau_2=0\,.
\end{equation}
\begin{equation}
  \label{eq:BT-f2}
  \begin{split}
&  4(\tau_1\pd_t\tau_2-\tau_2\pd_t\tau_1)+3(\tau_1\pd_x\pd_y\tau_2
-\pd_x\tau_1\pd_y\tau_2-\pd_y\tau_1\pd_x\tau_2+\tau_2\pd_x\pd_y\tau_1)
\\ & \qquad
+\tau_1\pd_x^3\tau_2-3\pd_x\tau_1\pd_x^2\tau_2+3\pd_x^2\tau_1\pd_x\tau_2-\tau_2\pd_x^3\tau_1
=0\,.
  \end{split}
\end{equation}
Let $w_i=2\pd_x\log\tau_1$ and $\pd_x^{-1}\pd_yw_i=2\pd_y\tau_i/\tau_i$
for $i=1$ and $2$.
Then
\begin{align*}
\pd_x^{-1}\pd_y(w_2-w_1)=&
2\left(\frac{\pd_y\tau_2}{\tau_2}-\frac{\pd_y\tau_1}{\tau_1}\right)
\\=& 2\left\{\pd_x\left(\frac{\pd_x\tau_2}{\tau_2}
+\frac{\pd_x\tau_1}{\tau_1}\right)
+\left(\frac{\pd_x\tau_2}{\tau_2}-\frac{\pd_x\tau_1}{\tau_1}\right)^2\right\}
\\=& \pd_x(w_2+w_1)+\frac12(w_2-w_1)^2\,.
\end{align*}

Via the B\"acklund transformation \eqref{eq:BT-f},
$N$-soliton solutions can be connected to $N-1$-soliton solutions.
See \cite{Fr-Nim}. In the following examples, $\tau_1$ and $\tau_2$ satisfy
\eqref{eq:BT-f}.
\begin{example}
\label{ex:1}
Suppose that  $u_i=2\pd_x^2\log\tau_i$ for $i=1$ and $2$.
  \begin{enumerate}[(i)]
\item (a $1$-line soliton and the null solution)
Let 
  \begin{equation}
    \label{eq:tau-1-sol}
    \tau_1=1\,,\quad\tau_2=e^{\theta_1}+e^{\theta_2}\,,\quad
u_1=0\,,\quad u_2=2c_{12}\sech^2\frac{\theta_2-\theta_1}{2}\,.
\end{equation}
Then $u_1=0$ and $u_2$ is a $1$-line soliton.
\item (a resonant line soliton and the null solution)
  Let $\tau_1=1$ and $\tau_2$ be as
\begin{equation}
  \label{eq:tauY}
e^{\theta_1}+ae^{\theta_2}+be^{\theta_3}\,,\quad a>0\,,\;b>0\,.  
\end{equation}
Then $u_1=0$ and $u_2$ is a resonant line soliton.

\item (a generic $2$-line soliton and a resonant line soliton)
  Let $\Theta$ and $K$ be as \eqref{eq:def-tau2line}.
Generic $2$-line solitons can be expressed as $u_2=2\pd_x^2\log\tau_2$,
where
\begin{equation}
  \label{eq:defA}
\tau_1=e^{\theta_2}+ae^{\theta_3}+be^{\theta_4}\,,\enskip
  \tau_2=\det(K\Theta A)\,,\enskip
  A^{\mathrm{T}}=
\begin{pmatrix}
  1 & 0 & -c & -d\\ 0 & 1 & a & b
\end{pmatrix}\,,
\end{equation}
$a$, $b$, $c$, $d>0$ and $ad-bc>0$.
\item (a $2$-line soliton of P-type and a $1$-line soliton)
Let $b=c=0$ and let  $\tau_1$ and $\tau_2$ be as \eqref{eq:defA}. Then
  \begin{equation}
    \label{eq:tau-P}
\tau_1=f_1\,,\quad \tau_2=\operatorname{Wr}(f_1,f_2)\,,\quad
f_1=e^{\theta_2}+ae^{\theta_3}\,,\quad f_2=de^{\theta_4}-e^{\theta_1}\,,
\end{equation}
and $u_1$ is a $1$-line soliton and 
$u_2$ is a $2$-line soliton of P-type.
\item (a $2$-line soliton of O-type and a $1$-soliton)
Let $a>0$, $b>0$ and $c\ge0$. Let  $\tau_1=e^{\theta_1}+ae^{\theta_2}$ and
  \begin{equation}
    \label{eq:degenerateA}
\tau_2=\det(K\Theta A)\,,\quad 
    A^{\textrm{T}}=
\begin{pmatrix}
  1 & b & 0 & -c\\ 0 & 0& 1 & a
\end{pmatrix}\,.
\end{equation}
Then $u_1$ is a $1$-line soliton and $u_2$ is a $2$-line soliton.
Especially if $c=0$,
  \begin{equation}
  \label{eq:tau-O}
 \tau_2=\operatorname{Wr}(f_1,f_2)\,,\quad
f_1=e^{\theta_1}+be^{\theta_2}\,,\quad f_2=e^{\theta_3}+ae^{\theta_4}\,,
\end{equation}
and $u_2$ is an elastic $2$-line soliton of O-type.
If we set $\tau_1=f_1$ or $\tau_1=f_2$, then $u_1$ is a $1$-line soliton.
\item (a resonant soliton and a $1$-line soliton)
Let $a>0$, $b>0$ and
  \begin{equation}
    \label{eq:tauY'}
\tau_2=K\diag(e^{\theta_1},e^{\theta_2},e^{\theta_3})A\,,\quad K=
\begin{pmatrix}
  1 & 1 & 1\\ \k_1 & \k_2 & \k_3
\end{pmatrix}\,,
\\ 
A^{\mathrm{T}}=
\begin{pmatrix}
  1 & 0 & -b \\ 0 & 1& a
\end{pmatrix}\,.
  \end{equation}
Since $\tau_2=\operatorname{Wr}(f_1,f_2)$ with 
$f_1=e^{\theta_1}-be^{\theta_3}$ or $f_1 =e^{\theta_1}+\tfrac{b}{a}e^{\theta_2}$
and $f_2=e^{\theta_2}+ae^{\theta_3}$,
the case \eqref{eq:tauY'} can be regarded as a degenerate case of \eqref{eq:tau-P} and \eqref{eq:tau-O}.

\item (a resonant soliton and $1$-line soliton)
  Let $\tau_1$ be as \eqref{eq:tauY}, $\tau_2=\operatorname{Wr}(f_1,f_2)$ with
$f_1=e^{\theta_1}+ae^{\theta_2}$ and $f_2=\tau_1$
or  $f_1=\tau_1$ and  $f_2=ae^{\theta_2}+be^{\theta_3}$.
  Then
  \begin{equation*}
    u_2=\frac{(\k_2-\k_1)^2}{2}   \sech^2
      \left(\frac{\theta_2-\theta_1}{2}+\mu\right)\,,
      \quad\mu=\frac12\log a\left(\frac{\k_3-\k_2}{\k_3-\k_1}\right)\,,
    \end{equation*}
 or 
  \begin{equation*}
    u_2=\frac{(\k_3-\k_2)^2}{2}\sech^2\left(\frac{\theta_3-\theta_2}{2}+\mu\right)\,,\quad\mu=\frac12\log\frac{b}{a}\left(\frac{\k_3-\k_1}{\k_2-\k_1}\right)\,.
    \end{equation*}
\end{enumerate}
  \end{example}
\bigskip

\subsection{Relations with the Miura transformation and the modified KP-II}
Now we recall the Miura transformations of the KP-II equation.
Let 
\begin{gather*}
e_{KP}(u)=4\pd_tu+\pd_x^3u+3\pd_x(u^2)+3\pd_x^{-1}\pd_y^2u\,,\\
e_{MKP}(v)=4\pd_tv+\pd_x^3v+3\pd_x^{-1}\pd_y^2v
-6v^2\pd_xv+6(\pd_xv)(\pd_x^{-1}\pd_yv)\,,\\
M_\pm(v)=\pm\pd_xv+\pd_x^{-1}\pd_yv-v^2\,,\\
\widetilde{M}_\pm(u,v)=(4\pd_x^{-1}\pd_t+\pd_x^2)v+6uv\mp3\pd_x(\pd_x^{-1}\pd_yv+v^2)+4v^3\,.
\end{gather*}
Formally, we have by a straightforward computation,
\begin{equation}
  \label{eq:kp-mkp}
  \begin{split}
& e_{KP}(u)=\nabla M_\pm(v)e_{MKP}(v)\,,\quad u=M_\pm(v)\,.
  \end{split}
\end{equation}
Thus if $v$ is a solution to the modified KP-II equation
\begin{equation}
  \label{eq:mkp}
4\pd_tv+\pd_x^3v+3\pd_x^{-1}\pd_y^2v-6v^2\pd_xv+6(\pd_xv)(\pd_x^{-1}\pd_yv)=0\,,
\end{equation}
then $u=M_\pm(v)$ is a solution of \eqref{eq:KPII}.
\par
Let $u_1=\pd_xw_1$, $u_2=\pd_xw_2$ and $v=(w_2-w_1)/2$.
Then \eqref{eq:BT1} can be read as the Miura transformations
\begin{gather}
\label{eq:Miura1}
u_1=M_-(v)\,,\quad u_2=M_+(v)\,,
\end{gather}
and \eqref{eq:BT2} can be read as
\begin{gather}
\label{eq:Miura-t1}
3\pd_x^{-1}\pd_yu_1+\widetilde{M}_-(u_1,v)=0\,,\quad
3\pd_x^{-1}\pd_yu_2+\widetilde{M}_+(u_2,v)=0\,.
\end{gather}
Suppose that $v$ is a solution of \eqref{eq:mkp}
and that $u_1$ and $u_2$ satisfy \eqref{eq:Miura1}.
Then it follows from \eqref{eq:kp-mkp} that
$u_1$ and $u_2$ are solutions of \eqref{eq:KPII}.
\par
On the other hand, if $w_1$ and $w_2$ are solutions of \eqref{eq:pkp}
satisfying \eqref{eq:BT1}, 
then $v=(w_2-w_1)/2$ is a solution of \eqref{eq:mkp}.
Indeed, by \eqref{eq:pkp} and \eqref{eq:BT1}, 
\begin{align*}
 4\pd_tv+\pd_x^3v+3\pd_x^{-1}\pd_y^2v
=& \frac{3}{2}\left\{(\pd_xw_1)^2-(\pd_xw_2)^2\right\}
  \\=& -3(\pd_xv)\pd_x(w_2+w_1)
\\=& 6v^2\pd_xv-6(\pd_xv)(\pd_x^{-1}\pd_yv)\,.
\end{align*}  
\par
Let $\tau_1$ and $\tau_2$ be as 
\eqref{eq:tau-P} or \eqref{eq:tau-O} in Example~\ref{ex:1} and
let
\begin{gather}
  \label{eq:u1u2v2}
u_1=2\pd_x^2\log\tau_1\,,\quad u_2=2\pd_x^2\log\tau_2\,,
\quad v_2=\pd_x\log\frac{\tau_2}{\tau_1}\,,
\quad \pd_x^{-1}\pd_yv_2=\pd_y\log\frac{\tau_2}{\tau_1}\,,
\\ \label{eq:v1}
v_1=\pd_x\log\tau_1\,, \quad \pd_x^{-1}\pd_yv_1=\pd_y\log\tau_1\,.  
\end{gather}
Then $u_2=M_+(v_2)$, $u_1=M_-(v_2)=M_+(v_1)$ and $ M_-(v_1)=0$.
\bigskip

\subsection{Linearized Miura transformations}
Let $e_{LKP}$ be a linearized operator of \eqref{eq:KPII} and 
$e_{LMKP}$ be a linearized operator of \eqref{eq:mkp}.
More precisely,
\begin{align*}
e_{LKP}(u):=&4\pd_t+\pd_x^3+6\pd_x(u\cdot)+3\pd_x^{-1}\pd_y^2\,,\\
e_{LMKP}(v):=&4\pd_t+\pd_x^3+3\pd_x^{-1}\pd_y^2-6\pd_x(v^2\cdot)
+6(\pd_xv)(\pd_x^{-1}\pd_y\cdot)+6(\pd_x\cdot)(\pd_x^{-1}\pd_yv)\,.
\end{align*}
Suppose that $u=M_\pm(v)$. Then by \eqref{eq:kp-mkp},
\begin{align*}
e_{LKP}(u)\nabla M_\pm(v)w=& \frac{d}{d\eta}e_{KP}(M_\pm(v+\eta w))\bigr|_{\eta=0}
\\=& \frac{d}{d\eta}\left\{\nabla M_\pm(v+\eta w)e_{MKP}(v+\eta w)\right\}\bigr|_{\eta=0}
\\=& \nabla M_\pm(v)e_{LMKP}(v)w-2we_{MKP}(v)\,.
\end{align*}
If $v$ is a solution of \eqref{eq:mkp},
\begin{equation}
  \label{eq:bH}
e_{LKP}(M_\pm(v))\nabla M_\pm(v)=\nabla M_\pm(v)e_{LMKP}(v)\,.
\end{equation}
\par
Especially, if  $\tau_1$ and $\tau_2$ are $\tau$-functions satisfying
\eqref{eq:BT-f} as in Example~\ref{ex:1} and  $u_1$, $u_2$ and $v_2$
are as \eqref{eq:u1u2v2}, then
\begin{gather}
  \label{eq:bH+}
e_{LKP}(u_2)\nabla M_+(v_2)=\nabla M_+(v_2)e_{LMKP}(v_2)\,,\\
\label{eq:bH-}
e_{LKP}(u_1)\nabla M_-(v_2)=\nabla M_-(v_2)e_{LMKP}(v_2)\,.  
\end{gather}
That is, if $W$ is a solution of a linearized MKP-II equation around $v_2$,
then $U_+=\nabla M_+(v_2)W$ is a solution of a linearized KP-II equation
around $u_2$,
and $U_-=\nabla M_-(v_2)W$ is a solution of a linearized KP-II equation around $u_1$.
\par
Roughly speaking, we will use the fact that
$\nabla M_+(v_2)(\nabla M_+(v_2))^{-1}$ gives an isomorphism between
solution spaces of $e_{LKP}(u_1)U=0$ and $e_{LKP}(u_2)U=0$
to prove linear stability properties of line soliton solutions.
\bigskip

\subsection{The Lax operator and the linearized Miura transformations}
As is well known, the KP-II equation \eqref{eq:KPII} can be expressed as
$[B,L]=0$
by using a Lax pair
\begin{equation}
  \label{eq:Laxpair}
L=-\pd_y+\pd_x^2+u\,, \quad
B=4\pd_t+4\pd_x^3+6u\pd_x+3\pd_xu+3\pd_x^{-1}\pd_yu\,.
\end{equation}
Let us introduce the relationship between the Lax operator of the
KP-II equation and the linearized Miura transformation around
solitons.
\begin{lemma}
  \label{lem:Lax-Miura}
Let $\tau_1$ and $\tau_2$ be $\tau$-functions satisfying \eqref{eq:BT-f}.
Let $u_1$, $u_2$, $v_1$ and $v_2$ be as \eqref{eq:u1u2v2}, \eqref{eq:v1}
and let $h=\tau_2/\tau_1$. Let
\begin{gather*}
L_i=-\pd_y+\pd_x^2+u_i\,,\quad B_i=4\pd_t+4\pd_x^3+6u_i\pd_x+3\pd_xu_i+3\pd_x^{-1}\pd_yu_i\,.
\end{gather*}
Then we have the following.
\begin{align}
&  \label{eq:Miura-Lax1}
\nabla M_+(v_2)=hL_2^*h^{-1}\pd_x^{-1}
=\pd_x^{-1}hL_1^*h^{-1}\,,
\\ & \label{eq:Miura-Lax2}
\nabla M_-(v_2)=-h^{-1}L_1h\pd_x^{-1}
=-\pd_x^{-1}h^{-1}L_2h\,,
\end{align}
\begin{align}
  \label{eq:Miura-Laxt1}
\nabla_v\widetilde{M}_+(u_2,v_2)=& 4(\pd_x^{-1}\pd_t+\pd_x^2)+6u_1-12v_2\pd_x+12v_2^2
  \\=& \notag -\pd_x^{-1}hB_1^*h^{-1}=-hB_2^*h^{-1}\pd_x^{-1}\,,
\end{align}
\begin{align}
  \label{eq:Miura-Laxt2}
\nabla_v\widetilde{M}_-(u_1,v_2)=& 4(\pd_x^{-1}\pd_t+\pd_x^2)+6u_2+12v_2\pd_x+12v_2^2
  \\=& \notag \pd_x^{-1}h^{-1}B_2h=h^{-1}B_1h\pd_x^{-1}\,.
\end{align}
\end{lemma}
\begin{proof}
By \eqref{eq:BT-f},
\begin{gather}
  \label{eq:h}
(\pd_y-\pd_x^2)h=
2h\left\{\frac{\pd_x^2\tau_1}{\tau_1}-\left(\frac{\pd_x\tau_1}{\tau_1}
\right)^2\right\}=hu_1\,,
\\
 \label{eq:th}
(\pd_x^2+\pd_y)h^{-1}=
2h^{-1}\left\{\left(\frac{\pd_x\tau_2}{\tau_2}\right)^2
-\frac{\pd_x^2\tau_2}{\tau_2}\right\}=-h^{-1} u_2\,.
\end{gather}
Let  for $i=1$ and $2$.
Using \eqref{eq:h}, \eqref{eq:th} and the fact that $\pd_xh=hv_2$,
we have \eqref{eq:Miura-Lax1} and \eqref{eq:Miura-Lax2}.
\par
By the definitions of $v_2$, $u_2$ and $u_1$, we have
$u_2-u_1=2\pd_xv_2$, $\pd_x^{-1}\pd_yv_2=(u_1+u_2)/2+v_2^2$ and
  \begin{multline*}
 \frac{\pd_x\pd_y\tau_2}{\tau_2}
      -\frac{\pd_x\tau_1\pd_y\tau_2}{\tau_1\tau_2}
      -\frac{\pd_y\tau_1\pd_x\tau_2}{\tau_1\tau_2}
      +\frac{\pd_x\pd_y\tau_1}{\tau_1}
=   \pd_yv_2+v_2\pd_x^{-1}\pd_yv_2+\pd_x^{-1}\pd_yu_1
    \\ =\frac{1}{2}\pd_x(u_2+u_1)+\frac{3}{2}u_2v_2-\frac12u_1v_2+v_2^3
         +\pd_x^{-1}\pd_yu_1\,,
  \end{multline*}
\begin{align*}
 \frac{\pd_x^3\tau_2}{\tau_2}-3\frac{\pd_x\tau_1\pd_x^2\tau_2}{\tau_1\tau_2}
+3\frac{\pd_x^2\tau_1\pd_x\tau_2}{\tau_1\tau_2}
 -\frac{\pd_x^3\tau_1}{\tau_1}
   =&
\frac{1}{2}\pd_x(u_2-u_1)+\frac{3}{2}v_2(u_1+u_2)+v_2^3\,.
\end{align*}
Substituting the above into \eqref{eq:BT-f2}, we have
\begin{equation*}
4h^{-1}\pd_th+\pd_x(2u_2+u_1)+3\pd_x^{-1}\pd_yu_1+6v_2u_2+4v_2^3 =0\,.
\end{equation*}
Since $\pd_xh=hv_2$ and
  \begin{align*}
& h^{-1}\pd_x^2h=v_2^2+\frac{1}{2}(u_2-u_1)\,,
\quad
h^{-1}\pd_x^3h=v_2^3+\frac{3}{2}v_2(u_2-u_1)+\frac12\pd_x(u_2-u_1)\,.
  \end{align*}
Combining the above, we have
  \begin{align*}
h^{-1}B_2h=& B_2+ 
\frac{1}{h}\left([\pd_t,h]+4[\pd_x^3,h]+6u_2[\pd_x,h]\right)
\\=& 4(\pd_t+\pd_x^3)+ 12v_2\pd_x^2+6(2v_2^2+2u_2-u_1)\pd_x
\\ & +6v_2(u_2-u_1)+3(\pd_x+\pd_x^{-1}\pd_y)(u_2-u_1)\,.
  \end{align*}
Since $u_2-u_1=2\pd_xv_2$ and
$(\pd_x+\pd_x^{-1}\pd_y+2v_2)(u_2-u_1) = 2\pd_x(u_2+2v_2^2)$, we have
\begin{equation*}
h^{-1}B_2h=\pd_x\nabla_v\widetilde{M}_-(u_1,v_2)\,.
\end{equation*}
We can prove the rest in the same way.  
\end{proof}
\begin{remark}
Let $L_0=-\pd_y+\pd_x^2$ and $B_0=4\pd_t+4\pd_x^3$.
Since  $\pd_y\tau_1=\pd_x^2\tau_1$,
\begin{align}
&  \label{eq:Miura-Lax3}
\nabla M_+(v_1)=\tau_1L_1^*\tau_1^{-1}\pd_x^{-1}
=\pd_x^{-1}\tau_1L_0^*\tau_1^{-1}\,,
\\ & \label{eq:Miura-Lax4}
\nabla M_-(v_1)=-\tau_1^{-1}L_0\tau_1\pd_x^{-1}
 =-\pd_x^{-1}\tau_1^{-1}L_1\tau_1\,,
\end{align}
\begin{align}
\label{eq:Miura-Laxt3}
\nabla_v\widetilde{M}_+(u_1,v_1)=& 4(\pd_x^{-1}\pd_t+\pd_x^2)-12v_1\pd_x+12v_1^2
\\=& \notag -\pd_x^{-1}\tau_1B_0^*\tau_1^{-1}=-\tau_1B_1^*\tau_1^{-1}\pd_x^{-1}\,,
\end{align}
\begin{align}
  \label{eq:Miura-Laxt4}
\nabla_v\widetilde{M}_-(0,v_1)=&
4(\pd_x^{-1}\pd_t+\pd_x^2)+6u_1+12v_1\pd_x+12v_1^2
  \\=& \notag \pd_x^{-1}\tau_1^{-1}B_1\tau_1=\tau_1^{-1}B_0\tau_1\pd_x^{-1}\,.     
\end{align}  
\end{remark}

\bigskip

\subsection{Jost solutions and continuous eigenfunctions}
Let us recall the definitions of Jost solutions and
dual Jost solutions.
We say that $\Phi$ is a Jost solution if $\Phi$ satisfies
\begin{equation}
\label{eq:Jost-eqs}
L\Phi=0\,,\quad B\Phi=0\,,
\end{equation}
and $\Phi^*$ is a dual Jost solution if $\Phi^*$ satisfies
\begin{gather}
\label{eq:dualJost-eqs}
L^*\Phi^*=0\,,\quad B^*\Phi^*=0\,,
\\ \notag
L^*=\pd_y+\pd_x^2+u\,,\quad
B^*=-4\pd_t-4\pd_x^3-6u\pd_x-3\pd_xu+3\pd_x^{-1}\pd_yu\,,
\end{gather}
where $\pd_x^{-1}\pd_yu=2\pd_x\pd_y\log\tau$.

The $x$ derivative of a product of a Jost solution and a dual Jost solution
is a solution of the linearized KP-II equation (see e.g. \cite{Lipov}).
\begin{lemma}
\label{lem:prodJdJ}
Let $\Phi$ be a solution of \eqref{eq:Jost-eqs}
and $\Phi^*$ be a solution of \eqref{eq:dualJost-eqs}.
Then $v=\pd_x(\Phi\Phi^*)$ is a solution of 
\begin{equation}
  \label{eq:linearKP}
\pd_x(4\pd_tv+6\pd_x(uv)+\pd_x^3v)+3\pd_y^2v=0\,,
\end{equation}
and $w=\Phi\Phi^*$ satisfies the adjoint equation
\begin{equation}
  \label{eq:adjLKP}
\pd_x(4\pd_tw+6u\pd_xw+\pd_x^3w)+3\pd_y^2w=0\,.  
\end{equation}
\end{lemma}
\begin{proof}
Following the proof of \cite[Lemma~3.4]{Sachs},   
we will prove Lemma~\ref{lem:prodJdJ} for the sake of self-containedness.
By the definitions of $\Phi$ and $\Phi^*$,
\begin{align}
  \label{eq:JT'1}
& 4\pd_t\Phi+\pd_x^3\Phi+3u\pd_x\Phi+3\pd_x\pd_y\Phi
+3\left(\pd_x^{-1}\pd_yu\right)\Phi=0\,,
\\
\label{eq:JT'2}
& 4\pd_t\Phi^*+\pd_x^3\Phi^*+3u\pd_x\Phi^*-3\pd_x\pd_y\Phi^*
-3\left(\pd_x^{-1}\pd_yu\right)\Phi^*=0\,,
\\ 
\label{eq:J-S}
& \pd_y\Phi=(\pd_x^2+u)\Phi\,,\quad \pd_y\Phi^*=-(\pd_x^2+u)\Phi^*\,.
\end{align}
By \eqref{eq:JT'1}--\eqref{eq:J-S},
$$\pd_yw=\pd_x(\Phi^*\pd_x\Phi-\Phi\pd_x\Phi^*)\,,$$
\begin{align*}
&  4\pd_tw+\pd_x^3w+6u\pd_xw
\\=& \Phi(4\pd_t\Phi^*+\pd_x^3\Phi^*+6u\pd_x\Phi^*)
+\Phi^*(4\pd_t\Phi+\pd_x^3\Phi+6u\pd_x\Phi)
+3(\pd_x^2\Phi\pd_x\Phi^*+\pd_x\Phi\pd_x^2\Phi^*)
\\=& 3\{(\pd_x^2\Phi)(\pd_x\Phi^*)+(\pd_x\Phi)(\pd_x^2\Phi^*)
+u(\Phi\pd_x\Phi^*+\Phi^*\pd_x\Phi)+\Phi\pd_x\pd_y\Phi^*-\Phi^*\pd_x\pd_y\Phi\}
\\=& 3\pd_y(\Phi\pd_x\Phi^*-\Phi^*\pd_x\Phi)\,.
\end{align*}
Combining the above, we see that $w$ is a solution of \eqref{eq:adjLKP} and that
$v=\pd_xw$ is a solution of \eqref{eq:linearKP}.
\end{proof}

Suppose that $u=2\pd_x^2\log\tau$ and that the $\tau$-function is represented as \eqref{eq:tau-linesol}.
Then the Jost solution and the dual Jost solution of the line soliton $u$
can be written as
\begin{equation}
  \label{eq:Jost-sol}
  \Phi(\bx,k)=e^{ikx-k^2y+ik^3t}
\frac{\tau_\Phi(\bx,k)}{\tau(\bx)}\,,\quad 
  \Phi^*(\bx,k)=e^{-ikx+k^2y-ik^3t}\frac{\tau_{\Phi^*}(\bx,k)}{\tau(\bx)}\,,
\end{equation}
where $\bx=(x,y,t)$ and the numerators $\tau_\Phi$ and $\tau_{\Phi^*}$ are the shifts of
the $\tau$-function:
\begin{gather*}
\tau_\Phi(\bx,k)=
\det\left(K\diag\left((ik-\k_j)e^{\theta_j}\right)_{1\le j\le M}A\right)\,,\\
\tau_{\Phi^*}(\bx,k)
=\det\left(K\diag\left((ik-\k_j)^{-1}e^{\theta_j}\right)_{1\le j\le M}A\right)\,.
\end{gather*}
Let
\begin{equation}
  \label{eq:Jost-zero-res}
\Phi_j(\bx):=\Phi(\bx,-i\k_j)\,,\quad
\Phi_j^*(\bx):=i\operatorname{Res}_{k=-i\k_j}\Phi^*(\bx,k)\,.
\end{equation}
We remark that 
$\sum_{j=1}^M \Phi_j(\bx)\Phi_j^*(\bxd)=0$. See \cite{CK,Dickey}. 
\bigskip

\section{Resonant modes for $2$-soliton solutions of P-type}
\label{sec:res2}
In this section, we express resonant modes of the linearized operator
around line soliton solutions by using the Jost solution and the dual
Jost solutions for the Lax pair \eqref{eq:Laxpair}.
\par

\subsection{Resonant modes for $1$-line soliton solutions}
\label{subsec:resonance-1}

First, we consider resonant continuous eigenfunctions of the linearized operator
around a $1$-line soliton
$$\varphi=2\pd_x^2\log\tau\,,\quad \tau=e^{\theta_1}+e^{\theta_2}\,.$$
\par
Let $a_{12}$, $c_{12}$ and $\omega_{12}$ be as \eqref{eqdef:a,c,omega} and 
$b_1$ and $b_2$ are real numbers satisfying $b_1+2a_{12}b_2=\omega_{12}$.
Let $X=x-b_1t$ and $Y=y-b_2t$. Then \eqref{eq:KPII} can be rewritten as
\begin{equation}
  \label{eq:KPII''}
  4\pd_tu+\pd_X^3u-4b_1\pd_Xu-4b_2\pd_Yu+3\pd_X(u^2)+3\pd_X^{-1}\pd_Y^2=0\,,
\end{equation}
and $\varphi$ is a stationary solution of \eqref{eq:KPII''}.
Linearizing \eqref{eq:KPII''} around $\varphi$, we have
$$\pd_t u=\mL u\,,\quad 
\mL=\frac14\{-\pd_X(\pd_X^2-4b_1+6\varphi)+4b_2\pd_Y+3\pd_X^{-1}\pd_Y^2\}\,.$$ 
\par
By \eqref{eq:tau-linesol} and \eqref{eq:Jost-sol} with $M=2$, $N=1$
and $A=(1,1)^{\mathrm{T}}$,
\begin{gather*}
\Phi(\bx,-i\beta)=e^{\beta x+\beta^2y -\beta^3t}
\left(\beta-\frac{\k_1e^{\theta_1}+\k_2e^{\theta_2}}{\tau}\right)\,,
\\
\Phi^*(\bx,-i\beta)=\frac{e^{-\beta x-\beta^2y +\beta^3t}}{\tau}
\left(\frac{e^{\theta_1}}{\beta-\k_1}+\frac{e^{\theta_2}}{\beta-\k_2}
\right)\,,
\\ \Phi(\bx)=\Phi(\bx,-i\k_1)=-\Phi(\bx,-i\k_2)
=(\k_1-\k_2)\frac{e^{\theta_1+\theta_2}}{\tau}\,,\quad
\Phi^*(\bx)=\frac{1}{\tau}\,.
\end{gather*}
Let
\begin{gather*}
\beta^\pm(\eta)=a_{12}\pm\gamma(\eta)\,,\quad
\gamma(\eta)=\sqrt{c_{12}+i\eta}\,,
\quad \lambda^\pm(\eta)=i\eta\{\gamma(\eta)\pm(3a_{12}-b_2)\}\,,
\\
Q^\pm(T,X,Y,\eta)=\Phi(\bx,-i\beta^\pm(\eta))\Phi^*(\bx)\,,\quad
\widetilde{Q}^\pm(T,X,Y,\eta)
=\eta\Phi(\bx)\Phi^*(\bx,-i\beta^\pm(\eta))\,.
\end{gather*}
Let $X_{12}=X+2a_{12}Y$. Since $\tau=2e^{(\theta_{1}+\theta_{2})/2}\cosh\sqrt{c_{12}}X_{12}$
and
\begin{equation*}
\beta^\pm(\eta)x+\beta^\pm(\eta)^2y-\beta^\pm(\eta)^3t-\frac{1}{2}(\theta_1+\theta_2)
= \pm \gamma(\eta)(X+2a_{12}Y) +i\eta Y\mp\lambda^\pm(\eta)T\,,
\end{equation*}
 \begin{gather*}
 \mL \pd_XQ^-(\eta)=\lambda^-(\eta)\pd_XQ^-(\eta)\,,\quad
 \mL \pd_X\widetilde{Q}^+(\eta)
=\lambda^+(\eta)\pd_X\widetilde{Q}^+(\eta)\,,\\
 \mL^*Q^+(\eta)=\lambda^+(\eta)Q^+(\eta)\,,\quad
 \mL^*\widetilde{Q}^-(\eta)=\lambda^-(\eta)\widetilde{Q}^-(\eta)\,.
\end{gather*}
Here $\beta^\pm(\eta)$ are chosen so that 
$e^{\a X_{12}}\pd_xQ$ and $e^{\a X_{12}}\pd_x\widetilde{Q}$ are bounded provided
$0<\a<2\sqrt{c_{12}}$ and $\eta\in\R$ is sufficiently close to $0$.
Indeed, $\Re\gamma(\eta)\simeq (\k_2-\k_1)/2+(\k_2-\k_1)^{-3}\eta^2$
around $\eta=0$,
\begin{gather*}
Q^-(\eta)=O(e^{-\Re\{\gamma(0)\pm\gamma(\eta)\}|X_{12}|})\,,
\quad \widetilde{Q}^+(\eta)=O(e^{-\Re\{\gamma(0)\pm\gamma(\eta)\}|X_{12}|})
\quad\text{as $X_{12}\to\pm\infty$,}\\
Q^+(\eta)=O(e^{-\Re\{\gamma(0)\mp\gamma(\eta)\}|X_{12}|})\,, \quad
\widetilde{Q}^-(\eta)=O(e^{-\Re\{\gamma(0)\mp\gamma(\eta)\}|X_{12}|})
\quad\text{as $X_{12}\to\pm\infty$.}
\end{gather*}
Suppose that $\mL$ is an operator on $L^2(\R^2;e^{2\a X_{12}}dXdY)$.
Then $\{\pd_xQ^-(\eta),\pd_x\widetilde{Q}^+(\eta)\}$
are continuous eigenfunctions of $\mL$ 
provided $\eta$ is sufficiently small.
We call them resonant modes because they grow exponentially in the rear of
the line solitons.
Those continuous eigenfunctions were found by \cite{Burtsev} and reflect dynamics of modulating $1$-line soliton
(see \cite{Miz15}).
\bigskip

\subsection{Resonant modes for $2$-line solitons of P-type}
\label{subsec:resonance-P}
Now we consider resonant modes of the linearized operator around
$u_2=2\pd_x^2\log\tau_2$ with $\tau_2$ given by \eqref{eq:tau-P} with $a=d=1$.
\par
Let $b_1$ and $b_2$ be constants satisfying \eqref{eq:b1-b2P}. In a moving coordinate
$X=x-b_1t$ and $Y=y-b_2t$, the $2$-line soliton solution $u_2$ is a stationary solution of
\eqref{eq:KPII''}.
Linearizing \eqref{eq:KPII''} around $u_2$, we have
\begin{gather}
  \label{eq:linequ2P}
  \pd_tu=\mL_2u\,,
\\ \label{eq:defmL0}
\mL_2=\mL_0-\frac{3}{2}\pd_x(u_2\cdot)\,,
\quad
\mL_0=\frac{1}{4}\left(-\pd_x^3+4b_1\pd_x+4b_2\pd_y-3\pd_x^{-1}\pd_y^2\right)\,.
\end{gather}
\par
Let $X_{ij}=X+2a_{ij}Y$.  We see from \eqref{eq:u2-asymptotics} that,
as $y\to\pm\infty$, $u_2$ is a superposition of $[2,3]$-soliton
$2c_{23}\sech^2\left(\frac{\k_3-\k_2}{2}X_{23}+\mu_{1,\pm}\right)$ and
$[1,4]$-soliton
$2c_{14}\sech^2\left(\frac{\k_4-\k_1}{2}X_{14}+\mu_{2,\pm}\right)$ .
Using Lemma~\ref{lem:prodJdJ}, we will find continuous eigenfunctions for
$\mL_2$ in exponentially weighted spaces
$\calX_1:=L^2(\R^2;e^{2\a X_{23}}dXdY)$ and
$\calX_2:=L^2(\R^2;e^{2\a X_{14}}dXdY)$.
Note that \eqref{eq:linequ2P} is well-posed for $t\ge0$
in $\calX_1$ and $\calX_2$
\par
The Jost solutions and the dual Jost solutions for $L_2$ are
\begin{align}
  \label{def:Phi2}
& \Phi^2(\bx,k)=\frac{e^{ik x-k^2y+ik^3t}}{\tau_2}
\operatorname{Wr}\left((ik-\pd_x)f_1,(ik-\pd_x)f_2\right)\,,
  \\ &
\label{def:Phi2*}       
 \Phi^{2,*}(\bx,k)=\frac{e^{-ik x+k^2y-ik^3t}}{\tau_2}
\operatorname{Wr}\left((ik-\pd_x)^{-1}f_1,(ik-\pd_x)^{-1}f_2\right)\,,
\end{align}
where $(ik-\pd_x)^{-1}e^{\theta_j}=(ik-\k_j)^{-1}e^{\theta_j}$.
For $n=1$, $2$, $3$, $4$, let
\begin{equation}
  \label{def:Phi2n}
  \Phi_n^2(\bx)=\Phi^2(\bx,-i\k_n)\,,\quad
  \Phi_n^{2,*}(\bx)=i\res_{k=-i\k_n}\Phi^{2,*}(\bx,k)\,.
\end{equation}
Then
\begin{gather*}
\Phi^2_1(\bx)=\Phi^2_4(\bx)=
\frac{e^{\theta_1+\theta_4}}{\tau_2}
(\k_1-\k_4)(\k_1-\pd_x)(\k_4-\pd_x)f_1\,,
\\
\Phi^2_2(\bx)=-\Phi^2_3(\bx)
=\frac{e^{\theta_2+\theta_3}}{\tau_2}
(\k_3-\k_2)(\k_3-\pd_x)(\k_2-\pd_x)f_2\,,
\\
\Phi_1^{2,*}(\bx)=-\Phi_4^{2,*}(\bx)=-\frac{f_1}{\tau_2}\,,
\quad
\Phi_2^{2,*}(\bx)=\Phi_3^{2,*}(\bx)=-\frac{f_2}{\tau_2}\,.
\end{gather*}
Let $\beta_{ij}^\pm(\eta)=a_{ij}\pm\gamma_{ij}(\eta)$,
$\gamma_{ij}(\eta)=\sqrt{c_{ij}+i\eta}$,
$\lambda_{ij}^\pm(\eta)=i\eta\{\gamma_{ij}(\eta)\pm(3a_{ij}-b_2)\}$ and
\begin{gather*}
Q_{ij}^\pm(T,X,Y,\eta)=\Phi^2(\bx,-i\beta_{ij}^\pm(\eta))\Phi_i^{2,*}(\bx)\,,\quad
\widetilde{Q}_{ij}^\pm(T,X,Y,\eta)
=i\eta\Phi^2_i(\bx)\Phi^{2,*}(\bx,-i\beta_{ij}^\pm(\eta))\,.
\end{gather*}
Lemma~\ref{lem:prodJdJ} implies that $\pd_XQ_{ij}^\pm$
and $\pd_X\widetilde{Q}_{ij}^\pm$ are solutions of \eqref{eq:linequ2P} and that
$Q_{ij}^\pm$ and $\widetilde{Q}_{ij}^\pm$ are solutions of
the adjoint equation of \eqref{eq:linequ2P}.
Let $0<\a<2\min\{\sqrt{c_{14}},\sqrt{c_{23}}\}$.
For $\{i,j\}=\{2,3\}$ and  $\{i,j\}=\{1,4\}$,
we have as  $X_{ij}\to\pm\infty$,
\begin{gather*}
Q_{ij}^-(\eta)=O(e^{-\Re\{\gamma_{ij}(0)\pm\gamma_{ij}(\eta)\}|X_{ij}|})\,,
\quad \widetilde{Q}_{ij}^+(\eta)=O(e^{-\Re\{\gamma_{ij}(0)\pm\gamma_{ij}(\eta)\}|X_{ij}|})\,,
\\
Q_{ij}^+(\eta)=O(e^{-\Re\{\gamma_{ij}(0)\mp\gamma_{ij}(\eta)\}|X_{ij}|})\,, \quad
\widetilde{Q}_{ij}^-(\eta)=O(e^{-\Re\{\gamma_{ij}(0)\mp\gamma_{ij}(\eta)\}|X_{ij}|})\,.
\end{gather*}
If $\eta\in\R$ is  small, then
$\{\pd_xQ_{ij}^-(\eta),\pd_x\widetilde{Q}_{ij}^+(\eta)\}$
and $\{Q_{ij}^+(\eta)\,,\,\widetilde{Q}_{ij}^-(\eta)\}$
are continuous eigenfunctions of $\mL_2$ in $L^2(\R^2;e^{2\a X_{ij}}dXdY)$
and those of $\mL_2^*$ in $L^2(\R^2;e^{-2\a X_{ij}}dXdY)$, respectively.
Those continuous eigenfunctions are associated with the modulation of
$[i,j]$-soliton.
\par
Using $\theta_i-\theta_j=(\k_i-\k_j)X_{ij}$ and the fact that
\begin{equation}
\label{eq:btitj}
\beta_{ij}^\pm(\eta)x+\beta_{ij}^\pm(\eta)^2y-\beta_{ij}^\pm(\eta)^3t-\frac{1}{2}(\theta_i+\theta_j)
= \pm \gamma_{ij}(\eta)X_{ij} +i\eta Y\mp\lambda_{ij}^\pm(\eta)T\,,
\end{equation}
we can show that
\begin{equation}
  \label{eq:LQ23}  
 \begin{split}
\mL_2 \pd_XQ_{23}^-(\eta)=\lambda_{23}^-(\eta)\pd_XQ_{23}^-(\eta)\,,
\quad \mL_2 \pd_X\widetilde{Q}_{23}^+(\eta)
=\lambda_{23}^+(\eta)\pd_X\widetilde{Q}_{23}^+(\eta)\,,
\\
\mL_2 \pd_XQ_{14}^-(\eta)=\lambda_{14}^-(\eta)\pd_XQ_{14}^-(\eta)\,,
\quad 
\mL_2 \pd_X\widetilde{Q}_{14}^+(\eta)
=\lambda_{14}^+(\eta)\pd_X\widetilde{Q}_{14}^+(\eta)\,,
\\
\mL_2^*Q_{23}^+(\eta)=\lambda_{23}^+(\eta)Q_{23}^+(\eta)\,,
\quad
\mL_2^*\widetilde{Q}_{23}^-(\eta)
=\lambda_{23}^-(\eta)\widetilde{Q}_{23}^-(\eta)\,,\\
\mL_2^*Q_{14}^+(\eta)=\lambda_{14}^+(\eta)Q_{14}^+(\eta)\,,
\quad
 \mL_2^*\widetilde{Q}_{14}^-(\eta)
=\lambda_{14}^-(\eta)\widetilde{Q}_{14}^-(\eta)\,.
 \end{split}
\end{equation}
Let $\lambda_{1,\pm}(\eta)=\lambda_{23}^\mp(\pm\eta)$,
$\lambda_{2,\pm}(\eta)=\lambda_{14}^\mp(\pm\eta)$, $(X,Y)$ be as
\eqref{eq:b1-b2P} and
\begin{align}
\label{eq:defg1+}
& g^2_{1,+}(X,Y,\eta)=(\k_3-\k_2)d_{1,+}(\eta)\pd_x
\left(\Phi^2(\bx,-i\beta_{23}^-(\eta))\Phi^{2,*}_2(\bx)\right)\,,
\\ & \label{eq:defg1-}
g^2_{1,-}(X,Y,\eta)=i\eta d_{1,-}(\eta)\pd_x\left(
\Phi^2_2(\bx)\Phi^{2,*}(\bx,-i\beta_{23}^+(-\eta))\right)\,,
\\ &
\label{eq:defg1*+}
     g^{2,*}_{1,+}(X,Y,\eta)=-\frac{i\bar{\eta}}{\k_3-\k_2}
     \Phi^2_2(\bx)\Phi^{2,*}(\bx,-i\beta_{23}^-(-\bar{\eta}))\,,
\\ & \label{eq:defg1*-}
g^{2,*}_{1,-}(X,Y,\eta)=\Phi^2(\bx,-i\beta_{23}^+(\bar{\eta}))\Phi^{2,*}_2(\bx)\,,
\\ &
\label{eq:defg2+}
  g^2_{2,+}(X,Y,\eta)=(\k_4-\k_1) d_{2,+}(\eta)
     \pd_x\left(\Phi^2(\bx,-i\beta_{14}^-(\eta)) \Phi^{2,*}_1(\bx)\right)\,,
\\ & \label{eq:defg2-}
     g^2_{2,-}(X,Y,\eta)=-i\eta d_{2,-}(\eta)\pd_x\left(
     \Phi^2_1(\bx)\Phi^{2,*}(\bx,-i\beta_{14}^+(-\eta)\right)\,,
\\ &
\label{eq:defg2*+}
     g^{2,*}_{2,+}(X,Y,\eta)=-\frac{i\bar{\eta}}{\k_4-\k_1}
     \Phi^2_1(\bx)\Phi^{2,*}(\bx,-i\beta_{14}^-(-\bar{\eta}))\,,
\\ & \label{eq:defg2*-}
g^{2,*}_{2,-}(X,Y,\eta)=-\Phi^2(\bx,-i\beta_{14}^+(\bar{\eta}))\Phi^{2,*}_1(\bx)\,,
\end{align}
where
$d_{1,\pm}(\eta)=\pm 1/\gamma_{23}(\pm\eta)$ and
$d_{2,\pm}(\eta)=\pm 1/\gamma_{14}(\pm\eta)$.
By  \eqref{eq:LQ23} and the definitions of $g^2_{j,\pm}$ and $g^{2,*}_{j,\pm}$,
we have the following.

\begin{lemma}
\label{lem:evmL2}
\begin{equation}
  \label{eq:lambdas}
  \begin{split}
\mL_2 g^2_{1,\pm}(x,y,\eta)& =\lambda_{1,\pm}(\eta)g^2_{1,\pm}(x,y,\eta)\,,\quad
\lambda_{1,\pm}(\eta)=i\eta\{b_2-3a_{23}\pm\gamma_{23}(\pm\eta)\}\,,
\\
\mL_2 g^2_{2,\pm}(x,y,\eta)& =\lambda_{2,\pm}(\eta)g^2_{2,\pm}(x,y,\eta)\,,\quad
\lambda_{2,\pm}(\eta)=i\eta\{b_2-3a_{14}\pm\gamma_{14}(\pm\eta)\}\,.
  \end{split}
\end{equation}
\begin{equation}
  \label{eq:adjoint}
\mL_2^* g^{2,*}_{j,\pm}(x,y,\eta)
=\overline{\lambda_{1,\pm}(\eta)}g^{2,*}_{j,\pm}(x,y,\eta)
\quad\text{for $j=1$, $2$.}
\end{equation}
\end{lemma}
\bigskip

\subsection{Linearized Miura transformation and resonant modes}
Let $\tau_1$ be as in \eqref{eq:tau-P}. Then the
$1$-line soliton solution $u_1=2\pd_x^2\log\tau_1$ is a stationary solution of
\eqref{eq:KPII''}.
If we linearize \eqref{eq:KPII''} around $u_1$, then
\begin{equation}
  \label{eq:linequ1P}
\pd_tu=\mL_1u\,,\quad
\mL_1=\mL_0-\frac{3}{2}\pd_x(u_1\cdot))\,.
\end{equation}
\par
Let $L_0=-\pd_y+\pd_x^2$, $L_1=-\pd_y+\pd_x^2+u_1$ and 
\begin{gather*}
\Phi^0(\bx,k)=e^{ikx-k^2y+ik^3t}\,,\quad
\Phi^{0,*}(\bx,k)=e^{-ikx+k^2y-ik^3t}\,,
\\
\Phi^1(\bx,k)=\frac{e^{ik x-k^2y+ik^3t}}{\tau_1}
\sum_{j=2,3}(ik-\k_j)e^{\theta_j}\,,
\quad \Phi^{1,*}(\bx,k)=\frac{e^{-ik x+k^2y-ik^3t}}{\tau_1}
\sum_{j=2,3}\frac{e^{\theta_j}}{ik-\k_j}\,,
\\
\Phi^1_2(\bx)=-\Phi^1_3(\bx)
=(\k_2-\k_3)\frac{e^{\theta_2+\theta_3}}{\tau_1}\,,
\quad \Phi^{1,*}_2(\bx)=\Phi^{1,*}_3(\bx)=\frac{1}{\tau_1}\,.
\end{gather*}
Note that $\Phi^j(\bx,k)$ and $\Phi^{j,*}(\bx,k)$ are
Jost solutions and the dual Jost solutions of $L_j$ for
$j=0$ and $1$, respectively and that
$\Phi_n^1(\bx)=\Phi^1(\bx,-i\k_n)$,
$\Phi_n^{1,*}(\bx)=i\res_{k=-i\k_n}\Phi^{1,*}(\bx,k)$
for $n=2$ and $3$.

Lemma~\ref{lem:prodJdJ} tells us that
$\pd_x(\Phi^i\Phi^{i,*})$
are  solutions of $\pd_tu=\mL_iu$ and that
$\Phi^i\Phi^{i,*}$ are solutions
of the adjoint equation $\pd_tu+\mL_i^*u=0$.
The Darboux transformations $\nabla M_\pm(v_i)$ give a correspondence
between $\pd_x(\Phi^i\Phi^{i,*})$
and $\pd_x(\Phi^{i-1}\Phi^{i-1,*})$,
and the Darboux transformations $\nabla M_\pm(v_i)^*$ give a correspondence
between $\Phi^i\Phi^{i,*}$ and $\Phi^{i-1}\Phi^{i-1,*}$.

\begin{lemma}
  Let $\beta\in \C$ and $\beta'\in \C\setminus\{\k_1,\k_2,\k_3,\k_4\}$. Then
  \label{lem:M-pdPhiPhi*}
  \begin{gather}
    \label{eq:M2+pdPhiPhi*}
\nabla M_+(v_2)\pd_x
\left(\Phi^1(\bx,-i\beta)\Phi^{2,*}(\bx,-i\beta')\right)
=2\pd_x\left(\Phi^2(\bx,-i\beta)\Phi^{2,*}(\bx,-i\beta')\right)\,,
\\
    \label{eq:M2-pdPhiPhi*}
\nabla M_-(v_2)\pd_x
\left(\Phi^1(\bx,-i\beta)\Phi^{2,*}(\bx,-i\beta')\right)
=2\pd_x\left(\Phi^1(\bx,-i\beta)\Phi^{1,*}(\bx,-i\beta')\right)\,,
\\
    \label{eq:M1+pdPhiPhi*}
\nabla M_+(v_1)\pd_x
\left(\Phi^0(\bx,-i\beta)\Phi^{1,*}(\bx,-i\beta')\right)
=2\pd_x\left(\Phi^1(\bx,-i\beta)\Phi^{1,*}(\bx,-i\beta')\right)\,,
\\    \label{eq:M1-pdPhiPhi*}
\nabla M_-(v_1)\pd_x
\left(\Phi^0(\bx,-i\beta)\Phi^{1,*}(\bx,-i\beta')\right)
=2\pd_x\left(\Phi^0(\bx,-i\beta)\Phi^{0,*}(\bx,-i\beta')\right)\,.
  \end{gather}
\end{lemma}

\begin{lemma}
  \label{lem:M-PhiPhi*}
Suppose that $\beta$, $\beta'\in \C\setminus\{\k_1,\k_2,\k_3,\k_4\}$. Then
  \begin{gather}
    \label{eq:M2+PhiPhi*}
\pd_x\nabla M_-(v_2)^*
\left(\Phi^1(\bx,-i\beta)\Phi^{1,*}(\bx,-i\beta')\right)
=2\pd_x\left(\Phi^2(\bx,-i\beta)\Phi^{1,*}(\bx,-i\beta')\right)\,,
\\
    \label{eq:M2-PhiPhi*}
\pd_x\nabla M_+(v_2)^*
\left(\Phi^2(\bx,-i\beta)\Phi^{2,*}(\bx,-i\beta')\right)
=2\pd_x\left(\Phi^2(\bx,-i\beta)\Phi^{1,*}(\bx,-i\beta')\right)\,,
\\
    \label{eq:M1+PhiPhi*}
\pd_x\nabla M_-(v_1)^*
\left(\Phi^0(\bx,-i\beta)\Phi^{0,*}(\bx,-i\beta')\right)
=2\pd_x\left(\Phi^1(\bx,-i\beta)\Phi^{0,*}(\bx,-i\beta')\right)\,,
\\    \label{eq:M1-PhiPhi*}
\pd_x\nabla M_+(v_1)^*
\left(\Phi^1(\bx,-i\beta)\Phi^{1,*}(\bx,-i\beta')\right)
=2\pd_x\left(\Phi^1(\bx,-i\beta)\Phi^{0,*}(\bx,-i\beta')\right)\,.
  \end{gather}
\end{lemma}

For the proof of Lemmas~\ref{lem:M-pdPhiPhi*} and \ref{lem:M-PhiPhi*},
see Appendix~\ref{sec:lem:M-pdPhiPhi*}.
\begin{corollary}
  \label{cl:kerM-P}
\begin{gather}
  \label{eq:Q-ker}
\nabla M_-(v_2)\left(\Phi^2(\bx,-i\beta)\Phi^{2,*}_1(\bx)\right)=0\,,
\\   \label{eq:Q-ker'}
\nabla M_-(v_1)\left(\Phi^1(\bx,-i\beta)\Phi^{1,*}_2(\bx)\right)=0\,.
\end{gather}
\end{corollary}

\begin{corollary}
  \label{cl:Mrelations-P}
Suppose that $\beta\in \C\setminus\{\k_1,\k_2,\k_3,\k_4\}$. Then
\begin{gather}
\label{eq:Mpmg1}
\nabla M_+(v_2)\left(\Phi^1(\bx,-i\beta)\Phi^{1,*}_2(\bx)\right)
=2\Phi^2(\bx,-i\beta)\Phi^{1,*}_2(\bx)\,,
\\  \label{eq:Mpmg2}
\nabla M_-(v_2)\left(\Phi^2(\bx,-i\beta)\Phi^{2,*}_2(\bx)\right)
=2\Phi^2(\bx,-i\beta)\Phi^{1,*}_2(\bx)\,,
\\ \label{eq:B2B3}
\nabla M_+(v_2)\left(\Phi^1_2(\bx)\Phi^{1,*}(\bx,-i\beta)\right)
=2\Phi^2_2(\bx)\Phi^{1,*}(\bx,-i\beta)\,,
\\ \label{eq:Mg2-}
\nabla M_-(v_2)\left(\Phi^2_2(\bx)\Phi^{2,*}(\bx,-i\beta)\right)
=2\Phi^2_2(\bx)\Phi^{1,*}(\bx,-i\beta)\,,
\\ \label{eq:Mg2-'}
\nabla M_-(v_2)\left(\Phi^2_1(\bx)\Phi^{2,*}(\bx,-i\beta)\right)
=2\Phi^2_1(\bx)\Phi^{1,*}(\bx,-i\beta)\,,
\\ \label{eq:Mg1-}
\nabla M_-(v_1)\left(\Phi^1_2(\bx)\Phi^{1,*}(\bx,-i\beta)\right)
=2\Phi^1_2(\bx)\Phi^{0,*}(\bx,-i\beta)\,.
\end{gather}
\end{corollary}
\begin{proof}[Proof of Corollaries~\ref{cl:kerM-P} and \ref{cl:Mrelations-P}]
 We formally have $\nabla M_+(v_2)^*=\nabla M_-(v_2)$. Since
$\Phi^{1,*}(\bx,-i\beta')$ is regular at $\beta'=\k_1$,
it follows from \eqref{eq:M2-PhiPhi*} that
\begin{align*}
&\pd_x\nabla  M_-(v_2)\left(\Phi^2(\bx,-i\beta)\Phi^{2,*}_1(\bx)\right)
\\=&  \lim_{\beta'\to \k_1}(\beta'-k_1)
\pd_x\nabla  M_-(v_2)\left(\Phi^2(\bx,-i\beta)\Phi^{2,*}(\bx,-i\beta')\right)
\\=& \lim_{\beta'\to \k_1}(\beta'-k_1)
\pd_x \left(\Phi^2(\bx,-i\beta)\Phi^{1,*}(\bx,-i\beta')\right)
= 0\,,
\end{align*}
and we have \eqref{eq:Q-ker} for sufficiently large $\beta>0$ because
$\nabla M_-(v_2)\Phi^2(\bx,-i\beta)\Phi^{2,*}_1(\bx)$ tends to $0$ as
$x\to-\infty$ for such $\beta$.  Since
$e^{-\beta x-\beta^2y+\beta^3t}\Phi^2(\bx,-i\beta)\Phi^{2,*}_1(\bx)$
is a polynomial in $\beta$, we have \eqref{eq:Q-ker} for any
$\beta\in\C$.

Using  Lemma~\ref{lem:M-PhiPhi*},
we can prove the rest in the similar manner.
Thus we complete the proof.  
\end{proof}
\par

Let $(X,Y)$ be as \eqref{eq:b1-b2P} and
\begin{gather}
\label{eq:defg11+}
g^1_{1,+}(X,Y,\eta)=(\k_3-\k_2)d_{1,+}(\eta)\pd_x\left\{
\Phi^1(\bx,-i\beta_{23}^-(\eta))\Phi^{1,*}_2(\bx)\right\}\,,  
\\
\label{eq:defg11-}
g^1_{1,-}(X,Y,\eta)=i\eta d_{1,-}(\eta)\pd_x\left\{
\Phi^1_2(\bx)\Phi^{1,*}(\bx,-i\beta_{23}^+(-\eta))\right\}\,,
\\
\label{eq:defg11*+}
g^{1,*}_{1,+}(X,Y,\eta)
=-\frac{i\bar{\eta}}{\k_3-\k_2}\Phi^1_2(\bx)
\Phi^{1,*}(\bx,-i\beta_{23}^-(-\bar{\eta}))\,,  
\\
\label{eq:defg11*-}
g^{1,*}_{1,-}(X,Y,\eta)=\Phi^1(\bx,-i\beta_{23}^+(\bar{\eta}))
\Phi^{1,*}_2(\bx)\,.
\end{gather}
Let $z_1=x+2a_{23}y$ and
$\mL_1(\eta)v(z_1):=e^{-iy\eta} \mL_1(e^{iy\eta}v(z_1))$.
We remark that $u_1$ depends only on $z_1$.
Let
\begin{gather}
  \label{eq:defg1z1'}
g(z_1,\eta)=\sqrt{c_{23}}d_{1,+}(\eta)\pd_{z_1}^2\left\{
e^{-\gamma_{23}(\eta)z_1}\sech\sqrt{c_{23}}z_1\right\}\,\\
  \label{eq:defg1z1'*}
g^*(z_1,\eta)=\frac12\pd_{z_1}\left\{
e^{\gamma_{23}(-\bar{\eta})z_1}\sech\sqrt{c_{23}}z_1\right\}\,.  
\end{gather}
By \eqref{eq:btitj}, we have for $\eta\in\C$,
\begin{gather}
  \label{eq:defg1z1}
g^1_{1,\pm}(x,y,\eta)=e^{iy\eta}g(z_1,\pm\eta)\,,\quad
g^{1,*}_{1,\pm}(x,y,\eta)=e^{iy\bar{\eta}}g^*(z_1,\pm\eta)\,,
\\  \label{eq:evmL1}
\mL_1(\eta)g(z_1,\pm\eta)=\lambda_{1,\pm}(\eta)g(z_1,\pm\eta)\,,
  \quad
\mL_1(\eta)^*g^*(z_1,\pm\eta)=\overline{\lambda_{1,\pm}(\eta)}g^*(z_1,\pm\eta)\,.
\end{gather}
\par

Darboux transformations $\nabla M_\pm(v_k)$ ($k=1$, $2$) connect
$g^k_{k,+}$ with $0$ and Darboux transformations $\nabla M_\pm(v_2)^*$
connect $g^{2,*}_{1,\pm}$ with $g^{1,*}_{1,\pm}$.  \begin{lemma}
\label{lem:g-tg} For $i=1$, $2$, \begin{gather} \label{eq:kernablaM}
\nabla M_-(v_i)\pd_x^{-1}g^i_{i,+}=0\,,\quad \nabla
M_+(v_i)^*g^{i,*}_{i,-}=0\,, \\ \label{eq:M+vi} \nabla
M_+(v_i)\pd_x^{-1}g^i_{i,+}=2g^i_{i,+}\,,\enskip \end{gather} and
\begin{equation} \label{eq:g-tg} \nabla M_-(v_2)\pd_x^{-1}g^2_{1,\pm}
=\nabla M_+(v_2)\pd_x^{-1}g^1_{1,\pm}\,, \quad \nabla
M_+(v_2)^*g^{2,*}_{1,\pm} =\nabla M_-(v_2)^*g^{1,*}_{1,\pm}\,.
\end{equation} \end{lemma} Since $\nabla M_+(v_i)^*=\nabla M_-(v_i)$
(formally), Lemma~\ref{lem:g-tg} follows immediately from
Corollaries~\ref{cl:kerM-P} and \ref{cl:Mrelations-P}.  \bigskip

\section{Projections to resonant modes
associated with $2$-line soliton of P-type}
\subsection{Orthogonality relations}
First, we will show the orthogonality relations between
$g^2_{j,\pm}(x,y,\eta)$ and $g^{2,*}_{j,\pm}(x,y,\eta)$.
\begin{lemma}
  \label{lem:orthogonality}
Let $\eta=\eta_R+i\eta_I$ and $\eta_1=\eta_{1,R}+i\eta_I$ and $\eta_0\in (0,4\sqrt{2}c_{23})$.
Suppose that $\eta_R\in(-\eta_0,\eta_0)$  and that $\eta_I\in (-c_{23},c_{23})$.
Then for $j=1$, $2$ and $\varphi\in C^1([-\eta_0,\eta_0])$,
\begin{gather}
\label{eq:gg*-1}
\lim_{M\to\infty}\int_{-\eta_0}^{\eta_0}\varphi(\eta_{1,R})\left(\int_{-M}^M \int_{-\infty}^\infty
g^2_{j,\pm}(x,y,\eta)\overline{g^{2,*}_{j,\pm}(x,y,\eta_1)}\,dxdy\right)\,d\eta_{1,R}
=\pm2\pi i\eta\varphi(\eta_R)\,,
\\
\label{eq:gg*-2}
\lim_{M\to\infty}\int_{-\eta_0}^{\eta_0}\varphi(\eta_{1,R})\left(\int_{-M}^M \int_{-\infty}^\infty
g^2_{j,\pm}(x,y,\eta)\overline{g^{2,*}_{j,\mp}(x,y,\eta_1)}\,dxdy\right)\,d\eta_{1,R}  =0\,.
\end{gather}
\end{lemma}
To prove Lemma~\ref{lem:orthogonality}, we use the following
identities.
\begin{claim}
  \label{cl:Phi22*s}
  \begin{gather}
\label{eq:pd12-22}
\Phi^2(\bx,-i\beta)\Phi^{2,*}_j(\bx)=
\pd_x\left(\Phi^1(\bx,-i\beta)\Phi^{2,*}_j(\bx)\right)
\quad\text{for $j=1$, $4$,}
\\  \label{eq:Phi11*-1}
\Phi^1(\bx,-i\beta)\Phi^{1,*}_j(\bx)=\pd_x\left(\Phi^0(\bx,-i\beta)\Phi^{1,*}_j(\bx)\right)
\quad\text{for $j=2$, $3$,}
\\  \label{eq:Phi22*-1}
\Phi^2(\bx,-i\beta)\Phi^{2,*}_2(\bx)=
-\pd_x\left\{\Phi^0(\bx,-i\beta)(\beta-\pd_x)f_2/\tau_2\right\}\,,
\\ \label{eq:Phi22*-3}
\Phi^2_1(\bx)\Phi^{2,*}(\bx,-i\beta)=
\frac{\k_1-\k_4}{(\beta-\k_1)(\beta-\k_4)}
\pd_xJ_1\,,
\\ \label{eq:Phi22*-2}
\Phi^2_2(\bx)\Phi^{2,*}(\bx,-i\beta)=
\frac{\k_2-\k_3}{(\beta-\k_2)(\beta-\k_3)}\pd_xJ_2\,,
\\ \label{eq:Phi11*}
\Phi^1_2(\bx)\Phi^{1,*}(\bx,-i\beta)
=\frac{-1}{(\beta-\k_2)(\beta-\k_3)}
\pd_x\left(\Phi^1_2(\bx)\Phi^{0,*}(\bx,-i\beta)\right)\,,
\end{gather}
where
\begin{align*}
& J_1=\Phi^{0,*}(\bx,-i\beta)\frac{e^{\theta_1+\theta_4}}{\tau_2}
\{(\k_2-\k_1)(\k_4-\k_2)\frac{e^{\theta_2}}{\beta-\k_2}
+(\k_3-\k_1)(\k_4-\k_3)\frac{e^{\theta_3}}{\beta-\k_3}\}\,,
\\ &
J_2=\Phi^{0,*}(\bx,-i\beta)\frac{e^{\theta_2+\theta_3}}{\tau_2}
\left\{(\k_4-\k_2)(\k_4-\k_3)\frac{e^{\theta_4}}{\beta-\k_4}
-(\k_3-\k_1)(\k_2-\k_1)\frac{e^{\theta_1}}{\beta-\k_1}\right\}\,.
\end{align*}
\end{claim}
Claim~\ref{cl:Phi22*s} follows from  straightforward computations.

\begin{proof}[Proof of Lemma~\ref{lem:orthogonality}]
By  \eqref{eq:lambdas} and \eqref{eq:adjoint}, we have for $j=1$ and $2$,
\begin{equation}
  \label{eq:g*Lg}
 \begin{split}
\{\lambda_{j,\pm}(\eta)-\lambda_{j,\pm}(\eta_1)\} g^2_{j,\pm}(\eta)
\overline{g^{2,*}_{j,\pm}(\eta_1)}
= & 
\overline{g^{2,*}_{j,\pm}(\eta_1)}\mL g^2_{j,\pm}(\eta)
-g^2_{j,\pm}(\eta)\mL^*\overline{g^{2,*}_{j,\pm}(\eta_1)}
\\= & \frac{1}{4}(\pd_xI_1+\pd_yI_2)\,,
 \end{split}  
\end{equation}
where we abbreviate $g^2_{j,\pm}(x,y,\eta)$ and $g^{2,*}_{j,\pm}(x,y,\eta_1)$
as  $g^2_{j,\pm}(\eta)$ and $g^{2,*}_{j,\pm}(\eta_1)$, respectively and
\begin{align*} 
I_1=& \left(-\pd_x^2+4b_1-6u_2\right)\{g^2_{j,\pm}(\eta)
\overline{g^{2,*}_{j,\pm}(\eta_1)}\}
+3\pd_xg^2_{j,\pm}(\eta)\pd_x\overline{g^{2,*}_{j,\pm}(\eta_1)}
\\ & +\pd_x^{-1}\pd_yg^2_{j,\pm}(\eta)\pd_x^{-1}
\pd_y\overline{g^{2,*}_{j,\pm}(\eta_1)}\,,
\end{align*}
$$
I_2=-3g^2_{j,\pm}(\eta)\pd_x^{-1}\pd_y\overline{g^{2,*}_{j,\pm}(\eta_1)}
-3\pd_x^{-1}\pd_yg^2_{j,\pm}(\eta)\overline{g^{2,*}_{j,\pm}(\eta_1)}
+4b_2g^2_{j,\pm}(\eta)\overline{g^{2,*}_{j,\pm}(\eta_1)}\,.$$
Suppose that $\psi_i$ and  $\tpsi_i$ ($i=1$, $2$) are solutions of
\eqref{eq:Jost-eqs} and \eqref{eq:dualJost-eqs}, respectively
and that $g^2_{j,\pm}(\eta)\overline{g^{2,*}_{j,\pm}}
=d_{j,\pm}(\eta)\pd_x(\psi_1\tpsi_1)\psi_2\tpsi_2$.
Then 
\begin{equation}
  \label{eq:psitpsi_y}
\pd_y(\psi_i\tpsi_i)=\pd_x(\pd_x\psi_i\tpsi_i-\psi_i\pd_x\tpsi_i)\,,
\end{equation}
and 
\begin{align}
\label{eq:I2}
I_2=& -3d_{j,\pm}(\eta)\pd_x\{(\pd_x\psi_1\tpsi_1-\psi_1\pd_x\tpsi_1)\psi_2\tpsi_2\}+d_{j,\pm}(\eta)I_3\,,
\\ \label{eq:I3}
I_3=& 
6(\psi_2\tpsi_1\pd_x\psi_1\pd_x\tpsi_2-\psi_1\tpsi_2\pd_x\psi_2\pd_x\tpsi_1)
+4b_2\pd_x(\psi_1\tpsi_1)\psi_2\tpsi_2\,.
\end{align}
\par
First, we will prove the orthogonal relation between $g^2_{1,+}$ and $g^{2,*}_{1,+}$.
Let $\beta=\beta_{23}^-(\eta)$, $\beta'=\beta_{23}^-(\eta_1)$ and
$$\psi_1=\Phi^2(\bx,-i\beta)|_{t=0}\,,
\enskip \tpsi_1=\Phi^{2,*}_2(\bx)|_{t=0}\,,
\enskip \psi_2=\Phi^2_2(\bx)|_{t=0}\,,
\enskip\tpsi_2=i\eta_1\Phi^{2,*}(\bx,-i\beta')|_{t=0}\,.$$
Let $z_1=x+2a_{23}y$,  $z_2=x+2a_{14}y$ and $Z_{1,\pm}$ be as
\eqref{eq:u2-asymptotics}. 
By \eqref{eq:t3-t2} and  \eqref{eq:btitj},
\begin{equation}
  \label{eq:tau2-23}
  \begin{split}
\tau_2e^{-(\theta_2+\theta_3)/2}=& 2\sqrt{(\k_4-\k_3)(\k_4-\k_2)}e^{\theta_4}\cosh Z_{1,+}
\\ & + 2\sqrt{(\k_3-\k_1)(\k_2-\k_1)}e^{\theta_1}\cosh Z_{1,-}\,,        
  \end{split}
\end{equation}
\begin{equation*}
\beta x+\beta^2y-\{\beta'x+(\beta')^2y\}
=(\gamma_{23}(\eta_1)-\gamma_{23}(\eta))z_1+iy(\eta-\eta_1)\,.
\end{equation*}
In view of  \eqref{eq:t3-t2}, we have as $(a_{14}-a_{23})y\to\pm\infty$,
\begin{align*}
& e^{-(\beta x+\beta^2y)}\psi_1= -\left(\beta-a_{14}\mp\sqrt{c_{14}}\right)
\left(\gamma_{23}(\eta)+\sqrt{c_{23}}\tanh Z_{1,\pm}\right)
+O(e^{-2\sqrt{c_{14}}|z_2|})\,,
\\ & 
e^{\beta' x+(\beta')^2y}\tpsi_2= \frac{-i}{\beta'-a_{14}\mp\sqrt{c_{14}}}
\left(\gamma_{23}(\eta_1)-\sqrt{c_{23}}\tanh Z_{1,\pm}\right)
+O(e^{-2\sqrt{c_{14}}|z_2|})\,,
\\ &
e^{(\theta_2+\theta_3)/2}|t_{t=0}\tpsi_1=
\frac{\mp1}{2}\frac{e^{\mu_{1,\pm}}}{\sqrt{c_{14}}\pm(a_{14}-\k_3)}\sech Z_{1,\pm}
(1+O(e^{-2\sqrt{c_{14}}|z_2|}))\,,
\\ &
e^{-(\theta_2+\theta_3)/2}|_{t=0}\psi_2 = \pm\sqrt{c_{23}}
e^{-\mu_{1,\pm}}\left\{\sqrt{c_{14}}\pm(a_{14}-\k_3)\right\}\sech Z_{1,\pm}
(1+O(e^{-2\sqrt{c_{14}}|z_2|}))\,,
\end{align*}
and for $i$, $j\ge0$,
\begin{gather*}
\pd_x^i\psi_2\pd_x^j\tpsi_1=O\left(e^{-(\k_3-\k_2)|z_1|}\right)\,,
\quad \pd_x^i\psi_1\pd_x^j\tpsi_2
=O\left(e^{(\gamma_{23}(\eta_1)-\gamma_{23}(\eta))z_1+iy(\eta-\eta_1)}\right)\,,
\\
  |\Re\left(\gamma_{23}(\eta)-\gamma_{23}(\eta_1)\right)|\le \Re\sqrt{c_{23}-\eta_I+i\eta_0}
<\k_3-\k_2\,.
\end{gather*}
Combining \eqref{eq:g*Lg}, \eqref{eq:I2}, \eqref{eq:I3} with the above,
we see that
\begin{align*}  
\int_{-M}^M\int_{-\infty}^\infty g^2_{1,+}(\eta)
\overline{g^{2,*}_{1,+}(\eta_1)}\,dxdy
=& \frac{d_{1,+}(\eta)}{4}
\left[\int_\R \frac{I_3}{\lambda_{1,+}(\eta)-\lambda_{1,+}(\eta_1)}\,dz_1
\right]_{y=-M}^{y=M}\,,
\end{align*}
and that 
$e^{iy(\eta_1-\eta)}I_3|_{\eta=\eta_1}$ converges exponentially to
$\sqrt{c_{23}}I_{3,\infty}(Z_{1,\pm})$
as $(a_{14}-a_{23})y\to\pm\infty$, where
\begin{align*}
I_{3,\infty}(x)= &
i\eta\left[\{3(\beta_{23}^-(\eta)^2-\k_2\k_3)+2b_2\gamma_{23}(\eta)\}\sech^2x
-6c_{23}\sech^4x\right]
\\ & 
+2b_2\sqrt{c_{23}}\sech^2x\tanh x(i\eta+2c_{23}\sech^2x)\,.
\end{align*}
We have
\begin{align*}
\int_{-\infty}^\infty I_{3,\infty}(x)\,dx
=& 2i\eta\{3\beta_{23}^-(\eta)^2+2b_2\gamma_{23}(\eta)-\omega_{23}\}\,.
\\ =& 4\eta\gamma_{23}(\eta)\pd_\eta\lambda_{1,+}(\eta)\,.
\end{align*}
It follows from Lemma~\ref{lem:lambda1-etaI} that $\pd_\eta\lambda_{1,+}\ne0$ and that
$$\frac{1}{\lambda_{1,+}(\eta)-\lambda_{1,+}(\eta_1)}
-\frac{1}{\pd_\eta\lambda_{1,+}(\eta)(\eta_R-\eta_{1,R})}$$
and its derivative are uniformly bounded for $\eta_{1,R}\in[-\eta_0,\eta_0]$.
Moreover,
$$\frac{1}{\lambda_{1,+}(\eta)-\lambda_{1,+}(\eta_1)}
\int_\R(e^{iy(\eta_1-\eta)}I_3-I_3|_{\eta=\eta_1})\,dz_1
$$
and its $\eta_1$-derivative are uniformly bounded for $y\in\R$ and
$\eta_{1,R}\in[-\eta_0,\eta_0]$.
Combining the above and using the Riemann-Lebesgue lemma,
we have for $\varphi\in C^1([-\eta_0,\eta_0])$,
\begin{align*}
  & \lim_{M\to\infty}
    \int_{-\eta_0}^{\eta_0}\varphi(\eta_{1,R}) \left(\int_{-M}^M\int_{-\infty}^\infty
    g^2_{1,+}(\eta)\overline{g^{2,*}_{1,+}(\eta_1)}\,dxdy\right)\,d\eta_{1,R}
  \\=&
\frac{id_{1,+}(\eta)}{2\pd_\eta\lambda_{1,+}(\eta)}\left(\int_\R I_{3,\infty}(x)\,dx\right)
       \lim_{M\to\infty}\int_{-\eta_0}^{\eta_0}\varphi(\eta_{1,R})
       \frac{\sin M(\eta_R-\eta_{1,R})}{\eta_R-\eta_{1,R}}\,d\eta_{1,R}
  \\=& 2\pi i\eta\varphi(\eta_R)\,.
\end{align*}
\par

Next, we will prove the orthogonal relation for $g^2_{2,+}$ and $g^{2,*}_{2,+}$.
Let
\begin{gather}
\label{eq:psi1,tpsi1}
\psi_1=\Phi^2(\bx,-i\beta)|_{t=0}\,,\enskip  
\tpsi_1=\Phi^{2,*}_1(\bx)|_{t=0}\,,\enskip \beta=\beta_{14}^-(\eta)\,,
\\ \label{eq:psi2,tpsi2}
\psi_2=\Phi^2_1(\bx)|_{t=0}\,,\enskip \tpsi_2=i\eta_1\Phi^{2,*}(\bx,-i\beta')|_{t=0}\,,
\enskip \beta'=\beta_{14}^-(\eta_1)\,.
\end{gather}
By \eqref{eq:t3-t2} and \eqref{eq:btitj},
\begin{equation}
  \label{eq:tau2-14}
  \begin{split}
    e^{-(\theta_1+\theta_4)/2}\tau_2=& 2\sqrt{(\k_4-\k_3)(\k_3-\k_1)}e^{\theta_3}\cosh Z_{2,+}
    \\ & +2\sqrt{(\k_4-\k_2)(\k_2-\k_1)}e^{\theta_2}\cosh Z_{2,-}\,,
  \end{split}
\end{equation}
where $Z_{2,\pm}$ be as \eqref{eq:u2-asymptotics} and 
\begin{equation*}
\beta x+\beta^2y-\{\beta'x+(\beta')^2y\}
=-(\gamma_{14}(\eta)-\gamma_{14}(\eta_1))z_2+iy(\eta-\eta_1)\,.
\end{equation*}
Combining the above with  the fact that $z_1=z_2+2(a_{23}-a_{14})y$,
we have as $(a_{23}-a_{14})y\to\pm\infty$,
\begin{align}
\label{eq:psi1-asymp}
& e^{-(\beta x+\beta^2y)}\psi_1= -(\beta-a_{23}\mp\sqrt{c_{23}})
\left(\gamma_{14}(\eta)+\sqrt{c_{14}}\tanh Z_{2,\pm}\right)
+O(e^{-2\sqrt{c_{23}}|z_1|})\,,
\\ & \label{eq:tpsi1-asymp}
e^{(\theta_1+\theta_4)/2}|_{t=0}\tpsi_1=
-\frac{1}{2}\frac{e^{\mu_{2,\pm}}}{\k_4-a_{23}\mp\sqrt{c_{23}}}\sech Z_{2,\pm}
(1+O(e^{-2\sqrt{c_{23}}|z_1|}))\,,
\\ & \label{eq:psi2-asymp}
e^{\beta' x+(\beta')^2y}\tpsi_2= \frac{-1}{\beta'-a_{23}\mp\sqrt{c_{23}}}
\left(\gamma_{14}(\eta_1)-\sqrt{c_{14}}\tanh Z_{2,\pm}\right)
+O(e^{-2\sqrt{c_{23}}|z_1|})\,,
\\ & \label{eq:tpsi2-asymp}
e^{-(\theta_1+\theta_4)/2}|_{t=0}\psi_2 = \sqrt{c_{14}}
e^{-\mu_{2,\pm}}\left(\k_4-a_{23}\mp\sqrt{c_{23}}\right)\sech Z_{2,\pm}
(1+O(e^{-2\sqrt{c_{23}}|z_1|}))\,.
\end{align}
Combining the above, we can prove the orthogonal relation of
$g^2_{2,+}$ and $g^{2,*}_{2,+}$ in the same way as that of $g^2_{1,+}$
and $g^{2,*}_{1,+}$.  We can prove the rest of \eqref{eq:gg*-1} in
exactly the same way.
\par
Finally, we will prove \eqref{eq:gg*-2}.
By \eqref{eq:defg1+}--\eqref{eq:defg1*-},
$$g^2_{1,+}(\eta)\overline{g^{2,*}_{1,-}(\eta_1)}
=(\k_3-\k_2)d_{1,+}(\eta)\pd_x(\psi_+\tpsi)\psi_-\tpsi\,,$$
where $\psi_+=\Phi^2(\bx,-i\beta)|_{t=0}$ with $\beta=\beta_{23}^-(\eta)$,
$\psi_-=\Phi^2(\bx,-i\beta')|_{t=0}$ with $\beta'=\beta_{23}^+(-\eta_1)$
and $\tpsi=\Phi^{2,*}_2(\bx)|_{t=0}$.
Since $\psi_+$, $\psi_-$ are solutions of \eqref{eq:Jost-eqs}
and $\tpsi$ is a solution of \eqref{eq:dualJost-eqs},
$$\pd_y(\psi_\pm\tpsi)=\pd_x(\pd_x\psi_\pm\tpsi-\psi_\pm\pd_x\tpsi)\,,$$
\begin{align*}
 2\pd_x(\psi_+\tpsi)\psi_-\tpsi
=& 
\pd_x(\psi_+\psi_-\tpsi^2)
+\psi_-\tpsi(\pd_x\psi_+\tpsi-\psi_+\pd_x\tpsi)
-\psi_+\tpsi(\pd_x\psi_-\tpsi-\psi_-\pd_x\tpsi)
\\=& 
\pd_x(\psi_+\psi_-\tpsi^2)
+\psi_-\tpsi\pd_x^{-1}\pd_y(\psi_+\tpsi)
+\psi_+\tpsi(\pd_x^{-1})^*\pd_y(\psi_-\tpsi)\,,
\end{align*}
where $\pd_x^{-1}f=-\int_x^\infty f$.
Hence it follows that 
$$2 \int_\R  \pd_x(\psi_+\tpsi)\psi_-\tpsi\,dz_1=\pd_y
\int_\R \pd_x^{-1}(\psi_+\tpsi)\psi_-\tpsi\,dz_1\,.$$
By \eqref{eq:Phi22*-1} with $j=2$ and the fact that
\begin{equation*}
\beta x+\beta^2y+\beta'x+(\beta')^2y
-(\theta_2+\theta_3)|_{t=0}=
\{\gamma_{23}(-\eta_1)-\gamma_{23}(\eta)\}z_1+iy(\eta-\eta_1)\,,
\end{equation*}
we see that $e^{-iy(\eta-\eta_1)}\pd_x^{-1}(\psi_+\tpsi)\psi_-\tpsi$ converges to
\begin{align*}
I_{4,\pm}:=&\frac{(\beta-a_{14}\mp\sqrt{c_{14}})
(\beta'-a_{14}\mp\sqrt{c_{14}})}
{4(\k_2-a_{14}\mp\sqrt{c_{14}})(\k_3-a_{14}\mp\sqrt{c_{14}})}
e^{-\{\gamma_{23}(\eta)-\gamma_{23}(-\eta_1)\}z_1}
\sech^2Z_{1,\pm}
\\ & \qquad \times
\left(\gamma_{23}(-\eta_1)+\frac{\k_3-\k_2}{2}\tanh Z_{1,\pm}\right)\
\end{align*}
exponentially as $(a_{14}-a_{23})y\to\pm\infty$.
Using the Riemann-Lebesgue lemma, we have for
$\varphi\in C^1([-\eta_0,\eta_0])$,
\begin{align*}
  & 2\lim_{M\to\infty}\int_{-\eta_0}^{\eta_0}\varphi(\eta_{1,R})
\left(\int_{-M}^M\int_{-\infty}^\infty \pd_x(\psi_+\tpsi)\psi_-\tpsi\,dz_1dy\right)\,d\eta_{1,R}
  \\ =& \sum_\pm \pm\lim_{M\to\infty} \int_{-\eta_0}^{\eta_0} e^{\pm iM(\eta-\eta_1)}
\varphi(\eta_{1,R})\left( \int_\R I_{4,\pm}\,dz_1\right)\,d\eta_{1,R}=0\,.
\end{align*}
We can prove the rest of \eqref{eq:gg*-2} in the same way by using
the identities in Claim~\ref{cl:Phi22*s}.
Thus we complete the proof.
\end{proof}

\par

Next, we will prove the orthogonal relations of
$g^1_{1,\pm}$ and $g^{1,*}_{1,\pm}$ by using
the Darboux transformations $\nabla M_\pm(v_2)$ and
Lemma~\ref{lem:orthogonality}.
Note that $g^1_{1,\pm}$, $g^{1,*}_{1,\pm}$ and
 $g^2_{1,\pm}$, $g^{2,*}_{1,\pm}$ satisfy the same orthogonal relations.
\begin{lemma}
  \label{lem:orthogonality[2,3]}
  Let $\eta=\eta_R+i\eta_I$ and $\eta_1=\eta_{1,R}+i\eta_I$ and
  $\eta_0\in (0,4\sqrt{2}c_{23})$.
Suppose that $\eta_R\in(-\eta_0,\eta_0)$  and that $\eta_I\in (-c_{23},c_{23})$.
Then for $j=1$, $2$ and $\varphi\in C^1([-\eta_0,\eta_0])$,
\begin{gather}
\label{eq:gg*-1'}
\lim_{M\to\infty}\int_{-\eta_0}^{\eta_0}\varphi(\eta_{1,R})\left(\int_{-M}^M \int_{-\infty}^\infty
g^1_{j,\pm}(x,y,\eta)\overline{g^{1,*}_{j,\pm}(x,y,\eta_1)}\,dxdy\right)\,d\eta_{1,R}
=\pm2\pi i\eta\varphi(\eta_R)\,,
\\
\label{eq:gg*-2'}
\lim_{M\to\infty}\int_{-\eta_0}^{\eta_0}\varphi(\eta_{1,R})\left(\int_{-M}^M \int_{-\infty}^\infty
g^1_{j,\pm}(x,y,\eta)\overline{g^{1,*}_{j,\mp}(x,y,\eta_1)}\,dxdy\right)\,d\eta_{1,R}  =0\,.
\end{gather}
\end{lemma}
\begin{proof}[Proof of \eqref{eq:gg*-1}  and \eqref{eq:gg*-2} for $j=1$]
 For $u=u(x,y)$ and $v=v(x,y)$, let
$$\la u,v\ra=\lim_{M\to\infty}\int_{-M}^M \int_{-\infty}^\infty 
u(x,y)\overline{v(x,y)}\,dxdy\,.$$ 
Using Claim~\ref{cl:Phi22*s} and the Riemann-Lebesgue lemma,
we have for $j=1$, $2$,$a$ and $b=\pm$,
\begin{equation*}
\lim_{M\to\infty}\int_{-\eta_0}^{\eta_0}\varphi(\eta_{1,R})\left(  
\int_{-\infty}^\infty\pd_x^{-2}g^j_{1,a}(x,\pm M,\eta)
\overline{g^{j,*}_{1,b}(x,\pm M,\eta_1)}\,dx\right)\,d\eta_{1,R}=0\,,
\end{equation*}
and
\begin{equation}
\label{eq:int-MM*}
\la \nabla M_\pm(v_2)\pd_x^{-1}g^j_{1,a}(\eta),g^{j,*}_{1,b}(\eta_1)\ra
=\la \pd_x^{-1}g^j_{1,a}(\eta), \nabla M_\pm(v_2)^*g^{j,*}_{1,b}(\eta_1)\ra\,.  
\end{equation}
By  \eqref{eq:g-tg} and \eqref{eq:int-MM*},
\begin{align*}
\la \nabla M_+(v_2)\pd_x^{-1}g^2_{1,a}(\cdot,\eta), g^{2,*}_{1,b}(\cdot,\eta_1)\ra
=& 
\la \pd_x^{-1}g^2_{1,a}(\cdot,\eta), \nabla M_+(v_2)^*g^{2,*}_{1,b}(\cdot,\eta_1)\ra
\\ =&
\la \pd_x^{-1}g^2_{1,a}(\cdot,\eta), \nabla M_-(v_2)^*g^{1,*}_{1,b}(\cdot,\eta_1)\ra
\\ =&
\la \nabla M_-(v_2)\pd_x^{-1}g^2_{1,a}(\cdot,\eta), g^{1,*}_{1,b}(\cdot,\eta_1)\ra
\\=& 
\la \nabla M_+(v_2)\pd_x^{-1}g^1_{1,a}(\cdot,\eta), g^{1,*}_{1,b}(\cdot,\eta_1)\ra\,.
\end{align*}
Similarly,
\begin{equation*}
\la \nabla M_-(v_2)\pd_x^{-1}g^2_{1,a}(\cdot,\eta), g^{2,*}_{1,b}(\cdot,\eta_1)\ra 
= \la \nabla M_-(v_2)\pd_x^{-1}g^1_{1,a}(\cdot,\eta),
g^{1,*}_{1,b}(\cdot,\eta_1)\ra\,.
\end{equation*}
Combining the above with the fact that
$\nabla M_+(v_2)-\nabla M_-(v_2)=2\pd_x$, we have
\begin{equation}
  \label{eq:13}
\begin{split}
\la g^2_{1,a}(\cdot,\eta),g^{2,*}_{1,b}(\cdot,\eta_1)\ra
=& \frac12\la \{\nabla M_+(v_2)-\nabla M_-(v_2)\}\pd_x^{-1}g^2_{1,a}(\cdot,\eta),
g^{2,*}_{1,b}(\cdot,\eta_1)\ra
\\=&  \frac12
\la \{\nabla M_+(v_2)-\nabla M_-(v_2)\}\pd_x^{-1}g^1_{1,a}(\cdot,\eta),
g^{1,*}_{1,b}(\cdot,\eta_1)\ra
\\ =&
  \la g^1_{1,a}(\cdot,\eta),g^{1,*}_{1,b}(\cdot,\eta_1)\ra\,.
\end{split}  
\end{equation}
Lemma~\ref{lem:orthogonality[2,3]} follows immediately
from Lemma~\ref{lem:orthogonality} and \eqref{eq:13}.
Thus we complete the proof.
\end{proof}

We can also prove Lemma~\ref{lem:orthogonality[2,3]} by using
Lemma~\ref{lem:orthgz} below.
\begin{lemma}[Lemma~2.1 in \cite{Miz15}]
  \label{lem:orthgz}
Suppose that $|\Im\eta|< c_{23}$. Then
  \begin{gather*}
\int_\R g(z,\eta)\overline{g^*(z,\eta)}\,dz=i\eta\,,\quad
       \int_\R g(z,\eta)\overline{g^*(z,-\eta)}\,dz=0\,.   
  \end{gather*}
\end{lemma}

By the definitions,
\begin{align*}
&    \lim_{\eta\to0}g^k_{1,+}(x,y,\eta) 
=\lim_{\eta\to0}g^k_{1,-}(x,y,\eta)
=2\pd_x\left(\Phi^k_2(\bx)\Phi^{k,*}_2(\bx)\right)\bigl|_{t=0}\,,
\\ &
\lim_{\eta\to0}g^{k,*}_{1,+}(x,y,\eta)
=\lim_{\eta\to0}g^{k,*}_{1,-}(x,y,\eta)
=-\Phi^k_2(\bx)\Phi^{k,*}_2(\bx)|_{t=0}
\end{align*}
for $k=1$, $2$ and
\begin{align*}
& \lim_{\eta\to0}g^2_{2,+}(x,y,\eta) 
=\lim_{\eta\to0}g^2_{2,-}(x,y,\eta)
=2\pd_x\left(\Phi^2_1(\bx)\Phi^{2,*}_1(\bx)\right)\bigl|_{t=0}\,,
\\ &
\lim_{\eta\to0}g^{2,*}_{2,+}(x,y,\eta) =
\lim_{\eta\to0}g^{2,*}_{2,-}(x,y,\eta)
=-\Phi^2_1(\bx)\Phi^{2,*}_1(\bx)|_{t=0}\,.
\end{align*}
To resolve these degeneracy of
$\spann\{g^{k}_{j,\pm}(\eta)\}$ and $\spann\{g^{k,*}_{j,\pm}(\eta)\}$
at $\eta=0$, let
\begin{equation}
  \label{def:gjk12}
\left\{  
\begin{aligned}
& g^k_{j,1}(x,y,\eta)=\frac{1}{2}\{g^k_{j,+}(x,y,\eta)+g^k_{j,-}(x,y,\eta)\}\,,
\\ &
g^k_{j,2}(x,y,\eta)=\frac{1}{2i\eta}  
\{g^k_{j,+}(x,y,\eta)-g^k_{j,-}(x,y,\eta)\}\,,
\\ & 
g^{k,*}_{j,1}(x,y,\eta)=\frac{i}{\bar{\eta}}\left\{g^{k,*}_{j,+}(x,y,\eta)
-g^{k,*}_{j,-}(x,y,\eta)\right\}\,,
\\ & 
g^{k,*}_{j,2}(x,y,\eta)=
g^{k,*}_{j,+}(x,y,\eta)+g^{k,*}_{j,-}(x,y,\eta)\,.
\end{aligned}\right.
\end{equation}
Then we have the following.
\begin{lemma}
  \label{lem:gg*ev}
  For $j$, $k=1$, $2$,
  \begin{gather*}
\mL_k g^k_{j,1}(\eta)=
\frac{\lambda_{j,+}(\eta)+\lambda_{j,-}(\eta)}{2}g^k_{j,1}(\eta)
+i\eta\frac{\lambda_{j,+}(\eta)-\lambda_{j,-}(\eta)}{2}g^k_{j,2}(\eta)\,,
\\
\mL_k g^k_{j,2}(\eta)=
\frac{\lambda_{j,+}(\eta)-\lambda_{j,-}(\eta)}{2i\eta}g^k_{j,1}(\eta)
+\frac{\lambda_{j,+}(\eta)+\lambda_{j,-}(\eta)}{2}g^k_{j,2}(\eta)\,,
\\
\mL_k^*g^{k,*}_{j,1}(\eta)=
\frac{\lambda_{j,+}(-\bar{\eta})+\lambda_{j,-}(-\bar{\eta})}{2}g^{k,*}_{j,1}(\eta)
-\frac{\lambda_{j,+}(-\bar{\eta})-\lambda_{j,-}(-\bar{\eta})}{2i\bar{\eta}}
g^{k,*}_{j,2}(\eta)\,,
\\
\mL_k^*g^*_{j,2}(\eta)=
\bar{\eta}\frac{\lambda_{j,+}(-\bar{\eta})-\lambda_{2,-}(-\bar{\eta})}{2i}
g^{k,*}_{j,1}(\eta)
+\frac{\lambda_{j,+}(-\bar{\eta})+\lambda_{j,-}(-\bar{\eta})}{2}
g^{k,*}_{j,2}(\eta)\,.    
  \end{gather*}
\end{lemma}
\begin{lemma}
  \label{lem:orthogonality2}
Let $\eta=\eta_R+i\eta_I$, $\eta_1=\eta_{1,R}+i\eta_I$ and $\eta_0\in (0,4\sqrt{2}c_{23})$.
Suppose that $\eta_R\in(-\eta_0,\eta_0)$ and that $\eta_I\in(-\sqrt{c_{23}},\sqrt{c_{23}})$.
Then for $j$, $k=1$, $2$ and  $\varphi\in C^1([-\eta_0,\eta_0])$,
\begin{gather*}
  \lim_{M\to\infty}\int_{-\eta_0}^{\eta_0}\varphi(\eta_{1,R})\left(
  \int_{-M}^M \int_{-\infty}^\infty
g^2_{j,k}(x,y,\eta)\overline{g^{2,*}_{j,k}(x,y,\eta_1)}\,dxdy\right)\,d\eta_{1,R}
=2\pi\delta_{jk}\varphi(\eta_R)\,,
\\
  \lim_{M\to\infty}\int_{-\eta_0}^{\eta_0}\varphi(\eta_{1,R})\left(
  \int_{-M}^M \int_{-\infty}^\infty
g^1_{j,k}(x,y,\eta)\overline{g^{1,*}_{j,k}(x,y,\eta_1)}\,dxdy\right)\,d\eta_{1,R}
=2\pi\delta_{jk}\varphi(\eta_R)\,.
\end{gather*}
\end{lemma}
Lemma~\ref{lem:orthogonality2} follows immediately from
Lemmas~\ref{lem:orthogonality} and \ref{lem:orthogonality[2,3]}.
\par
Finally, we introduce continuous eigenfunctions associated with
linearized modified KP-II equations. Let
\begin{gather}
  \label{eq:defg1M}
g^1_M(\tbx,\eta)=(\k_3-\k_2)d_{1,+}(\eta)
\Phi^1(\bx,-i\beta_{23}^-(\eta))\Phi^{1,*}_2(\bx)\,,
\\ \label{eq:defg1M*}
g^{1,*}_M(\tbx,\eta)
=\frac{1}{\k_3-\k_2}\Phi^1_2(\bx)\Phi^{0,*}(\bx,-i\beta_{23}^-(-\eta))\,,
\\   \label{eq:defg2M}
g^2_M(\tbx,\eta)=d_{2,+}(\eta)(\k_4-\k_1)
\Phi^2(\bx,-i\beta_{14}^-(\eta))\Phi^{2,*}_1(\bx)\,,
\\   \label{eq:defg2M*}
g^{2,*}_M(\tbx,\eta)=\frac{1}{\k_4-\k_1}
\Phi^2_1(\bx)\Phi^{1,*}(\bx,-i\beta_{14}^-(-\eta))\,.
\end{gather}
Then we have the following.
\begin{lemma}
  \label{lem:eigenfunctions-mKP}
  For $k=1$ and $2$, $  g^k_M(\cdot,\eta)=\pd_x^{-1}g^k_{k,+}(\cdot,\eta)$,
and
\begin{gather}
\label{eq:12}    
\nabla M_-(v_k) g^k_M(\cdot,\eta)=0\,, \quad
\nabla M_+(v_k) g^k_M(\cdot,\eta)=2g^k_{k,+}(\cdot,\eta)\,,\\
\label{eq:12-1}
g^{k,*}_M(\cdot,\eta)
=\frac{1}{2}\nabla M_+(v_k)^*g^{k,*}_{k,1}(\cdot,\eta)\,,
\\
\label{eq:orthgMgM*}
\la g^k_M(\cdot,\eta), g_M^{k,*}(\cdot,\eta_1)\ra
= 2\pi\delta(\eta-\eta_1)\,.
\end{gather}
\end{lemma}
\begin{proof}
  We see that \eqref{eq:12} follows from Lemma~\ref{lem:g-tg} and that
  \eqref{eq:12-1} follows from \eqref{eq:Mg2-'} and \eqref{eq:kernablaM}.
Combining  \eqref{eq:M+vi} with \eqref{def:gjk12} and
Lemma~\ref{lem:orthogonality2},
we can prove
\begin{equation*}
 \la g^k_M(\cdot,\eta), g_M^{k,*}(\cdot,\eta_1)\ra
 = \la g^k_{k,+}(\cdot,\eta), g^{k,*}_{k,1}(\cdot,\eta_1)\ra
= 2\pi\delta(\eta-\eta_1)\,,
\end{equation*}
in the same way as \eqref{eq:int-MM*}. 
\end{proof}

\bigskip

\subsection{Spectral projections}
\label{subsec:specproj}
Let $\calX_{1}=L^2(\R^2;e^{2\a z_1}dxdy)$ and
$\calX_{2}=L^2(\R^2;e^{2\a z_2}dxdy)$ for an $\a>0$.
 Let
\begin{align*}
  &  g_1(z_1,\eta)=\frac{1}{2}
\left\{g(z_1,\eta)+g(z_1,-\eta)\right\}\,,
\quad
g_2(z_1,\eta)=\frac{1}{2i\eta}
\left\{g(z_1,\eta)-g(z_1,-\eta)\right\}\,,
\\ &
g^*_1(z_1,\eta)=\frac{i}{\bar{\eta}}
\left\{g^*(z_1,\eta)-g^*(z_1,-\eta)\right\}\,,
\quad
g^*_2(z_1,\eta)=g^*(z_1,\eta)+g^*(z_1,-\eta)\,.
\end{align*}
Then $g^1_{1,j}(x,y,\eta)=g_j(z_1,\eta)e^{iy\eta}$
and $g^{1,*}_{1,j}(x,y,\eta)=g_j^*(z_1,\eta)e^{iy\bar{\eta}}$ for $j=1$ and $2$.

Let $\a\in(0,2\sqrt{c_{23}})$ and $\eta_{*,1}$  be a positive number satisfying
\begin{gather}
  \label{def:eta*1}
  \Re\gamma_{23}(\eta_{*,1})=\sqrt{c_{23}}+\a\,.
\end{gather}
Let $\eta_0\in(0,\eta_{*,1})$ and
\begin{gather*}
P_1(\eta_0)f(z_1,y)=\frac{1}{\sqrt{2\pi}}
\sum_{j=1,2}\int_{-\eta_0}^{\eta_0}
\tilde{a}_j(\eta)g_j(z_1,\eta)e^{iy\eta}\,d\eta\,,
\\
a_j(\eta)=\frac{1}{\sqrt{2\pi}}
\int_{\R^2}f(z_1,y)\overline{g^*_j(z_1,\eta)}e^{-iy\eta}\,dz_1dy\,.
  \end{gather*}
Then it follows from \cite[Lemma~3.1]{Miz15} that
$P_1(\eta_0):\calX_1\mapsto\calX_1$ is a spectral projection associated with
$\mL_1$.
\par

Next we introduce a spectral projection on $\calX_2$ associated with $\mL_2$.
Let $P_2(\eta_0)$ be an operator defined by
\begin{equation}
  \label{eq:defP2}
\left\{\begin{gathered}
P_2(\eta_0)f(x,y)=\frac{1}{\sqrt{2\pi}}\sum_{k=1,2}\int_{-\eta_0}^{\eta_0}
a_k(\eta)g^2_{2,k}(x,y,\eta)\,d\eta\,,
\\
a_k(\eta)=\frac{1}{\sqrt{2\pi}}\lim_{M\to\infty}\int_{-M}^M\int_\R f(x,y)
\overline{g^{2,*}_{2,k}(x,y,\eta)}\,dxdy\,,
\end{gathered}
\right.
\end{equation}
and let $\eta_{*,2}$ be a positive number satisfying
\begin{gather}
    \label{def:eta*2}
\Re\gamma_{14}(\eta_{*,2})=\sqrt{c_{14}}+\a\,. 
\end{gather}
We will show that
$P_2(\eta_0): \calX_2\to 
\spann\{g^2_{2,1}(\eta),\,g^2_{2,2}(\eta)\}_{|\eta|\le \eta_0}$
is a spectral projection associated with $\mL_2$
if $\eta_0\in(0,\eta_{*,2})$.
\begin{lemma}
  \label{lem:spectral-proj}
Let $\a\in(0,2\sqrt{c_{14}})$, $\eta_1\in(0,\eta_{*,2})$ and
let $\mL_2$ be a closed operator on $\calX_2$.
Suppose that $\eta_0\in(0,\eta_1]$.
Then
\begin{enumerate}
\item \label{it:p01}
$\|P_2(\eta_0)f\|_{\calX_2}
\le C\|f\|_{\calX_2}$  for any $f\in \calX_2$,
where $C$ is a positive constant depending only on $\a$ and $\eta_1$,
\item \label{it:p03}
$\mL_2 P_2(\eta_0)f=P_2(\eta_0)\mL_2 f$ for any $f\in D(\mL_2)$,
\item \label{it:p04}
$P_2(\eta_0)^2=P_2(\eta_0)$ on $\calX_2$,
\item \label{it:p05}
$e^{t\mL_2}P_2(\eta_0)=P_2(\eta_0)e^{t\mL_2}$ on $\calX_2$.
\end{enumerate}
\end{lemma}
Before we start to prove Lemma~\ref{lem:spectral-proj}, let us introduce
a linearized operator of \eqref{eq:KPII} around $[1,4]$-soliton and
its continuous eigenfunctions.
Let $\tilde{u}_1=2\pd_x^2\log(e^{\theta_1}+e^{\theta_4})
=2c_{14}\sech^2\sqrt{c_{14}}z_2$ and let
$$
\wmL_1= \frac14\{-\pd_x(\pd_x^2-4b_1+6\tu_1)
+4b_2\pd_y-3\pd_x^{-1}\pd_y^2\}$$
be an operator on $\calX_2$.
Then $\widetilde{\mL_1}$ has continuous eigenfunctions
associated with $\lambda_{2,\pm}(\eta)$. Let
$\wmL_1(\eta)u(z_2):=e^{-iy\eta}\widetilde{\mL_1}
(e^{iy\eta}u(z_2))$ and
\begin{equation}
  \label{eq:defg1z2}
\left\{  
\begin{aligned}
&  \tg(z_2,\eta)=\sqrt{c_{14}}d_{2,+}(\eta)\pd_x^2\left\{
    e^{-\gamma_{14}(\eta)z_2}\sech\sqrt{c_{14}}z_2\right\}\,,\\
&   \tg^*(z_2,\eta)=\frac12\pd_x\left\{
    e^{\gamma_{14}(-\bar{\eta})z_2}\sech\sqrt{c_{14}}z_2\right\}\,.
\end{aligned}
\right.
\end{equation}
Then we have the following.
\begin{lemma} 
  \label{lem:evtmL1}
  For $\eta\in\C$,
  \begin{equation*}
\wmL_1(\eta) \tg(z_2,\pm\eta)=
\lambda_{2,\pm}(\eta)\tg(z_2,\pm\eta)\,,\quad
\wmL_1(\eta)^* \tg^*(z_2,\pm\eta)
 =\overline{\lambda_{2,\pm}(\eta)}\tg^*(z_2,\pm\eta)\,.    
\end{equation*}
If $|\Im\eta|<\sqrt{c_{14}}$,
  \begin{gather*}
 \int_\R \tg(z_2,\eta)\overline{\tg^*(z_2,\eta)}\,dz_2=i\eta\,,
\quad
  \int_\R \tg(z_2,\eta)\overline{\tg^*(z_2,-\eta)}\,dz_2=0\,.  
\end{gather*}
\end{lemma}

Let
\begin{align*}
  & \tg_1(z_2,\eta)=\frac{1}{2}
\left\{\tg(z_2,\eta)+\tg(z_2,-\eta)\right\}\,,
\\ &
\tg_2(z_2,\eta)=\frac{1}{2i\eta}
\left\{\tg(z_2,\eta)-\tg(z_2,-\eta)\right\}\,,
\\ &
\tg^*_1(z_2,\eta)=\frac{i}{\bar{\eta}}
\left\{\tg^*(z_2,\eta)-\tg^*(z_2,-\eta)\right\}\,,
\\ &
\tg^*_2(z_2,\eta)=\tg^*(z_2,\eta)+\tg^*(z_2,-\eta)\,,
\end{align*}
and 
\begin{gather*}
  \widetilde{P}_1(\eta_0)f(z_2,y)=\frac{1}{\sqrt{2\pi}}
\sum_{j=1,2}\int_{-\eta_0}^{\eta_0}
\tilde{a}_j(\eta)\tg_j(z_2,\eta)e^{iy\eta}\,d\eta\,,
\\
\tilde{a}_j(\eta)=\frac{1}{\sqrt{2\pi}}
\int_{\R^2}
f(z_2,y)\overline{\tg^*_j(z_2,\eta)}e^{-iy\eta}\,dz_2dy
=\int_\R (\mathcal{F}_yf)(z_2,\eta)\overline{\tg^*_j(z_2,\eta)}\,dz_2\,.
\end{gather*}
We remark that $\widetilde{P}_1(\eta_0)$ is a projection from $\calX_2$
to $\spann\{\tg_1(\cdot,\eta)e^{iy\eta},
\tg_2(\cdot,\eta)e^{iy\eta}\}_{-\eta_0\le\eta\le\eta_0}\}$
if $\eta_0\in(0,\eta_{*,2})$.

To prove Lemma~\ref{lem:spectral-proj}, we need the asymptotic profiles
of $g^2_{j,k}$ and $g^{2,*}_{j,k}$ as $y\to\pm\infty$.

\begin{claim}  
  \label{cl:gg*-asymp} 
As $z_1\to\infty$,
  \begin{gather*}
g^2_{2,\pm}(x,y,\eta)=c_{\pm,3}(\eta)\tg(\frac{Z_{2,+}}{\sqrt{c_{14}}},\pm\eta)
e^{iy\eta}(1+O(e^{-2\sqrt{c_{23}}z_1}))\,,
\\
g^{2,*}_{2,\pm}(x,y,\eta)=\frac{1}{c_{\pm,3}(-\bar{\eta})}
\tg^*(\frac{Z_{2,+}}{\sqrt{c_{14}}},\pm\eta)e^{iy\bar{\eta}}
(1+O(e^{-2\sqrt{c_{23}}z_1}))\,,
\end{gather*}
and as $z_1\to-\infty$,
\begin{gather*}
g^2_{2,\pm}(x,y,\eta)=c_{\pm,2}(\eta)\tg(\frac{Z_{2,-}}{\sqrt{c_{14}}},\pm\eta)
e^{iy\eta}(1+O(e^{2\sqrt{c_{23}}z_1}))\,,
\\
 g^{2,*}_{2,\pm}(x,y,\eta)
 =\frac{1}{c_{\pm,2}(-\bar{\eta})}
 \tg^*(\frac{Z_{2,-}}{\sqrt{c_{14}}},\pm\eta)e^{iy\bar{\eta}}
(1+O(e^{2\sqrt{c_{23}}z_1}))\,,
\end{gather*}
where for $j=2$ and $3$,
  \begin{align*}
& c_{+,j}(\eta)=
(\kappa_j-\beta_{14}^-(\eta))
(\kappa_4-\kappa_j)^{-\frac12(1-\gamma_{14}(\eta)/\sqrt{c_{14}})}
(\kappa_j-\kappa_1)^{-\frac12(1+\gamma_{14}(\eta)/\sqrt{c_{14}})}\,,
\\
& c_{-,j}(\eta)=
\frac{(\kappa_4-\kappa_j)^{\frac12(1+\gamma_{14}(-\eta)/\sqrt{c_{14}})}
(\kappa_j-\kappa_1)^{\frac12(1-\gamma_{14}(-\eta)/\sqrt{c_{14}})}}
{\beta_{14}^+(-\eta)-\k_j}\,.
  \end{align*}
  \par
As $z_2\to\infty$,
  \begin{gather*}
g^2_{1,\pm}(x,y,\eta)=c_{\pm,4}(\eta)g(\frac{Z_{1,+}}{\sqrt{c_{23}}},\pm\eta)
e^{iy\eta}(1+O(e^{-2\sqrt{c_{14}}z_2}))\,,
\\
g^{2,*}_{1,\pm}(x,y,\eta)=\frac{1}{c_{\pm,4}(-\bar{\eta})}
g^*(\frac{Z_{1,+}}{\sqrt{c_{23}}},\pm\eta)e^{iy\bar{\eta}}
(1+O(e^{-2\sqrt{c_{14}}z_2}))\,,
\end{gather*}
and as $z_2\to-\infty$,
\begin{gather*}
g^2_{1,\pm}(x,y,\eta)=c_{\pm,1}(\eta)g(\frac{Z_{1,-}}{\sqrt{c_{23}}},\pm\eta)
e^{iy\eta}(1+O(e^{2\sqrt{c_{14}}z_2}))\,,
\\
g^{2,*}_{1,\pm}(x,y,\eta)
=\frac{1}{c_{\pm,1}(-\bar{\eta})}
g^*(\frac{Z_{1,-}}{\sqrt{c_{23}}},\pm\eta)e^{iy\bar{\eta}}
(1+O(e^{2\sqrt{c_{14}}z_2}))\,,
\end{gather*}
where for $j=1$ and $4$,
  \begin{align*}
& c_{+,j}(\eta)=(-1)^j(\k_j-\beta_{23}^-(\eta))
|\kappa_j-\kappa_3|^{-\frac12(1-\gamma_{23}(\eta)/\sqrt{c_{23}})}
|\kappa_j-\kappa_2|^{-\frac12(1+\gamma_{23}(\eta)/\sqrt{c_{23}})}\,,
\\
& c_{-,j}(\eta)=(-1)^{j-1}
\frac{|\kappa_j-\kappa_3|^{\frac12(1+\gamma_{23}(-\eta)/\sqrt{c_{23}})}
|\kappa_j-\kappa_2|^{\frac12(1-\gamma_{23}(-\eta)/\sqrt{c_{23}})}}
{\beta_{23}^+(-\eta)-\k_j}\,.
\end{align*}
 \end{claim}
\begin{proof}
By \eqref{eq:btitj}, \eqref{eq:pd12-22} and
\eqref{eq:Phi22*-1}--\eqref{eq:Phi22*-2},
\begin{gather}
 \label{eq:g21pmpm*}
g^2_{1,\pm}(x,y,\eta)=e^{iy\eta}\sum_{j=1,4}\pd_x^2g^2_{1,\pm,j}(x,y,\eta)\,,
\quad
g^{2,*}_{1,\pm}(x,y,\eta)=e^{iy\bar{\eta}}\sum_{j=1,4}
\pd_xg^{2,*}_{1,\pm,j}(x,y,\eta)\,,
\\  \label{eq:g2pmpm*}
g^2_{2,\pm}(x,y,\eta)=e^{iy\eta}\sum_{j=2,3}\pd_x^2g^2_{2,\pm,j}(x,y,\eta)\,,
\quad
g^{2,*}_{2,\pm}(x,y,\eta)=e^{iy\bar{\eta}}\sum_{j=2,3}\pd_xg^{2,*}_{2,\pm,j}(x,y,\eta)\,,
\end{gather}
where
\begin{align}
  & \label{eq:g21+j}
g^2_{1,+,j}=(-1)^j(\k_3-\k_2)d_{1,+}(\eta)(\k_j-\beta_{23}^-(\eta))
e^{-\gamma_{23}(\eta)z_1}
\frac{e^{\theta_j+(\theta_2+\theta_3)/2}}{\tau_2}\bigr|_{t=0}\,,
\\ & \label{eq:g21-j}
 g^2_{1,-,j}=(-1)^j(\k_3-\k_2)d_{1,-}(\eta)
\frac{(\k_j-\k_2)(\k_j-\k_3)}{\beta_{23}^+(-\eta)-\k_j}
e^{-\gamma_{23}(-\eta)z_1}
\frac{e^{\theta_j+(\theta_2+\theta_3)/2}}{\tau_2}\bigr|_{t=0}\,
\\ & \label{eq:g21*+j}
  g^{2,*}_{1,+,j}=(-1)^{j}\frac{(\k_j-\k_2)(\k_j-\k_3)}
  {\k_j-\beta_{23}^-(-\bar{\eta})}
e^{\gamma_{23}(-\bar{\eta})z_1}
    \frac{e^{\theta_j+(\theta_2+\theta_3)/2}}{\tau_2}\bigr|_{t=0}\,
\\ & \label{eq:g21*-j}
g^{2,*}_{1,-,j}=(-1)^j(\k_j-\beta_{23}^+(\bar{\eta}))e^{\gamma_{23}(\bar{\eta})z_1}
    \frac{e^{\theta_j+(\theta_2+\theta_3)/2}}{\tau_2}\bigr|_{t=0}\,,   
\end{align}
\begin{align}
  & \label{eq:g2+j}
g^2_{2,+,j}=(\k_4-\k_1)d_{2,+}(\eta)(\k_j-\beta_{14}^-(\eta))
e^{-\gamma_{14}(\eta)z_2}
\frac{e^{\theta_j+(\theta_1+\theta_4)/2}}{\tau_2}\bigr|_{t=0}\,,
\\ & \label{eq:g2-j}
 g^2_{2,-,j}=-(\k_4-\k_1)d_{2,-}(\eta)
\frac{(\k_j-\k_1)(\k_4-\k_j)}{\beta_{14}^+(-\eta)-\k_j}
e^{-\gamma_{14}(-\eta)z_2}
\frac{e^{\theta_j+(\theta_1+\theta_4)/2}}{\tau_2}\bigr|_{t=0}\,
\\ & \label{eq:g2*+j}
  g^{2,*}_{2,+,j}=\frac{(\k_j-\k_1)(\k_4-\k_j)}
  {\k_j-\beta_{14}^-(-\bar{\eta})}
e^{\gamma_{14}(-\bar{\eta})z_2}
    \frac{e^{\theta_j+(\theta_1+\theta_4)/2}}{\tau_2}\bigr|_{t=0}\,
\\ & \label{eq:g2*-j}
g^{2,*}_{2,-,j}=(\beta_{14}^+(\bar{\eta})-\k_j)e^{\gamma_{14}(\bar{\eta})z_2}
    \frac{e^{\theta_j+(\theta_1+\theta_4)/2}}{\tau_2}\bigr|_{t=0}\,.
\end{align}
Combining the above with \eqref{eq:tau2-23} and \eqref{eq:tau2-14},
we have Claim~\ref{cl:gg*-asymp}.
\end{proof}

\begin{lemma}
  \label{lem:gg*2-asymp}
As $z_1\to\pm\infty$,
\begin{align*}
  &  g^2_{2,1}= e^{iy\eta}\left\{(1+O(\eta))
\tg_1(\frac{Z_{2,\pm}}{\sqrt{c_{14}}},\eta)
+O(\eta^2)\tg_2(\frac{Z_{2,\pm}}{\sqrt{c_{14}}},\eta)\right\}
(1+O(e^{-2\sqrt{c_{23}}|z_1|}))\,,
\\ &              
     g^2_{2,2}= e^{iy\eta}\left\{
(1+O(\eta)) \tg_2(\frac{Z_{2,\pm}}{\sqrt{c_{14}}},\eta)
-\left(c_\pm+O(\eta)\right)\tg_1(\frac{Z_{2,\pm}}{\sqrt{c_{14}}},\eta)\right\}
(1+O(e^{-2\sqrt{c_{23}}|z_1|}))\,,
\\ &              
g^{2,*}_{2,1}= e^{iy\bar{\eta}}\left\{
(1+O(\eta))  \tg^*_1(\frac{Z_{2,\pm}}{\sqrt{c_{14}}},\eta)
+\left(c_\pm+O(\eta)\right)\tg_2^*(\frac{Z_{2,\pm}}{\sqrt{c_{14}}},\eta)\right\}
(1+O(e^{-2\sqrt{c_{23}}|z_1|}))\,,
\\ &              
g^{2,*}_{2,2}= e^{iy\bar{\eta}}\left\{
(1+O(\eta))\tg^*_2(\frac{Z_{2,\pm}}{\sqrt{c_{14}}},\eta)
+O(\eta^2)\tg^*_2(\frac{Z_{2,\pm}}{\sqrt{c_{14}}},\eta)\right\}
(1+O(e^{-2\sqrt{c_{23}}|z_1|}))\,,
\end{align*}
where $c_+=i\{\pd_\eta c_{+,3}(0)-\pd_\eta c_{-,3}(0)\}/2$ and
$c_-=i\{\pd_\eta c_{+,2}(0)-\pd_\eta c_{-,2}(0)\}/2$ are real numbers.
Moreover, for $k=1$, $2$ and $\eta_0\in(0,\eta_{*,2})$,
  \begin{equation*}
\sup_{|\eta|\le \eta_0}\left(\|e^{\a z_2}g^2_{2,k}(\eta)\|_{L^\infty_yL^2_{z_2}}+
\|e^{-\a z_2}g^{2,*}_{2,k}(\eta)\|_{L^\infty_yL^2_{z_2}}\right)<\infty\,.
  \end{equation*}
\end{lemma}

\begin{lemma}
  \label{lem:gg*1-asymp}
As $z_2\to\pm\infty$,
\begin{align*}
&  g^2_{1,1}= e^{iy\eta}\left\{
\left(1+O(\eta)\right)g_1(\frac{Z_{1,\pm}}{\sqrt{c_{23}}},\eta)
+ O(\eta^2) g_2(\frac{Z_{1,\pm}}{\sqrt{c_{23}}},\eta)
\right\} (1+O(e^{-2\sqrt{c_{14}}|z_2|}))\,,
\\ &              
g^2_{1,2}= e^{iy\eta}\left\{
\left(1+O(\eta)\right)g_2(\frac{Z_{1,\pm}}{\sqrt{c_{23}}},\eta)     
-\left(c_\pm+O(\eta)\right) g_1(\frac{Z_{1,\pm}}{\sqrt{c_{23}}},\eta)
\right\} (1+O(e^{-2\sqrt{c_{14}}|z_2|}))\,,
\\ &              
g^{2,*}_{1,1}= e^{iy\bar{\eta}}\left\{
\left(1+O(\eta)\right) g^*_1(\frac{Z_{1,\pm}}{\sqrt{c_{23}}},\eta)
+\left(c_\pm+O(\eta)\right) g^*_2(\frac{Z_{1,\pm}}{\sqrt{c_{23}}},\eta)\right\}
(1+O(e^{-2\sqrt{c_{14}}|z_2|}))\,,
\\ &              
g^{2,*}_{1,2}= e^{iy\bar{\eta}}\left\{
\left(1+O(\eta)\right) g^*_2(\frac{Z_{1,\pm}}{\sqrt{c_{23}}},\eta)
+O(\eta^2) g^*_1(\frac{Z_{1,\pm}}{\sqrt{c_{23}}},\eta)\right\}
     (1+O(e^{-2\sqrt{c_{14}}|z_2|}))\,,
\end{align*}
where $c_+=i\{\pd_\eta c_{+,4}(0)-\pd_\eta c_{-,4}(0)\}/2$ and
$c_-=i\{\pd_\eta c_{+,1}(0)-\pd_\eta c_{-,1}(0)\}/2$ are real numbers.  
Moreover, for $k=1$, $2$  and $\eta_0\in(0,\eta_{*,1})$,,
  \begin{equation*}
\sup_{|\eta|\le \eta_0}\left(
\|e^{\a z_1}g^2_{1,k}(\eta)\|_{L^\infty_yL^2_{z_1}}+
\|e^{-\a z_1}g^{2,*}_{1,k}(\eta)\|_{L^\infty_yL^2_{z_1}}\right)<\infty\,.
\end{equation*}
\end{lemma}

Lemmas~\ref{lem:gg*2-asymp} and \ref{lem:gg*1-asymp} immediately follow
from Claim~\ref{cl:gg*-asymp} below. Note that $c_{\pm,j}(0)=1$.
Now we are in position to prove Lemma~\ref{lem:spectral-proj}.
\begin{proof}[Proof of Lemma~\ref{lem:spectral-proj}]
  Since mappings
  \begin{gather*}
    \calX_2 \ni f\mapsto \int_{\R^2}f(x,y)\overline{\tg^*_k(z_2)}e^{-iy\eta}\in
    L^2(-\eta_0,\eta_0)\,,\\
     L^2(-\eta_0,\eta_0)\ni a\mapsto \int_{-\eta_0}^{\eta_0}a(\eta)\tg_k(z_2,\eta)
    e^{iy\eta}\,d\eta\in\calX_2
  \end{gather*}
  are bounded on $\calX_2$, it follows from Lemma~\ref{lem:gg*2-asymp} that
  $P_2(\eta_0)$ is bounded on $\calX_2$.
  We have \eqref{it:p03} from Lemma~\ref{lem:gg*ev} and
\eqref{it:p05} from \eqref{it:p03} since $\mL_2$ generates a $C^0$-semigroup
on $\calX_2$.
\par

Finally, we will prove $P_2(\eta_0)^2=P_2(\eta_0)$.
Let $\varphi$, $\psi\in C_0^\infty(\R^2)$, $\tbx=(x,y)$,
$\tbx'=(x',y')$, $\tbx''=(x'',y'')$ and $z_2'=x'+2a_{14}y'$. Then 
\begin{align*}
  & \left\la P_2(\eta_0)^2\varphi,\psi\right\ra
\\ =&
 \frac{1}{(2\pi)^2}\sum_{j,k=1,2}
      \int_{\R^2}d\tbxdd\overline{\psi(\tbxdd)}  \int_{-\eta_0}^{\eta_0}d\eta'
g^2_{2,k}(\tbxdd,\eta')  \lim_{M\to\infty}
\int_{\R^2}d\tbxd \varphi(\tbxd)I_{j,k}(\bxd,\eta',M)\,,
\end{align*}
where
\begin{equation*}
  I_{j,k}(\bxd,\eta',M)=\int_{-\eta_0}^{\eta_0}\overline{g^{2,*}_{2,j}(\tbxd,\eta)}
 \left( \int_{-M}^M\int_{-\infty}^\infty
g^2_{2,j}(\tbx,\eta)\overline{g^{2,*}_{2,k}(\tbx,\eta')}\,d\tbx\right)
\, d\eta\,.
\end{equation*}
It follows from Lemmas~~\ref{lem:orthogonality2} and \ref{lem:gg*2-asymp}
that $e^{-\a z_2'}I_{j,k}(\bxd,\eta,M)$ is uniformly bounded and
\begin{equation*}
  \lim_{M\to\infty}I_{j,k}(\bxd,\eta,M)=
  \begin{cases}
     2\pi \overline{g^{2,*}_{2,k}(\bxd,\eta')} & \text{if $j=k$,}
      \\  0 &  \text{if $j\ne k$.}
  \end{cases}
\end{equation*}
Hence it follows that 
\begin{align*}
  & \left\la P_2(\eta_0)^2\varphi,\psi\right\ra
  \\=& \frac{1}{2\pi}\sum_{k=1,2}
\int_{\R^2}d\tbxdd\, \overline{\psi(\tbxdd)}
\int_{-\eta_0}^{\eta_0}d\eta' g^2_{2,k}(\tbxdd,\eta')
\left(\int_{\R^2}d\tbxd\, \varphi(\tbxd)\overline{g^{2,*}_{2,j}(\tbxd,\eta')}\right)
\\=& \la P_2(\eta_0)\varphi,\psi\ra\,.
\end{align*}
Thus we complete the proof.
\end{proof}

Finally, we introduce a spectral projection on $\calX_1$.
Let $\eta=\eta_R+i\eta_I$ with $\eta_I=2(a_{23}-a_{14})\a$ and
$P_2'(\eta_{0,1},\tilde{\eta}_{0,2})=P_{21}'(\eta_{0,1})+P_{22}'(\tilde{\eta}_{0,2})$,
where
\begin{gather*}
P_{21}'(\eta_{0,1})f(x,y)=
\frac{1}{\sqrt{2\pi}}\sum_{k=1,2}\int_{-\eta_{0,1}}^{\eta_{0,1}}
a_{1,k}(\eta_R)g^2_{1,k}(\tbx,\eta_R)\,d\eta_R\,,
  \\ 
P_{22}'(\tilde{\eta}_{0,2})f(x,y)=       
\frac{1}{\sqrt{2\pi}}\sum_{k=1,2}\int_{-\tilde{\eta}_{0,2}}^{\tilde{\eta}_{0,2}}
a_{2,k}(\eta)g^2_{2,k}(\tbx,\eta)\,d\eta_R\,,
\\
 a_{j,k}(\eta)= \frac{1}{\sqrt{2\pi}}\int_{\R^2}
f(x,y)\overline{g^{2,*}_{j,k}(\tbx,\eta)}\,d\tbx\,.
\end{gather*}
\begin{lemma}
  \label{lem:projP21}
Let
  \begin{gather}
    \label{ass-alpha'}
    0<\a<\min\{\sqrt{c_{23}},\,\frac{c_{14}}{3|a_{14}-a_{23}|}\,,
    \,\k_2-\k_1\,,\,\k_4-\k_3\}\,,
\\  \label{eq:l2ne0}
\a\ne \pm\frac{(b_2-3a_{14})^2-c_{14}}{2(a_{14}-a_{23})}\,.  
  \end{gather}
Let $\eta_{0,1}$ be a small positive number,
$\tilde{\eta}_{0,2}\in (\tilde{\eta}_{*,+},\eta_+')$ if $a_{14}>a_{23}$ and
$\tilde{\eta}_{0,2}\in (\tilde{\eta}_{*,-},\eta_-')$ if $a_{14}<a_{23}$.
Let $\mL_2$ be a closed operator on $\calX_1$.
Then we have the following.
\begin{enumerate}
\item \label{it:p011}
$\|P_2'(\eta_{0,1},\tilde{\eta}_{0,2})f\|_{\calX_1} \le C\|f\|_{\calX_1}$ for any $f\in \calX_1$,
where $C=C(\eta_{0,1},\tilde{\eta}_{0,2})$ is a positive constant depending only on $\a$, $\eta_{0,1}$ and $\tilde{\eta}_{0,2}$ and
$C(\eta_{0,1},\tilde{\eta}_{0,2})$ is bounded as
$(\eta_{0,1},\tilde{\eta}_{0,2})\to(0,0)$.
\item \label{it:p031}
  $\mL_2P_2'(\eta_{0,1},\tilde{\eta}_{0,2})f
  =P_2'(\eta_{0,1},\tilde{\eta}_{0,2})\mL_2 f$ for any $f\in D(\mL_1)$.
\item \label{it:p041}
  $P_2'(\eta_{0,1},\tilde{\eta}_{0,2})^2
  =P_2'(\eta_{0,1},\tilde{\eta}_{0,2})$ on $\calX_1$.
\item \label{it:p051}
  $e^{t\mL_2}P_2'(\eta_{0,1},\tilde{\eta}_{0,2})
  =P_2'(\eta_{0,1},\tilde{\eta}_{0,2})e^{t\mL_2}$ on $\calX_1$.
\end{enumerate}
\end{lemma}
By Lemma~\ref{lem:lambda2-etaI}, curves $\{\lambda_{1,\pm}(\eta)\}_{|\eta|\le \eta_{0,1}}$ and
$\{\lambda_{2,\pm}(\eta_R+i\eta_I)\}_{|\eta_R|\le \tilde{\eta}_{0,2}}$ do not intersect
provided $\eta_{0,1}$ is sufficiently small, and
we can prove Lemma~\ref{lem:projP21} in the same way as Lemma~\ref{lem:spectral-proj}
by using Lemma~\ref{lem:gg*1-asymp}.

\bigskip

\section{Green functions of Darboux transformations
  for $2$-line soliton of P-type}
In this section, we will construct Green functions of
the Darboux transformations $\nabla M_\pm$ using an intertwining property 
of $\nabla M_\pm(v_i)$ and the Lax operators $L_i$
(see \eqref{eq:Miura-Lax1}--\eqref{eq:Miura-Lax4}).

\subsection{Green functions of Lax operators}
To begin with, we will introduce Green functions of $L_0$ and $L_1$
and their adjoints. 
Using the Jost solutions and the dual Jost solutions of the Lax pair,
Boiti, Pempinelli and Pogrebkov constructed Green functions
of $L=-\pd_y+\pd_x^2 +u_N$, where $u_N$ is a soliton potential
(see e.g.~\cite{BPP, BPP12}). They regard the operator $L$ as a heat operator
and use Green functions of $L$ which has discontinuity in $y$.
For example, their Green function for $L_0$ is 
$-(4\pi y)^{-1/2}e^{-x^2/4y}$ for $y>0$ and $0$ for $y<0$.

To have $L^2$-boundedness in the transverse variable,
we will use Green functions oscillating in $y$
whose $x$-derivative has a jump discontinuity
along the crest of the largest line soliton.
\par
Let  $\tbx=(x,y)$, $\tbxd=(x',y')$, $\eta=\eta_R+i\eta_I$ and
$(\eta_R,\eta_I)\in\R^2$.
Let  $z_1'=x'+2a_{23}y'$, $z_2'=x'+a_{14}y'$ and
\begin{gather}
\mathfrak{g}^{0,-}_\pm(\tbx,\tbxd,\eta)=-\frac{1}{2\gamma_{23}(\eta)}
\Phi^0(x,y,0,-i\beta_{23}^\pm(\eta))\Phi^{0,*}(x',y',0,-i\beta_{23}^\pm(\eta))\,,
\\  \label{eq:def-g1pm}
\mathfrak{g}^{1,-}_\pm(\tbx,\tbxd,\eta)=
-\frac{1}{2\gamma_{14}(\eta)}\Phi^1(x,y,0,-i\beta_{14}^\pm(\eta))
\Phi^{1,*}(x',y',0,-i\beta_{14}^\pm(\eta))\,,
\end{gather}
and for $k=0$, $1$,
\begin{align}
\label{eq:def-fg}
& \mathfrak{g}^{k,-}(\tbx,\tbxd,\eta)=\left\{
  \begin{aligned}
& \mathfrak{g}^{k,-}_-(\tbx,\tbxd,\eta)\quad\text{if $z_{k+1}>z_{k+1}'$,}
\\ &
\mathfrak{g}^{k,-}_+(\tbx,\tbxd,\eta)\quad\text{if $z_{k+1}<z_{k+1}'$,}
  \end{aligned}\right.
\\ &
\label{eq:def-fg*}
\mathfrak{g}^{k,+}(\tbx,\tbxd,\eta)=\left\{
  \begin{aligned}
& \mathfrak{g}^{k,+}_+(\tbx,\tbxd,\eta)\quad\text{if $z_{k+1}>z_{k+1}'$,}
\\ & 
\mathfrak{g}^{k,+}_-(\tbx,\tbxd,\eta)\quad\text{if $z_{k+1}<z_{k+1}'$,}
  \end{aligned}\right.
\\ & \label{eq:def-fg*2}
\mathfrak{g}^{k,+}_\pm(\tbx,\tbxd,\eta)
=\mathfrak{g}^{k,-}_\pm(\tbxd,\tbx,-\eta)\,,
\end{align}
\begin{gather*}
\mathcal{G}^{k,-}(\tbx,\tbxd,\eta_I)=\frac{1}{2\pi}
\int_{-\infty}^\infty \mathfrak{g}^{k,-}(\tbx,\tbxd,\eta)\,d\eta_R\,,
\\
\mathcal{G}^{k,+}(\tbx,\tbxd,\eta_I)=
\frac{1}{2\pi}\int_{-\infty}^\infty\mathfrak{g}^{k,+}(\tbx,\tbxd,\eta)\,d\eta_R\,.
\end{gather*}
\begin{lemma}
  \label{lem:L-Green}
Let $k=0$, $1$ and $\eta_I\in\R$. Then
\begin{equation*}
L_k\mathcal{G}^{k,-}(\cdot,\tbxd,\eta_I)=\delta(\cdot-\tbxd)\,,\quad
L_k^*\mathcal{G}^{k,+}(\cdot,\tbxd,\eta_I)=\delta(\cdot-\tbxd)\,.  
\end{equation*}
\end{lemma}
\begin{proof}
  Let $\theta_j=\k_jx+\k_j^2y-\k_j^3t$, $\theta_j'=\k_jx'+\k_j^2y'-\k_j^3t'$,
  $\theta_{j,0}=\theta_j|_{t=0}$ and $\theta'_{j,0}=\theta_j'|_{t=0}$.
By \eqref{eq:btitj}
\begin{gather*}
  \mathfrak{g}^{0,-}(\tbx,\tbxd,\eta)
  =-\frac{e^{\frac12\sum_{j=2,3}(\theta_{j,0}-\theta'_{j,0})}}
{2\gamma_{23}(\eta)}e^{-\gamma_{23}(\eta)|z_1-z_1'|}e^{i\eta(y-y')}\,,
\\
|\mathfrak{g}^{0,-}(\tbx,\tbxd,\eta)|\lesssim
|\eta|^{-1/2}e^{\frac12\sum_{j=2,3}(\theta_{j,0}-\theta_{j',0})}
e^{-\eta_I(y-y')-|\eta|^{1/2}|z_1-z_1'|/\sqrt{2}}
\quad\text{as $\eta_R\to\pm\infty$,}  
\end{gather*}
and $\mathfrak{g}^{0,-}(\tbx,\tbxd,\eta)$ is integrable with respect to $\eta_R$.
Moreover, 
\begin{gather}
\label{eq:dg0-jump}
\left\{
  \begin{aligned}
& \mathfrak{g}^{0,-}(z_1'-2a_{23}y+0,y,\tbx',\eta)
=\mathfrak{g}^{0,-}(z_1'-2a_{23}y-0,y,\tbx',\eta)\,,
\\ &
\pd_x\mathfrak{g}^{0,-}(z_1'-2a_{23}y+0,y,\tbx',\eta)
-\pd_x\mathfrak{g}^{0,-}(z_1'-2a_{23}y-0,y,\tbx',\eta)=e^{(i\eta-\k_2\k_3)(y-y')}\,.
  \end{aligned}\right.
\end{gather}
Since $L_0\Phi^0=0$, 
it follows from \eqref{eq:dg0-jump} that for $\varphi\in C_0^\infty$,
\begin{align*}
& \left\la L_0\mathcal{G}^{0,-}(\cdot,\tbxd,\eta_I),\varphi(\cdot)\right\ra
 =\frac{1}{2\pi}
\int_{\R^3}\mathfrak{g}^{0,-}(\tbx,\tbxd,\eta)L_0^*\varphi(\tbx)\,d\tbx d\eta
\\ =&  \frac{1}{2\pi}\int_{\R^2}
\left[\pd_x\mathfrak{g}^{0,-}(\tbx,\tbx',\eta)\varphi(\tbx)-
\mathfrak{g}^{0,-}(\tbx,\tbx',\eta)(\pd_x+2a_{23})\varphi(\tbx)\right]_{z_1=z_1'-0}^{z_1=z_1'+0}\,
       dyd\eta   \\=&
\frac{1}{2\pi} \int_\R \varphi(x'+2a_{23}(y'-y),y)e^{(i\eta-\k_2\k_3)(y-y')}\,dyd\eta
\\=& \varphi(\tbxd)\,.
\end{align*}
Thus we have
$L_0\mathcal{G}^{0,-}(\cdot,\tbxd,\eta_I)=\delta(\cdot-\tbxd)$.
\par
Using a formula
\begin{gather}
\label{eq:phi-phi*1}
\Phi^1(\bx,-i\beta)\Phi^{1,*}(\bxd,-i\beta)=\Phi^0(\bx,-i\beta)\Phi^{0,*}(\bxd,-i\beta)
\left(1+\sum_{j=2,3}e^{\theta_j'-\theta_j}
\frac{\Phi^1_j(\bx)\Phi^{1,*}_j(\bxd)}{\beta-\k_j}\right)\,.
\end{gather}
we have
\begin{equation}
  \label{eq:frakg1}
\left\{
  \begin{aligned}
& \mathfrak{g}^{1,-}(\tbx,\tbxd,\eta)= \mathfrak{g}^{1,-}_1(\tbx,\tbxd,\eta)
+\mathfrak{g}^{1,-}_2(\tbx,\tbxd,\eta)\,,
\\ &   
\mathfrak{g}^{1,-}_1(\tbx,\tbxd,\eta)
= -\frac{1}{2\gamma_{14}(\eta)}
e^{\frac{1}{2}\sum_{j=1,4}(\theta_{j,0}-\theta_{j,0}')
-\gamma_{14}(\eta)|z_2-z_2'|+i\eta(y-y')}\,,
\\ &
\mathfrak{g}^{1,-}_2(\tbx,\tbxd,\eta)= 
\mathfrak{g}^{1,-}_1(\tbx,\tbxd,\eta)
\sum_{j=2,3}
\left(e^{\theta_{j,0}'-\theta_{j,0}}
\frac{\Phi^1_j(\bx)\Phi^{1,*}_j(\bxd)}{\beta_{14}^\mp(\eta)-\k_j}\right)
\quad\text{if  $\pm(z_2-z_2')>0$.}
\end{aligned}\right.
\end{equation}
Since $L_1\Phi^1(\bx,-i\beta)=0$ and 
\begin{gather*}
\mathfrak{g}^{1,-}_1(z_2'-2a_{14}y+0,y,\tbx',\eta)
=\mathfrak{g}^{1,-}_1(z_2'-2a_{14}y-0,y,\tbx',\eta)\,,
\\
\pd_x\mathfrak{g}^{1,-}_1(z_2'-2a_{14}y+0,y,\tbx',\eta)
-\pd_x\mathfrak{g}^{1,-}_1(z_2'-2a_{14}y-0,y,\tbx',\eta)=e^{(i\eta-\k_1\k_4)(y-y')}\,,
\end{gather*}
we have for $\varphi\in C_0^\infty$,
\begin{align*}
&\int_{\R^3}\mathfrak{g}^{1,-}(\tbx,\tbxd,\eta)L_1^*\varphi(\tbx)\,d\tbx d\eta
\\ =& \int_{\R^2} \left[\pd_x\mathfrak{g}^{1,-}(\tbx,\tbx',\eta)\varphi(\tbx)-
\mathfrak{g}^{1,-}(\tbx,\tbx',\eta)(\pd_x+2a_{14})\varphi(\tbx)
\right]_{z_2=z_2'-0}^{z_2=z_2'+0}\,dyd\eta
\\=&  2\pi\varphi(\tbxd)
+
\int_{\R^2} \left[\pd_x\mathfrak{g}^{1,-}_2(\tbx,\tbx',\eta)\varphi(\tbx)-
\mathfrak{g}^{1,-}_2(\tbx,\tbx',\eta)(\pd_x+2a_{14})\varphi(\tbx)
\right]_{z_2=z_2'-0}^{z_2=z_2'+0}\,dyd\eta\,.
\end{align*}
Using the fact that $\sum_{j=2,3}\Phi^1_j(\bx)\Phi^{1,*}_j(\bxd)=0$ and
\begin{equation*}
\int_{-\infty}^\infty\frac{e^{i(y-y')\eta}}{i\eta+a}\,d\eta
=\pi\{1+\sgn(y-y')\}e^{-a(y-y')}\quad\text{for an $a>0$,}  
\end{equation*}
we have
\begin{align*}
& \int_{\R^2} \left[\mathfrak{g}^{1,-}_2(\tbx,\tbx',\eta)
(\pd_x+2a_{14})\varphi(\tbx)\right]_{z_2=z_2'-0}^{z_2=z_2'+0} \,dyd\eta
\\=&
\sum_{j=2,3}\int_{\R^2}
(\pd_x+2a_{14})\varphi(\tbx)\Phi^1_j(\tbx,0)\bigr|_{x=z_2'-2a_{14}y}
\Phi^{1,*}_j(\tbxd,0)
\frac{e^{\{i\eta+(\k_4-\k_j)(\k_j-\k_1)\}(y-y')}}
{(\k_4-\k_j)(\k_j-\k_1)+i\eta}\,dyd\eta
\\=&
2\pi\sum_{j=2,3}\int_{y'}^\infty
(\pd_x+2a_{14})\varphi(\tbx)\Phi^1_j(\tbx,0)\bigr|_{x=z_2'-2a_{14}y}\Phi^{1,*}_j(\tbxd,0)\,
dy=0\,.
\end{align*}
Since
$$\frac{\pd_x\left(e^{-\theta_{j,0}}\mathfrak{g}^{1,-}_1(\tbx,\tbxd,\eta)\right)}{\beta_{14}^-(\eta)-\k_j}
\biggr|_{z_2=z_2'+0}=
\frac{\pd_x\left(e^{-\theta_{j,0}}\mathfrak{g}^{1,-}_1(\tbx,\tbxd,\eta)\right)}{\beta_{14}^+(\eta)-\k_j}
\biggr|_{z_2=z_2'-0}\quad\text{for $j=2$, $3$,}$$
\begin{align*}
& \int_{\R^2} \left[\pd_x\mathfrak{g}^{1,-}_2(\tbx,\tbx',\eta)
  \varphi(\tbx)\right]_{z_2=z_2'-0}^{z_2=z_2'+0} \,dyd\eta
\\  =&
\sum_{j=2,3}\int_{\R^2}
\varphi(\tbx)\pd_x\Phi^1_j(\tbx,0)\bigr|_{x=z_2'-2a_{14}y}
\Phi^{1,*}_j(\tbxd,0)
\frac{e^{\{i\eta+(\k_4-\k_j)(\k_j-\k_1)\}(y-y')}}
{(\k_4-\k_j)(\k_j-\k_1)+i\eta}\,dyd\eta
\\=&
2\pi\sum_{j=2,3}\int_{y'}^\infty
\varphi(\tbx)\pd_x\Phi^1_j(\tbx,0)\bigr|_{x=z_2'-2a_{14}y}\Phi^{1,*}_j(\tbxd,0)\,
dy=0\,.
\end{align*}
Thus we have
$L_1\mathcal{G}^{1,-}(\cdot,\tbxd,\eta_I)=\delta(\cdot-\tbxd)$.
We can prove 
$L_k^*\mathcal{G}^{k,+}(\cdot,\tbxd,\eta_I)=\delta(\cdot-\tbxd)$
in exactly the same way.
\end{proof}
\bigskip

\subsection{The inverse of Darboux transformations.}
\label{subsec:invMpm2}
Now we introduce Green functions of $\nabla M_\pm(v_k)$.
Let $h_1=\tau_1|_{t=0}$ and $h_2=\tau_2/\tau_1|_{t=0}$. For $k=1$, $2$, let
\begin{align}
\label{eq:deftgk-pm}
\tfg^{k,-}_\pm(\tbx,\tbxd,\eta)=-
\frac{h_k(\tbxd)}{h_k(\tbx)}\mathfrak{g}^{k-1,-}_\pm(\tbx,\tbxd,\eta)\,,
\\
\tfg^{k,+}_\pm(\tbx,\tbxd,\eta)
=\frac{h_k(\tbx)}{h_k(\tbxd)}\mathfrak{g}^{k-1,+}_\pm(\tbx,\tbxd,\eta)\,,
\\
\tfg^{k,-}=\tfg^{k,-}_\mp\,,\quad
\tfg^{k,+}=\tfg^{k,+}_\pm\quad\text{if $\pm(z_k-z_k')>0$,}
\end{align}
and
\begin{gather*}
\widetilde{\mathcal{G}}^{k,-}(\tbx,\tbxd,\eta_I)
=\frac{1}{2\pi}\int_\R \tfg^{k,-}(\tbx,\tbxd,\eta)\,d\eta_R
=-\frac{h_k(\tbxd)}{h_k(\tbx)}
\mathcal{G}^{k-1,-}(\tbx,\tbxd,\eta_I)\,,
\\
\widetilde{\mathcal{G}}^{k,+}(\tbx,\tbxd,\eta_I)
=\frac{1}{2\pi}\int_\R \tfg^{k,+}(\tbx,\tbxd,\eta)\,d\eta_R
=\frac{h_k(\tbx)}{h_k(\tbxd)}\mathcal{G}^{k-1,+}(\tbx,\tbxd,\eta_I)\,.
\end{gather*}
\begin{remark}
By \eqref{eq:pd12-22} and \eqref{eq:Phi11*-1},
\begin{align}
\label{eq:tg1--}
& \pd_x\tfg^{1,-}_-(\tbx,\tbxd,\eta)=\frac{1}{2(\k_3-\k_2)}g^1_M(\tbx,\eta)
h_1(\tbxd)\Phi^{0,*}(\bxd,-i\beta_{23}^-(\eta))\,,
\\ & \label{eq:tg2--}
\pd_x\tfg^{2,-}_-(\tbx,\tbxd,\eta)=-\frac{1}{2(\k_4-\k_1)}
g^2_M(\tbx,\eta)h_2(\tbxd)\Phi^{1,*}(\bxd,-i\beta_{14}^-(\eta))\,,
\\ \label{eq:tfg1-+}
&  \pd_{x'}\tfg^{1,+}_+(\tbx,\tbxd,\eta) =
-\frac{h_1(\bx)}{2\gamma_{23}(-\eta)}
\Phi^{0,*}(\tbx,0,-i\beta_{23}^+(-\eta))g^{1,*}_{1,-}(\tbxd,-\bar{\eta})\,,
 \\  & \label{eq:tfg2-+}
\pd_{x'} \tfg^{2,+}_+(\tbx,\tbxd,\eta)=
-\frac{h_2(\tbx)}{2\gamma_{14}(-\eta)}
\Phi^{1,*}(\tbx,0, -i\beta_{14}^+(-\eta))g^{2,*}_{2,-}(\tbxd,-\bar{\eta})\,.
\end{align}  
\end{remark}
It follows from \eqref{eq:Miura-Lax1}--\eqref{eq:Miura-Lax4}
and Lemma~\ref{lem:L-Green} that for $k=1$ and $2$,
\begin{gather}
  \label{eq:22}
\nabla M_-(v_k)\pd_x\widetilde{\mathcal{G}}^{k,-}(\cdot,\tbxd,\eta_I)
=\pd_x\widetilde{\mathcal{G}}^{k,-}(\cdot,\tbxd,\eta_I)\nabla M_-(v_k)
=\delta(\cdot-\tbxd)\,,\\
  \label{eq:22'}
\nabla M_+(v_k)\widetilde{\mathcal{G}}^{k,+}(\cdot,\tbxd,\eta_I)\pd_{x'}
=\widetilde{\mathcal{G}}^{k,+}(\cdot,\tbxd,\eta_I)\pd_{x'}\nabla M_+(v_k)
=\delta(\cdot-\tbxd)\,.  
\end{gather}
\par

Let
$\mathcal{Y}_k=\{v\in \calX_k\mid \pd_xv\in \calX_k\,,\,
\pd_x^{-1}\pd_yv\in\calX_k\}$ be a normed space equipped with a norm
$\|v\|_{\mathcal{Y}_k}=\|v\|_{\calX_k}+\|\pd_xv\|_{\calX_k}
+\|\pd_x^{-1}\pd_yv\|_{\calX_k}$.  Let $\eta_{*,1}$ and $\eta_{*,2}$
be as in \eqref{def:eta*1} and \eqref{def:eta*2}.
Since $\sup_{y\in\R}\|e^{\a z_k}\pd_x^{-1}g^k_{k,+}(\cdot,y,\eta)\|_{H^1(R_x)}<\infty$
for $\eta\in (-\eta_{*,k},\eta_{*,k})$, it follows from \eqref{eq:kernablaM} that
$\nabla M_-(v_k)$ does not have bounded inverse on $\calX_k$.
However, for sufficiently large $\eta_R$,
\begin{equation*}
\int_{\R^2}\tfg^{k,\pm}(\tbx,\tbxd,\eta)f(\tbxd)\,d\tbx  
\end{equation*}
is bounded on $\calX_k$. 
\par
For $\eta_1$ and $\eta_2$ satisfying $\eta_2\ge\eta_1\ge0$ and
$\eta_R$, $\eta_I\in\R$, let $\eta=\eta_R+i\eta_I$ and 
\begin{align}
\label{eq:defTk-}
& T^{k,-}(\eta_1,\eta_2,\eta_I)f=\frac{1}{2\pi}\pd_x
\int_{\eta_1\le|\eta_R|\le \eta_2}\int_{\R^2}
\tfg^{k,-}(\tbx,\tbxd,\eta)f(\tbxd)\,d\tbxd d\eta_R\,,
\\ & \label{eq:defTk+}
T^{k,+}(\eta_1,\eta_2,\eta_I)f=\frac{1}{2\pi}
\int_{\eta_1\le|\eta_R|\le \eta_2}\int_{\R^2}
\tfg^{k,+}(\tbx,\tbxd,\eta)\pd_{x'}f(\tbxd)\,d\tbxd d\eta_R\,.
\end{align}
Let $T^{k,\pm}(\eta_1,\eta_2):=T^{k,\pm}(\eta_1,\eta_2,0)$.
Then we have the following.
\begin{lemma}
  \label{lem:T2pm}
  Let $\a\in(0,2\sqrt{c_{14}})$ and $\eta_0>\eta_{*,2}$.
  If $\eta_2\ge\eta_1\ge\eta_0$, then
\begin{gather}
  \label{eq:T2-bound}
\|T^{2,-}(\eta_1,\eta_2)f\|_{\calX_2}
+\|\pd_x^{-1}\pd_yT^{2,-}(\eta_1,\eta_2)f\|_{\calX_2}  \le C\|f\|_{\calX_2}\,,
\\
  \label{eq:T2+bound}
\|T^{2,+}(\eta_1,\eta_2)f\|_{\calX_2}
+\|\pd_x^{-1}\pd_yT^{2,+}(\eta_1,\eta_2)f\|_{\calX_2} \le C\|f\|_{\calX_2}\,,
\end{gather}
where $C$ is  a positive constant depending only on $\a$ and $\eta_0$.
\end{lemma}
\begin{lemma}
  \label{lem:T2pm'}
  Suppose that $\a\in(0,2\sqrt{c_{14}})$ and that
$\Re\gamma_{14}(\eta_0+i\eta_I)>\sqrt{c_{14}}+\a$ for an $\eta_0\ge0$.
If $\eta_2\ge\eta_1\ge\eta_0$, then
\begin{equation*}
\|e^{\eta_Iy}T^{2,\pm}(\eta_1,\eta_2,\eta_I)f\|_{\calX_2}
\le C\|e^{\eta_Iy}f\|_{\calX_2}\,,
\end{equation*}
where $C$ is  a positive constant depending only on $\a$ and $\eta_0$.
\end{lemma}

\begin{proof}[Proof of Lemmas~\ref{lem:T2pm} and \ref{lem:T2pm'}]
Let $\eta=\eta_R+i\eta_I$ and
\begin{gather*}
\tfg^{2,-}_0(\tbx,\tbxd,\eta)=
-\frac{e^{\theta_{1,0}'}+e^{\theta_{4,0}'}}{e^{\theta_{1,0}}+e^{\theta_{4,0}}}
\mathfrak{g}^{1,-}_1(\tbx,\tbxd,\eta)\,,
\\  
\mathcal{T}_0^{2,-}(\eta_1,\eta_2,\eta_I)f=\frac{1}{2\pi}
\int_{\eta_1\le|\eta_R|\le \eta_2}\int_{\R^2}
\tfg^{2,-}_0(\tbx,\tbxd,\eta)f(\tbxd)\,d\tbxd d\eta_R\,,
\\
\mathcal{T}_j^{2,-}(\eta_1,\eta_2,\eta_I)f=\frac{1}{2\pi}
\int_{\eta_1\le|\eta_R|\le \eta_2}\int_{\R^2}
\frac{\tfg^{2,-}_0(\tbx,\tbxd,\eta)}{\sum_\pm\chi_\pm(z_2-z_2')
  \beta_{14}^\mp(\eta)-\k_{j+1}}f(\tbxd)\,d\tbxd d\eta_R\,,
\end{gather*}
where $\chi_+(x)$ and $\chi_-(x)$ are indicator functions of
$[0,\infty)$ and $(-\infty,0]$, respectively.
Note that $\Re\beta_{14}^+(\eta)>\k_4+\a$ and $\Re\beta_{14}^-(\beta)<
\k_1-\a$ by the assumption. Let
$$
\mathcal{W}_0(\tbx)=\frac{e^{\theta_{1,0}}+e^{\theta_{4,0}}}{h_2(\tbx)}\,,
\quad \mathcal{W}_1(\tbx)=(\k_2-\k_3)\frac{e^{\theta_{3,0}}}{h_1}\,,
\quad \mathcal{W}_2(\tbx)=\frac{e^{\theta_{2,0}}}{h_1}\,.
$$
Then it follow from  \eqref{eq:frakg1}, \eqref{eq:deftgk-pm} and \eqref{eq:defTk-} that
\begin{equation}
  \label{eq:T2-decomp}
  T^{2,-}=\pd_x\mathcal{W}_0\left\{\mathcal{T}_0^{2,-}+\mathcal{W}_1\mathcal{T}_1^{2,-}\mathcal{W}_2
    -\mathcal{W}_2\mathcal{T}_2^{2,-}\mathcal{W}_1\right\}\mathcal{W}_0^{-1}\,.
\end{equation}
By \eqref{eq:btitj},
\begin{equation}
  \label{eq:T2decomp2}
\left\{
  \begin{aligned}
& \tfg^{2,-}_0(\tbx,\tbxd,\eta)=e^{i\eta(y-y')}k_2(z_2,z_2',\eta)\,,
\\ &
k_2(z_2,z_2',\eta)=\frac{1}{2\gamma_{14}(\eta)}
e^{-\gamma_{14}(\eta)|z_2-z_2'|}
\sech\sqrt{c_{14}}z_2\cosh\sqrt{c_{14}}z_2'\,,
  \end{aligned}
\right.
\end{equation}
\begin{gather*}
\end{gather*}
\begin{align*}
\mathcal{T}^{2,-}_0(\eta_1,\eta_2,\eta_I)f(\tbx)
=\frac{1}{\sqrt{2\pi}}\int_{\eta_1\le|\eta_R|\le\eta_2}e^{iy\eta}\left(
  \int_\R k_2(z_2,z_2',\eta)\mathcal{F}_{y'}f^\#(z_2',\eta)\,dz_2'\right)
  \,d\eta_R\,,
\end{align*}
where $f^\#(z_2,y)=f(x,y)$ and
\begin{equation}
  \label{eq:T2-decomp1}
  \begin{split}
\pd_x\mathcal{T}_j^{2,-}(\eta_1,\eta_2,\eta_I)f=& \frac{1}{2\pi}
\int_{\eta_1\le|\eta_R|\le \eta_2}\int_{\R^2}
\frac{\pd_x\tfg^{2,-}_0(\tbx,\tbxd,\eta)}{\sum_\pm\chi_\pm(z_2-z_2')
\beta_{14}^\mp(\eta)-\k_{j+1}}f(\tbxd)\,d\tbxd d\eta_R
\\ &
-\frac{1}{2\pi}\int_{\eta_1\le|\eta_R|\le \eta_2}\int_\R
\frac{f^\#(z_2,y')}{i\eta+(\k_4-\k_{j+1})(\k_{j+1}-\k_1)}e^{i\eta(y-y')}\,dy'
d\eta_R\,.
  \end{split}
\end{equation}
Since
$\cosh \sqrt{c_{14}}z_2'\sech \sqrt{c_{14}}z_2\le 2e^{\sqrt{c_{14}}|z_2-z_2'|}$,
\begin{equation}
\label{eq:k20-est}
|\pd_x^jk_2(z_2,z_2',\eta)| \lesssim
\la\eta\ra^{(j-1)/2}e^{(-\Re\gamma_{14}(\eta)+\sqrt{c_{14}})|z_2-z_2'|}
\quad\text{for $j\ge0$.}
\end{equation}
\par
Let $\mathcal{T}_j^{2,-}(\eta_1,\eta_2)=\mathcal{T}_j^{2,-}(\eta_1,\eta_2,0)$
for $j=0$, $1$ and $2$. If $\pm\eta\ge \eta_0$, then 
$\Re\gamma_{14}(\eta)\ge \Re\gamma_{14}(\eta_0)>\sqrt{c_{14}}+\a$  and 
for $j=0$ and $1$,
$$\|e^{\a(z_2-z_2')}\pd_x^jk_2(z_2,z_2',\eta)\|_{L^\infty_{z_2}L^1_{z_2'}\cap L^\infty_{z_2'}L^1_{z_2}}
\lesssim  \la\eta\ra^{-(2-j)/2}\|f\|_{\calX_2}\,.$$
Let $I_{\eta_1,\eta_2}=[-\eta_2,-\eta_1]\cup[\eta_1,\eta_2]$.
By the Plancherel theorem,
\begin{align*}
  & \|\mathcal{T}^{2,-}_0(\eta_1,\eta_2)f\|_{\calX_2}
    =
\left\|\int_\R e^{\a z_2}k_2(z_2,z_2',\eta)\mathcal{F}_yf^\#(z_2',\eta)\,dz_2'
\right\|_{L^2(\R_{z_2}\times I_{\eta_1,\eta_2})}
\\ & \le 
\sup_{\eta\in I_{\eta_1,\eta_2}}
\left\|e^{\a(z_2-z_2')}k_2(z_2,z_2',\eta)\right\|_{L^\infty_{z_2}L^1_{z_2'}\cap L^\infty_{z_2'}L^1_{z_2}}
\|e^{\a z_2}\mathcal{F}_yf^\#(z_2,\eta)\|_{L^2(\R_{z_2}\times I_{\eta_1,\eta_2})}
\\ & \lesssim \|f\|_{\calX_2}\,.
\end{align*}
Similarly, we have
\begin{gather*}
\|\pd_x\mathcal{T}^{2,-}_0(\eta_1,\eta_2)f\|_{\calX_2}
+\|\pd_y\mathcal{T}^{2,-}_0(\eta_1,\eta_2)f\|_{\calX_2}
  \lesssim \|f\|_{\calX_2}\,,
\\  
  \sum_{j=2,3} (\|\mathcal{T}_j^{2,-}(\eta_1,\eta_2)f\|_{\calX_2}+
  \|\pd_x\mathcal{T}_j^{2,-}(\eta_1,\eta_2)f\|_{\calX_2}
+\|\pd_y\mathcal{T}_j^{2,-}(\eta_1,\eta_2)f\|_{\calX_2})
  \lesssim \|f\|_{\calX_2}\,.
\end{gather*}
Note that  $\beta_{14}^\pm(\eta)$ are bounded away from $\k_2$ and $\k_3$
for $\eta\in\R$.

Combining the above with \eqref{eq:T2-decomp} and the fact that $\mathcal{W}_j$
$(j=0$, $1$, $2$), $\mathcal{W}_0^{-1}$ and their derivatives are uniformly
bounded, we have \eqref{eq:T2-bound}.  We can prove
\eqref{eq:T2+bound} and Lemma~\ref{lem:T2pm'} in the same way.
Thus we complete the proof.
\end{proof}

If $\eta_1$ is small, 
then $T^{k,\pm}(\eta_1,\eta_2)$ are not bounded on $\calX_k$
because $e^{\a(z_k-z_k')}\tfg^{k,\pm}(\tbx,\tbxd,\eta)$ grow exponentially as
$z_k'\to-\infty$.
We will construct $T^{2,-}_{low}(\eta_1,\eta_2)$ such that
$T^{2,-}_{low}(\eta_1,\eta_2)$ is bounded on $\calX_2$ and 
$\nabla M_-(v_2))\{T^{2,-}(\eta_1,\eta_2)-T^{2,-}_{low}(\eta_1,\eta_2)\}=0$
provided $\eta_2$ is sufficiently small.
Let
\begin{gather*}
T^{2,-}_{low}(\eta_1,\eta_2,\eta_I)=T^{2,-}_1(\eta_1,\eta_2,\eta_I)
+T^{2,-}_{2,low}(\eta_1,\eta_2,\eta_I)+T^{2,-}_3(\eta_1,\eta_2,\eta_I)\,,
\\
T^{2,-}_1(\eta_1,\eta_2,\eta_I)f(\tbx)=
\frac{1}{2\pi}\int_{I_{\eta_1,\eta_2}}\int_{-\infty}^\infty
\int_{z_2}^\infty \pd_x\tfg^{2,-}_+(\tbx,\tbxd,\eta)f(\tbxd)\,dz_2'dy'd\eta_R\,,
\\
T^{2,-}_{2,low}(\eta_1,\eta_2,\eta_I)f(\tbx)=
\frac{1}{2\pi}\int_{I_{\eta_1,\eta_2}}\int_{-\infty}^\infty
\int_0^{z_2}\pd_x\tfg^{2,-}_-(\tbx,\tbxd,\eta)f(\tbxd)\,dz_2'dy'd\eta_R\,,
\end{gather*}
\begin{align*}
T^{2,-}_3(\eta_1,\eta_2,\eta_I)f(\tbx)
=  \frac{1}{2\pi}\sum_{j=1,2}&\int_{I_{\eta_1,\eta_2}}\int_\R
\frac{\mathcal{W}_0(\tbx)\mathcal{W}_j(\tbx)\mathcal{W}_{3-j}(\tbxd)}{\mathcal{W}_0(\tbxd)}\Bigr|_{z_2'=z_2}
\\ & \quad\times      \frac{(-1)^je^{i(y-y')\eta}}
      {i\eta+(\k_4-\k_{j+1})(\k_{j+1}-\k_1)}f^\#(z_2,y')\,dy'd\eta_R\,.  
\end{align*}
\begin{lemma}
  \label{lem:20}
Let $\eta=\eta_R+i\eta_I$ with $(\eta_R,\eta_I)\in\R^2$ and $\eta_0\ge0$.
Suppose that $f\in C_0^\infty(\R^2)$. Then
\begin{gather}
  \label{eq:20-}
 \nabla M_-(v_k)T^{2,-}_{low}(0,\eta_0,\eta_I)f
 +\nabla M_-(v_k)T^{2,-}(\eta_0,\infty,\eta_I)f=f\,,\\
 \label{eq:20}
 \nabla M_+(v_k)T^{k,+}(0,\infty,\eta_I)f=f
 \quad\text{for $k=1$ and $2$.} 
\end{gather}
\end{lemma}
\begin{proof}[Proof of Lemma~\ref{lem:20}]
We have \eqref{eq:20} from \eqref{eq:22'}. 
By \eqref{eq:tg2--}
\begin{align*}
&  T^{2,-}(\eta_1,\eta_2,\eta_I)f(\tbx)
                 - T^{2,-}_{low}(\eta_1,\eta_2,\eta_I)f(\tbx)
\\  =&
\frac{1}{2\pi}\pd_x\int_{I_{\eta_1,\eta_2}}\int_{-\infty}^\infty
\int^0_{-\infty}\tfg^{2,-}_-(\tbx,\tbxd,\eta)f(\tbxd)\,dz_2'dy'd\eta_R
  \\=&
-\frac{1}{4(\k_4-\k_1)\pi}
\int_{I_{\eta_1,\eta_2}} d\eta_R g^2_M(\tbx,\eta)
\int_{-\infty}^\infty dy' \int^0_{-\infty}dz_2'
h_2(\tbxd)\Phi^{1,*}(\bxd,-i\beta_{14}^-(\eta))f(\tbxd)\,.
\end{align*}
Since $\nabla M_-(v_2)g^2_M(\cdot,\eta)=0$ by
Lemma~\ref{lem:eigenfunctions-mKP}, we have \eqref{eq:20-}.  Thus we
complete the proof.
\end{proof}

Next, we will restrict the domain of $T^{2,+}$ so that
$T^{2,+}$ is bounded on $\calX_2$ provided $\eta_2<\eta_{*,2}$.
By integration by parts, we have for $f\in C_0^\infty(\R^2)$,
$$T^{2,+}(\eta_1,\eta_2,\eta_I)=\sum_{1\le j\le 3}
T^{2,+}_j(\eta_1,\eta_2,\eta_I)\,,$$
where 
\begin{align*}
& T^{2,+}_1(\eta_1,\eta_2,\eta_I)f(\tbx)  
= -\frac{1}{2\pi}\int_{I_{\eta_1,\eta_2}}\left(
\int_{-\infty}^\infty\int^{z_2}_{-\infty}\pd_{x'}\tfg^{2,+}_+(\tbx,\tbxd,\eta)
f(\tbxd)\,dz_2'dy'\right)\,d\eta_R\,,
  \\ &
T^{2,+}_2(\eta_1,\eta_2,\eta_I)f(\tbx)  
= -\frac{1}{2\pi}\int_{I_{\eta_1,\eta_2}}\left(
\int_{-\infty}^\infty\int_{z_2}^\infty\pd_{x'}\tfg^{2,+}_-(\tbx,\tbxd,\eta)
f(\tbxd)\,dz_2'dy'\right)\,d\eta_R\,,
\end{align*}
\begin{align*}
T^{2,+}_3(\eta_1,\eta_2,\eta_I)f(\tbx)
  = \frac{1}{2\pi}\sum_{j=1,2}
  & \int_{I_{\eta_1,\eta_2}}\int_\R
\frac{\mathcal{W}_0(\tbxd)\mathcal{W}_{3-j}(\tbx)\mathcal{W}_j(\tbxd)}{\mathcal{W}_0(\tbx)}\Bigr|_{z_2'=z_2}
\\ & \quad \times     \frac{(-1)^{j-1}e^{i(y-y')\eta}}
      {i\eta-(\k_4-\k_{j+1})(\k_{j+1}-\k_1)}f^\#(z_2,y')\,dy'd\eta_R\,.
\end{align*}
Let
$T^{2,+}_{low}(\eta_1,\eta_2,\eta_I)=
T^{2,+}_{1,low}(\eta_1,\eta_2,\eta_I)+\sum_{j=2,3}T^{2,+}_j(\eta_1,\eta_2,\eta_I)$
and
\begin{align*}
& T^{2,+}_{1,low}(\eta_1,\eta_2,\eta_I)f(\tbx)  
\\=& \chi_+(z_2)\frac{1}{2\pi}\int_{I_{\eta_1,\eta_2}}\left(
\int_{-\infty}^\infty\int^\infty_{z_2}\pd_{x'}\tfg^{2,+}_+(\tbx,\tbxd,\eta)
f(\tbxd)\,dz_2'dy'\right)\,d\eta_R
\\ & -\chi_-(z_2)\frac{1}{2\pi}\int_{I_{\eta_1,\eta_2}}\left(
\int_{-\infty}^\infty\int_{-\infty}^{z_2}\pd_{x'}\tfg^{2,+}_+(\tbx,\tbxd,\eta)
f(\tbxd)\,dz_2'dy'\right)\,d\eta_R\,.
\end{align*}
Then
\begin{align*}
& T^{2,+}_{low}(\eta_1,\eta_2,\eta_I)f(\tbx)
    -T^{2,+}(\eta_1,\eta_2,\eta_I)f(\tbx)  
\\ =&
\frac{1}{2\pi}\chi_+(z_2)\int_{I_{\eta_1,\eta_2}}\int_{\R^2}
\pd_{x'}\tfg^{2,+}_+(\tbx,\tbxd,\eta)f(\tbxd)\,dz_2'dy'd\eta_R\,,
\end{align*}
and it follows from \eqref{eq:tfg2-+} that 
$T^{2,+}(\eta_1,\eta_2,\eta_I)f
=T^{2,+}_{low}(\eta_1,\eta_2,\eta_I)f$ for $f\in C_0^\infty(\R^2)$ satisfying
\begin{equation}
  \label{eq:secular}
  \la f, g^{2,*}_{2,-}(\cdot,\eta)\ra=0\quad\text{for
$\eta=\eta_R+i\eta_I$ with $\eta_R\in I_{\eta_1,\eta_2}$.}
\end{equation}
\par

Let $T^{2,\pm}_{low}(\eta_1,\eta_2):=T^{2,\pm}_{low}(\eta_1,\eta_2,0)$.
Then we have the following.
\begin{lemma}
  \label{lem:T2lowpm}
  Let $\a\in(0,2\sqrt{c_{14}})$, $\eta_{*,2}$ be as \eqref{def:eta*2} and
  $0\le\eta_1\le\eta_2\le\eta_0<\eta_{*,2}$.
  \begin{enumerate}
  \item There exists a positive constant $C$ depending only on $\a$ and $\eta_0$
such that \allowbreak
$\|T^{2,\pm}_{low}(\eta_1,\eta_2)f\|_{\calX_2} \le C\|f\|_{\calX_2}$ for $f\in \calX_2$.
\item
 Suppose \eqref{eq:secular}. Then 
 there exists a positive constant $C$ depending only on $\a$ and $\eta_0$
 such that $\|T^{2,+}(\eta_1,\eta_2)f\|_{\calX_2} \le C\|f\|_{\calX_2}$ for $f\in \calX_2$.
  \end{enumerate}
\end{lemma}
\begin{lemma}
  \label{lem:T2lowpm'}
  Assume that $0<\a<\min\{\sqrt{c_{14}},\k_4-\k_3,\k_2-\k_1\}$ and that
  $0<|\eta_I|<c_{14}$.
  \begin{enumerate}
  \item Suppose that $\gamma_{14}(i\eta_I)>\sqrt{c_{14}}-\a$ and that
    $\eta_0$ is a nonnegative number satisfying
   $\Re\gamma_{14}(\eta_0+i\eta_I)<\sqrt{c_{14}}+\a$.
 If $0\le\eta_1\le\eta_2\le \eta_0$, then   
\begin{equation*}
  \|e^{\eta_Iy}T^{2,-}_{low}(\eta_1,\eta_2,\eta_I)f\|_{\calX_2}
  \le C\|e^{\eta_Iy}f\|_{\calX_2}\,,
\end{equation*}
where $C$ is a positive constant depending only on $\a$ and $\eta_0$.
\item
Suppose that $\gamma_{14}(-i\eta_I)>\sqrt{c_{14}}-\a$ and that
$\eta_0$ is a nonnegative number satisfying
   $\Re\gamma_{14}(\eta_0-i\eta_I)<\sqrt{c_{14}}+\a$.
 If $0\le\eta_1\le\eta_2\le \eta_0$, then   
\begin{equation*}
  \|e^{\eta_Iy}T^{2,+}_{low}(\eta_1,\eta_2,\eta_I)f\|_{\calX_2}
  \le C\|e^{\eta_Iy}f\|_{\calX_2}\,,
\end{equation*}
where $C$ is a positive constant depending only on $\a$ and $\eta_0$. 
  \end{enumerate}
\end{lemma}

\begin{proof}[Proof of Lemmas~\ref{lem:T2lowpm} and \ref{lem:T2lowpm'}]
We can estimate $T^{2,-}_1$, $T^{2,-}_3$, $T^{2,+}_2$ and $T^{2,+}_3$ in
exactly the same way as Lemma~\ref{lem:T2pm}.
\par

  By \eqref{eq:frakg1}, \eqref{eq:deftgk-pm} and \eqref{eq:T2decomp2},
\begin{equation}
  \label{eq:tfg2-}
  \begin{split}
    \tfg^{2,-}_\pm(\tbx,\tbxd,\eta)
    =& \frac{e^{\pm\gamma_{14}(\eta)(z_2-z_2')+i(y-y')\eta}}{2\gamma_{14}(\eta)}
    \frac{\cosh\sqrt{c_{14}}z_2'}{\cosh\sqrt{c_{14}}z_2}
\\   &  \quad\times \frac{\mathcal{W}_0(\tbx)}{\mathcal{W}_0(\tbxd)}    
\left(1-\sum_{j=1,2}(-1)^j\frac{\mathcal{W}_j(\tbx)\mathcal{W}_{3-j}(\tbxd)}
{\beta_{14}^\pm(\eta)-\k_{j+1}}\right)\,.
  \end{split}
\end{equation}
We have 
\begin{equation}
  \label{eq:26}
\cosh b \sech a\le
2e^{-|a-b|}\quad\text{if  $\pm a\ge \pm b\ge0$.}
\end{equation}
By \eqref{eq:tfg2-} and \eqref{eq:26},
\begin{align*}
  \|T^{2,-}_{2,low}(\eta_1,\eta_2)f(\tbx)\|_{\calX_2}
  \lesssim &
  \|e^{\{\a-\gamma_{14}(\eta)\}z_2}\sech\sqrt{c_{14}}z_2\|_{L^1(\R)}\|f\|_{\calX_2}
  \\ \lesssim &
 \frac{1}{\a+\sqrt{c_{14}}-\Re\gamma_{14}(\eta_0)}\|f\|_{\calX_2}\,.
\end{align*}
Since $\tfg^{2,+}_\pm(\tbx,\tbxd,\eta)=\tfg^{2,-}_\pm(\tbx,\tbxd,-\eta)$,
\begin{align*}
\|\chi_+(z_2)T^{2,+}_{1,low}(\eta_1,\eta_2)f\|_{\calX_2}
\lesssim & 
\|e^{(\gamma_{14}(-\eta)-\sqrt{c_{14}}-\a)z_2}\|_{L^1(0,\infty)}\|f\|_{\calX_2}
  \\ \lesssim &
\frac{1}{\sqrt{c_{14}}+\a-\Re\gamma_{14}(\eta_0)}\|f\|_{\calX_2}\,,
\end{align*}
\begin{align*}
 \|\chi_-(z_2)T^{2,+}_{1,low}(\eta_1,\eta_2)f\|_{\calX_2}
\lesssim &
\|e^{-(\gamma_{14}(-\eta)+\sqrt{c_{14}}-\a)z_2}\|_{L^1(0,\infty)}\|f\|_{\calX_2}
  \\ \lesssim &
\frac{1}{\sqrt{c_{14}}+\Re\gamma_{14}(\eta_0)-\a}\|f\|_{\calX_2}\,.
\end{align*}
\par
If $f\in\calX_2$ satisfies \eqref{eq:secular}, then
$T^{2,+}_{1,low}(\eta_1,\eta_2)f=T^{2,+}_1(\eta_1,\eta_2)f$ and
$ \|T^{2,+}_1(\eta_1,\eta_2)f\|_{\calX_2}\lesssim   \|f\|_{\calX_2}$.
Thus we prove  Lemma~\ref{lem:T2lowpm}.
\par
We can prove Lemma~\ref{lem:T2lowpm'} in the same way.
Note that $\gamma_{14}(\pm i\eta_I)>\sqrt{c_{14}}-\a$ and that
$$\beta_{14}^+(\pm i\eta_I)>\k_4-\a>\k_3\,,\quad
\beta_{14}^-(\pm i\eta_I)<\k_1+\a<\k_2\,.$$
\end{proof}

Since $\|(T^{2,\pm})^*\|_{B(\calX_2^*)}=\|T^{2,\pm}\|_{B(\calX_2)}$ and
$\|(T^{2,\pm}_{low})^*\|_{B(\calX_2^*)}=
\|T^{2,\pm}_{low}\|_{B(\calX_2)}$, we have the following.
\begin{lemma}
  \label{lem:T2-*}
  Let $\a\in(0,2\sqrt{c_{14}})$ and $\eta_{*,2}$ be as \eqref{def:eta*2}.
  \begin{enumerate}
  \item   If $\eta_2\ge\eta_1\ge\eta_0>\eta_{*,2}$, then 
\begin{equation*}
  \|T^{2,\pm}(\eta_1,\eta_2)^*f\|_{\calX_2^*} \le C\|f\|_{\calX_2^*}\,,
\end{equation*}
where $C$ is a positive constant depending only on $\a$ and $\eta_0$.
\item
  If $0\le\eta_1\le\eta_2\le\eta_0<\eta_{*,2}$, then
\begin{equation*}
  \|T^{2,\pm}_{low}(\eta_1,\eta_2)^*f\|_{\calX_2^*} \le C\|f\|_{\calX_2^*}\,,
\end{equation*}
where $C$ is a positive constant depending only on $\a$ and $\eta_0$.
  \end{enumerate}
\end{lemma}

\par
In view of the proof of Lemma~\ref{lem:20}, the range of
$T^{2,-}_2(\eta_1,\eta_2,\eta_I)-T^{2,-}_{2,low}(\eta_1,\eta_2,\eta_I)$
is spanned by $\{g^2_M(\cdot,\eta)\}_{\eta_1\le|\eta_R|\le \eta_2}$.
Thus for $f\in C_0^\infty(\R^2)$ satisfying
  \begin{equation}
    \label{eq:ortht2*}
 \int_{\R^2}f(\tbx)\overline{g^2_M(\tbx,\eta)}\,dxdy=0
\quad\text{for $\eta=\eta_R+i\eta_I$ and $\eta_R\in I_{\eta_1,\eta_2}$,}    
  \end{equation}
we have
$$T^{2,-,*}_2(\eta_1,\eta_2,\eta_If=T^{2,-,*}_{2,low}(\eta_1,\eta_2,\eta_I)f\,.$$
\par
Roughly speaking, $\nabla M_-(v_2)^*$ has a bounded inverse
on $(\ker\nabla M_-(v_2))^\perp$ for low frequencies in $y$.
\begin{lemma}
  \label{lem:T2-*low}
  Let $\a\in(0,2\sqrt{c_{14}})$, $\eta_{*,2}$ be as \eqref{def:eta*2}.
If $0\le \eta_1\le \eta_2\le \eta_0<\eta_{*,2}$
and $f\in\calX_2^*$ satisfies \eqref{eq:ortht2*} with $\eta_I=0$,
\begin{equation*}
 \|T^{2,-}(\eta_1,\eta_2)^*f\|_{\calX_2^*} \le C\|f\|_{\calX_2^*}\,,
\end{equation*}
where $C$ is a positive constant depending only on $\a$ and $\eta_0$.
\end{lemma}

Since
\begin{align*}
  &  T^{2,+}_{low}(\eta_1,\eta_2,\eta_I)^*f(\tbx)
    -T^{2,+}(\eta_1,\eta_2,\eta_I)^*f(\tbx)
\\=& \frac{1}{2\pi}\int_{I_{\eta_1,\eta_2}}\left(
\int_{-\infty}^\infty\int_0^\infty\pd_x\tfg^{2,+}_+(\tbxd,\tbx,-\bar{\eta})
f(\tbxd)\,dz_2'dy'\right)\,d\eta_R\,,
\end{align*}
it follows from Lemma~\ref{lem:g-tg} and \eqref{eq:tfg2-+} that
$$\nabla M_+(v_2)^*T^{2,+}_{low}(\eta_1,\eta_2,\eta_I)^*
=\nabla M_+(v_2)^*T^{2,+}(\eta_1,\eta_2,\eta_I)^*\,.$$
We can prove the following in exactly the same way as Lemma~\ref{lem:20}.
\begin{lemma}
  \label{lem:21}
Suppose that $\eta_0\ge0$, $\eta_I\in\R$ and that $f\in C_0^\infty(\R^2)$.
Then \begin{gather*}
 \nabla M_-(v_2)^*T^{2,-,*}(0,\infty,\eta_I)f=f\,,\\
 \nabla M_+(v_2)^*T^{2,+}_{low}(0,\eta_0,\eta_I)^*f
+\nabla M_+(v_2)^*T^{2,+}(\eta_0,\infty,\eta_I)^*f=f\,.
\end{gather*}
\end{lemma}

\bigskip

\section{Linear stability of $2$-line solitons of P-type}
\label{sec:LS-P}
\subsection{Linear stability of the null solution}
In this subsection, we will investigate linear stability of the null solution
in $\calX_1$ and $\calX_2$.
\par
By the Plancherel theorem,
\begin{equation}
  \label{eq:Plancherel-X}
\|f\|_{\calX_k}=\|\hat{f}(\xi+i\a,\eta+i\a a_i)\|_{L^2(\R^2_{\xi,\eta})}  
\quad\text{for $k=1$, $2$ and $f\in C_0(\R^2)$,}
\end{equation}
where $a_1=2a_{23}$ and $a_2=2a_{14}$.
If $m(\xi,\eta)$ is bounded and analytic on
$$\{(\xi,\eta)\in\C^2\mid 0\le \Im \xi\le\a\,,\; 0\le (\sgn a_k)\Im \eta
\le \a|a_k|\}\,,$$
 a multiplier operator $m(D)$ is bounded on $\calX_k$
and satisfies
\begin{equation}
  \label{eq:m(D)-bound}
\|m(D)\|_{B(\calX_k)}=\sup_{(\xi,\eta)\in\R^2}|m(\xi+i\a,\eta+i\a a_k)|<\infty\,.
\end{equation}
\par

First, we investigate the spectrum of $\mL_0$ in $\calX_1$ and in $\calX_2$.
\begin{lemma}
\label{lem:specL0}
Let $c_1=c_{23}$, $c_2=c_{14}$ and $\a>0$.
Suppose that $\mL_0$ is a closed operator on $\calX_k$. Then 
$\sigma(\mL_0)\subset \{\lambda\in\C\mid \Re\lambda\le -\a(c_k-\a^2/4)\}$.
\end{lemma}
\begin{corollary}
  \label{cor:specL0*}
Let $c_1=c_{23}$, $c_2=c_{14}$ and $\a>0$
Suppose that $\mL_0$ is a closed operator on $\calX_k^*$. Then 
$\sigma(\mL_0)\subset \{\lambda\in\C\mid \Re\lambda\ge \a(c_k-\a^2/4)\}$. 
\end{corollary}
\begin{proof}
  Let $v(x,y)=u(-x,-y)$. Then $v\in \calX_k^*$ if and only if $u\in \calX_k$.
  Since $\mL_0$ is antisymmetric with respect to the change of variables
  $(x,y)\mapsto (-x,-y)$, Corollary~\ref{cor:specL0*} follows immediately
  from Lemma~\ref{lem:specL0}.
\end{proof}
\begin{proof}[Proof of Lemma~\ref{lem:specL0}]
Since
$$
\widehat{\mL_0f}(\xi,\eta)=ip(\xi,\eta)\hat{f}(\xi,\eta)\,,\quad
p(\xi,\eta)=\frac{1}{4}\left(\xi^3-3\frac{\eta^2}{\xi}
  +4b_1\xi+4b_2\eta\right)\,,$$
it follows from \eqref{eq:Plancherel-X} and \eqref{eq:m(D)-bound} that
\begin{align*}
\|(\lambda-\mL_0)^{-1}f\|_{\calX_k}
=& \left\|\frac{\hat{f}(\xi+i\a,\eta+i\a a_k)}
{\lambda-ip(\xi+i\a,\eta+i\a a_k}\right\|_{L^2(\R^2_{\xi,\eta})}
\\ \le &
\sup_{(\xi,\eta)\in\R^2}
\frac{1}{|\Re \lambda+\Im p(\xi+i\a,\eta+i\a a_k)|}\|f\|_{\calX_k}\,.
\end{align*}
By \eqref{eq:b1-b2P}, we have
$4c_k=4(b_1+a_kb_2)-3a_k^2$ for $k=1$, $2$ and
for $\xi$, $\eta$,
\begin{equation}
\label{eq:Imp}
\begin{split}
\Im p(\xi+i\a,\eta+ia_k\a)  =&\frac{\a}{4}
\left(3\xi^2-\a^2+3\frac{(\eta-a_k\xi)^2}{\xi^2+\a^2}+4c_k\right)
\\ \ge & \a(c_k-\a^2/4)\,.
\end{split}
\end{equation}
Thus we prove Lemma~\ref{lem:specL0}.
\end{proof}
\begin{remark}
  \label{rem:PlX}
  In view of the proof of Lemma~\ref{lem:specL0},
$$\sigma(\mL_0)=  
\overline{\{ip(\xi+i\a,\eta+i\a a_k)\mid (\xi,\eta)\in\R^2\}}\,.$$
Moreover, it follows from
\eqref{eq:Plancherel-X} and \eqref{eq:Imp} that
for $\a>0$, $k=1$, $2$ and $t\ge0$,
\begin{align*}
 \|e^{t\mL_0}f\|_{\calX_k}
=& \|e^{-t\Im p(\xi+i\a,\eta+i\a a_k)}\hat{f}(\xi+i\a,\eta+i\a a_k)\|_{L^2(\R^2)}
\le  e^{-\a(c_k-\a^2/4)t}\|f\|_{\calX_k}\,.
\end{align*}  
\end{remark}
\bigskip
  
\subsection{Spectral stability of $[2,3]$-soliton}
\label{subsec:LS-[2,3]P}
Darboux transformations was used to prove
linear stability of $1$-line solitons in an exponentially weighted space.
See \cite[Proposition~3.2]{Miz15}.
In this subsection, we will investigate the spectrum of a linearized operator
around a $1$-line soliton by using Darboux transformations.

First, we consider the case where  weight functions
increase in directions of the motion of $1$-line solitons.
 In the moving coordinate $(t,z_1,y-b_2t)$, \eqref{eq:linequ1P} can be read as
\begin{gather*}
\pd_tu= \mL_1u\,,\quad \mL_1=\mL_0-\frac{3}{2}\pd_{z_1}(u_1\cdot)\,,
\\
\mL_0=
\frac{1}{4}\left\{-\pd_{z_1}(\pd_{z_1}^2-4c_{23})+4(b_2-3a_{23})\pd_y
-3\pd_{z_1}^{-1}\pd_y^2\right\}\,.
\end{gather*}

\begin{proposition}
  \label{lem:wmL1-decay}
\begin{enumerate}
\item Let $\a\in(0,\sqrt{c_{23}})$, $\eta_{*,1}$ be as \eqref{def:eta*1} and
  let $\mL_0$ and $\mL_1$ be closed operators on $\calX_1$. Then
  $\sigma(\mL_1)=\sigma(\mL_0)\cup\{\lambda_{1,\pm}(\eta)\mid
  \eta\in[-\eta_{*,1},\eta_{*,1}]\}$. Moreover,
 $\lambda_{1,\pm}(\eta)-\mL_1$ is invertible on
 $\left(I-P_1(\eta_0)\right)\calX_1$ if $\eta_0\in(0,\eta_{*,1})$ and
 $\eta\in(-\eta_0,\eta_0)$.

\item
Let $\a\in(0,\sqrt{c_{14}})$, $\eta_{*,2}$ be as \eqref{def:eta*2} and
let $\mL_0$ and $\wmL_1$ be closed operators on $\calX_2$. Then
  $\sigma(\wmL_1)=\sigma(\mL_0)\cup\{\lambda_{2,\pm}(\eta)\mid
  \eta\in[-\eta_{*,2},\eta_{*,2}]\}$. Moreover,
 $\lambda_{2,\pm}(\eta)-\mL_1$ is invertible on
 $\left(I-\widetilde{P}_1(\eta_0)\right)\calX_2$ if $\eta_0\in(0,\eta_{*,2})$
 and $\eta\in(-\eta_0,\eta_0)$.
\end{enumerate}
\end{proposition}

Linearizing \eqref{eq:mkp} around
$v_1=a_{23}+\psi$ with $\psi=\sqrt{c_{23}}\tanh\sqrt{c_{23}}z_1$,
we have
\begin{gather*}
  \pd_tv=\mL_{M1}v\,,\\
  \mL_{M1}=\mL_0-\frac{3}{4}\left\{\pd_{z_1}(u_1\cdot)
    +u_1\pd_{z_1}^{-1}\pd_y\right\}\,.
\end{gather*}
Let $\mL_0(\eta)$, $\mL_1(\eta)$, $\mL_{M1}(\eta)$ and $\mM_\pm(\eta)$
be restrictions of $\mL_0$, $\mL_1$, $\mL_{M1}$ and $\nabla M_\pm(v_1)$
on $e^{iy\eta}L^2(\R;e^{2\a z_1}dz_1)$, respectively. More precisely,
\begin{gather*}
\mL_0(\eta)=e^{-iy\eta} \mL_0e^{iy\eta}\,,\quad
\mL_1(\eta)=e^{-iy\eta} \mL_1e^{iy\eta}\,,\quad
\mL_{M1}(\eta)=e^{-iy\eta}\mL_{M1}e^{iy\eta}\,,
\\
\mL_1(\eta)=-\frac{1}{4}\pd_{z_1}(\pd_{z_1}^2-4c_{23}+6u_1)
+i(b_2-3a_{23})\eta+\frac{3}{4}\eta^2\pd_{z_1}^{-1}\,,
\\
\mM_\pm(\eta)=e^{-iy\eta}\nabla M_\pm(v_1)e^{iy\eta}
=\pm\pd_{z_1}+i\eta\pd_{z_1}^{-1}-2\psi\,.
\end{gather*}
Note that $u_1$ and $\psi$ depend only on $z_1$.
By \eqref{eq:bH},
\begin{equation}
  \label{eq:1}
\mL_1(\eta)\mM_+(\eta)=\mM_+(\eta)\mL_{M1}(\eta)\,,\quad
\mL_0(\eta)\mM_-(\eta)=\mM_-(\eta)\mL_{M1}(\eta)\,.
\end{equation}
Using \eqref{eq:Miura-Lax3} and \eqref{eq:Miura-Lax4}, we have
\begin{equation}
\label{eq:mM+}
\begin{split}
\mM_+(\eta)=& \pd_{z_1}^{-1}\cosh(\sqrt{c_{23}}z_1)
\left(\pd_{z_1}^2-\gamma_{23}(-\eta)^2\right)\sech(\sqrt{c_{23}}z_1)
\\=& \cosh(\sqrt{c_{23}}z_1)
\left(\pd_{z_1}^2+u_1-\gamma_{23}(-\eta)^2\right)
\sech(\sqrt{c_{23}}z_1)\pd_{z_1}^{-1}\,,  
\end{split}
\end{equation}
\begin{equation}
 \label{eq:mM-}
 \begin{split}
\mM_-(\eta)=& \sech(\sqrt{c_{23}}z_1)(\gamma_{23}(\eta)^2-\pd_{z_1}^2)
\cosh(\sqrt{c_{23}}z_1)\pd_{z_1}^{-1}
\\=& \pd_{z_1}^{-1}\sech(\sqrt{c_{23}}z_1)(\gamma_{23}(\eta)^2-\pd_{z_1}^2-u_1)
\cosh(\sqrt{c_{23}}z_1)\,.
 \end{split}
\end{equation}
\par

Let
\begin{gather*}
T^{1,-}(\eta)f(z_1)
=\int_\R k_1(z_1,z_1',\eta)f(z_1')\,dz_1'\,,
\\
T^{1,+}(\eta)f(z_1)
=\int_\R k_1(z_1',z_1,-\eta)f(z_1')\,dz_1'\,,
\end{gather*}
where
\begin{gather*}
k_1(z_1,z_1',\eta)=\frac{1}{2\gamma_{23}(\eta)}
\pd_{z_1}\left(e^{-\gamma_{23}(\eta)|z_1-z_1'|}
\sech\sqrt{c_{23}}z_1\cosh\sqrt{c_{23}}z_1'\right)\,.
\end{gather*}
We remark that 
$\tfg^{1,-}(\tbx,\tbxd,\eta)=e^{i\eta(y-y')}\pd_{z_1}^{-1}k_1(z_1,z_1',\eta)$
and that, roughly speaking,
$T^{1,\pm}(\eta)$ are the inverse of $\mM_\pm(\eta)$.
\begin{lemma}\protect{(\cite[Lemmas~2.8]{Miz15})}
  \label{lem:T1+-old}
  Let $\a\in(0,2\sqrt{c_{23}})$ and $\eta_0>\eta_{*,1}$.
If $\eta\in\R$ satisfies $|\eta|\ge \eta_0$,
\begin{gather*}
  \mM_\pm(\eta)T^{1,\pm}(\eta)=T^{1,\pm}(\eta)\mM_\pm(\eta)=I
  \quad\text{on $L^2(\R;e^{2\a z_1}dz_1)$,} 
\\ \
\sum_{j=0,1,2}\la \eta\ra^{(2-j)/2}
\|\pd_{z_1}^{j-1}T^{1,\pm}(\eta)u\|_{L^2(\R;e^{2\a z_1}dz_1)}
\le  C\|u\|_{L^2(\R;e^{2\a z_1}dz_1)}\,,  
\end{gather*}
where $C$ is a positive constant depending only on $\a$ and $\eta_0$.
\end{lemma}
 Let
 \begin{gather}
\label{eq:defgM(eta)}
g_M(z_1,\eta)=\pd_x^{-1}g(z_1,\eta)
= \sqrt{c_{23}}d_{1,+}(\eta)
\pd_{z_1}\left(e^{-\gamma_{23}(\eta)z_1}\sech\sqrt{c_{23}}z_1\right)\,,
\\ \label{eq:mM+gM*}
g_M^*(z_1,\eta)=\frac{1}{2}\mM_+(\eta)^*g_1^*(z_1,\eta)
=-\frac{1}{2} e^{\gamma_{23}(-\bar{\eta})z_1}\sech\sqrt{c_{23}}z_1\,.
\end{gather}
In view of \eqref{eq:defg1z1'*} and \eqref{eq:defgM(eta)},
we have $T^{1,\pm}(\eta)=T^{1,\pm}_+(\eta)+T^{1,\pm}_-(\eta)$, where
\begin{align}  
(T^{1,+}_+(\eta)f)(z_1)=& \frac{\cosh\sqrt{c_{23}}z_1}{2\gamma_{23}(-\eta)}
\int_{z_1}^\infty\pd_{z_1'}\left(e^{\gamma_{23}(-\eta)(z_1-z_1')}
\sech\sqrt{c_{23}}z_1'\right)f(z_1')\,dz_1'
\notag  \\ =& \label{eq:expT1++low}
\frac{e^{\gamma_{23}(-\eta)z_1}\cosh\sqrt{c_{23}}z}{\gamma_{23}(-\eta)}
\int_{z_1}^\infty f(z_1')\overline{g^*(-z_1',-\eta)}\,dz_1'\,,
\end{align}
\begin{align}
(T^{1,+}_-(\eta)f)(z_1)=& \frac{\cosh\sqrt{c_{23}}z_1}{2\gamma_{23}(-\eta)}
\int^{z_1}_{-\infty}\pd_{z_1'}\left(e^{\gamma_{23}(-\eta)(z_1'-z_1)}
\sech\sqrt{c_{23}}z_1'\right)f(z_1')\,dz_1'
\notag \\
\label{eq:expT1+-low}
=&
\frac{e^{-\gamma_{23}(-\eta)z_1}\cosh\sqrt{c_{23}}z}{\gamma_{23}(-\eta)}
\int^{z_1}_{-\infty}f(z_1')\overline{g^*(z_1',-\eta)}\,dz_1'\,,
\end{align}
\begin{align}
\notag
T^{1,-}_+(\eta)f(z_1)=& \frac{1}{2\gamma_{23}(\eta)}
\int^\infty_{z_1}\pd_{z_1}\left(e^{\gamma_{23}(\eta)(z_1-z_1')}\sech\sqrt{c_{23}}z_1\right)
\cosh\sqrt{c_{23}}z_1'f(z_1')\,dz_1'
  \\=& \label{eq:expT1-+low}
\frac{g_M(-z_1,\eta)}{2\sqrt{c_{23}}}
\int^\infty_{z_1} e^{-\gamma_{23}(\eta)z_1'}\cosh\sqrt{c_{23}}z_1'f(z_1')\,dz_1'\,,
\end{align}
\begin{align}
\notag
 T^{1,-}_-(\eta)f(z_1)=&
\frac{1}{2\gamma_{23}(\eta)}
\int_{-\infty}^{z_1}\pd_{z_1}\left(e^{-\gamma_{23}(\eta)(z_1-z_1')}
\sech\sqrt{c_{23}}z_1\right)\cosh\sqrt{c_{23}}z_1'f(z_1')\,dz_1'
  \\ =& \label{eq:expT1--low}
\frac{g_M(z_1,\eta)}{2\sqrt{c_{23}}}
\int_{-\infty}^{z_1} e^{\gamma_{23}(\eta)z_1'}\cosh\sqrt{c_{23}}z_1'f(z_1')\,dz_1'\,.
\end{align}

For small $\eta$, we replace $T^{1,+}_-(\eta)$ and $T^{1,-}_-(\eta)$ by
\begin{multline*}
(T^{1,+}_{-,low})(\eta)f(z_1)=
\frac{\cosh\sqrt{c_{23}}z_1}{2\gamma_{23}(-\eta)}\chi_-(z_1)
\int^{z_1}_{-\infty}\pd_{z_1'}\left(e^{\gamma_{23}(-\eta)(z_1'-z_1)}\sech\sqrt{c_{23}z_1'}
\right)f(z_1')\,dz_1'
\\
-\frac{\cosh\sqrt{c_{23}}z_1}{2\gamma_{23}(-\eta)}\chi_+(z_1)
\int_{z_1}^\infty\pd_{z_1'}\left(e^{\gamma_{23}(-\eta)(z_1'-z_1)}\sech\sqrt{c_{23}z_1'}
\right)f(z_1')\,dz_1'\,,
\end{multline*}
\begin{align}
(T^{1,-}_{-,low})(\eta)f(z_1)=&
\frac{1}{2\gamma_{23}(\eta)}
\int_0^{z_1}\pd_{z_1}\left(e^{-\gamma_{23}(\eta)(z_1-z_1')}
  \sech\sqrt{c_{23}}z_1\right)\cosh\sqrt{c_{23}}z_1'f(z_1')\,dz_1'
\notag  \\ =& \label{eq:expT1--low'}
\frac{1}{2\sqrt{c_{23}}}g_M(z_1,\eta)
\int_0^{z_1} e^{\gamma_{23}(\eta)z_1'}\cosh\sqrt{c_{23}}z_1'f(z_1')\,dz_1'\,,
\end{align}
and let $T^{1,\pm}_{low}(\eta)=T^{1,\pm}_+(\eta)+T^{1,\pm}_{-,low}(\eta)$. Then we have the following.
\begin{lemma}\protect{(\cite[Lemmas~2.6 and 2.7]{Miz15})}
  \label{lem:T1+-'old}
Let $\a\in(0,2\sqrt{c_{23}})$ and $\eta\in(-\eta_{*,1},\eta_{*,1})$.
\begin{enumerate}
\item 
$\ker(\mM_+(\eta))=\{0\}$ and $\operatorname{Range}(\mM_+(\eta))={}^\perp\spann\{g^*(\cdot,-\eta)\}$.
If $|\eta|\le \eta_0<\eta_{*,1}$ and $f\in L^2(\R;e^{2\a z_1}dz_1)$ satisfies
\begin{equation}
  \label{eq:orthT1+}
  \int_\R f(z_1)\overline{g^*(z_1,-\eta)}\,dz_1=0\,,
\end{equation}
then $v=T^{1,+}_{low}(\eta)f$ is a solution of $\mM_+(\eta)v=f$ satisfying
\begin{equation*}
\|v\|_{H^1(\R;e^{2\a z_1}dz_1)}+\|\pd_{z_1}^{-1}v\|_{L^2(\R;e^{2\a z_1}dz_1)}
\le  C\|f\|_{L^2(\R;e^{2\a z_1}dz_1)}\,,
\end{equation*}
where $C$ is a positive constant depending only on $\eta_0$ and $\a$.
Moreover, if $f\in  L^2(\R;e^{2\a z_1}dz_1)$ satisfies $\int_\R f(z_1)\overline{g^*_k(z_1,\eta)}\,dz_1=0$
for $k=1$ and $2$, then a solution $v\in L^2(\R;e^{2\a z_1}dz_1)$ of $\mM_+(\eta)v=f$ satisfies
$$\int_\R v(z_1)\overline{g^*_M(z_1,\eta)}\,dz_1=0\,.$$
\item
$\ker(\mM_-(\eta))=\spann\{g_M(\cdot,\eta)\}$ and
$\operatorname{Range}(\mM_-(\eta))=L^2(\R;e^{2\a z_1}dz_1)$.
Moreover, if $|\eta|\le\eta_0<\eta_{*,1}$, then
$v=T^{1,-}_{low}(\eta)f$ is a solution of $\mM_-(\eta)v=f$ satisfying
\begin{equation*}
\|v\|_{H^1(\R;e^{2\a z_1}dz_1)}+\|\pd_{z_1}^{-1}v\|_{L^2(\R;e^{2\a z_1}dz_1)}\le
 C\|f\|_{L^2(\R;e^{2\a z_1}dz_1)}\,,
\end{equation*}
where $C$ is a positive constant $C$ depending only on $\eta_0$ and $\a$.
\end{enumerate}
  \end{lemma}

Now we are in position to prove Proposition~\ref{lem:wmL1-decay}.
\begin{proof}[Proof of Proposition~\ref{lem:wmL1-decay}]
  Since $e^{\a z_1}\pd_{z_1}^jg(z_1,\eta)$ are bounded for every
  $j\ge0$ and $\eta\in[-\eta_{*,1},\eta_{*,1}]$, it follows from
  \eqref{eq:evmL1} that $\lambda_{1,\pm}(\eta)\in\sigma(\mL_1)$ for
  $\eta\in[-\eta_{*,1},\eta_{*,1}]$.
  \par
  Let $(\xi,\eta)\in\R^2$, $\varphi=\pd_x(\Phi^1(\tbx,k)\Phi^{1,*}(\tbx,\ell))$ and
  $$k=\frac{1}{2}\left(\xi+i\a-\frac{i\eta-2\a a_{23}}{\xi+i\a}\right)\,,\quad
  \ell=\frac12\left(\xi+i\a+\frac{i\eta-2\a a_{23}}{\xi+i\a}\right)\,.$$
Then $\varphi\simeq e^{-\a z_1+i(x\xi+y\eta)}$ as $z_1\to\pm\infty$ and 
$\mL_1\varphi=ip(\xi+i\a,\eta+2i\a a_{23})\varphi$ follows from Lemma~\ref{lem:prodJdJ}.
Since $\varphi$ and its derivatives are bounded in $e^{-\a z_1}L^\infty(\R^2)$ and
we have $ip(\xi+i\a,\eta+2i\a a_{23})\in \sigma(\mL_1)$ and  $\sigma(\mL_0)\subset \sigma(\mL_1)$ by Lemma~\ref{lem:specL0}.
\par
Next, we will show that $\sigma(\mL_1)\subset \sigma(\mL_0)
\cup\{\lambda_{1,\pm}(\eta)\mid\eta\in[-\eta_{*,1},\eta_{*,1}]\}
=:\tilde{\sigma}(\mL_1)$.
Suppose, on the contrary, that
$\lambda\in \sigma(\mL_1)\setminus \tilde{\sigma}(\mL_1)$.
If $\lambda\in \sigma_p(\mL_1)\cup\sigma_c(\mL_1)$, then there exists $\{u_n\}$ and $\{f_n\}$ such that
$$(\lambda-\mL_1)u_n=f_n\,,\quad \|u_n\|_{\calX_1}=1\,,\quad \lim_{n\to\infty}\|f_n\|_{\calX_1}=0\,.$$
Let $V_0=-\frac{3}{2}\pd_x(u_1\cdot)$ and
$$u_{1,n}=(\lambda-\mL_0)^{-1}V_0u_n\,,\quad u_{2,n}=(\lambda-\mL_0)^{-1}f_n\,,
\quad f_{1,n}=V_0u_{2,n}\,.$$
Then $\mL_1=\mL_0+V_0$, $u_n=u_{1,n}+u_{2,n}$, $\lim_{n\to\infty}\|u_{2,n}\|_{\calX_1}=0$ and
\begin{equation*}
  (\lambda-\mL_1)u_{1,n}=f_{1,n}\,,\quad\lim_{n\to\infty}\|u_{1,n}\|_{\calX_1}=1\,,\quad
  \lim_{n\to\infty}\|f_{1,n}\|_{\calX_1}=0\,.
\end{equation*}
Moreover, by Lemma~\ref{lem:imev-5} and the fact that $u_1=O(e^{-2\sqrt{c_{23}}|z_1|})$,
there exists an $\eps>0$ such that
\begin{equation*}
\sup_{n\ge1}\|e^{\eps|z_1|}u_{1,n}\|_{\calX_1}<\infty\,,\quad
  \lim_{n\to\infty}\|e^{\eps|z_1|}f_{1,n}\|_{\calX_1}=0\,.
\end{equation*}
Let $\eps'\in(0,\eps)$ and let $\eta_1$ and $\eta_2$ be positive
numbers such that $0<\eta_1<\eta_{*,1}<\eta_2$,
\begin{equation}
  \label{eq:epseps'}
\sqrt{c_{23}}+\a-\eps<\Re\gamma_{23}(\eta_1)<\sqrt{c_{23}}+\a-\eps'\,,\quad
\sqrt{c_{23}}+\a+\eps'<\Re\gamma_{23}(\eta_2)<\sqrt{c_{23}}+\a+\eps\,,
\end{equation}
and that the inverse of 
$$
\begin{pmatrix}
 \la g_1(\eta_1),g^*_1(\eta)\ra &  \la g_1(\eta_1),g^*_2(\eta)\ra \\
  \la g_2(\eta_1),g^*_1(\eta)\ra &  \la g_2(\eta_1),g^*_2(\eta)\ra 
\end{pmatrix}=I+O(\eta-\eta_1)$$
is uniformly bounded for $\eta\in[-\eta_2,\eta_2]$.
Let $u_{1,n}(\eta)(z_1)=\mathcal{F}_yu_{1,n}(z_1,\eta)$
and $f_{1,n}(\eta)(z_1)=\mathcal{F}_yf_{1,n}(z_1,\eta)$.
By \eqref{eq:evmL1},
\begin{align*}
  \la u_{1,n}(\eta), g^*_1(\eta)\ra=&
\sum_\pm\frac{1}{2\left(\lambda-\lambda_{1,\pm}(\eta)\right)}
\la f_{1,n}(\eta),g^*_1(\eta)\ra
  \\ & + \frac{\Re\gamma_{23}(\eta)}
{\left(\lambda-\lambda_{1,+}(\eta)\right)
\left(\lambda-\lambda_{1,-}(\eta)\right)}\la f_{1,n}(\eta),g^*_2(\eta)\ra\,,
\end{align*}
\begin{align*}
\la u_{1,n}(\eta), g^*_2(\eta)\ra=&
-\frac{\eta^2\Re\gamma_{23}(\eta)}
{\left(\lambda-\lambda_{1,+}(\eta)\right)
\left(\lambda-\lambda_{1,-}(\eta)\right)}\la f_{1,n}(\eta),g^*_1(\eta)\ra
\\ & + \sum_\pm\frac{1}{2\left(\lambda-\lambda_{1,\pm}(\eta)\right)}
\la f_{1,n}(\eta),g^*_2(\eta)\ra \,,
\end{align*}
and $\sum_{k=1,2}|\la u_{1,n}, g^*_k(\eta)\ra|=o(1)$
in $L^2(-\eta_2,\eta_2)$ as $n\to\infty$.
Let $u_{3,n}(\eta)=u_{1,n}(\eta)$ if $|\eta|\ge \eta_2$ and 
$$
u_{3,n}(\eta)=\left\{
\begin{aligned}
 &  u_{1,n}(\eta)-\sum_{k=1,2}b_{k,n}(\eta)g_k(z_1,\eta)
    \quad\text{if $|\eta|\le\eta_1$,}
\\ &
u_{1,n}(\eta)-\sum_{k=1,2}b_{k,n}(\eta)g_k(z_1,\pm\eta_1)
\quad\text{if $\pm\eta\in(\eta_1,\eta_2)$.}
  \end{aligned}\right.
$$
We choose $b_{1,n}(\eta)$ and $b_{2,n}(\eta)$ such that
\begin{equation}
  \label{eq:42}
  \la u_{3,n}(\eta),g^*_k(\eta)\ra=0\quad\text{for $k=1$, $2$ and
    $\eta\in[-\eta_2,\eta_2]$.}
\end{equation}
Then there exists a $C>0$ such that for $\eta\in[-\eta_2,\eta_2]$ and $n\in\N$,
\begin{equation*}
\sum_{k=1,2}|b_{k,n}(\eta)|
\lesssim \sum_{k=1,2}\left|\la u_{1,n}(\eta),g^*_k(\eta)\ra\right|
\le C \|e^{\eps z_1}f_{1,n}(\eta)\|_{L^2(\R;e^{2\a z_1}dz_1)}\,.
\end{equation*}
By \eqref{eq:42}, $f_{3,n}(\eta)=(\lambda-\mL_1(\eta))u_{3,n}(\eta)$ satisfies
\begin{equation}
  \label{eq:42'}
  \la f_{3,n}(\eta),g^*_k(\eta)\ra=0\quad\text{for $k=1$, $2$ and
    $\eta\in[-\eta_2,\eta_2]$.}  
\end{equation}
Let
$$u_{3,n}(z_1,y)=\frac{1}{\sqrt{2\pi}}\int_\R u_{3,n}(\eta)(z_1)
e^{iy\eta}\,d\eta\,,\quad
f_{3,n}(z_1,y)=\frac{1}{\sqrt{2\pi}}\int_\R
f_{3,n}(\eta)(z_1)e^{iy\eta}\,d\eta\,.$$
Then 
$$\|e^{\eps'|z_1|}\pd_x^i\pd_y^j(u_{1,n}-u_{3,n})\|_{\calX_1}
\lesssim \sum_{k=1,2}\|b_{k,n}\|_{L^2(-\eta_2,\eta_2)}\to0
\quad\text{as $n\to\infty$,}$$
\begin{gather}
\label{eq:u3n}  
  (\lambda-\mL_1)u_{3,n}=f_{3,n}\,,\\
\label{eq:f3nlim}
  \lim_{n\to\infty}\|u_{3,n}\|_{\calX_1}=1\,,
\quad \sup_{n\ge1} \|e^{\eps'|z_1|}u_{3,n}\|_{\calX_1}<\infty\,,
\quad \lim_{n\to\infty}\|e^{\eps'|z_1|}f_{3,n}\|_{\calX_1}=0\,.
\end{gather}
Let $u_{M,n}(\eta)=T^{1,+}(\eta)u_{3,n}(\eta)$ and $f_{M,n}(\eta)=T^{1,+}(\eta)f_{3,n}(\eta)$.
Let  $\eps_1\in(0,\eps')$, $\eta_3\in(\eta_1,\eta_{*,1})$ and
$\eta_4\in(\eta_{*,1},\eta_2)$ be numbers such that
$$\sqrt{c_{23}}+\a-\eps'<\Re\gamma_{23}(\eta_3)<\sqrt{c_{23}}+\a-\eps_1\,,\quad
\sqrt{c_{23}}+\a+\eps_1<\Re\gamma_{23}(\eta_4)<\sqrt{c_{23}}+\a+\eps'\,.$$
By Lemma~\ref{lem:T1+-old} and \eqref{eq:epseps'}, there exists a
$C>0$ such that for $\eta\in\R$ satisfying $|\eta|\ge\eta_3$,
\begin{gather*}
\|e^{(\a-\eps')z_1}u_{M,n}(\eta)\|_{H^1(\R)}
+|\eta|\|e^{(\a-\eps')z_1}\pd_{z_1}^{-1}u_{M,n}(\eta)\|_{L^2(\R)}
\le C\|e^{(\a-\eps')z_1}u_{3,n}(\eta)\|_{L^2(\R)}\,,
\\
\|e^{(\a-\eps')z_1}f_{M,n}(\eta)\|_{H^1(\R)}
+|\eta|\|e^{(\a-\eps')z_1}\pd_{z_1}^{-1}f_{M,n}(\eta)\|_{L^2(\R)}
\le C\|e^{(\a-\eps')z_1}f_{3,n}(\eta)\|_{L^2(\R)}\,,
\end{gather*}
and that for $\eta\in\R$ satisfying $|\eta|\ge \eta_4$,
\begin{gather*}
\|e^{(\a+\eps_1)z_1}u_{M,n}(\eta)\|_{H^1(\R)}
+|\eta|\|e^{(\a+\eps_1)z_1}\pd_{z_1}^{-1}u_{M,n}(\eta)\|_{L^2(\R)}
\le C\|e^{(\a+\eps_1)z_1}u_{3,n}(\eta)\|_{L^2(\R)}\,,
\\
\|e^{(\a+\eps_1)z_1}f_{M,n}(\eta)\|_{H^1(\R)}
+|\eta|\|e^{(\a+\eps_1)z_1}\pd_{z_1}^{-1}f_{M,n}(\eta)\|_{L^2(\R)}
\le C\|e^{(\a+\eps_1)z_1}f_{3,n}(\eta)\|_{L^2(\R)}\,.
\end{gather*}
By \eqref{eq:42} and \eqref{eq:42'},
$u_{M,n}(\eta)=T^{1,+}_{low}(\eta)u_{3,n}(\eta)$
and $f_{M,n}(\eta)=T^{1,+}_{low}(\eta)f_{3,n}(\eta)$
for $\eta\in[-\eta_2,\eta_2]$.
Hence it follows from Lemma~\ref{lem:T1+-'old} and \eqref{eq:epseps'} that
there exists a $C>0$ such that for $\eta\in[-\eta_3,\eta_3]$,
\begin{gather*}
\|e^{(\a-\eps_1)z_1}u_{M,n}(\eta)\|_{H^1(\R)}
\le C\|e^{(\a-\eps_1)z_1}u_{3,n}(\eta)\|_{L^2(\R)}\,,
\\
\|e^{(\a-\eps_1)z_1}f_{M,n}(\eta)\|_{H^1(\R)}
\le C\|e^{(\a-\eps_1)z_1}f_{3,n}(\eta)\|_{L^2(\R)}\,,
\end{gather*}
and that for $\eta\in[-\eta_4,\eta_4]$,
\begin{gather*}
\|e^{(\a+\eps') z_1}u_{M,n}(\eta)\|_{H^1(\R)}
\le C\|e^{(\a+\eps') z_1}u_{3,n}(\eta)\|_{L^2(\R)}\,,
\\
\|e^{(\a+\eps') z_1}f_{M,n}(\eta)\|_{H^1(\R)}
\le C\|e^{(\a+\eps') z_1}f_{3,n}(\eta)\|_{L^2(\R)}\,.
\end{gather*}
Combining the above, we have
\begin{equation}
  \label{eq:uMnest-1}
  \begin{split}
& \sum_{j=0,1}\|e^{\a z_1+\eps_1|z_1|}\pd_{z_1}^ju_{M,n}(\eta)\|_{L^2(\R)}
+\la \eta\ra\|e^{\a z_1+\eps_1|z_1|}\pd_{z_1}^{-1}u_{M,n}(\eta)\|_{L^2(\R)}
\\ \le & C\|e^{\a z_1+\eps'|z_1|}u_{3,n}(\eta)\|_{L^2(\R)}\,, 
  \end{split}
\end{equation}
\begin{equation}
 \label{eq:fMnest-1}
  \begin{split}
& \sum_{j=0,1}\|e^{\a z_1+\eps_1|z_1|}\pd_{z_1}^jf_{M,n}(\eta)\|_{L^2(\R)}
+\la \eta\ra\|e^{\a z_1+\eps_1|z_1|}\pd_{z_1}^{-1}f_{M,n}(\eta)\|_{L^2(\R)}
\\ \le & C\|e^{\a z_1+\eps'|z_1|}f_{3,n}(\eta)\|_{L^2(\R)}\,,
  \end{split}
\end{equation}
where $C$ is a constant independent of $n\in\N$ and $\eta\in\R$.
\par
It follows from Lemmas~\ref{lem:T1+-old} and \ref{lem:T1+-'old} that
\begin{equation*}
\mM_+(\eta)u_{M,n}(\eta)=u_{3,n}(\eta)\,,\quad\mM_+(\eta)f_{M,n}(\eta)=f_{3,n}(\eta)\,.
\end{equation*}
Hence by  \eqref{eq:1},
\begin{align*}
&   \mM_+(\eta)\left\{\left(\lambda-\mL_{M1}(\eta)\right)u_{M,n}(\eta)-f_{M,n}(\eta)\right\}
\\  =& \left(\lambda-\mL_1(\eta)\right)u_{3,n}(\eta)-f_{3,n}(\eta)
= 0\,.
\end{align*}
Since $\ker\mM_+(\eta)=\{0\}$ on $L^2(\R;e^{2\a z_1}dz_1)$ by \eqref{eq:mM+},
\begin{equation}
  \label{eq:40}
\left(\lambda-\mL_{M1}(\eta)\right)u_{M,n}(\eta)=f_{M,n}(\eta)\,.  
\end{equation}
Let $u_n'(\eta)=\mM_-(\eta)u_{M,n}(\eta)$, $f_n'(\eta)=\mM_-(\eta)f_{M,n}(\eta)$ and
\begin{gather*}
u_n'(z_1,y)=\frac{1}{\sqrt{2\pi}}\int_\R u_n'(z_1,\eta)e^{iy\eta}\,d\eta\,,\quad
f_n'(z_1,y)=\frac{1}{\sqrt{2\pi}}\int_\R f_n'(z_1,\eta)e^{iy\eta}\,d\eta\,.
\end{gather*}
We have
$u_n'\in\calX_1$ and $\lim_{n\to\infty}\|f_n'\|_{\calX_1}=0$
from \eqref{eq:f3nlim}--\eqref{eq:fMnest-1}.
By \eqref{eq:1} and \eqref{eq:40}, we have
$\left(\lambda-\mL_0(\eta)\right)u_n'(\eta)-f_n'(\eta)=0$
and $$(\lambda-\mL_0)u_n'=f_n'\,.$$
Since $\lambda\in\rho(\mL_0)$,
\begin{equation}
  \label{eq:41}
  \lim_{n\to\infty}\|u_n'\|_{\calX_1}=0\,.
\end{equation}
Let $\eps_2\in(0,\eps_1)$, $\eta_5\in(\eta_3,\eta_{*,1})$
and $\eta_6\in(\eta_{*,1},\eta_4)$ be numbers such that
$$\sqrt{c_{23}}+\a-\eps_1<\Re\gamma_{23}(\eta_5)<\sqrt{c_{23}}+\a-\eps_2\,,\quad
\sqrt{c_{23}}+\a+\eps_2<\Re\gamma_{23}(\eta_6)<\sqrt{c_{23}}+\a+\eps_1\,.$$
Let $u_{M,n}^*(\eta)=e^{2\a z_1}u_{M,n}(\eta)$ for $\eta$ satisfying
$|\eta|\ge\eta_6$ and
\begin{equation*}
u_{M,n}^*(\eta)=\left\{
  \begin{aligned}
&   e^{2\a z_1}u_{M,n}(\eta)-a_n(\eta)g_M^*(z_1,\eta)
    \quad\text{if $\eta\in[-\eta_5,\eta_5]$,}
\\ &
e^{2\a z_1}u_{M,n}(\eta)-a_n(\eta)g_M^*(z_1,\pm\eta_3)
\quad\text{if $\pm\eta\in(\eta_5,\eta_6)$,}
  \end{aligned}\right.
\end{equation*}
\begin{equation*}
  u_{M,n}^*(z_1,y)=\frac{1}{\sqrt{2\pi}}
  \int_\R u_{M,n}^*(\eta)e^{iy\eta}\,d\eta\,.
\end{equation*}
Here we choose $a_n(\eta)$ such that
\begin{equation}
  \label{eq:uMn*orth-1}
\la g_M(\eta),u_{M,n}^*(\eta)\ra=0\quad\text{for $\eta\in[-\eta_6,\eta_6]$.}
\end{equation}
By \eqref{eq:uMn*orth-1}, we have for $\eta\in[-\eta_6,\eta_6]$,
\begin{gather*}
|a_n(\eta)|\le C\|e^{(\a-\eps_1) z_1}u_{M,n}(\eta)\|_{L^2(\R)}\,,
\end{gather*} 
where $C$ is a constant independent of $n$ and $\eta\in[-\eta_4,\eta_4]$.
Thus we have
\begin{equation*}
  \|e^{\eps_2|z_1|}u_{M,n}^*(\eta)\|_{\calX_1^*}
  \lesssim \|e^{\eps_1|z_1|}u_{M,n}\|_{\calX_1}\,.
\end{equation*}
Let $\eps_3\in(0,\eps_2)$, $\eta_7\in(\eta_5,\eta_{*,1})$
and $\eta_8\in(\eta_{*,1},\eta_6)$ be numbers such that
$$\sqrt{c_{23}}+\a-\eps_2<\Re\gamma_{23}(\eta_7)<\sqrt{c_{23}}+\a-\eps_3\,,\quad
\sqrt{c_{23}}+\a+\eps_3<\Re\gamma_{23}(\eta_8)<\sqrt{c_{23}}+\a+\eps_2\,.$$
Let $w_n(\eta)=T^{1,-}(\eta)^*u_{M,n}^*(\eta)$ and
$$w_n(z_1,y)=\frac{1}{\sqrt{2\pi}}\int_\R u_{M,n}^*(\eta)e^{iy\eta}\,d\eta\,.$$ 
It follows from Lemma~\ref{lem:T1+-old} that there exists a $C>0$ such that
\begin{align*}
  & \|e^{-(\a-\eps_2)z_1}w_n(\eta)\|_{L^2(\R)} \le
    C\|e^{-(\a-\eps_2)z_1}u_{M,n}^*(\eta)\|_{L^2(\R)}
    \quad\text{if $|\eta|\ge \eta_7$,}
\\ & \|e^{-(\a+\eps_3)z_1}w_n(\eta)\|_{L^2(\R)} \le
    C\|e^{-(\a+\eps_3)z_1}u_{M,n}^*(\eta)\|_{L^2(\R)}
    \quad\text{if $|\eta|\ge\eta_8$.}
\end{align*}
In view of \eqref{eq:expT1--low} and \eqref{eq:expT1--low'},
\begin{align*}
  T^{1,-}(\eta)^*f-T^{1,-}_{low}(\eta)^*f=&
\frac{\chi_-(z_1)}{2\sqrt{c_{23}}}e^{\gamma_{23}(-\eta)z_1}\cosh\sqrt{c_{23}}z_1
\la f, g_M(\eta)\ra\,,
\end{align*}
and $T^{1,-}(\eta)^*u_{M,n}^*(\eta)=T^{1,-}(\eta)_{low}^*u_{M,n}^*(\eta)$
for $\eta\in[-\eta_6,\eta_6]$.
Thus by Lemma~\ref{lem:T1+-'old}, there exists a $C>0$ such that
\begin{align*}
  & \|e^{-(\a-\eps_3)z_1}w_n(\eta)\|_{L^2(\R)} \le
    C\|e^{-(\a-\eps_3)z_1}u_{M,n}^*(\eta)\|_{L^2(\R)}
    \quad\text{if $\eta\in[-\eta_7, \eta_7]$,}
\\ & \|e^{-(\a+\eps_2)z_1}w_n(\eta)\|_{L^2(\R)} \le
    C\|e^{-(\a+\eps_2)z_1}u_{M,n}^*(\eta)\|_{L^2(\R)}
    \quad\text{if $\eta\in[-\eta_8,\eta_8]$.}
\end{align*}
Combining the above, we have
\begin{equation}
  \label{eq:44}
  \|w_n\|_{\calX_1^*}\lesssim \|e^{\eps_1|z_1}u_{M,n}^*\|_{\calX_1^*}\,.
\end{equation}
By \eqref{eq:mM+gM*}, \eqref{eq:42} and \eqref{eq:42'},
\begin{equation*}
\la u_{M,n}(\eta),g_M^*(\eta)\ra=\la f_{M,n}(\eta),g_M^*(\eta)\ra=0\quad
  \text{for $\eta\in[-\eta_2,\eta_2]$,}
\end{equation*}
and
\begin{align*}
  \|u_{M,n}\|_{\calX_1}^2=& \int_\R \|e^{\a z_1}u_{M,n}(\eta)\|_{L^2(\R)}^2\,d\eta
  \\=& \left(1+O(\eta_6-\eta_5)\right)
       \int_\R \la u_{M,n}(\eta),u_{M,n}^*(\eta)\ra\,d\eta\,.  
\end{align*}
Since $\mM_-(\eta)^*w_n=u_{M,n}^*(\eta)$ and
$\mM_-(\eta)u_{M,n}(\eta)=u_n'(\eta)$, it follows from \eqref{eq:41},
\eqref{eq:44} and the above that as $n\to\infty$,
\begin{align*}
 \|u_{M,n}\|_{\calX_1}^2
  = & \left(1+O(\eta_6-\eta_5)\right)
      \int_\R \la u_n'(\eta),w_n(\eta)\ra\,d\eta
\\  \lesssim &  \|u_n'\|_{\calX_1}\|w_n\|_{\calX_1^*}\to0\,.
\end{align*}
Combining the above with
\eqref{eq:f3nlim}, \eqref{eq:41} and the fact that
$$u_{3,n}(\eta)-u_n'(\eta)=(\mM_+(\eta)-\mM_-(\eta))u_{M,n}(\eta)
=2\pd_xu_{M,n}(\eta)\,,$$
we have as $n\to\infty$,
$$
\|u_{M,n}\|_{\calX_1}+\|(-\pd_x+\pd_x^{-1}\pd_y)u_{M,n}\|_{\calX_1}\to0\,,
\quad \|\pd_xu_{M,n}\|_{\calX_1}\to\frac12\,,$$
which contradicts to Lemma~\ref{cl:m0-bound}.
Thus we prove $\lambda\not\in\sigma_p(\mL_1)\cup\sigma_c(\mL_1)$
if $\lambda\in\rho(\mL_0)$ and $\lambda\ne\lambda_{1,\pm}(\eta)$
for $\eta\in[-\eta_{*,1},\eta_{*,1}]$.
We see that $\lambda\not\in\sigma_r(\mL_1)$ as well because
$\bar{\lambda}\in\sigma_p(\mL_1^*)$ if $\lambda\in\sigma_r(\mL_1)$.
Suppose, on the contrary,
that $\mL_1^*\varphi=\bar{\lambda}\varphi$ for $\varphi\in\calX_1^*$.
Since $\pd_x\mL_1^*=-\mL_1\pd_x$ (formally) and
$u_1$ is invariant under a change of variables $(x,y)\mapsto (-x,-y)$,
$$\mL_1\tilde{\varphi}=\bar{\lambda}\tilde{\varphi}\,,\quad
\tilde{\varphi}(x,y)=\pd_x\varphi(-x,-y)\in\calX_1\,,$$
and $\bar{\lambda}\in\sigma_p(\mL_1)$. Since $\mL_1$ is real valued and
$\bar{\lambda}\in\sigma_p(\mL_1)$, we have $\lambda\in\sigma_p(\mL_1)$,
which is a contradiction.
\par
Finally, we will prove simplicity of continuous eigenvalues
$\{\lambda_{1,\pm}(\eta)\}$. 
Suppose the contrary. Then there exist
$\lambda\in\{\lambda_{1,\pm}(\eta)\mid -\eta_0<\eta<\eta_0\}$,
$u_n\in\calX_1$ and $f_n\in\calX_1$ such that for $n\in\N$,
  $$(\lambda-\mL_1)u_n=f_n\,,\enskip \|u_n\|_{\calX_1}=1\,,
  \enskip \|f_n\|_{\calX_1}\le \frac1n\,,\enskip
P_1(\eta_0)u_n=0\,,\enskip P_1(\eta_0)f_n=0\,.$$
Let $u_{1,n}$ and $u_{2,n}$ be as above.
Since $P_1(\eta_0)u_n=0$, we have for $k=1$, $2$,
\begin{align*}
\lim_{n\to\infty} \la u_{1,n},g^*_k(\eta)\ra
  =& -\lim_{n\to\infty}\la u_{2,n},g^*_k(\eta)\ra
 = 0\quad\text{in $L^2(-\eta_0,\eta_0)$,}
\end{align*}
and the rest of the proof is exactly the same as the case where
$\lambda\ne\lambda_{1,\pm}(\eta)$ for any $\eta\in[-\eta_{*,1},\eta_{*,1}]$.

We can prove the latter part of Proposition~\ref{lem:wmL1-decay}
in exactly the same way.
This completes the proof of Proposition~\ref{lem:wmL1-decay}.
\end{proof}
\par
Next, we will investigate spectral stability of $\mL_1$
in $\calX_2$ whose weight function increases in a direction of motion
of $[1,4]$-soliton. Let
\begin{gather*}
\eta_{*,+}=2(\sqrt{c_{23}}+\a)\sqrt{\a(\a+2a_{14}-2\k_2)}\,,\quad
\eta_{*,-}=2(\sqrt{c_{23}}+\a)\sqrt{\a(\a+2\k_3-2a_{14})}\,,
\\
\eta_{\#,+}=2(\sqrt{c_{23}}-\a)\sqrt{\a(\a+2a_{14}-2\k_3)}\,,\quad
\eta_{\#,-}=2(\sqrt{c_{23}}-\a)\sqrt{\a(\a+2\k_2-2a_{14})}\,.
\end{gather*}

\begin{proposition}
  \label{prop:specmL1}
Assume
\begin{equation}
\label{eq:ass-alpha}
0<\a<\min\left\{\sqrt{c_{23}}\,,\,
  \frac{(b_2-3a_{23})^2-c_{23}}{4|a_{14}-a_{23}|}\,,\,
\frac{c_{23}}{3|a_{14}-a_{23}|}\right\}\,.
\end{equation}
Let $\mL_0$ and $\mL_1$ be  closed operators on $\calX_2$.
\begin{enumerate}
\item If $\k_2<a_{14}<\k_3$, 
$\sigma(\mL_1)=\sigma(\mL_0)\cup\{\lambda_{1,\pm}(\eta)\mid \eta\in[-\eta_{*,\pm},\eta_{*,\pm}]\}$.
\item If $\k_2>a_{14}$ and $\a<2(\k_2-a_{14})$,
  $$\sigma(\mL_1)=\sigma(\mL_0) \cup\{\lambda_{1,-}(\eta)\mid \eta\in[-\eta_{*,-},\eta_{*,-}]\}
   \cup\{\tilde{\lambda}_{1,-}(\eta)\mid \eta\in[-\eta_{\#,-},\eta_{\#,-}]\}\,.$$
\item If $a_{14}>\k_3$ and $\a<2(a_{14}-\k_3)$,
  $$\sigma(\mL_1)=\sigma(\mL_0)\cup\{\lambda_{1,+}(\eta)\mid \eta\in[-\eta_{*,+},\eta_{*,+}]\}
  \cup\{\tilde{\lambda}_{1,+}(\eta)\mid \eta\in[-\eta_{\#,+},\eta_{\#,+}]\} \,.$$
\item If $\k_2=a_{14}$,
  $$\sigma(\mL_1)=\sigma(\mL_0)\cup\{\lambda_{1,\pm}(\eta)\mid \eta\in[-\eta_{*,\pm},\eta_{*,\pm}]\}
\cup\{\tilde{\lambda}_{1,-}(\eta)\mid \eta\in[-\eta_{\#,-},\eta_{\#,-}]\}  \,.$$
\item If $\k_3=a_{14}$,
  $$\sigma(\mL_1)=\sigma(\mL_0)\cup\{\lambda_{1,\pm}(\eta)\mid \eta\in[-\eta_{*,\pm},\eta_{*,\pm}]\}
\cup\{\tilde{\lambda}_{1,+}(\eta)\mid \eta\in[-\eta_{\#,+},\eta_{\#,+}]\}  \,.$$
\end{enumerate}
\end{proposition}

\begin{corollary}
  \label{lem:mL-1inX2}
Assume that \eqref{eq:ass-alpha}.  
Suppose that $\mL_1$ is a closed operator on $\calX_2$.
Then there exists $b\in(0,\a(c_{14}-\a^2/4))$ such that
\begin{equation}
  \label{eq:boundmL1}
\sup_{\Re\lambda\ge -b}\|(\lambda-\mL_1)^{-1}\|_{B(\calX_2)}<\infty\,.  
\end{equation}
\end{corollary}
\par

We can prove that $\sigma(\mL_0)\subset \sigma(\mL_1)$ in exactly the
same way as Proposition~\ref{lem:wmL1-decay}.
\par
Let us investigate the continuous eigenvalues of $\mL_1$ which have to do with
modulations of $[2,3]$-soliton. By \eqref{eq:evmL1},
\begin{gather*}
\mL_1(\eta)g(z_1,\pm\eta)=\lambda_{1,\pm}(\eta)g(z_1,\pm\eta)\,,\quad
\mL_1(\eta)g(-z_1,\pm\eta)=\tilde{\lambda}_{1,\pm}(\eta)g(-z_1,\pm\eta)\,,
\end{gather*}
where $\tilde{\lambda}_{1,\pm}(\eta)=i\eta(b_2-3a_{23}\mp\gamma_{23}(\pm\eta))$.
Let $\eta=\eta_R+i\eta_I$, $\eta_I=2(a_{14}-a_{23})\a$ and $\eta_R\in\R$. Then
$\a z_2+i\eta y=\a z_1+i\eta_R y$ and  it follows from \eqref{eq:defg1z1'}
and \eqref{eq:defg1z1} that
\begin{gather*}
\left\|e^{\a z_2} g^1_{1,\pm}(\tbx,\eta)\right\|_{L^2(\R_{z_1})}
\lesssim 
\left\|e^{\a z_1}g(z_1,\pm\eta)\right\|_{L^2(\R)}\,.
\end{gather*}
We remark that for $\eta_I=2\a(a_{14}-a_{23})$,
$$\calX_2=e^{-y\eta_I}\calX_1=\bigoplus_{\eta=\eta_R+i\eta_I,\,\eta_R\in\R}
e^{iy\eta}L^2(\R;e^{2\a z_1}dz_1)\,.$$
We remark that if $\a-\sqrt{c_{23}}<\Re \gamma_{23}(\pm\eta)<\sqrt{c_{23}}+\a$,
then $g(z_1,\pm\eta)\in L^2(\R;e^{2\a z_1}dz_1)$ and
$\lambda^1_{1,\pm}(\eta)$ are continuous eigenvalues of  $\mL_1$.
\par
Suppose that $0<\a<2(\k_2-a_{14})$. Then for every $\eta_R\in\R$,
\begin{equation}
\Re \gamma_{23}(\eta)\ge \gamma_{23}(i\eta_I)>\sqrt{c_{23}}+\a\,,
\quad g(z_1,\eta)\not\in L^2(\R;e^{2\a z_1}dz_1)\,.
\end{equation}
We have
$\Re\gamma_{23}(\pm\eta_{\#,-}-i\eta_I)=\sqrt{c_{23}}-\a$
and
\begin{equation*}
  \{\tilde{\lambda}_{1,-}(\eta_R+i\eta_I)\mid
  \eta_R\in[-\eta_{\#,-},\eta_{\#,-}]\}  \subset \sigma(\mL_1)
\end{equation*}
because $g(-z_1,-\eta)\in L^2(\R;e^{2\a z_1}dz_1)$ provided
$|\eta_R|<\eta_{\#,-}$.
\par
Suppose that $0<\a<2(a_{14}-\k_3)$. Then for every $\eta_R\in\R$,
\begin{equation}
\Re \gamma_{23}(-\eta)>\sqrt{c_{23}}+\a\,, \quad
\quad g(z_1,-\eta)\not\in L^2(\R;e^{2\a z_1}dz_1)\,.
\end{equation}
We have
$\Re\gamma_{23}(\pm\eta_{\#,+}+i\eta_I)=\sqrt{c_{23}}-\a$ 
and
\begin{equation*}
  \{\tilde{\lambda}_{1,+}(\eta_R+i\eta_I)\mid
  \eta_R\in[-\eta_{\#,+},\eta_{\#,+}]\}\subset\sigma(\mL_1)
\end{equation*}
because $g(-z_1,\eta)\in L^2(\R;e^{2\a z_1}dz_1)$ provided $|\eta_R|<\eta_{\#,+}$.
\par
If $\k_2\le a_{14}$, then
$\Re \gamma_{23}(\pm\eta_{*,+}+i\eta_I)=\sqrt{c_{23}}+\a$ and
$$\{\lambda_{1,+}(\eta_R+i\eta_I)\mid \eta_R\in[-\eta_{*,+},\eta_{*,+}]\}
\subset \sigma(\mL_1)$$
because $g(z_1,\eta)\in L^2(\R;e^{2\a z_1}dz_1)$ for 
$\eta_R\in(-\eta_{*,+},\eta_{*,+})$.
\par
If $\k_3\ge a_{14}$, then
$\Re \gamma_{23}(\pm\eta_{*,-}-i\eta_I)=\sqrt{c_{23}}+\a$
 and
$$\{\lambda_{1,-}(\eta_R+i\eta_I)\mid \eta_R\in[-\eta_{*,-},\eta_{*,-}]\}
\subset \sigma(\mL_1)$$
because $g(z_1,-\eta)\in L^2(\R;e^{2\a z_1}dz_1)$
for $\eta_R\in(-\eta_{*,-},\eta_{*,-})$.
\par

For $\eta=\eta_R+i\eta_I$ with $\eta_R\in\R$ and $\eta_I=2(a_{14}-a_{23})\a$,
we see that $\lambda_{1,\pm}(\eta)$ belong to the left half plane and that
$\lambda_{1,\pm}(\eta)$ are bounded away from the imaginary axis.

\begin{lemma}
  \label{lem:lambda1-etaI}
Assume  \eqref{eq:ass-alpha}.
Let $\eta=\eta_R+i\eta_I$ and $\eta_I=2\a(a_{14}-a_{23})$. Then 
\begin{gather}
\label{eq:l1-monotone}
\eta_R\pd_{\eta_R}\Re \lambda_{1,\pm}(\eta)<0
\quad\text{for any $\eta_R\in\R\setminus\{0\}$,}
\\
\Re\lambda_{1,\pm}(\eta_R+i\eta_I)\le \Re \lambda_{1,\pm}(i\eta_I)
\le -C\a\quad\text{for any $\eta_R\in\R$,}
\end{gather}
where $C$ is a positive constant that does not depend on $\a$.
\end{lemma}

\begin{lemma}
  \label{lem:lambda1-etaI'}
Assume  \eqref{eq:ass-alpha}.
Let $\eta=\eta_R+i\eta_I$ and $\eta_I=2\a(a_{14}-a_{23})$. Then 
\begin{equation}
\label{eq:l1-monotone'}
\eta_R\pd_{\eta_R}\Re \tilde{\lambda}_{1,\pm}(\eta)>0
\quad\text{for any $\eta_R\in\R\setminus\{0\}$.}
\end{equation}
Moreover, we have the following.
\begin{enumerate}
\item If $0<\a<2(a_{14}-\k_3)$ and $\eta_R\in[-\eta_{\#,+},\eta_{\#,+}]$,
$$\Re\tilde{\lambda}_{1,+}(\eta_R+i\eta_I)
\le \Re\tilde{\lambda}_{1,+}(\eta_{\#,+}+i\eta_I)<0\,.$$
\item
 If $0<\a<2(\k_2-a_{14})$ and $\eta_R\in[-\eta_{\#,-},\eta_{\#,-}]$,
$$\Re\tilde{\lambda}_{1,-}(\eta_R+i\eta_I)
\le \Re\tilde{\lambda}_{1,-}(\eta_{\#,-}+i\eta_I)<0\,.$$
\end{enumerate}
\end{lemma}

To prove Lemmas~\ref{lem:lambda1-etaI} and \ref{lem:lambda1-etaI'},
we need to see that in which direction continuous eigenmodes propagate.
\begin{lemma}
  \label{cl:sgnpdlambda1}
  Let $a_{23}$, $a_{14}$ and $b_2$ be as \eqref{eqdef:a,c,omega} and
  \eqref{eq:b1-b2P}. 
  Then 
\begin{gather}
  \label{eq:pdlambda-11}
  b_2>\k_2+2\k_3\quad\text{and}\quad
  b_2>2\k_1+\k_4 \quad\text{if $a_{14}>a_{23}$,}
\\ \label{eq:pdlambda-12}
b_2<2\k_2+\k_3\quad\text{and}\quad b_2<\k_1+2\k_4
\quad\text{if $a_{14}<a_{23}$,}  
\end{gather}
$i^{-1}(a_{14}-a_{23})\pd_\eta\lambda_{1,\pm}(0)>0$. Moreover,
$i^{-1}\pd_\eta\lambda_{2,+}(0)>0$ if $a_{14}>a_{23}$
and $i^{-1}\pd_\eta\lambda_{2,-}(0)<0$ if $a_{14}<a_{23}$.
\end{lemma}
\begin{proof}
By \eqref{eq:b1-b2P'},
\begin{align*}
b_2-3a_{23}-\sqrt{c_{23}}=&
b_2-(\k_2+2\k_3)
\\=& 2(a_{14}-a_{23})+\frac{(\k_2-\k_1)(\k_4-\k_2)}
{2(a_{14}-a_{23})}>0\quad\text{if $a_{14}>a_{23}$,}
\end{align*}
\begin{align*}
b_2-3a_{23}+\sqrt{c_{23}}=&
b_2-(2\k_2+\k_3)
\\=& 2(a_{14}-a_{23})+\frac{(\k_3-\k_1)(\k_4-\k_3)}
{2(a_{14}-a_{23})}<0\quad\text{if $a_{14}<a_{23}$,}
\end{align*}
\begin{align*}
b_2-3a_{14}+\sqrt{c_{14}}=&
b_2-(2\k_1+\k_4)
  \\=& a_{23}-\k_1+\frac{(\k_4-a_{23})(a_{23}-\k_1)-c_{23}}{2(a_{14}-a_{23})}
       >0\quad\text{if $a_{14}>a_{23}$,}
\end{align*}
\begin{align*}
b_2-3a_{14}-\sqrt{c_{14}}=&
b_2-(\k_1+2\k_4)
  \\=& a_{23}-\k_4-\frac{(\k_4-a_{23})(a_{23}-\k_1)-c_{23}}{2(a_{23}-a_{14})}
       <0\quad\text{if $a_{14}<a_{23}$.}
\end{align*}
Note that $\k_1<a_{23}<\k_4$ and that
$(\k_4-a_{23})(a_{23}-\k_1)>(\k_3-a_{23})(a_{23}-\k_2)=c_{23}$.
\par
By \eqref{eq:lambdas},
\begin{equation}
  \label{eq:pdlambdas-0}
    \pd_\eta\lambda_{1,\pm}(0)=i(b_2-3a_{23}\pm\sqrt{c_{23}})\,,\quad
    \pd_\eta\lambda_{2,\pm}(0)=i(b_2-3a_{14}\pm\sqrt{c_{14}})\,.
  \end{equation}
Thus we prove Lemma~\ref{cl:sgnpdlambda1}.
\end{proof}

\begin{proof}[Proof of Lemma~\ref{lem:lambda1-etaI}]
Let $\gamma_{R,\pm}(\eta_R)=\Re\gamma_{23}(\eta_R\pm i\eta_I)$ and
$\gamma_{I,\pm}(\eta_R)=\Im\gamma_{23}(\eta_R\pm i\eta_I)$. Then 
for $\eta_R\in\R$,
\begin{gather}
  \label{eq:g23-reim1}
  \gamma_{R,\pm}(\eta_R)>0\,,\quad \eta_R\gamma_{I,\pm}(\eta_R)\ge0\,,\quad
  \pd_{\eta_R}\gamma_{I,\pm}(\eta_R)\ge0\,,
\\   \label{eq:g23-reim2}
\pd_{\eta_R}\gamma_{R,\pm}(\eta_R)=\frac{\gamma_{I,\pm}(\eta_R)}{2\sqrt{(c_{23}\mp\eta_I)^2+\eta_R^2}}\,,
\quad
\pd_{\eta_R}\gamma_{I,\pm}(\eta_R)=\frac{\gamma_{R,\pm}(\eta_R)}{2\sqrt{(c_{23}\mp\eta_I)^2+\eta_R^2}}
\,.
\end{gather}
By \eqref{eq:ass-alpha}, \eqref{eq:g23-reim1} and \eqref{eq:g23-reim2},
\begin{align*}
\eta_R\pd_{\eta_R}\Re \lambda_{1,\pm}(\eta)=&
-\eta_R\eta_I\pd_{\eta_R}\gamma_{R,\pm}(\pm\eta_R)\mp \eta_R\gamma_{I,\pm}(\pm\eta_R)-\eta_R^2\pd_{\eta_R}\gamma_{I,\pm}(\pm\eta_R)
\\=&
 \mp \eta_R\gamma_{I,\pm}(\pm\eta_R)
\left(1\pm\frac{\eta_I}{2\sqrt{(c_{23}\mp\eta_I)^2+\eta_R^2}}\right)
-\eta_R^2\pd_{\eta_R}\gamma_{I,\pm}(\pm\eta_R)
\\ \le & -\eta_R^2\pd_{\eta_R}\gamma_{I,\pm}(\pm\eta_R)  \le  0\,,
\end{align*}
and
$\Re\lambda_{1,\pm}(\eta)\le  \Re\lambda_{1,\pm}(i\eta_I)
=-\eta_I\{b_2-3a_{23}\pm\gamma_{23}(\pm i\eta_I)\}$.

If  $a_{14}>a_{23}$, then $\eta_I>0$, $b_2-3a_{23}>\sqrt{c_{23}}$
by \eqref{eq:pdlambda-11} and
\begin{align*}
b_2-3a_{23}+\gamma_{23}(i\eta_I) \ge &  b_2-3a_{23}-\gamma_{23}(-i\eta_I)
\\ \ge &
\frac{(b_2-3a_{23})^2-c_{23}-\eta_I}{b_2-3a_{23}+\gamma_{23}(-i\eta_I)}>0
\end{align*}
by \eqref{eq:ass-alpha}.
Thus we have
$\Re\lambda_{1,+}(i\eta_I)\le \Re\lambda_{1,-}(i\eta_I)\lesssim
-\a$.
\par
If $a_{23}>a_{14}$, then $\eta_I<0$, $3a_{23}-b_2>\sqrt{c_{23}}$
by \eqref{eq:pdlambda-12} and
\begin{align*}
b_2-3a_{23}-\gamma_{23}(-i\eta_I) \le &  b_2-3a_{23}+\gamma_{23}(i\eta_I)
\\ \le &
\frac{c_{23}-(3a_{23}-b_2)^2-\eta_I}{3a_{23}-b_2+\gamma_{23}(i\eta_I)}<0\,,
\end{align*}
and
$\Re\lambda_{1,-}(i\eta_I)\le \Re\lambda_{1,+}(i\eta_I)\lesssim
-\a$. 
Thus we complete the proof.
\end{proof}

\begin{proof}[Proof of Lemma~\ref{lem:lambda1-etaI'}]
Since $\pd_{\eta_R}\Re\tilde{\lambda}_{1,\pm}(\eta)
=-\pd_{\eta_R}\Re\lambda_{1,\pm}(\eta)$, we have \eqref{eq:l1-monotone'}.
\par
Since $\gamma_{R,\pm}(\eta_{\#,\pm})=\sqrt{c_{23}}-\a$,
$\gamma_{I,+}(\eta_{\#,+})=\sqrt{\a(\a+2a_{14}-2\k_3)}$,
$\gamma_{I,-}(\eta_{\#,-})=\sqrt{\a(\a+2\k_2-2a_{14})}$ and
\begin{equation*}
\Re \tilde{\lambda}_{1,\pm}(\eta)=
\eta_I(\pm\gamma_{R,\pm}(\pm\eta_R)-b_2+3a_{23})
\pm \eta_R\gamma_{I,\pm}(\pm\eta_R)\,,
\end{equation*}
we have for $\eta=\eta_{\#,+}+i\eta_I$,
\begin{gather*}
\Re \tilde{\lambda}_{1,+}(\eta)=
-\a\{3(\k_3-a_{14})^2+c_{14}\}-6(a_{14}-\k_3)\a^2-2\a^3\,,
\\
\Re \tilde{\lambda}_{1,-}(\eta)=
-\a\{3(\k_2-a_{14})^2+c_{14}\}-6(\k_2-a_{14})\a^2-2\a^3\,.
\end{gather*}
Thus we prove Lemma~\ref{lem:lambda1-etaI'}.
\end{proof}

Next, we will estimate Green functions of Darboux transformations.
Using properties of the Green functions
(Lemmas~\ref{lem:T1+}--\ref{lem:T1-'k3small} below),
we can prove Proposition~\ref{prop:specmL1} in the same way
  as Proposition~\ref{lem:wmL1-decay}.
Let $\eta_I(\a)=2\a(a_{14}-a_{23})$ and
\begin{gather*}
  \eta_{*,\pm}(\a,\beta)=2(\sqrt{c_{23}}+\beta)
  \sqrt{\beta^2+2\sqrt{c_{23}}\beta\pm \eta_I(\a)}\,,
  \\
\eta_{\#,\pm}(\a,\beta)=2(\sqrt{c_{23}}-\beta)
  \sqrt{\beta^2-2\sqrt{c_{23}}\beta\pm \eta_I(\a)}\,.
\end{gather*}
Note that
 $\eta_{*,\pm}(\a,\a)=\eta_{*,\pm}$, $\eta_{\#,\pm}(\a,\a)=\eta_{\#,\pm}$
and that
$$\Re\gamma_{23}(\eta_{*,\pm}(\a,\beta)\pm i\eta_I(\a))=\sqrt{c_{23}}+\beta\,,
\quad
\Re\gamma_{23}(\eta_{\#,\pm}(\a,\beta)\pm i\eta_I(\a))=\sqrt{c_{23}}-\beta\,.$$

If $\pm\eta$ are sufficiently large, then
$\mM_\pm(\eta):H^1(\R;e^{2\a z_1}dz_1)\to L^2(\R;e^{2az_1}dz_1)$
are bijective.
\begin{lemma}
  \label{lem:T1+}
Let $\a$, $\beta\in(0,\sqrt{c_{23}})$ and $\beta$ be sufficiently
close to $\a$.  Let $\eta=\eta_R+i\eta_I$, $\eta_R\in\R$ and
$\eta_I=2\a(a_{14}-a_{23})$.
Suppose that  $c_{23}+\eta_I>0$ and that
$\eta_0=0$ if $0<\a<2(a_{14}-\k_3)$ and $\eta_0>\eta_{*,-}$
otherwise.  Then for every $\eta_R$ satisfying $|\eta_R|\ge \eta_0$,
\begin{equation}
  \label{eq:mM+inv}
\mM_+(\eta)T^{1,+}(\eta)=T^{1,+}(\eta)\mM_+(\eta)=I\quad\text{on $L^2(\R;e^{2\a z_1}dz_1)$.}
\end{equation}
Moreover, if $|\eta_R|\ge\eta_{*,-}(\a,\beta)$, then
\begin{gather}
  \\ \label{eq:T1+est}
\sum_{j=0,1,2}\la \eta\ra^{(2-j)/2}
\|\pd_{z_1}^{j-1}T^{1,+}(\eta)f\|_{L^2(\R;e^{2\beta z_1}dz_1)}
\le  C\|f\|_{L^2(\R;e^{2\beta z_1}dz_1)}\,,  
\end{gather}
where $C$ is a positive constant depending only on $\a$, $\beta$ and $\eta_0$.
\end{lemma}

\begin{lemma}
  \label{lem:T1-}
Let $\a\in(0,\sqrt{c_{23}})$, $\eta=\eta_R+i\eta_I$, $\eta_R\in\R$ and
$\eta_I=2\a(a_{14}-a_{23})$.
Suppose that $c_{23}>\eta_I$ and that
$\eta_0=0$ if $0<\a<2(\k_2-a_{14})$ and  $\eta_0>\eta_{*,+}$ otherwise.
Then for every  $\eta_R$ satisfying $|\eta_R|\ge \eta_0$,
\begin{gather}
  \label{eq:mM-inv}
  \mM_-(\eta)T^{1,-}(\eta)=T^{1,-}(\eta)\mM_-(\eta)=I
  \quad\text{on $L^2(\R;e^{2\a z_1}dz_1)$,} 
\end{gather}
Moreover, if $\eta_0>\eta_{*,+}(\a,\beta)$, then
\begin{gather} \label{eq:T1-est}
\sum_{j=0,1,2}\la \eta\ra^{(2-j)/2}\|\pd_{z_1}^{j-1}T^{1,-}(\eta)f
\|_{L^2(\R;e^{2\beta z_1}dz_1)}
\le C\|f\|_{L^2(\R;e^{2\beta z_1}dz_1)}\,,
\end{gather}
where $C$ is a positive constant depending only on $\a$, $\beta$ and $\eta_0$.
\end{lemma}

\begin{proof}[Proof of Lemmas~\ref{lem:T1+} and \ref{lem:T1-}]
  Equations \eqref{eq:mM+inv} and \eqref{eq:mM-inv} follow immediately from
\eqref{eq:mM+} and \eqref{eq:mM-}.
As \eqref{eq:k20-est}, we have
\begin{equation}
  \label{eq:k1-est}
|k_1(z_1,z_1',\eta)|+\la \eta\ra^{1/2}|\pd_x^{-1}k_1(z_1,z_1',\eta)|
\lesssim e^{(-\Re\gamma_{23}(\eta)+\sqrt{c_{23}})|z_1-z_1'|}\,,
\end{equation}
By the assumption assumption,
$\Re\gamma_{23}(\eta)\ge \Re\gamma_{23}(\eta_0+i\eta_I)>\sqrt{c_{23}}+\beta$
and
\begin{equation}
\label{eq:T1-bound}
\begin{split}
& \la\eta\ra^{1/2}\|T^{1,-}(\eta)f\|_{L^2(\R;e^{2\beta z_1}dz_1)} 
+ \la\eta\ra\|\pd_x^{-1}T^{1,-}(\eta)f\|_{L^2(\R;e^{2\beta z_1}dz_1)} 
\\ & \le 
\sum_{j=0,1}
\la\eta\ra^{(2-j)/2}\left\|e^{\beta(z_1-z_1')}\pd_{z_1}^{j-1}k_1(z_1,z_1',\eta)
\right\|_{L^\infty_{z_1}L^1_{z_1'}\cap L^\infty_{z_1'}L^1_{z_1}}
 \|f\|_{L^2(\R;e^{2\beta z_1}dz_1)} 
 \\ & \lesssim
 \frac{1}{\Re\gamma_{23}(\eta_0+i\eta_I)-\sqrt{c_{23}}-\beta}
\|f\|_{L^2(\R;e^{2\beta z_1}dz_1)}\,,
\end{split}
\end{equation}
and \eqref{eq:T1-est} follows from \eqref{eq:mM-inv} and \eqref{eq:T1-bound}.
We can prove \eqref{eq:T1+est} in exactly the same way.
\end{proof}

Next, we will consider the cases where $\eta_R$ is small and
$g(z_1,\pm\eta)$ are eigenfunctions of $\mL_1(\eta)$.
\begin{lemma}
  \label{lem:T1+'}
Assume that $\a\in(0,\sqrt{c_{23}})$, $c_{23}+\eta_I>0$ and that $\beta$ is
sufficiently close to $\a$.
Let $\eta=\eta_R+i\eta_I$,  $\eta_R\in\R$ and $\eta_I=2\a(a_{14}-a_{23})$.
  \begin{enumerate}
  \item   Suppose that $\k_2< a_{14}\le \k_3$.
If $\a<  2(a_{14}-\k_2)$ and $|\eta_R|<\eta_{*,-}$, then
$\ker(\mM_+(\eta))=\{0\}$,
$\operatorname{Range}(\mM_+(\eta))={}^\perp\spann\{g^*(\cdot,-\eta)\}$
and for every $f\in L^2(\R;e^{2\a z_1}dz_1)$ satisfying \eqref{eq:orthT1+},
$v=T^{1,+}_{low}(\eta)f$ is a solution of $\mM_+(\eta)v=f$.
Moreover, if $|\eta_R|\le \eta_0<\eta_{*,-}(\a,\beta)$, then
\begin{equation}
  \label{eq:T1+est'}
\|v\|_{H^1(\R;e^{2\beta z_1}dz_1)}+\|\pd_{z_1}^{-1}v\|_{L^2(\R;e^{2\beta z_1}dz_1)}
\le  C\|f\|_{L^2(\R;e^{2\beta z_1}dz_1)}\,,
\end{equation}
where $C$ is a positive constant depending only on $\eta_0$, $\a$ and $\beta$.
\item Suppose that $\k_2\ge a_{14}$ and $\eta_{\#,-}<|\eta_R|<\eta_{*,-}$.
Then $\ker(\mM_+(\eta))=\{0\}$,
$\operatorname{Range}(\mM_+(\eta))={}^\perp\spann\{g^*(\cdot,-\eta)\}$
and  for every $f\in L^2(\R;e^{2\a z_1}dz_1)$ satisfying \eqref{eq:orthT1+},
$v=T^{1,+}_{low}(\eta)f$ is a solution of $\mM_+(\eta)v=f$.
Moreover, if $\eta_{\#,-}(\a,\beta)<\eta_1\le|\eta_R|\le
\eta_2<\eta_{*,-}(\a,\beta)$, then $v$ satisfies \eqref{eq:T1+est'}, where
$C$ is a positive constant depending only on $\eta_1$, $\eta_2$, $\a$ and
$\beta$.
\item
  If $f\in  L^2(\R;e^{2\a z_1}dz_1)$ satisfies $\int_\R f(z_1)\overline{g^*_k(z_1,\eta)}\,dz_1=0$
  for $k=1$ and $2$, then a solution $v\in L^2(\R;e^{2\a z_1}dz_1)$ of $\mM_+(\eta)v=f$ satisfies
$$\int_\R v(z_1)\overline{g^*_M(z_1,\eta)}\,dz_1=0\,.$$    
\end{enumerate}
  \end{lemma}

\begin{lemma}
  \label{lem:T1-'}
  Assume that $0<\a<\sqrt{c_{23}}$, $c_{23}>\eta_I$  and that
  $\beta$ is sufficiently close to $\a$.
Let $\eta=\eta_R+i\eta_I$,  $\eta_R\in\R$ and $\eta_I=2\a(a_{14}-a_{23})$.

\begin{enumerate}
\item  Assume that $\k_2\le a_{14}<\k_3$ and that $\a<2(\k_3-a_{14})$.
If $|\eta_R|<\eta_{*,+}$, then
  \begin{equation}
    \label{eq:lemT1-'}
\ker(\mM_-(\eta))=\spann\{g_M(\cdot,\eta)\}\,,\quad
\operatorname{Range}(\mM_-(\eta))=L^2(\R;e^{2\a z_1}dz_1)\,,
\end{equation}
and for $f\in L^2(\R;e^{2\a z_1}dz_1)$,
$v=T^{1,-}_{low}(\eta)f$ is a solution of $\mM_-(\eta)v=f$ in $ L^2(\R;e^{2\a z_1}dz_1)$.
Moreover, if $|\eta_R|\le \eta_0<\eta_{*,+}(\a,\beta)$,
\begin{gather}
  \label{eq:T1-est'}
\|v\|_{H^1(\R;e^{2\beta z_1}dz_1)}+\|\pd_{z_1}^{-1}v\|_{L^2(\R;e^{2\beta z_1}dz_1)}\le
 C\|f\|_{L^2(\R;e^{2\beta z_1}dz_1)}\,,
\end{gather}
where $C$ is a positive constant $C$ depending only on $\eta_0$, $\a$ and
$\beta$.
\item Assume that $\k_3\le a_{14}$. If $\eta_{\#,+}<|\eta_R|<\eta_{*,+}$,
then \eqref{eq:lemT1-'} holds and for $f\in L^2(\R;e^{2\a z_1}dz_1)$,
$v=T^{1,-}_{low}(\eta)f$ is a solution of $\mM_-(\eta)v=f$ in $ L^2(\R;e^{2\a z_1}dz_1)$.  Moreover, if
  $\eta_{\#,+}(\a,\beta)<\eta_1\le|\eta_R|\le\eta_2<\eta_{*,+}(\a,\beta)$,
  then $v$ satisfies \eqref{eq:T1-est'}, where $C$ is a positive
  constant $C$ depending only on $\eta_1$, $\eta_2$, $\a$ and $\beta$.
\end{enumerate}
\end{lemma}
\begin{proof}[Proof of Lemmas~\ref{lem:T1+'} and \ref{lem:T1-'}]
Since
\begin{gather}
\Re\gamma_{23}(\eta)\ge \gamma_{23}(i\eta_I)>\sqrt{c_{23}}-\a
\quad\text{if $0<\a<2(\k_3-a_{14})$,}\\
\label{eq:100}
\Re\gamma_{23}(-\eta)\ge \gamma_{23}(-i\eta_I)>\sqrt{c_{23}}-\a
  \quad\text{if $0<\a<2(a_{14}-\k_2)$,}
\end{gather}
it follows from \eqref{eq:mM+} and \eqref{eq:mM-} that
\begin{gather*}
  \ker(\mM_+(\eta))=\{0\}\,,\quad
\ker(\mM_-(\eta))=\spann\{g_M(\cdot,\eta)\}\,,
\\
\ker(\mM_-(\eta)^*)=\{0\}\,,
\quad \ker(\mM_+(\eta)^*)=\spann\{g^*(\cdot,-\eta)\}\,,
\end{gather*}
\begin{equation}
  \label{eq:2}
\mM_\pm(\eta)T^{1,\pm}(\eta)f=f\quad\text{for $f\in C_0^\infty(\R)$.}
\end{equation}
\par
Next, we will show \eqref{eq:T1+est'} for $f$ satisfying \eqref{eq:orthT1+}.
By \eqref{eq:k1-est} and \eqref{eq:100},
\begin{equation}
\label{eq:T1++est}
\sum_{j=0,1}  \|\pd_{z_1}^{j-1}T^{1,+}_+(\eta)u\|_{L^2(\R;e^{2\beta z_1}dz_1)}
 \lesssim \|f\|_{L^2(\R;e^{2\beta z_1}dz_1)}\,,
\end{equation}
provided $\a-\beta$ is small.
In view of \eqref{eq:expT1+-low}, we have 
$T^{1,+}_-(\eta)f=T^{1,+}_{-,low}(\eta)f$ for $f$ satisfying \eqref{eq:orthT1+}.
Applying \eqref{eq:26} with $a=\sqrt{c_{23}}z_1$ and $b=\sqrt{c_{23}}z_1'$,
we have 
\begin{align*}
\|T^{1,+}_{-,low}f\|_{L^2(\R;e^{2\beta z_1}dz_1)}
 \lesssim &
                \|e^{(\beta-\gamma_{23}(-\eta))z_1-\sqrt{c_{23}}|z_1|}\|_{L^1(\R)}
                \|f\|_{L^2(\R;e^{2\beta z_1}dz_1)}
\\ \lesssim &  \|f\|_{L^2(\R;e^{2\beta z_1}dz_1)}\,.
\end{align*}
Combining the above with \eqref{eq:2} and 
$\pd_{z_1}^{-1}\in B(L^2(\R;e^{2\beta z_1}dz_1)$, we have \eqref{eq:T1+est'}.
By \eqref{eq:mM+gM*},
\begin{align*}
\int_\R \mM_+(\eta)v(z_1)\overline{g_1^*(z_1,\eta)}\,dz_1
=& 2\int_\R v(z_1)\overline{g^*_M(z_1,\eta)}\,dz_1=0\,.
\end{align*}
\par
Next, we will prove Lemma~\ref{lem:T1-'}.
We can prove
$$\|T^{1,-}_+(\eta)u\|_{L^2(\R;e^{2\beta z_1}dz_1)}
 \lesssim \|f\|_{L^2(\R;e^{2\beta z_1}dz_1)}$$ in exactly the same way as \eqref{eq:T1++est}.
Applying \eqref{eq:26} with $a=\sqrt{c_{23}}z_1$ and $b=\sqrt{c_{23}}z_1'$,
we have
\begin{align*}
\|T^{1,-}_{-,low}(\eta)f\|_{L^2(\R;e^{2\beta z_1}dz_1)}\lesssim &
\|e^{(\beta-\Re\gamma_{23}(\eta))z_1}e^{-\sqrt{c_{23}}|z_1|}\|_{L^1(\R)}
\|f\|_{L^2(\R;e^{2\beta z_1}dz_1)}
\\ \lesssim & \|f\|_{L^2(\R;e^{2\beta z_1}dz_1)}\,.
\end{align*}
Since $T^{1,-}(\eta)f-T^{1,-}_{-,low}(\eta)f$ is a constant multiple of $g_M(\cdot,\eta)$
and belongs to $\ker(\mM_-(\eta))$,
$$\mM_-(\eta)T^{1,-}_{-,low}(\eta)=I\quad\text{on $L^2(\R;e^{2\a z_1}dz_1)$.}$$
Combining the above, we have \eqref{eq:T1-est'}.
Thus we complete the proof.
\end{proof}

Secondly, we investigate Green functions of $\mM_+(\eta)$ in the case
where $g(\pm z_1,-\eta)$ are eigenfunctions of $\mL_1(\eta)$.
Let $T^{1,+}_{low\#}(\eta)=T^{1,+}_{+,low}(\eta)+T^{1,+}_{-,low}(\eta)$, where
\begin{multline*}
 (T^{1,+}_{+,low}(\eta)f)(z_1)=
\frac{\cosh\sqrt{c_{23}}z_1}{2\gamma_{23}(-\eta)}\chi_+(z_1)
\int_{z_1}^\infty\pd_{z_1'}\left(e^{\gamma_{23}(-\eta)(z_1-z_1')}\sech\sqrt{c_{23}z_1'}
\right)f(z_1')\,dz_1'
\\ 
-\frac{\cosh\sqrt{c_{23}}z_1}{2\gamma_{23}(-\eta)}\chi_-(z_1)
\int^{z_1}_{-\infty}\pd_{z_1'}\left(e^{\gamma_{23}(-\eta)(z_1-z_1')}\sech\sqrt{c_{23}z_1'}\right)f(z_1')\,dz_1'\,.
\end{multline*}

\begin{lemma}
  \label{lem:T1+'k2large}
Assume that $\k_2\ge a_{14}$ and
$0<\a<\min\{\sqrt{c_{23}},\frac{c_{23}}{2(a_{23}-a_{14})}\}$.
Let $\eta=\eta_R+i\eta_I$, $\eta_R\in\R$ and $\eta_I=2\a(a_{14}-a_{23})$.
If $|\eta_R|<\eta_{\#,-}$, then $\ker\mM_+(\eta)=\{0\}$ and
$\operatorname{Range}\mM_+(\eta)={}^\perp\spann\{g^*(\pm\cdot,-\eta)\}$.
\par
If $|\eta_R|<\eta_{\#,-}$ and $f\in L^2(\R;e^{2\a z_1}dz_1)$ satisfies
  \begin{equation}
    \label{eq:T1+orth-pm}
 \int_\R f(z_1)\overline{g^*(\pm z_1,-\eta)}\,dz_1=0\,,
\end{equation}
then $v=T^{1,+}_{low\#}(\eta)f$ is a solution of $\mM_+(\eta)v=f$  in $ L^2(\R;e^{2\a z_1}dz_1)$.
Moreover, if $|\eta_R|\le\eta_0<\eta_{\#,-}(\a,\beta)$,
\begin{equation}
  \label{eq:T1+sh-est}
\|v\|_{H^1(\R;e^{2\beta z_1}dz_1)}+\|\pd_{z_1}^{-1}v\|_{L^2(\R;e^{2\beta z_1}dz_1)}
\le  C\|f\|_{L^2(\R;e^{2\beta z_1}dz_1)}\,, 
\end{equation}
where $C$ is a positive constant depending only on $\eta_0$, $\a$ and $\beta$.
\end{lemma}

Finally, we will investigate Green functions of $\mM_-(\eta)$ in the case
where $g(\pm z_1,\eta)$ are eigenfunctions of $\mL_1(\eta)$.
Let $T^{1,-}_{low\#}(\eta)=T^{1,-}_{+,low}(\eta)+T^{1,-}_{-,low}(\eta)$ and
\begin{align*}
T^{1,-}_{+,low}(\eta)f(z_1)=&
-\frac{1}{2\gamma_{23}(\eta)}
\int_0^{z_1}\pd_{z_1}\left(e^{\gamma_{23}(\eta)(z_1-z_1')}
\sech\sqrt{c_{23}}z_1\right)\cosh\sqrt{c_{23}}z_1'f(z_1')\,dz_1'\,.
\end{align*}

\begin{lemma}
  \label{lem:T1-'k3small}
Assume that $\k_3\le a_{14}$  and that $0<\a<\min\{\sqrt{c_{23}},\frac{c_{23}}{2(a_{14}-a_{23})}\}$.
Let $\eta=\eta_R+i\eta_I$, $\eta_R\in\R$ and $\eta_I=2\a(a_{14}-a_{23})$.
If $|\eta_R|<\eta_{\#,+}$, then
\begin{equation*}
\ker(\mM_-(\eta))=\spann\{g_M(\pm\cdot,\eta)\}\,,\quad
\operatorname{Range}(\mM_-(\eta))=L^2(\R;e^{2\a z_1}dz_1)\,,
\end{equation*}
and for every $f\in L^2(\R;e^{2\a z_1}dz_1)$,
$v=T^{1,-}_{low\#}(\eta)f$ is a solution of $\mM_-(\eta)v=f$.
Moreover, if $|\eta_R|\le \eta_0<\eta_{\#,+}(\a,\beta)$, then
\begin{equation}
\label{eq:T1-sh-est}
\|v\|_{H^1(\R;e^{2\beta z_1}dz_1)}+\|\pd_{z_1}^{-1}v\|_{L^2(\R;e^{2\beta z_1}dz_1)}\le
 C\|f\|_{L^2(\R;e^{2\beta z_1}dz_1)}\,,
\end{equation}
where $C$ is a positive constant $C$ depending only on $\eta_0$, $\a$ and
$\beta$.
\end{lemma}
\begin{proof}[Proof of Lemma~\ref{lem:T1+'k2large}]
Since $0<\Re\gamma_{23}(-\eta)<\sqrt{c_{23}}-\a$ by the assumption,
it follows from \eqref{eq:mM+} that
$\ker(\mM_+(\eta))=\{0\}$ and
$\ker(\mM_+(\eta)^*)=\spann\{g^*(\pm\cdot,-\eta)\}$.
\par
Next, we will estimate $v=T^{1,+}_{low\#}(\eta)f$.
Since $|\eta_R|<\eta_{\#,-}(\a,\beta)<\eta_{*,-}(\a,\beta)$,
we can estimate $T^{1,+}_{-,low}(\eta)$ in exactly the same way as Lemma~\ref{lem:T1+'}.
By \eqref{eq:26},
$$
\left|T^{1,+}_{+,low}(\eta)f(z_1)\right|
\lesssim \int_\R e^{\Re\gamma_{23}(-\eta)(z_1-z_1')-\sqrt{c_{23}}|z_1-z_1'|}
|f(z_1')|\,dz_1'\,.$$
Since $0<\Re\gamma_{23}(-\eta)\le \Re\gamma_{23}(-\eta_0-i\eta_I)
<\sqrt{c_{23}}-\beta$,
\begin{align*}
\|T^{1,+}_{+,low}(\eta)f\|_{L^2(\R;e^{2\beta z_1}dz_1)}
 \lesssim &
\|e^{(\beta-\gamma_{23}(-\eta))z_1-\sqrt{c_{23}}|z_1|}\|_{L^1(\R)}
                \|f\|_{L^2(\R;e^{2\beta z_1}dz_1)}
\\ \le &  C\|f\|_{L^2(\R;e^{2\beta z_1}dz_1)}\,,
\end{align*}
where $C$ is a constant depending only on $\eta_0$, $\a$ and $\beta$.
By \eqref{eq:expT1++low} and \eqref{eq:expT1+-low}, we have
$T^{1,+}_{\pm,low}(\eta)f=T^{1,+}_\pm(\eta)f$ for $f$ satisfying
\eqref{eq:T1+orth-pm}, and $\mM_+(\eta)T^{1,+}_{low\#}(\eta)f=f$.
Combining the above with the fact that
$\pd_{z_1}^{-1}\in B(L^2(\R;e^{2\a z_1}dz_1)$, we have
\eqref{eq:T1+sh-est}.  Thus we complete the proof.
\end{proof}

\begin{proof}[Proof of Lemma~\ref{lem:T1-'k3small}]
Since $0<\Re\gamma_{23}(\eta)<\sqrt{c_{23}}-\a$ by the assumption,
we have $\ker(\mM_-(\eta))=\spann\{g_M(\pm\cdot,\eta)\}$ and
$\ker(\mM_-(\eta)^*)=\{0\}$.

Since $|\eta_R|<\eta_{\#,+}(\a,\beta)<\eta_{*,+}(\a,\beta)$, we can estimate $T^{1,-}_{-,low}(\eta)f$
in exactly the same way as the proof of Lemma~\ref{lem:T1-'}.
By \eqref{eq:26} and the fact that
$\Re\gamma_{23}(\eta)\le \Re\gamma_{23}(\eta_0+i\eta_I)<\sqrt{c_{23}}-\beta$, 
\begin{align*}  
 \|T^{1,-}_{+,low}(\eta)f\|_{L^2(\R;e^{2\beta z_1}dz_1)}
  \lesssim &
\|e^{(\beta+\Re\gamma_{23}(\eta))z_1}e^{-\sqrt{c_{23}}|z_1|}\|_{L^1(\R)}
\|f\|_{L^2(\R;e^{2\beta z_1}dz_1)}
\\ \le & C\|f\|_{L^2(\R;e^{2\beta z_1}dz_1)}\,,
\end{align*}
where $C$ is a constant depending only on $\eta_0$, $\a$ and $\beta$.
Thus we prove $T^{1,-}_{low\#}(\eta)$ is bounded on $L^2(\R;e^{2\beta z_1}dz_1)$
provided $\beta$ is sufficiently close to $\a$.
Since $T^{1,-}(\eta)f-T^{1,-}_{-,low\#}(\eta)f\in\spann\{g_M(\pm\cdot,\eta)\}$
and $g_M(\pm\cdot,\eta)\in\ker(\mM_-(\eta))$, it follows that
$\mM_-(\eta)T^{1,-}_{-,low\#}(\eta)=I$ on $L^2(\R;e^{2\a z_1}dz_1)$.
Combining the above with the definition of $\mM_-(\eta)$,
we have \eqref{eq:T1-sh-est}.
Thus we complete the proof.
\end{proof}

\bigskip

\subsection{Spectral stability of $[1,4]$-soliton in $\calX_1$}
In this subsection, we will investigate spectral stability of $\wmL_1$
in $\calX_1$ whose weight function increases in a direction of a motion
of $[2,3]$-soliton.
To begin with, we investigate the upper bound of the growth rate of
continuous eigenfunctions $e^{t\lambda_{2,\pm}(\eta)+iy\eta}\tg(z_1,\pm\eta)$
in $\calX_1$.
Let $\eta=\eta_R+i\eta_I$, $\eta_I=2(a_{23}-a_{14})\a$ and $\eta_R\in\R$. Then
$$\calX_1=e^{-y\eta_I}\calX_2=\bigoplus_{\eta=\eta_R+i\eta_I,\,\eta_R\in\R}
e^{iy\eta}L^2(\R;e^{2\a z_2}dz_2)\,.$$
We have $\Re\gamma_{14}(\eta'_\pm\pm i\eta_I)=\sqrt{c_{14}}+\a$ for
\begin{gather}
  \label{def:eta+'}
  \eta_+'=2(\sqrt{c_{14}}+\a)\sqrt{\a(\a+2a_{23}-2\k_1)}\,,\\
  \label{def:eta-'}
\eta_-'=2(\sqrt{c_{14}}+\a)\sqrt{\a(\a+2\k_4-2a_{23})}\,.
\end{gather}
\begin{lemma}
\label{lem:lambda2-etaI}  
Assume \eqref{ass-alpha'} and \eqref{eq:l2ne0}.
Let $\eta=\eta_R+i\eta_I$, $\eta_R\in\R$ and $\eta_I=2(a_{23}-a_{14})\a$.
Then $\Re\lambda_{2,\pm}(\eta_R+i\eta_I)$ are even in $\eta_R$,
$\pd_{\eta_R}\Re \lambda_{2,\pm}(\eta_R+i\eta_I)<0$ if $\eta_R>0$ and
$ \Re\lambda_{2,\pm}(\eta'_\pm+i\eta_I)<0$.
If $\pm(a_{14}-a_{23})>0$, then
there exists $\tilde{\eta}_{*,\pm}\in(0,\eta_\pm')$ satisfying
$\Re\lambda_{2,\pm}(\tilde{\eta}_{*,\pm})=0$.

For every $\eta_R\in\R$, $\lambda_{2,\pm}(\eta_R+i\eta_I)\ne0$ and
\begin{equation}
  \label{eq:23}
  (a_{14}-a_{23})\left\{\Re\lambda_{2,+}(\eta)-\Re\lambda_{2,-}(\eta)\right\}>0\,.
\end{equation}
\end{lemma}
\begin{proof}
 Let
$\tilde{\gamma}_{R,\pm}(\eta_R)=\Re \gamma_{14}(\eta_R\pm i\eta_I)$ and
$\tilde{\gamma}_{I,\pm}(\eta_R)=\Im \gamma_{14}(\eta_R\pm i\eta_I)$.
We can prove
\begin{gather*}
  \tilde{\gamma}_{R,\pm}(\eta_R)>0\,,\quad
  \eta_R\tilde{\gamma}_{I,\pm}(\eta_R)\ge 0\,,
  \quad \pd_{\eta_R}\tilde{\gamma}_{I,\pm}(\eta_R)>0
  \quad \text{for $\eta_R\in\R$,}
  \\  \eta_R\pd_{\eta_R}\tilde{\gamma}_{R,\pm}(\eta_R)>0\,,\quad
  \eta_R\pd_{\eta_R}\Re\lambda_{2,\pm}(\eta)<0
  \quad \text{if $\eta_R\ne0$,}
\end{gather*}
in the same way as the proof of Lemma~\ref{lem:lambda1-etaI}.
Since $\tilde{\gamma}_{R,\pm}$ is even and $\tilde{\gamma}_{I,\pm}$
is odd, $\Re\lambda_{2,\pm}$ are even.
\par
Since $\tilde{\gamma}_{R,\pm}(\eta'_\pm)=\sqrt{c_{14}}+\a$,
$\eta'_\pm=2(\sqrt{c_{14}}+\a)\tilde{\gamma}_{I,\pm}(\eta'_\pm)$
and $\tilde{\gamma}_{I,\pm}(\eta'_\pm)^2=\a^2+2\sqrt{c_{14}}\a\pm\eta_I$,
  \begin{align*}
 \Re\lambda_{2,\pm}(\eta'_\pm+i\eta_I)=&
-\eta_I(b_2-3a_{14}\pm \tilde{\gamma}_{R,\pm}(\pm\eta'_\pm))
\mp\eta'_\pm\tilde{\gamma}_{I,\pm}(\pm\eta'_\pm)
 \\ =& -\eta_I\left\{b_2-3a_{14}\pm 3(\sqrt{c_{14}}+\a)\right\}
       -2\a(\sqrt{c_{14}}+\a)(2\sqrt{c_{14}}+\a)\,.
  \end{align*}
  Substituting \eqref{eq:b1-b2P'} into the above, we have
  \begin{align*}
& \Re\lambda_{2,+}(\eta'_++i\eta_I)=
-\{3(a_{23}-\k_1)^2+c_{23}\}\a-6(a_{23}-\k_1)\a^2-2\a^3<0\,,
\\ &                                       
\Re\lambda_{2,-}(\eta'_-+i\eta_I)=
-\{3(\k_4-a_{23})^2+c_{23}\}\a-6(\k_4-a_{23})\a^2-2\a^3<0\,.
  \end{align*}
\par

Suppose that  $a_{14}>a_{23}$. Then $\eta_I<0$ and it follows from
Lemma~\ref{cl:sgnpdlambda1} that
\begin{align*}
  \Re\lambda_{2,+}(i\eta_I)
  =& -\eta_I\left\{b_2-3a_{14}+\gamma_{14}(i\eta_I)\right\}
 \ge  -\eta_I(b_2-3a_{14}+\sqrt{c_{14}})>0\,.
\end{align*}
Since $\Re\lambda_{2,+}(\cdot+i\eta_I)$ is monotone decreasing for $\eta_R>0$,
$\lambda_{2,+}(i\eta_I)>0$
and $\Re\lambda_{2,+}(\eta_+'+i\eta_I)<0$,
there exists $\tilde{\eta}_{*,+}\in(0,\eta_+')$ such that
$\Re\lambda_{2,+}(\eta)>0$ if $|\eta_R|<\tilde{\eta}_{*,+}$
and that $\Re\lambda_{2,+}(\eta)<0$ if $|\eta_R|>\tilde{\eta}_{*,+}$.
\par
Suppose that  $a_{14}<a_{23}$.  Then $\eta_I>0$ and by
Lemma~\ref{cl:sgnpdlambda1},
\begin{align*}
\Re\lambda_{2,-}(i\eta_I)
=& -\eta_I\left\{b_2-3a_{14}-\gamma_{14}(-i\eta_I)\right\}
\ge -\eta_I(b_2-3a_{14}-\sqrt{c_{14}})>0\,,
\end{align*}
and there exists $\tilde{\eta}_{*,-}>0$ such that
$\Re\lambda_{2,-}(\eta)>0$ if $|\eta_R|<\tilde{\eta}_{*,-}$
and that $\Re\lambda_{2,-}(\eta)<0$ if $|\eta_R|>\tilde{\eta}_{*,-}$.
\par
Now we will prove \eqref{eq:23}. We have $\tilde{\gamma}_{R,\pm}(\eta_R)>0$
and
\begin{align*}
  \Re\lambda_{2,+}(\eta)-\Re\lambda_{2,-}(\eta)=&
-\eta_I\left\{\tilde{\gamma}_{R,+}(\eta_R)+\tilde{\gamma}_{R,-}(-\eta_R)\right\}
-\eta_R\left\{\tilde{\gamma}_{I,+}(\eta_R)
+\tilde{\gamma}_{I,-}(-\eta_R)\right\}\,.
\\=& 
-\eta_I\left\{\tilde{\gamma}_{R,+}(\eta_R)+\tilde{\gamma}_{R,-}(\eta_R)\right\}
-\eta_R\left\{\tilde{\gamma}_{I,+}(\eta_R)
-\tilde{\gamma}_{I,-}(\eta_R)\right\}\,.
\end{align*}
We see that $\eta_I$ and
$\tilde{\gamma}_{I,+}(\eta_R)-\tilde{\gamma}_{I,-}(\eta_R)$ have
the same sign if $\eta_R>0$, and have the opposite sign if $\eta_R<0$.
Thus we have \eqref{eq:23}.
\par

If $\lambda_{2,\pm}(\eta)=0$, then $b_2-3a_{14}\pm\gamma_{14}(\pm\eta)=0$
and we have $\eta_R=0$ and 
$$\a=\pm\frac{(b_2-3a_{14})^2-c_{14}}{2(a_{14}-a_{23})}\,.$$
Thus we have $\lambda_{2,\pm}(\eta)\ne0$.
This completes the proof of Lemma~\ref{lem:lambda2-etaI}.
\end{proof}

\begin{remark}
  \label{rem:ins2}
Let $\eta=\eta_R+i\eta_I$ and $\eta_I=2(a_{23}-a_{14})\a$.
Then
$\Re\lambda_{2,+}(\eta)>\Re\lambda_{2,-}(\eta)$ and   
$\eta_+'<\eta_-'$ if $a_{14}>a_{23}$ and $\Re\lambda_{2,+}(\eta)<\Re\lambda_{2,-}(\eta)$ and
$\eta_+'>\eta_-'$ if $a_{14}<a_{23}$.
The endpoints of the curves of continuous eigenvalues
$\{\lambda_{2,\pm}(\eta_R+i\eta_I)\}_{|\eta_R|\le\eta'_\pm}$ lie in the stable
half plane.
\end{remark}

\begin{proposition}
\label{prop:spec[1-4]}
Assume \eqref{ass-alpha'}. Let $\mL_0$ and $\wmL_1$ be closed operators on $\calX_1$. Then
$$\sigma(\wmL_1)=\sigma(\mL_0)\cup
\{\lambda_{2,\pm}(\eta_R+i\eta_I)\mid \eta_R\in[-\eta_\pm',\eta_\pm']\}\,.$$
\end{proposition}
\begin{remark}
It follows from Lemma~\ref{lem:lambda2-etaI} that  $\Re \lambda_{2,\pm}(\eta)>0$
if $\pm(a_{14}-a_{23})>0$ and $|\eta_R|<\tilde{\eta}_{*,\pm}$ and that
$$\sigma(\wmL_1)\cap \{\lambda\in\C\mid \Re\lambda>0\}\ne\emptyset\,.$$
Continuous eigenmodes $\{g^1_{2,\pm}(\eta_R+i\eta_I)\}$ have to do with
modulations of $[1,4]$-soliton.
Proposition~\ref{prop:spec[1-4]} suggests that some modulations
of $[1,4]$-soliton go ahead of $[2,3]$-soliton.
\end{remark}

To prove Proposition~\ref{prop:spec[1-4]}, we will use Darboux transformations.
In the moving coordinate
$(t,z_2,y)$, the operator $\wmL_1$ can be written as
\begin{gather*}
\wmL_1=\mL_0-\frac{3}{2}\pd_{z_2}(\tu_1\cdot)\,,
\\
\mL_0=
\frac{1}{4}\left\{-\pd_{z_2}(\pd_{z_2}^2-4c_{14})+4(b_2-3a_{14})\pd_y
-3\pd_{z_2}^{-1}\pd_y^2\right\}\,.
\end{gather*}
Linearizing \eqref{eq:mkp} around
$\tv_1=a_{14}+\tpsi$ with $\tpsi=\sqrt{c_{14}}\tanh\sqrt{c_{14}}z_2$,
we have
\begin{gather*}
  \pd_tv=\wmL_{M1}v\,,\\
  \wmL_{M1}=\mL_0-\frac{3}{4}\left\{\pd_{z_2}(\tu_1\cdot)
    +\tu_1\pd_{z_2}^{-1}\pd_y\right\}\,.
\end{gather*}
Let $\wmL_0(\eta)$, $\wmL_1(\eta)$, $\wmL_{M1}(\eta)$ and $\wmM_\pm(\eta)$
be operators on $L^2(\R;e^{2\a z_2}dz_2)$ defined by
\begin{gather*}
\wmL_0(\eta)=e^{-iy\eta} \wmL_0e^{iy\eta}\,,\quad
\wmL_1(\eta)=e^{-iy\eta} \wmL_1e^{iy\eta}\,,\quad
\wmL_{M1}(\eta)=e^{-iy\eta}\wmL_{M1}e^{iy\eta}\,,
\\
\wmM_\pm(\eta)=e^{-iy\eta}\nabla M_\pm(\tv_1)e^{iy\eta}
=\pm\pd_{z_2}+i\eta\pd_{z_2}^{-1}-2\tpsi\,.
\end{gather*}
Note that $\tu_1$ and $\tpsi$ depend only on $z_2$.
By \eqref{eq:bH},
\begin{equation}
  \label{eq:1t}
\wmL_1(\eta)\wmM_+(\eta)=\wmM_+(\eta)\wmL_{M1}(\eta)\,,\quad
\wmL_0(\eta)\wmM_-(\eta)=\wmM_-(\eta)\wmL_{M1}(\eta)\,.
\end{equation}
\par

Let $k_2(z_2,z_2',\eta)$ be as \eqref{eq:T2decomp2},
$\tk_1(z_2,z_2',\eta)=\pd_{z_2}k_2(z_2,z_2',\eta)$ and
\begin{gather*}
\wT^{1,-}(\eta)f(z_2)
=\int_\R \tk_1(z_2,z_2',\eta)f(z_2')\,dz_2'\,,
\\
\wT^{1,+}(\eta)f(z_2)
=\int_\R \tk_1(z_2',z_2,-\eta)f(z_2')\,dz_2'\,.
\end{gather*}
Let $\eta_I(\a)=2\a(a_{23}-a_{14})$ and $\eta_\pm'$ be as in \eqref{def:eta+'}
and \eqref{def:eta-'}. Let
$$\eta_\pm'(\a,\beta)=2(\sqrt{c_{14}}+\beta)
\sqrt{\beta^2+2\sqrt{c_{14}}\beta\pm\eta_I(\a)}\,.$$
Note that  $\eta_\pm'(\a,\a)=\eta_\pm'$ and $\eta_\pm'(\a,\beta)>0$
if $\a>0$ and $\beta$ is sufficiently close to $\a$.

If $|\eta_R|$ is sufficiently large, Darboux transformations $\mM_\pm(\eta)$ are bijective.
\begin{lemma}
  \label{lem:wT1pm}
  Let $\a\in(0,\sqrt{c_{14}})$, $\eta=\eta_R+i\eta_I$, $\eta_R\in\R$ and
  $\eta_I=2\a(a_{23}-a_{14})$. 
 Suppose that $|\eta_R|>\eta_\mp'$. Then
 $$\wmM_\pm(\eta)\wT^{1,\pm}(\eta)=\wT^{1,\pm}(\eta)\wmM_\pm(\eta)=I\quad\text{on $L^2(\R;e^{2\a z_2}dz_2)$.}$$
 Moreover, if $\beta$ is sufficiently close to $\a$ and $|\eta_R|\ge\eta_0>\eta_\mp'(\a,\beta)$,
\begin{equation}
 \label{eq:wT1-est}
  \sum_{j=0,1,2}\la \eta\ra^{(2-j)/2}
  \|\pd_{z_2}^{j-1}\wT^{1,\pm}(\eta)f\|_{L^2(\R;e^{2\beta z_2}dz_2)}
  \le C\|f\|_{L^2(\R;e^{2\beta z_2}dz_2)}\,,
\end{equation}
where $C$ is a positive constant depending only on $\a$, $\beta$ and $\eta_0$.
\end{lemma}

We have $\wT^{1,\pm}(\eta)=\wT^{1,\pm}_+(\eta)+\wT^{1,\pm}_-(\eta)$, where
\begin{align*}  
(\wT^{1,+}_+(\eta)f)(z_2)=& \frac{\cosh\sqrt{c_{14}}z_2}{2\gamma_{14}(-\eta)}
\int_{z_2}^\infty\pd_{z_2'}\left(e^{\gamma_{14}(-\eta)(z_2-z_2')}
\sech\sqrt{c_{14}}z_2'\right)f(z_2')\,dz_2'\,,
  \\
(\wT^{1,+}_-(\eta)f)(z_2)=& \frac{\cosh\sqrt{c_{14}}z_2}{2\gamma_{14}(-\eta)}
\int^{z_2}_{-\infty}\pd_{z_2'}\left(e^{\gamma_{14}(-\eta)(z_2'-z_2)}
\sech\sqrt{c_{14}}z_2'\right)f(z_2')\,dz_2'
 \\
\wT^{1,-}_+(\eta)f(z_2)=& \frac{1}{2\gamma_{14}(\eta)}
\int^\infty_{z_2}\pd_{z_2}\left(e^{\gamma_{14}(\eta)(z_2-z_2')}\sech\sqrt{c_{14}}z_2\right)
\cosh\sqrt{c_{14}}z_2'f(z_2')\,dz_2'
  \\
\wT^{1,-}_-(\eta)f(z_2)=&
\frac{1}{2\gamma_{14}(\eta)}
\int_{-\infty}^{z_2}\pd_{z_2}\left(e^{-\gamma_{14}(\eta)(z_2-z_2')}
\sech\sqrt{c_{14}}z_2\right)\cosh\sqrt{c_{14}}z_2'f(z_2')\,dz_2',.
\end{align*}

For small $\eta$, we replace $\wT^{1,+}_-(\eta)$ and $\wT^{1,-}_-(\eta)$ by
\begin{multline*}
(\wT^{1,+}_{-,low})(\eta)f(z_2)=
\frac{\cosh\sqrt{c_{14}}z_2}{2\gamma_{14}(-\eta)}\chi_-(z_2)
\int^{z_2}_{-\infty}\pd_{z_2'}\left(e^{\gamma_{14}(-\eta)(z_2'-z_2)}\sech\sqrt{c_{14}z_2'}
\right)f(z_2')\,dz_2'
\\
-\frac{\cosh\sqrt{c_{14}}z_2}{2\gamma_{14}(-\eta)}\chi_+(z_2)
\int_{z_2}^\infty\pd_{z_2'}\left(e^{\gamma_{14}(-\eta)(z_2'-z_2)}\sech\sqrt{c_{14}z_2'}
\right)f(z_2')\,dz_2'\,,
\end{multline*}
\begin{align*}
(\wT^{1,-}_{-,low})(\eta)f(z_2)=&
\frac{1}{2\gamma_{14}(\eta)}
\int_0^{z_2}\pd_{z_2}\left(e^{-\gamma_{14}(\eta)(z_2-z_2')}
  \sech\sqrt{c_{14}}z_2\right)\cosh\sqrt{c_{14}}z_2'f(z_2')\,dz_2'\,,
\end{align*}
and let $\wT^\pm(\eta)=\wT^\pm_+(\eta)+\wT^\pm_{-,low}(\eta)$.

In the case where $\tg(z_2,\pm\eta)$ are eigenfunctions of $\wmL_1(\eta)$, we have the following.
\begin{lemma}
  \label{lem:wT1+'}
  Assume \eqref{ass-alpha'}.
Let $\eta=\eta_R+i\eta_I$, $\eta_R\in\R$ and $\eta_I=2\a(a_{23}-a_{14})$.
Suppose that $|\eta_R|<\eta_-'$. Then $\ker(\wmM_+(\eta))=\{0\}$,
$\operatorname{Range}(\wmM_+(\eta))={}^\perp\spann\{\tg^*(\cdot,-\eta)\}$ and
for $f\in L^2(\R;e^{2\a z_2}dz_2)$ satisfying
\begin{equation*}
  \int_\R f(z_2)\overline{\tg^*(z_2,-\eta)}\,dz_2=0\,,
\end{equation*}
$v=\wT^{1,+}_{low}(\eta)f$ is a solution of $\wmM_+(\eta)v=f$ in $L^2(\R;e^{2\a z_2}dz_2)$.
Moreover, if $\beta$ is sufficiently close to $\a$ and $|\eta_R|\le\eta_0<\eta_-'(\a,\beta)$,
\begin{equation*}
\|v\|_{H^1(\R;e^{2\beta z_2}dz_2)}+\|\pd_{z_2}^{-1}v\|_{L^2(\R;e^{2\beta z_2}dz_2)}
\le  C\|f\|_{L^2(\R;e^{2\a z_2}dz_2)}\,, 
\end{equation*}
where $C$ is a positive constant depending only on $\eta_0$, $\a$ and $\beta$.
\end{lemma}
Let
\begin{gather*}
\tg_M(z_2,\eta)=\pd_x^{-1}\tg(z_2,\eta)
= \sqrt{c_{14}}d_{2,+}(\eta)
\pd_{z_2}\left(e^{-\gamma_{14}(\eta)z_2}\sech\sqrt{c_{14}}z_2\right)\,,\\
\tg_M^*(z_2,\eta)=\frac{1}{2\bar{\eta}}\wmM_+(\eta)^*\tg_1^*(z_2,\eta)
=-\frac{1}{2} e^{\gamma_{14}(-\bar{\eta})z_2}\sech\sqrt{c_{14}}z_2\,.
\end{gather*}
\begin{lemma}
  \label{lem:wT1-'}
Assume \eqref{ass-alpha'}.  
Let $\eta=\eta_R+i\eta_I$ $\eta_R\in\R$ and $\eta_I=2\a(a_{23}-a_{14})$.
Suppose that  $|\eta_R|<\eta_+'$. Then
\begin{equation*}
\ker(\wmM_-(\eta))=\spann\{\tg_M(\cdot,\eta)\}\,,\quad
\operatorname{Range}(\wmM_-(\eta))=L^2(\R;e^{2\a z_2}dz_2)\,,
\end{equation*}
and for any $f\in L^2(\R;e^{2\a z_2}dz_2)$ satisfying
$$\int_\R v(z_2)\overline{\tg^*_M(z_2,\eta)}\,dz_2=0\,,$$
$v=\wT^-_{low}(\eta)f$ is a solution of $\wmM_-(\eta)v=f$ in $L^2(\R;e^{2\a z_2}dz_2)$.
Moreover, if $\beta$ is sufficiently close to $\a$ and $|\eta_R|\le\eta_0<\eta_+'(\a,\beta)$,
\begin{gather*}
\|v\|_{H^1(\R;e^{2\beta z_2}dz_2)}+\|\pd_{z_2}^{-1}v\|_{L^2(\R;e^{2\beta z_2}dz_2)}\le
 C\|f\|_{L^2(\R;e^{2\beta z_2}dz_2)}\,,
\end{gather*}
where $C$ is a positive constant $C$ depending only on $\eta_0$, $\a$ and $\beta$.
\end{lemma}
We can prove Lemmas~\ref{lem:wT1pm}--\ref{lem:wT1-'} 
in exactly the same way as Lemmas~\ref{lem:T1+}--\ref{lem:T1-'k3small}.
\par
Using Lemmas~\ref{lem:wT1pm}--\ref{lem:wT1-'},
we can prove Proposition~\ref{prop:spec[1-4]} in the same
way as Propositions~\ref{lem:wmL1-decay} and \ref{prop:specmL1}.
\begin{proof}[Proof of Proposition~\ref{prop:spec[1-4]}]
For $\eta=\eta_R+i\eta_I$ satisfying  $|\eta_R|< \eta_\pm'$, we have
$\sqrt{c_{14}}-\a<\Re\gamma_{14}(\pm\eta)<\sqrt{c_{14}}+\a$.
Hence it follows from \eqref{eq:defg1z2} and Lemma~\ref{lem:evtmL1} that
$\lambda_{2,\pm}(\eta_R+i\eta_I)\in \sigma(\wmL_1)$
for $\eta_R\in[-\eta_\pm',\eta_\pm']$.
As in the proof of Proposition~\ref{lem:wmL1-decay},
we see that $ip(\xi+i\a,\eta+2i\a a_{23})$
with $(\xi,\eta)\in\R^2$ are continuous eigenvalues of $\wmL_1$.
Thus we have $$\tilde{\sigma}(\wmL_1):=
\sigma(\mL_0)\cup\{\lambda_{2,\pm}(\eta_R+i\eta_I)\mid\eta_R\in[-\eta_\pm',\eta_\pm']\}
\subset\sigma(\wmL_1)\,.$$
Using Lemmas~\ref{lem:wT1pm}--\ref{lem:wT1-'}, we can prove $\sigma(\wmL_1)=\tilde{\sigma}(\wmL_1)$
in the same way as Proposition~\ref{lem:wmL1-decay}.
Thus we complete the proof.
\end{proof}

\bigskip

\subsection{Linear stability of $2$-line solitons of P-type in $\calX_2$}
\label{subsec:LSP-2}
In this subsection, we will prove that $2$-line solitons of
P-type are linearly stable in $\calX_2$.
\begin{proposition}
  \label{prop:linearstability-1}
Assume \eqref{eq:ass-alpha}.
Let $P_2(\eta_0)$ be a projection defined by \eqref{eq:defP2},
$\eta'\in(0,\eta_{*,2})$ and $\eta_{*,2}$ be as \eqref{def:eta*2}.
There exist positive constants $K$ and $b$ such that if
$\eta_0\in(0,\eta']$
$$\left\|e^{t\mL_2}\left(I-P_2(\eta_0)\right)\right\|_{B(\calX_2)}\le Ke^{-bt}\quad\text{for $t\ge0$.}$$
\end{proposition}

To begin with, we will see that
$\{\lambda_{2,\pm}(\eta)\}_{\eta\in\R}$ lie in the left half plane of $\C$.
\begin{lemma}
  \label{lem:lambda2-est}
  If $\eta\in\R\setminus\{0\}$,
$\eta\pd_\eta\Re\lambda_{2,\pm}(\eta)<0$ and
$\Re\lambda_{2,\pm}(\eta)<\Re\lambda_{2,\pm}(0)=0$.
\par
Let $\mL_0$ and $\mL_2$ be operators on $\calX_2$. Then
$\lambda_{2,+}(\pm\eta_{*,2})$, $\lambda_{2,-}(\pm\eta_{*,2})\in\sigma(\mL_0)$
and 
\begin{equation}
  \label{eq:25}
\{\lambda_{2,\pm}(\eta)\mid -\eta_{*,2}\le\eta\le\eta_{*,2}\}
\subset \sigma(\mL_2)\,.
\end{equation}
\end{lemma}
\begin{proof}
  We can prove that $\Re\lambda_{2,\pm}(\eta)<0$ for
  $\eta\in\R\setminus\{0\}$ in exactly the same way as
  Lemma~\ref{lem:lambda1-etaI}.  We have \eqref{eq:25} from
  Lemma~\ref{lem:evmL2} and the fact that
  $\left\|e^{\a z_2}g^2_{2,\pm}(\cdot,\eta)\right\|_{L^2_x}$ is
  bounded in $y$.
\par
For $\xi_{*,2}=\sqrt{\a(\a+\k_4-\k_1)}$,
we have $\eta_{*,2}=(2\a+\k_4-\k_1)\xi_{*,2}$ and
$\lambda_{2,\pm}(\eta_{*,2})=ip(\mp\xi_{*,2}+i\a,
\eta_{*,2}+2a_{14}(\mp\xi_{*,2}+i\a))$,
$\lambda_{2,\pm}(-\eta_{*,2})=ip(\pm\xi_{*,2}+i\a,
-\eta_{*,2}+2a_{14}(\pm\xi_{*,2}+i\a))
\in\sigma(\mL_0)$.
Thus we complete the proof.
\end{proof}
\begin{remark}  
  \label{rem:lambda1-est}
Lemma~\ref{lem:lambda2-est} tells us that endpoints of a curve of
continuous eigenvalues $\{\lambda_{2,+}(\eta)\}_{-\eta_{*,2}\le\eta\le\eta_{*,2}}$
belongs to $\sigma(\mL_0)$.
\par
Similarly, $\eta_{*,1}=(2\a+\k_3-\k_2)\xi_{*,1}$, $\xi_{*,1}=\sqrt{\a(\a+\k_3-\k_2)}$ and
$\lambda_{1,\pm}(\eta_{*,1})=ip(\mp\xi_{*,1}+i\a,
\eta_{*,1}+2a_{23}(\mp\xi_{*,1}+i\a))$,
$\lambda_{1,\pm}(-\eta_{*,1})=ip(\pm\xi_{*,1}+i\a,
-\eta_{*,1}+2a_{23}(\pm\xi_{*,1}+i\a))$, and we have 
$\lambda_{1,+}(\pm\eta_{*,1})$, $\lambda_{1,-}(\pm\eta_{*,1})\in\sigma(\mL_0)$
in $\calX_1$.
\end{remark}

If $\lambda\in\sigma(\mL_2)$, then either $\lambda$ is a point
spectrum, a continuous spectrum or a residual spectrum.
We will show that $\lambda\in\sigma_p(\mL_2)$ if
$\lambda\in\sigma(\mL_2)$ and
$\lambda$ is a resolvent point of $\mL_1$ and $\wmL_1$
and show simplicity of continuous eigenvalues $\lambda_{2,\pm}(\eta)$.
\begin{lemma}
  \label{lem:spec1}
  Assume $\a\in(0,2\sqrt{c_{23}})$. Let $\mL_0$, $\mL_1$, $\wmL_1$ and
  $\mL_2$ be closed operators on $\calX_2$.  If
  $\lambda\in \sigma(\mL_2)\cap\rho(\mL_0)\cap\rho(\mL_1)\cap\rho(\wmL_1)$, then
  $\lambda$ is an isolated eigenvalue of $\mL_2$ with finite
  multiplicity.
\end{lemma}

\begin{lemma}
  \label{lem:anal-F1}
Assume \eqref{eq:ass-alpha}.
Let $\mL_2$ be a closed operator on $\calX_2$. Let
$\eta_{0,2}\in(0,\eta_{*,2})$ and
$\widetilde{\mathcal{C}}=
\{\lambda_{2,\pm}(\eta)\mid \eta\in(-\eta_0,\eta_0)\}$.
Then $\widetilde{\mathcal{C}}\cap\sigma_p(\mL_2)$ is a discrete set and
$\lambda-\mL_2$ is invertible on $(I-P_2(\eta_0))\calX_2$
if $\lambda\in\widetilde{\mathcal{C}}\setminus\sigma_p(\mL_2)$.
\end{lemma}

\begin{proof}[Proof of Lemmas~\ref{lem:spec1} and \ref{lem:anal-F1}]
Let $\chi_0(r)$ and $\chi_1(r)$ be nonnegative smooth even functions
satisfying $\chi_0(r)+\chi_1(r)=1$ for $r\in\R$,
$\chi_1(r)=1$ for $r\ge 3$ and $\chi_1(r)=0$ for $r\in[-1,1]$.
Let $\tilde{\chi}_0(r)$ and $\tilde{\chi}_1(r)$ be smooth even
functions satisfying 
$$\tilde{\chi}_0(r)=
\begin{cases}
1 & \quad\text{for $r\in[-3,3]$,}\\
0 & \quad\text{for $r\ge 4$,}
\end{cases}
\qquad
\tilde{\chi}_1(r)=
\begin{cases}
1 & \quad\text{for $r\ge1$,}\\
0 & \quad\text{for $r\in[-\frac12,\frac12]$.}
\end{cases}
$$
Note that $\chi_0(r)\tilde{\chi}_0(r)=\chi_0(r)$,
$\chi_1(r)\tilde{\chi}_1(r)=\chi_1(r)$ and
\begin{gather*}
\chi_0(z_1)\chi_0(z_2)+\chi_1(z_1)\chi_1(z_2)+
\chi_0(z_1)\chi_1(z_2)+\chi_0(z_2)\chi_1(z_1)=1\,.
\end{gather*}
Let $\chi_{i,R}(r)=\chi_i(r/R)$ and $\tilde{\chi}_{i,R}(r)=\tilde{\chi}_i(r/R)$
for $i=0$, $1$. Let $\chi_{1,R}^\pm(r)=\chi_{1,R}(r)$ and
$\tilde{\chi}_{1,R}^\pm(r)=\tilde{\chi}_{1,R}(r)$
if $\pm r\ge0$ and $\chi_{1,R}^\pm(r)=\tilde{\chi}_{1,R}^\pm(r)=0$ if $\pm r\le0$.
Let $u_1^\pm(z_1)=u_1(Z_{1,\pm}/\sqrt{c_{23}})$,
$\tu^\pm_1(z_2)=\tu_1(Z_{2,\pm}/\sqrt{c_{14}})$,
\begin{gather*}
\mL_1^\pm=\mL_0-\frac{3}{2}\pd_x(u_1^\pm\cdot)\,,\quad
\widetilde{\mL_1}^\pm=\mL_0-\frac{3}{2}\pd_x(\tu_1^\pm\cdot)\,,
\end{gather*}
and $\widetilde{S}_\pm$ be shift operators defined by
$(\widetilde{S}_\pm f)(z_1,y)=f(Z_{2,\pm}/\sqrt{c_{14}},y)$ so that
$\widetilde{S}_\pm\tu_1=\tu_1^\pm$.

Let $Q_2(\eta_0)=I-P_2(\eta_0)$ and
$\widetilde{Q}_1(\eta_0)=I-\widetilde{P}_1(\eta_0)$.
Note that $Q_2(\eta_0)$ and $\widetilde{Q}_1(\eta_0)$ are spectral projections
associated with $\mL_2$ and $\wmL_1$, respectively.
\par
Let 
$W(\lambda):Q_2(\eta_0)\calX_2\to Q_2(\eta_0)\calX_2$ be a mapping
defined by
\begin{align*}
W(\lambda)=& Q_2(\eta_0)\tilde{\chi}_{0,R}(z_1)\tilde{\chi}_{0,R}(z_2)
         (\lambda-\mL_0)^{-1}\chi_{0,R}(z_1)\chi_{0,R}(z_2)Q_2(\eta_0)
\\ & +
Q_2(\eta_0)\sum_{a=\pm}\tilde{\chi}_{0,R}(z_1)\tilde{\chi}_{1,R}^a(z_2)
(\lambda-\mL_1^a)^{-1}\chi_{0,R}(z_1)\chi_{1,R}^a(z_2)Q_2(\eta_0)
  \\ & +
       Q_2(\eta_0)\sum_{a,b=\pm}\tilde{\chi}_{1,R}^a(z_1)\widetilde{S}_a
       (\lambda-\wmL_1)^{-1}\widetilde{Q}_1(\eta_0)\widetilde{S}_b^{-1}
       \chi_{1,R}^b(z_1)Q_2(\eta_0)\,.
\end{align*}
Then
$W(\lambda)(\lambda-\mL_2)= Q_2(\eta_0)\left(II+B_1^1+B_1^2+K_1^1\right)
Q_2(\eta_0)$,
where
\begin{align*}
  II=& \chi_{0,R}(z_1)+
 \sum_{a,b=\pm}\tilde{\chi}_{1,R}^a(z_1)\widetilde{S}_a\widetilde{Q}_1(\eta_0)
       \widetilde{S}_b^{-1}\chi_{1,R}^b(z_1)
  \\ =& I-
\sum_{a,b=\pm}\tilde{\chi}_{1,R}^a(z_1)\widetilde{S}_a\widetilde{P}_1(\eta_0)
      \widetilde{S}_b^{-1}\chi_{1,R}^b(z_1)
\end{align*}
for sufficiently large $R$, and
\begin{multline*}
B_1^1= \tilde{\chi}_{0,R}(z_1)\tilde{\chi}_{0,R}(z_2)
(\lambda-\mL_0)^{-1}[\mL_2,\chi_{0,R}(z_1)\chi_{0,R}(z_2)]
\\  +
\sum_{a=\pm}\tilde{\chi}_{0,R}(z_1)\tilde{\chi}_{1,R}^a(z_2)
(\lambda-\mL_1^a)^{-1}[\mL_2,\chi_{0,R}(z_1)\chi_{1,R}^a(z_2)]
\\ +\sum_{a,b=\pm}\tilde{\chi}_{1,R}^a(z_1)\widetilde{S}_a(\lambda-\wmL_1)^{-1}
\widetilde{Q}_1(\eta_0)(\widetilde{S}_b)^{-1}[\mL_2,\chi_{1,R}^b(z_1)]\,,
\end{multline*}
\begin{align*}  
  B_1^2=& \frac{3}{2}\sum_{a=\pm}\tilde{\chi}_{0,R}(z_1)\tilde{\chi}_{1,R}^a(z_2)
(\lambda-\mL_1^a)^{-1}\pd_x\{(u_2-u_1^a)\chi_{0,R}(z_1)\chi_{1,R}^a(z_2)\cdot\}
\\ &+\frac{3}{2}\sum_{a,b=\pm}\tilde{\chi}_{1,R}^a(z_1)\widetilde{S}_a(\lambda-\wmL_1)^{-1}
     \widetilde{Q}_1(\eta_0)(\widetilde{S}_b)^{-1}
     \pd_x\{\chi_{1,R}^b(z_1)(u_2-\tu_1^b)\cdot\}\,,
\end{align*}
\begin{align*}  
  K_1^1=& \frac{3}{2}\tilde{\chi}_{0,R}(z_1)\tilde{\chi}_{0,R}(z_2)
(\lambda-\mL_0)^{-1}\pd_x\{\chi_{0,R}(z_1)\chi_{0,R}(z_2)u_2\cdot\}\,.
\end{align*}
Since $\chi_{0,R}$ and $\tilde{\chi}_{0,R}$ have compact supports, it follows
from \eqref{eq:com-1} that $K^1_1$ is compact.
By \eqref{eq:u2-asymptotics},
\begin{equation}
  \label{eq:u2-u1,tu1}
u_2-u_1^\pm=O(e^{-2\sqrt{c_{14}}|z_2|})\quad\text{and}\quad
u_2-\tu_1^\pm=O(e^{-2\sqrt{c_{23}}|z_1|})\text{ as $z_1\to\pm\infty$.}
\end{equation}
It follows from \eqref{eq:u2-u1,tu1} and \eqref{eq:com-1}--\eqref{eq:com-2''}
that
\begin{equation}
  \label{eq:27a}
\|B_1^1\|_{B(\calX_2)}+\|B_1^2\|_{B(\calX_2)}=o(1)\quad\text{as $R\to\infty$.}
\end{equation}
Let $B_1=B_1^1+B_1^2+B_1^3$, $K_1=K_1^1-K_1^2$ and
\begin{gather*}
B_1^3=  \sum_{a,b=\pm} \tilde{\chi}_{1,R}^a(z_1)
\{P_2(\eta_0)-\widetilde{S}_a\widetilde{P}_1(\eta_0)\widetilde{S}_b^{-1}\}
\chi_{1,R}^b(z_1)\,,\\
K_1^2=(1-\tilde{\chi}_{1,R}(z_1))P_2(\eta_0)\chi_{0,R}(z_1)\,. 
\end{gather*}
Since $P_2(\eta_0)Q_2(\eta_0)=Q_2(\eta_0)P_2(\eta_0)=O$,
$$Q_2(\eta_0)\sum_{a,b=\pm}\tilde{\chi}_{1,R}^a(z_1)P_2(\eta_0)
\chi_{1,R}^b(z_1)Q_2(\eta_0)
=Q_2(\eta_0)K_1^2Q_2(\eta_0)\,,$$
and
\begin{equation}
  \label{eq:27b}
W_2(\lambda)(\lambda-\mL_2)= Q_2(\eta_0)\left(I+B_1+K_1\right)
Q_2(\eta_0)\,.  
\end{equation}
By the definition,
\begin{gather*}
 e^{\a z_2}K_1^2f(\bx)=\frac{1}{2\pi}\sum_{k=1,2}\int_{|\eta|\le \eta_0}
    \mathcal{K}_k(\tbx,\tbxd,\eta)e^{\a z_2'}f(\tbxd)\,d\tbxd d\eta\,,
    \\
\mathcal{K}_k(\tbx,\tbxd,\eta)=e^{\a(z_2-z_2')}
(1-\tilde{\chi}_{1,R}(z_1))g^2_{2,k}(\tbx,\eta)\chi_{0,R}(z_1')g^{2,*}_{2,k}(\tbxd,\eta)\,.
\end{gather*}
By Lemma~\ref{lem:gg*1-asymp},
\begin{equation}
  \label{eq:27c}
\|B_1^3\|_{B(\calX_2)}\lesssim \eta_{0,1}+e^{-2\sqrt{c_{23}}R}\,,
\end{equation}
and $K_1^2$ is a Hilbert Schmidt operator.
Since $I+B_1$ is invertible and $K_1$ is compact,
Lemma~\ref{lem:anal-F} follows from the analytic Fredholm theorem
(\cite[Theorem~VI.14]{RS1}).
\par
If $\lambda\in\rho(\mL_1)$, then
$\lambda-\mL_2$ is invertible on $P_2(\eta_0)\calX_2$
and Lemma~\ref{lem:spec1} follows from the above argument.
Thus we complete the proof.
\end{proof}

To prove nonexistence of unstable eigenvalues, 
we will show that unstable eigenfunctions should be strongly localized
in $z_2$ provided $\lambda\not\in\sigma(\mL_1)$.
\begin{lemma}
  \label{lem:imev}
Assume that $\a\in(0,2\sqrt{c_{14}})$.
Let $\mL_0$, $\mL_1$ and $\mL_2$ be closed operators on $\calX_2$.
If $\lambda\in\rho(\mL_0)\cap\rho(\mL_1)$ and $\mL_2u=\lambda u$ for 
$u\in\calX_2$, then there exists an $\a'>0$ such that $e^{\a'|z_2|}u\in \calX_2$.
\end{lemma}

\begin{proof}[Proof of Lemma~\ref{lem:imev}]
Let 
  \begin{align*}
W_1(\lambda)=& \chi_{0,R}(z_2)(\lambda-\mL_0)^{-1}\tilde{\chi}_{0,R}(z_2)
+\sum_\pm \chi_{1,R}^\pm(z_2)(\lambda-\mL_1^\pm)^{-1}    
\tilde{\chi}_{1,R}^\pm(z_2)\,.
  \end{align*}
Then
\begin{equation}
  \label{eq:imev-6}
W_1(\lambda)(\lambda-\mL_2)=I+B_2+V_1\,,
\end{equation}
where 
\begin{align*}
  B_2=\chi_{0,R}(z_2)(\lambda-\mL_0)^{-1}[\mL_2,\tilde{\chi}_{0,R}(z_2)]
  + \sum_\pm \chi_{1,R}^\pm(z_2) (\lambda-\mL_1^\pm)^{-1}
  [\mL_2, \tilde{\chi}_{1,R}^\pm(z_2)]\,,
\end{align*}
\begin{align*}
  V_1=&
\frac32\chi_{0,R}(z_2)(\lambda-\mL_0)^{-1}\pd_x(u_2\tilde{\chi}_{0,R}(z_2)\cdot)
\\ &  +      
\frac32\sum_\pm\chi_{1,R}^\pm(z_2)(\lambda-\mL_1^\pm)^{-1}
\pd_x\{(u_2-u_1^\pm)\tilde{\chi}_{1,R}^\pm(z_2)\cdot\} \,.
\end{align*}
Using Lemma~\ref{lem:imev-5} and following the proof of Lemmas~\ref{cl:com-L0} and \ref{cl:com-L1},
we see that there exists an $\a'>0$ such that for $i=1$ and $2$
and $\beta\in[-\a',\a']$,
\begin{equation}
  \label{eq:imev-1}
\|e^{\beta z_2}B_2e^{-\beta z_2}\|_{B(\calX_2)}<1\,,
\end{equation}
provided $R$ is sufficiently large. Moreover, 
\begin{equation}
  \label{eq:imev-2}
e^{\a'|z_2|}V_1\in B(\calX_2)\quad\text{for  an $\a'>0$,}
\end{equation}
since $\chi_{0,R}(z_2)$ has a compact support and 
$u_2-u_1^\pm=O(e^{-2\sqrt{c_{14}}|z_2|})$.
\par

Suppose that $\mL_2u=\lambda u$ for $u\in\calX_2$. Then
\begin{equation}
\label{eq:imev-3}
  W_1(\lambda)(\lambda-\mL_2)u= (I+B_2)u+V_1u=0\,,
\end{equation}
and it follows from \eqref{eq:imev-1}--\eqref{eq:imev-3} that
for $\beta\in[-\a',\a']$,
\begin{align*}
\|e^{\beta z_2}u\|_{\calX_2}=&
\|e^{\beta z_2}(I+B_2)^{-1}V_1u\|_{\calX_2}
 \lesssim  \|u\|_{\calX_2}\,.
\end{align*}
This completes the proof of Lemma~\ref{lem:imev}.
\end{proof}

\begin{corollary}
  \label{cor:imev}
  Assume that $\a\in(0,2\sqrt{c_{14}})$.  Let $\mL_0$, $\mL_1$ and
  $\mL_2$ be closed operators on $\calX_2$.  If $\mL_2^*u=\lambda u$
  for $\lambda\in\rho(\mL_0^*)\cap\rho(\mL_1^*)$ and $u\in\calX_2^*$,
  then there exists an $\a'>0$ such that $e^{\a'|z_2|}u\in \calX_2^*$.
\end{corollary}
We can prove Corollary~\ref{cor:imev} in exactly the same way as
Lemma~\ref{lem:imev}.
\par

Darboux transformations $\nabla M_\pm(v_2)$ tell us that if $\lambda$ is a unstable eigenvalue of $\mL_2$,
then $\lambda\in\sigma_p(\mL_1)$,
which contradicts to Corollary~\ref{lem:mL-1inX2}.
The first step is to show that $\nabla M_+(v_2)$ is injective.
\begin{lemma}
  \label{lem:functional}
Assume that $\a\in(0,\sqrt{c_{14}})$.
\begin{enumerate}
\item If $\varphi\in\calX_2$ and $\nabla M_+(v_2)\varphi=0$, then $\varphi=0$.
\item If $\varphi\in\calX_2^*$ and $\nabla M_-(v_2)^*\varphi=0$, then
  $\varphi=0$.
\end{enumerate}
\end{lemma}

\begin{proof}
Let $\varphi(x,y)=\pd_x\tilde{\varphi}(x,y)\in\calX_2$. Then
\begin{equation}
  \label{eq:tphi}
  \pd_y\tilde{\varphi}+\pd_x^2\tilde{\varphi}=2v_2\pd_x\tilde{\varphi}\,,
\end{equation}
Note that $\tilde{\varphi}\in \calX_2$ since $\pd_x$ is invertible on $\calX_2$.
\par
Let $a(x,y)=-\Phi^2_1(x,y,0)\Phi^{2,*}_1(x,y,0)$. Then $a(x,y)$ is positive and
    \begin{equation}
\label{eq:a-asymp}
a(x,y)=O\left(\sech^2\frac{\theta_4-\theta_1}{2}\right)\,.
\end{equation}
Since $L_2\Phi^2_1=0$, it follows from \eqref{eq:Miura-Lax2}  that 
\begin{equation}
\label{eq:a}
\pd_ya-\pd_x^2a-2\pd_x(av_2)=\pd_x\nabla M_-(v_2)a=0\,.
\end{equation}
Combining \eqref{eq:tphi} and \eqref{eq:a}, we have
$\pd_y(a\tilde{\varphi}^2)=\pd_x\{(\pd_xa)\tilde{\varphi}^2
  -2a\tilde{\varphi}\varphi+2av_2\tilde{\varphi}^2\}+2a\varphi^2$
and
\begin{align}
  \label{eq:aphi^2}
  \frac{d}{dy}\int_\R a(x,y)\tilde{\varphi}(x,y)^2\,dx
=  2\int_\R a\varphi^2\,dx\,.
\end{align}
By \eqref{eq:a-asymp}, we have $a\tilde{\varphi}^2\in L^1(\R^2)$.
Since $\pd_x^i\pd_y^j\tilde{\varphi}\in\calX_2$ for every $i$, $j\ge0$
by Lemma~\ref{cl:m0-bound}, we see that \eqref{eq:aphi^2} is uniformly
bounded for $y\in\R$ and that $\int_R a\tilde{\varphi}^2\,dx$ tends to
$0$ as $|y|\to\infty$.  Integrating \eqref{eq:aphi^2} over $\R$ with
respect to $y$, we have $ \int_{\R^2}a(x,y)\varphi(x,y)^2\,dxdy=0$.
Since $a(x,y)$ is positive, $\varphi=\pd_x\tilde{\varphi}=0$ on
$\R^2$.  We can prove the latter part of Lemma~\ref{lem:functional} in
exactly the same way.  This completes the proof of
Lemma~\ref{lem:functional}.
\end{proof}

If $u$ is a solution of $\nabla M_+(v_2)u=f$ and $f$ satisfies
a secular term condition for the linearized KP-II equation around $u_2$,
then $u$ satisfies a secular term condition for the linearized
modified KP-II equation around $v_2$.
Let $\calX_2(\eta_0)=\{u\in\calX_2\mid P_2(\eta_0)u=0\}$ and
\begin{gather*}
\mathcal{Y}_2(\eta_0)=\{v\in \mathcal{Y}_2\mid 
\la v, g_M^{2,*}(\cdot,\eta)\ra=0
\quad\text{for $\eta\in[-\eta_0,\eta_0]$}\}\,.
\end{gather*}
\begin{lemma}
  \label{lem:M+orth}
  Let $\eta_0\in(0,\eta_{*,2})$, 
  $f\in\calX_2(\eta_0)$ and $u\in\calX_2$. If
  $\nabla M_+(v_2)u=f$, then $u\in\mathcal{Y}_2(\eta_0)$.
\end{lemma}
\begin{proof}
We have $u\in\mathcal{Y}_2$ from Lemma~\ref{cl:m0-bound}.
By \eqref{eq:tau2-14} and \eqref{eq:defg2M*},
\begin{align*}
g^{2,*}_M(\tbx,\eta)=e^{\gamma_{14}(-\eta)z_2+iy\eta}\sech Z_{2,\pm} \left\{
  C_\pm(\eta)+O(e^{-2\sqrt{c_{23}}|z_1|})\right\}
  \quad\text{as $(a_{23}-a_{14})y\to\pm\infty$,}
\end{align*}
where $C_\pm(\eta)$ are smooth functions of $\eta$.
Note that $e^{\gamma_{14}(-\eta)z_2}\sech Z_{2,\pm}$
are bounded in $L^2(\R;e^{-2\a z_2}dz_2)$ for $\eta\in[-\eta_0,\eta_0]$ and that
$e^{\gamma_{14}(-\eta)z_2}\sech Z_{2,\pm}e^{-2\sqrt{c_{23}}|z_1|}\in \calX_2^*$.
By \eqref{eq:12-1},
\begin{align*}
  \la u, g^{2,*}_M\ra =& \frac{1}{2}\lim_{R\to\infty}
\la \chi_{0,R}(y)u, \nabla M_+(v_2)^*g^{2,*}_{2,1}\ra
  \\=& \frac{1}{2}\lim_{R\to\infty}\la \chi_{0,R}(y)\nabla M_+(v_2)u,
       g^{2,*}_{2,1}\ra
       = \frac12\la f,g^{2,*}_{2,1}\ra=0\,.
\end{align*}
\end{proof}

Now we are in position to prove that there is no point spectrum
in a stable half plane.
\begin{lemma}
  \label{lem:Liouville}
Assume \eqref{eq:ass-alpha}. Let $\mL_1$ and $\mL_2$ be closed operators on $\calX_2$.
If $\mL_2\varphi=\lambda\varphi$ for $\lambda\in\rho(\mL_1)$ and $\varphi\in\calX_2$, then $\varphi=0$.
\end{lemma}
To prove Lemma~\ref{lem:Liouville}, we need the following.
\begin{claim}
\label{cl:Liouville-1}
Let $\eps>0$ and let $\eta_1$ and $\eta_2$ be positive numbers satisfying
$0<\eta_1<\eta_{*,2}<\eta_2$ and $\Re\gamma_{14}(\eta_2)<\sqrt{c_{14}}+\a+\eps$.
Let 
\begin{equation*}  
g^{2,*}_{Mr}(\tbx,\eta)
= g^{2,*}_M(\tbx,\eta)e^{-iy\eta}-g^{2,*}_M(\tbx,\pm\eta_1)e^{\mp iy\eta_1}
\quad \text{if $\pm\eta\in [\eta_1,\eta_2]$,}
\end{equation*}
  $$(Sa)(\eta)=a(\eta)-\frac{1}{2\pi}\int_{\R^2}
\left(\int_{I_{\eta_1,\eta_2}}a(\eta')g^{2,*}_{Mr}(\tbx,\eta')e^{iy\eta'}\,d\eta'\right)
\overline{g^2_M(\tbx,\eta)}\,d\tbx\,.$$
If $\eta_2-\eta_1$ is sufficiently small, then
$S:L^2(-\eta_2,\eta_2)\to L^2(-\eta_2,\eta_2)$ has a bounded inverse.
\end{claim}
\begin{proof}[Proof of Claim~\ref{cl:Liouville-1}]
Let $a(\eta)=0$ outside of $[-\eta_2,\eta_2]$ and let
$g^{2,*}_{Mr}(\tbx,\eta)=0$ if $\left||\eta|-\eta_2\right|\ge \eta_2-\eta_1$ and
\begin{equation*}
g^{2,*}_{Mr}(\tbx,\eta)=g^{2,*}_{Mr}(\tbx,\pm\eta_2)
\frac{2\eta_2-\eta_1\mp\eta}{\eta_2-\eta_1}
\quad\text{for  $\pm\eta\in[\eta_2,2\eta_2-\eta_1]$.}
\end{equation*}
Let $S_1=I-S$.
By Parseval's identity,
\begin{equation}
  \label{eq:clLpf1}
(S_1a)(\eta)=\frac{1}{2\pi}\int_{\R^2}\overline{g^2_M(\tbx,\eta)}
\left(\int_\R \mathcal{F}_\eta^{-1}g^{2,*}_{Mr}(\tbx,y-y_1)
\mathcal{F}^{-1}_\eta a(y_1)\,dy_1\right)\,d\tbx\,.
\end{equation}
By  \eqref{eq:btitj},
\begin{align*}
& g^{2,*}_{Mr}(\tbx,\eta)
=\sum_\pm b_{1,\pm}(z_2,y)b_{2,\pm}(z_2,\eta)\,,\\
& b_{1,\pm}(z_2,y)=\frac{1}{2(\k_4-\k_1)}\Phi^2_1(\tbx)e^{-(\theta_1+\theta_4)/2}
e^{\pm\sqrt{c_{23}}z_1}\sech\sqrt{c_{23}}z_1\,,
\\ &
b_{2,\pm}(z_2,\eta)=0
     \quad\text{if $|\eta|\le \eta_1$ or $|\eta|\ge 2\eta_2-\eta_1$,}      
\\ &
b_{2,\pm}(z_2,\eta)=\frac{e^{\gamma_{14}(-\eta)z_2}}
{\beta_{14}^-(-\eta)-a_{23}\mp\sqrt{c_{23}}}
-\frac{e^{\gamma_{14}(-\eta_1\sgn\eta)z_2}}{\beta_{14}(-\eta_1\sgn\eta)-a_{23}\mp\sqrt{c_{23}}}
     \quad\text{if $\eta_1\le |\eta|\le\eta_2$,}
  \\ & b_{2,\pm}(z_2,\eta)=
       b_{2,\pm}(z_2,\eta_2\sgn\eta)\frac{2\eta_2-\eta_1-|\eta|}{\eta_2-\eta_1}
\quad\text{if $\eta_2\le |\eta|\le 2\eta_2-\eta_1$.}
\end{align*}
It follows from \eqref{eq:tau2-14} and \eqref{eq:defg2M*} that
$$\|b_{1,\pm}\|_{L^\infty_y}\lesssim \sech\sqrt{c_{14}}z_2\,,\quad
\|b_{2,\pm}\|_{H^1(\R_\eta)}
\lesssim (\eta_2-\eta_1)^{1/2}\la z_2\ra
 \max(e^{\sqrt{c_{14}}z_2},e^{\Re\gamma_{14}(\eta_2)z_2})\,,$$
and that
\begin{align*}
\|\mathcal{F}_\eta^{-1}g^{2,*}_{Mr}(z_2,y,z)\|_{L^1_zL^\infty_y}
\lesssim & 
\sum_\pm \|b_{1,\pm}(z_2,y)\|_{L^\infty_y}
\|b_{2,\pm}(z_2,\eta)\|_{H^1(\R_\eta)}
\\ \lesssim & 
(\eta_2-\eta_1)^{1/2}\la z_2\ra\sech\sqrt{c_{14}}z_2
\max(e^{\sqrt{c_{14}}z_2},e^{\Re\gamma_{14}(\eta_2)z_2})\,,
\end{align*}
\begin{equation}
  \label{eq:clLpf2}
\begin{split}
 \left\|e^{-\eps z_2}\int_\R \mathcal{F}_\eta^{-1}g^{2,*}_{Mr}(z_2,y,y-y_1)
\mathcal{F}^{-1}_\eta a(y_1)\,dy_1\right\|_{\calX_2^*}
 \lesssim &
(\eta_2-\eta_1)^{1/2}\|a\|_{L^2(-\eta_2,\eta_2)}\,.
\end{split}  
\end{equation}
Note that $\operatorname{supp}\mathcal{F}^{-1}a\subset [-\eta_2,\eta_2]$.

It follows from \eqref{eq:psi1,tpsi1}, \eqref{eq:psi1-asymp},
\eqref{eq:tpsi1-asymp} and \eqref{eq:defg2M} that, as $y\to\pm\infty$,
\begin{align}
  \label{eq:clLpf3}
g^2_M(\tbx,\eta)=e^{-\gamma_{14}(\eta)z_2+iy\eta}
\{C_\pm(\eta)\sech Z_{2,\pm} & (\gamma_{14}(\eta)+\sqrt{c_{14}}\tanh Z_{2,\pm})
  \\ & \notag
+O(e^{-(\sqrt{c_{14}}|z_2|+2\sqrt{c_{23}}|z_1|)})\}\,.
\end{align}
where $C_\pm(\eta)$ are smooth functions of $\eta$.
Combining \eqref{eq:clLpf1}--\eqref{eq:clLpf3}, we have
$$\|S_1a\|_{L^2(-\eta_2,\eta_2)}\lesssim (\eta_2-\eta_1)^{1/2}
\|a\|_{L^2(-\eta_2,\eta)}\,.$$
Thus we complete the proof.
\end{proof}

\begin{proof}[Proof of Lemma~\ref{lem:Liouville}]
By Lemma~\ref{lem:imev} and Proposition~\ref{lem:mLu=f}, there exists an $\eps>0$ such that
$$
\sum_{0\le i+2j\le2}\|e^{\eps |z_2|}\pd_x^i\pd_y^j\varphi\|_{\calX_2}
\lesssim \|\varphi\|_{\calX_2}\,.
$$
Let $\eta_{*,\eps}$, $\eta_{*,-\eps}$, $\eta_1$ and $\eta_2$ be
positive numbers satisfying 
\begin{gather}
  \label{eq:etapmeps}
\Re\gamma_{14}(\eta_{*,\eps})=\sqrt{c_{14}}+\a+\eps\,,\quad
\Re\gamma_{14}(\eta_{*,-\eps})=\sqrt{c_{14}}+\a-\eps\,,\\
  \label{eq:30}
\eta_{*,-\eps}<\eta_1<\eta_{*,2}<\eta_2<\eta_{*,\eps}\,.  
\end{gather}
Since $\mL_2^*g^{2,*}_{2,\pm}(\eta)=\lambda_{2,\pm}(-\eta)g^{2,*}_{2,\pm}(\eta)$,
we may assume
\begin{equation}
  \label{eq:21}
\la \varphi, g^{2,*}_{2,k}(\eta)\ra=0
\quad\text{for $k=1$, $2$ and $\eta\in(-\eta_{*,\eps},\eta_{*,\eps})$}  
\end{equation}
by taking $\eps$ smaller if necessary.
\par

Let $\varphi_M=T^{2,+}(0,\infty)\varphi$.
Since
\begin{gather*}
\sqrt{c_{14}}+\a-\eps<\Re\gamma_{14}(\eta_1)<\sqrt{c_{14}}+\a<
\Re\gamma_{14}(\eta_2)<\sqrt{c_{14}}+\a+\eps\,,
\end{gather*}
it follows from Lemma~\ref{lem:T2pm} that
\begin{gather*}
\| e^{-\eps z_2}T^{2,+}(\eta_1,\eta_2)\varphi\|_{\calX_2}
\lesssim \|e^{-\eps z_2}\varphi\|_{\calX_2}\,,
\quad
\|T^{2,+}(\eta_2,\infty)\varphi\|_{\calX_2}
\lesssim \|\varphi\|_{\calX_2}\,.
\end{gather*}
Moreover, it follows from Lemma~\ref{lem:T2lowpm} and \eqref{eq:21} that
\begin{gather*}
\|T^{2,+}_{low}(0,\eta_1)\varphi\|_{\calX_2}
\lesssim \|\varphi\|_{\calX_2}\,,
\quad
\|e^{\eps z_2}T^{2,+}_{low}(\eta_1,\eta_2)\varphi\|_{\calX_2}
\lesssim \|e^{\eps z_2}\varphi\|_{\calX_2}\,,
\\
T^{2,+}(0,\eta_1)\varphi=T^{2,+}_{low}(0,\eta_1)\varphi\,,\quad
T^{2,+}(\eta_1,\eta_2)\varphi=T^{2,+}_{low}(\eta_1,\eta_2)\varphi\,.
\end{gather*}
Thus we have
$\varphi_M=T^{2,+}(0,\eta_1)\varphi+T^{2,+}(\eta_1,\eta_2)\varphi+T^{2,+}(\eta_2,\infty)\varphi
\in \calX_2$. By Lemmas~~\ref{lem:20} and \ref{cl:m0-bound},
$\nabla M_+(v_2)\varphi_M=\varphi$ and $\pd_x^i\pd_y^j\varphi_M\in\calX_2$ for any $i\in\Z$ and $j\ge0$.
\par
Let
$$
\mL_M=\mL_0+\frac32\pd_x(v_2^2\cdot)-
\frac{3}{2}(\pd_xv_2\pd_x^{-1}\pd_y+\pd_x^{-1}\pd_yv_2\pd_x)\,,$$
where $\pd_x^{-1}\pd_yv_2=\pd_y\tau_2/\tau_2-\pd_y\tau_1/\tau_1$.
By \eqref{eq:bH+} and the assumption,
\begin{align*}
\nabla M_+(v_2)(\lambda-\mL_M)\varphi_M=& (\lambda-\mL_2)\nabla M_+(v_2)\varphi_M=
(\lambda-\mL_2)\varphi=0\,,
\end{align*}
and it follows from Lemma~\ref{lem:functional} that
$\mL_M\varphi_M=\lambda \varphi_M$.
Thus by \eqref{eq:bH-},
\begin{align*}
  (\lambda-\mL_1)\nabla M_-(v_2)\varphi_M=\nabla M_-(v_2)(\lambda-\mL_M)\varphi_M=0\,.
  \quad 
\end{align*}
Since $\lambda\in\rho(\mL_1)$, it follows that $\nabla M_-(v_2)\varphi_M=0$ and that
$$\varphi=\nabla M_+(v_2)\varphi_M-\nabla M_-(v_2)\varphi_M=2\pd_x\varphi_M\,.$$
Using \eqref{eq:imev-6} and the fact that $(\lambda-\mL_2)\pd_x\varphi_M=0$, we have
\begin{align*}
  (I+B_2+\pd_x^{-1}[B_2,\pd_x])\varphi_M+(V_1+\pd_x^{-1}[V_1,\pd_x])\varphi_M=0\,,
\end{align*}
and we can prove $\|e^{\eps|z_2|}\varphi_M\|_{\calX_2}\lesssim \|\varphi_M\|_{\calX_2}$
following the proof of Lemma~\ref{lem:imev}.
By Lemma~\ref{lem:M+orth} and \eqref{eq:21},
\begin{align}
\label{eq:orthphi}
\int_{\R^2} \varphi_M(\tbx)
\overline{g^{2,*}_M(\tbx,\eta)}\,d\tbx=0
\quad\text{for $\eta\in(-\eta_{*,\eps},\eta_{*,\eps})$.}  
\end{align}
Let $b(\eta)= \int_{\R^2} e^{2\a z_2}\varphi_M(\tbx)\overline{g^2_M(\tbx,\eta)}\,d\tbx$.
By \eqref{eq:clLpf3},
$$\|b\|_{L^2(-\eta_2,\eta_2)}\lesssim \|e^{-\eps z_2}\varphi_M\|_{\calX_2}\lesssim\|\varphi_M\|_{\calX_2}\,.$$
By Claim~\ref{cl:Liouville-1}, there exists
an $a\in L^2(-\eta_2,\eta_2)$ satisfying $Sa=b$ and $\|a\|_{L^2}=O(\|b\|_{L^2})$.
Combining the above with \eqref{eq:orthgMgM*}, we see that
\begin{gather*}
\varphi^*_M:=e^{2\a z_2}\varphi_M-\frac{1}{2\pi}\int_{-\eta_2}^{\eta_2}
a(\eta)\left(g^{2,*}_M(\tbx,\eta)-g^{2,*}_{Mr}(\tbx,\eta)e^{iy\eta}\right)\,d\eta
\end{gather*}
satisfies
\begin{equation}
  \label{eq:orthphi*}
\lim_{R\to\infty}\int_{-R}^R\int_\R g^2_M(\tbx,\eta)
\overline{\varphi_M^*(\tbx)}\,dz_2dy=0
\quad\text{for $\eta\in[-\eta_2,\eta_2]$.}
\end{equation}
By \eqref{eq:orthphi} and \eqref{eq:clLpf2},
\begin{align*}
\int_{\R^2}\varphi_M^*(\tbx)\overline{\varphi_M(\tbx)}\,d\tbx=&
\|\varphi_M\|_{\calX_2}^2+\frac{1}{2\pi}\int_{I_{\eta_1,\eta_2}}
a(\eta)\left(\int_{\R^2}\overline{\varphi_M(\tbx)}g^{2,*}_{Mr}(\tbx,\eta)
e^{iy\eta}\,d\tbx\right)
  \\=& \|\varphi_M\|_{\calX_2}^2\left(1+O\left((\eta_2-\eta_1)^{1/2}\right)\right)
\lesssim \|\varphi_M\|_{\calX_2}^2\,,
\end{align*}
provided $\eta_2-\eta_1$ is small.
Note that $e^{-\eps z_2}\varphi_M^*\in\calX_2^*$ and $e^{\eps' z_2}\varphi_M^*\in\calX_2^*$
for $\eps'$ satisfying $0<\eps'<\min\{\a+\sqrt{c_{14}}
-\Re\gamma_{14}(\eta_1),\eps\}$.
\par
Suppose that $\eta_1'\in(\eta_1,\eta_{*,2})$ and $\Re\gamma_{14}(\eta_1')>\sqrt{c_{14}}+\a-\eps'$.
It follows from Lemma~\ref{lem:T2-*low} and  \eqref{eq:orthphi*} that
\begin{gather*}
\| T^{2,-,*}(0,\eta_1')\varphi_M^*\|_{\calX_2^*}
\lesssim \|\varphi_M^*\|_{\calX_2^*}\,,\quad
\|e^{-\eps z_2}T^{2,-,*}(\eta_1',\eta_2)\varphi_M^*\|_{\calX_2^*}
\lesssim \|e^{-\eps z_2}\varphi_M^*\|_{\calX_2^*}\,.
\end{gather*}
By Lemma~\ref{lem:T2-*},
\begin{gather*}
  \| e^{\eps' z_2}T^{2,-,*}(\eta_1',\eta_2)\varphi_M^*\|_{\calX_2^*}
\lesssim \|e^{\eps' z_2}\varphi_M^*\|_{\calX_2^*}\,,\quad
\|T^{2,-,*}(\eta_2,\infty)\varphi_M^*\|_{\calX_2^*}
\lesssim \|\varphi_M^*\|_{\calX_2^*}\,.
\end{gather*}
Combining the above, we have
$$v:=
T^{2,-,*}(0,\eta_1')\varphi_M^*+T^{2,-,*}(\eta_1',\eta_2)\varphi_M^*+
T^{2,-,*}(\eta_2,\infty)\varphi_M^* \in\mathcal{Y}_2^*\,,$$
and  $\nabla M_-(v_2)^*v=\varphi_M^*$ by Lemma~\ref{lem:21}.
 Since $\nabla M_-(v_2)^*v=\varphi_M^*$ and $\nabla M_-(v_2)\varphi_M=0$,
\begin{align*}
\la \varphi_M,\varphi_M^*\ra=&  \la \nabla M_-(v_2)\varphi,v\ra=0\,,
\end{align*}
and $\varphi=\pd_x\varphi_M=0$ provided $\eta_2-\eta_1$ is sufficiently small.
 Thus we complete the proof.
\end{proof}

Now we are in position to prove Proposition~\ref{prop:linearstability-1}.
\begin{proof}[Proof of Proposition~\ref{prop:linearstability-1}]
  By Lemmas~\ref{lem:spec1}, \ref{lem:anal-F1} and \ref{lem:Liouville}
  and Corollary~\ref{lem:mL-1inX2}, there exists a $b'>0$ such that
  $\lambda\in \rho(\mL_2|_{\calX_2(\eta_0)})$ if $\Re\lambda>-b'$.
  Thus by the Gearhart-Pr\"{u}ss theorem (\cite{Gh,Pruss}), we have
  Proposition~\ref{prop:linearstability-1}.  Thus we complete the
  proof.
\end{proof}
\bigskip

\subsection{Linear stability of $2$-line solitons of P-type in  $\calX_1$}
\label{sec:LSX1}
 Lemma~\ref{cl:sgnpdlambda1} tells us that a $2$-line soliton solution
$u_2$ is spectrally unstable in a weighted space $\calX_1$.
\begin{proposition}
  \label{prop:instabilityX1}
Suppose that $0<\a<\min\{\sqrt{c_{14}},\,\frac{c_{14}}{3|a_{14}-a_{23}|}\}$
and that $\mL_2$ be a closed operator on $\calX_1$. Then
$\sigma(\mL_2)\cap \{\lambda\in\C\mid \Re\lambda>0\}\ne\emptyset$.
\end{proposition}
\begin{proof}
  Let $\eta_I=2(a_{23}-a_{14})\a$ and $\eta_R\in(-\eta_\pm',\eta_\pm')$.
Then $|\Re\gamma_{14}(\eta_R\pm i\eta_I)-\a|<\sqrt{c_{14}}$,
\begin{gather*}
\sup_{y\in\R}\|e^{\a z_1}g^2_{2,\pm}(\tbx,\eta)\|_{L^2(\R_x)}
  \lesssim \left\|e^{\a z_2}e^{-\gamma_{14}(\pm\eta)z_2}\sech\sqrt{c_{14}}z_2
  \right\|_{L^2(\R)}<\infty\,,
  \end{gather*}
and  $\lambda_{2,\pm}(\eta_R+i\eta_I)\in \sigma(\mL_2)$.
By Lemma~\ref{lem:lambda2-etaI}, we have $\Re\lambda_{2,\pm}(\eta)>0$
if $\pm(a_{14}-a_{23})>0$ and $|\eta_R|<\tilde{\eta}_{*,\pm}$.
This completes the proof of Proposition~\ref{prop:instabilityX1}.
\end{proof}
Next, we will show spectral stability of $\mL_2$ on a subspace of $\calX_1$
that does not contain continuous eigenfunctions corresponding to
modulations of $[1,4]$-soliton.

\begin{lemma}
  \label{lem:wT2pm}
Assume \eqref{ass-alpha'} and that $\eps$ is sufficiently close
  to $0$. 
  \begin{enumerate}
  \item If $\eta_2\ge \eta_1\ge\eta_0>\eta_+'(\a,\a+\eps)$, then
\begin{gather*}
  \|e^{\eps z_2}T^{2,-}(\eta_1,\eta_2,\eta_I(\a))f\|_{\calX_1}
  \le C\|e^{\eps z_2}f\|_{\calX_1}\,,
\\
\|e^{-\eps z_2}T^{2,-}(\eta_1,\eta_2,\eta_I(\a))^*f\|_{\calX_1^*}
\le C\|e^{-\eps z_2}f\|_{\calX_1^*}\,,
\end{gather*}
where $C$ is  a positive constant depending only on $\a$, $\eps$ and $\eta_0$.

\item  
If $\eta_2\ge \eta_1\ge\eta_0>\eta_-'(\a,\a+\eps)$, then
\begin{gather*}
  \|e^{\eps z_2}T^{2,+}(\eta_1,\eta_2,\eta_I(\a))f\|_{\calX_1}
  \le C\|e^{\eps z_2}f\|_{\calX_1}\,,
\\
\|e^{-\eps z_2}T^{2,+}(\eta_1,\eta_2,\eta_I(\a))^*f\|_{\calX_1^*}
\le C\|e^{-\eps z_2}f\|_{\calX_1^*}\,,
\end{gather*}
where $C$ is  a positive constant depending only on $\a$, $\eps$ and $\eta_0$.
  \end{enumerate}
\end{lemma}

\begin{lemma}
  \label{lem:wT2pm-low}
Assume \eqref{ass-alpha'} and that $\eps$ is sufficiently close to $0$.
  \begin{enumerate}
  \item If $0\le\eta_1\le \eta_2\le\eta_0<\eta_+'(\a,\a+\eps)$,
\begin{gather*}
 \|e^{\eps z_2}T^{2,-}_{low}(\eta_1,\eta_2,\eta_I(\a))f\|_{\calX_1}
  \le C\|e^{\eps z_2}f\|_{\calX_1}\,,
\\
\|e^{-\eps z_2}T^{2,-}_{low}(\eta_1,\eta_2,\eta_I(\a))^*f\|_{\calX_1^*}
\le C\|e^{-\eps z_2}f\|_{\calX_1^*}\,,
\end{gather*}
where $C$ is  a positive constant depending only on $\a$, $\eps$ and $\eta_0$.

\item  
If $0\le \eta_1\le \eta_2\le\eta_0<\eta_-'(\a,\a+\eps)$,
\begin{gather*}
 \|e^{\eps z_2}T^{2,+}_{low}(\eta_1,\eta_2,\eta_I(\a))f\|_{\calX_1}
  \le C\|e^{\eps z_2}f\|_{\calX_1}\,,
\\
\|e^{-\eps z_2}T^{2,+}_{low}(\eta_1,\eta_2,\eta_I(\a))^*f\|_{\calX_1^*}
\le C\|e^{-\eps z_2}f\|_{\calX_1^*}\,,
\end{gather*}
where $C$ is  a positive constant depending only on $\a$, $\eps$ and $\eta_0$.
  \end{enumerate}
\end{lemma}

\begin{proof}[Proof of Lemmas~\ref{lem:wT2pm} and \ref{lem:wT2pm-low}]
Lemma~\ref{lem:wT2pm} follows immediately from Lemma~\ref{lem:T2pm'}.
Since $\gamma_{14}(i\eta_I(\a))>\sqrt{c_{14}}-\a$ for $0<\a<2(\k_4-a_{23})$
and $\gamma_{14}(-i\eta_I(\a))>\sqrt{c_{14}}-\a$ for
$0<\a<2(a_{23}-\k_1)$, Lemmas~\ref{lem:wT2pm-low}
follows from Lemmas~\ref{lem:T2lowpm'}.  
\end{proof}

  \begin{lemma}
    \label{lem:imev2}
  Assume $\a\in(0,2\sqrt{c_{23}})$. Let $\mL_0$, $\mL_1$, $\wmL_1$ and
  $\mL_2$ be closed operators on $\calX_1$.
\begin{enumerate}
\item
  If $\lambda\in\sigma(\mL_2)\cap\rho(\mL_0)\cap\rho(\mL_1)\cap\rho(\wmL_1)$,
  then $\lambda$ is an isolated eigenvalue of $\mL_2$ with finite multiplicity.
\item
If $\lambda\in\rho(\mL_0)\cap\rho(\mL_1)$ and $\mL_2u=\lambda u$ for 
$u\in\calX_1$, then there exists an $\a'>0$ such that $e^{\a'|z_2|}u\in \calX_1$.
\end{enumerate}
\end{lemma}
\begin{proof}
  We can prove the former part in the same way as
  Lemma~\ref{lem:spec1} and the latter part in the same way as Lemma~\ref{lem:imev}.
\end{proof}

\begin{lemma}
  \label{lem:functional'}
  Assume  \eqref{ass-alpha'}. If $\nabla M_+(v_2)\varphi=0$ and
  $e^{\eps|z_2|}\varphi\in\calX_1$ for an $\eps>0$, then $\varphi=0$.
\end{lemma}
\begin{proof}
  Since $\nabla M_+(v_2)\varphi=0$ and $e^{\pm\eps z_2}\varphi\in\calX_1$,
  it follows from Lemma~\ref{cl:m0-bound} that
  $e^{\eps|z_2|}\pd_x\varphi$, $e^{\eps|z_2|}\pd_x^{-1}\pd_y\varphi\in\calX_1$.
\par
  Let $u=\left(T^{2,+}_{low}(0,\eta_-',\eta_I)^*+
    T^{2,+}(\eta_-',\infty,\eta_I)^*\right)  e^{2\a z_1}\varphi$.
  Since $\eta_-'(\a,\a-\eps)<\eta_-'<\eta_-'(\a,\a+\eps)$, it follows from
  Lemmas~\ref{lem:wT2pm} and \ref{lem:wT2pm-low} that
  \begin{equation*}
e^{-\eps z_2}T^{2,+}_{low}(0,\eta_-',\eta_I)^*e^{\eps z_2}\in B(\calX_1^*)\,,
\quad e^{\eps z_2}T^{2,+}(\eta_-',\infty,\eta_I)^*e^{-\eps z_2}\in B(\calX_1^*)\,,
\end{equation*}
and $e^{-\eps |z_2|}u\in \calX_1^*$.
By Lemma~\ref{lem:21}, we have 
$\nabla M_+(v_2)^*u=e^{2\a z_1}\varphi$ and
\begin{align*}
  \|\varphi\|_{\calX_1}^2=& \la \varphi,\nabla M_+(v_2)^*u\ra
= \la M_+(v_2)\varphi,u\ra=0\,.
\end{align*}
This completes the proof of Lemma~\ref{lem:functional'}.
\end{proof}

\begin{lemma}
  \label{lem:Liouville'}
  Assume \eqref{ass-alpha'}.
 Let $\mL_0$, $\mL_1$, $\wmL_1$ and $\mL_2$ be closed operators on $\calX_1$.
 If $\lambda\in\rho(\mL_0)\cap\rho(\mL_1)\cap\rho(\wmL_1)$, then
$\lambda\in \rho(\mL_2)$.
\end{lemma}
\begin{proof}
  Suppose on the contrary that $\lambda\in\sigma(\mL_2)$. Then
  it follows from Lemma~\ref{lem:imev2} that there exist
  $\varphi\in\calX_1$ and an $\eps>0$ such that
  $\mL_2\varphi=\lambda\varphi$ and that $e^{\eps|z_2|}\varphi\in\calX_1$.
By Proposition~\ref{lem:mLu=f}, we have
$\pd_x^j\pd_y^k\varphi\in  e^{-\eps|z_2|} \calX_1$ for every $j$, $k\ge0$.
Suppose that $\eps$, $\eps'$, $\eta_1$ and $\eta_2$ are positive numbers
satisfying $\eps>\eps'$ and 
$$\eta_-'(\a,\a-\eps)<\eta_1<\eta_-'(\a,\a-\eps')
<\eta_-'<\eta_-'(\a,\a+\eps')<\eta_2<\eta_-'(\a,\a+\eps)\,.$$
Since $\lambda\in\rho(\wmL_1)$, it follows from Proposition~\ref{prop:spec[1-4]}
that if $\eps$ is sufficiently small,
$$\la \varphi,g^{2,*}_{2,-}(\eta_R+i\eta_I)\ra=0
\quad\text{for $\eta_R$ satisfying $|\eta_R|<\eta_-'(\a,\a+\eps)$}\,.$$
Thus by using Lemmas~\ref{lem:wT2pm} and \ref{lem:wT2pm-low},
we can prove
  \begin{gather*}
 e^{\eps'|z_2|} T^{2,+}(0,\eta_1,\eta_I(\a))\varphi\,,\quad
e^{\eps'|z_2|} T^{2,+}(\eta_2,\infty,\eta_I(\a))\varphi\in \calX_1\,,\\
    e^{\eps|z_2|}T^{2,+}(\eta_1,\eta_2,\eta_I(\a))\varphi\in\calX_1\,,\\
 \varphi_M:=T^{2,+}(0,\infty,\eta_I(\a))\varphi\in e^{-\eps'|z_2|}\calX_1\,,
  \end{gather*}
in the same way as the proof of Lemma~\ref{lem:Liouville}.
By Lemmas~\ref{lem:20} and \ref{cl:m0-bound}, we have
$\nabla M_+(v_2)\varphi_M=\varphi$ and
$\pd_x^j\pd_y^k\varphi_M\in e^{-\eps'|z_2|}\calX_1$  for every $j\in\Z$
and $k\ge0$. Using an intertwining property \eqref{eq:bH+}, we have
$$\nabla M_+(v_2)(\lambda-\mL_M)\varphi_M=(\lambda-\mL_2)\varphi=0\,.$$
Hence it follows from Lemma~\ref{lem:functional'} that
$(\lambda-\mL_M)\varphi_M=0$. By \eqref{eq:bH-},
\begin{gather}
  \label{eq:L1ef}
(\lambda-\mL_1)\nabla M_-(v_2)\varphi_M
=\nabla M_-(v_2)(\lambda-\mL_M)\varphi_M=0\,.
\end{gather}
Since $\lambda\in\rho(\mL_1)$,
it follows from \eqref{eq:L1ef} that $\nabla M_-(v_2)\varphi_M=0$ and that
$\varphi_M=\frac12\pd_x^{-1}\varphi\in e^{-\eps|z_2|}\calX_1$.
Following the proof of Claim~\ref{cl:Liouville-1} and
Lemma~\ref{lem:Liouville}, we have  $\varphi=2\pd_x\varphi_M=0$.
Thus we complete the proof.
\end{proof}
Lemma~\ref{lem:Liouville'} tells us that
$\sigma(\mL_2)\subset \sigma(\mL_1)\cup\sigma(\wmL_1)$.
More precisely, we have the following.
\begin{proposition}
  \label{prop:specL2-X1}
 Assume \eqref{ass-alpha'}. 
Let $\mL_0$ and $\mL_2$ be closed operators on $\calX_1$ and
let $\eta_I=2(a_{23}-a_{14})\a$.
Then
\begin{align*}
\sigma(\mL_2)=\sigma(\mL_0)\cup &
\{\lambda_{1,\pm}(\eta)\mid \eta\in[-\eta_{*,1},\eta_{*,1}]\}
  \\ &
\cup\{\lambda_{2,\pm}(\eta_R+i\eta_I)\mid \eta\in[-\eta_\pm',\eta_\pm']\}\,.
\end{align*}
\end{proposition}
\begin{remark}
  \label{rem:6.5}
If $\lambda\in \sigma(\mL_2)\cap\rho(\mL_1)$,
we can show that the associated continuous eigenspace is spanned by
$\{g^2_{2,\pm}(\eta_R+i\eta_I)\}$ in the same way as
Lemmas~\ref{lem:anal-F1} and \ref{lem:Liouville} by using
Lemmas~\ref{lem:wT2pm} and \ref{lem:wT2pm-low}.
\end{remark}

Next, we will investigate continuous eigenspaces that belong to
$\lambda_{1,\pm}(\eta)$.
\begin{lemma}
  \label{lem:anal-F}
Assume \eqref{ass-alpha'} and \eqref{eq:l2ne0}.  
Let $\mL_2$ be a closed operator on $\calX_1$. Let
$\eta_{0,1}>0$ be a small number and
$\mathcal{C}=\{\lambda_{1,\pm}(\eta)\mid \eta\in(-\eta_{0,1},\eta_{0,1})\}$.
Then $\mathcal{C}\cap\sigma_p(\mL_2)$ is a discrete set and
$\lambda-\mL_2$ is invertible on $(I-P_{21}'(\eta_{0,1}))\calX_1$
if $\lambda\in\mathcal{C}\setminus\sigma_p(\mL_2)$.
\end{lemma}
We can prove Lemma~\ref{lem:anal-F} in exactly the same way as
Lemma~\ref{lem:anal-F1}.
By Lemma~\ref{lem:anal-F},
continuous eigenspaces associated with $\lambda_{1,\pm}(\eta)$ are
generically spanned by $\{g^2_{1,k}(\eta)\}_{k=1,2}$ if $\eta\simeq0$,
and in the exceptional cases, $\lambda_{1,\pm}(\eta)$ have corresponding eigenfunctions.
We assume smallness of $\eta$ so that $\lambda_{1,\pm}(\eta)\not\in
\{\lambda_{2,\pm}(\eta_R+i\eta_I)\mid |\eta_R|\le \eta_\pm'\}$.

The argument of the proof of Lemma~\ref{lem:Liouville} fails
because it requires that $\lambda\in\rho(\mL_1)$.
Instead of showing the completeness of every continuous eigenspace
belonging to $\lambda_{1,\pm}(\eta)$, we will prove that 
there are no eigenvalues of $\mL_2$ near the imaginary axis.
For the purpose, we will use the Darboux transformation
$\nabla \widetilde{M}_\pm$ as well as $\nabla M_\pm$ and prove
$\lambda=0$ is not an eigenvalue of $\mL_2$.
\par
Suppose that $\mL_2U_2=0$. Then in a moving coordinate \eqref{eq:b1-b2P},
\begin{equation}
  \label{eq:0ef-u2}
3\pd_Y^2U_2-4b_2\pd_X\pd_YU_2+\pd_X^2(\pd_X^2-4b_1+6u_2)U_2=0\,.
\end{equation}
\begin{proposition}
  \label{prop:0notev}
 Let $0<\a<\min\{\sqrt{c_{23}}, \left||b_2-3a_{14}|-\sqrt{c_{14}}\right|\}$.
If $U_2\in \calX_1$ is a solution of \eqref{eq:0ef-u2}, then $U_2=0$.
\end{proposition}
To prove Proposition~\ref{prop:0notev}, we will use the Darboux transformations
\eqref{eq:Miura-Lax1}--\eqref{eq:Miura-Laxt2}.
If we linearize \eqref{eq:Miura1} and \eqref{eq:Miura-t1} around $v=v_2$,
\begin{gather}
  \label{eq:LMiura-1}
  U_1=\nabla M_-(v_2)V\,,\quad U_2=\nabla M_+(v_2)V\,,
\\
\label{eq:LMiura-t1}
3(\pd_x+\pd_x^{-1}\pd_y+2v_2)U_1+\nabla_v\widetilde{M}_-(u_1,v_2)V=0\,,
\\
\label{eq:LMiura-t2}
3(-\pd_x+\pd_x^{-1}\pd_y+2v_2)U_2+\nabla_v\widetilde{M}_+(u_2,v_2)V=0\,,
\end{gather}
$$
\nabla_v\widetilde{M}_-(u_1,v_2)V=
4\left\{\pd_X^2+3v_2\pd_X-b_1+3v_2^2+\frac{3}{2}u_2\right\}V-4b_2\pd_X^{-1}\pd_YV\,,
$$
and formally, $\nabla_v\widetilde{M}_-(u_1,v_2)^*=\nabla_v\widetilde{M}_+(u_2,v_2)$.
\par
The Darboux transformations \eqref{eq:LMiura-1}--\eqref{eq:LMiura-t2}
map a solution of \eqref{eq:0ef-u2} to a solution of
\begin{equation}
  \label{eq:0ef-u1}
3\pd_Y^2U_1-4b_2\pd_X\pd_YU_1+\pd_X^2(\pd_X^2-4b_1+6u_1)U_1=0\,.
\end{equation}
\begin{lemma}
  \label{lem:U1-U2}
If $U_1$ and $U_2$ satisfy \eqref{eq:LMiura-1}--\eqref{eq:LMiura-t2}
and $U_2$ is a solution of \eqref{eq:0ef-u2}, then $U_1$ is a solution of
\eqref{eq:0ef-u1}.
\end{lemma}
\begin{proof}
By \eqref{eq:0ef-u2}, $e_{LKP}(u_2)U_2=0$.
Using \eqref{eq:LMiura-1} and the fact that $u_2=M_+(v_2)$ and $u_1=M_-(v_2)$, we have
\begin{align*}
  2\pd_x^{-1}\pd_y^2V=& \pd_y(U_1+U_2+4v_2V)
  \\=& \pd_y(U_1+U_2)+\pd_x\{8v_2^2V+2(u_1+u_2)V+2v_2(U_1+U_2)\}
       \\ \qquad & -2(u_1+u_2)\pd_xV-(U_1+U_2)(u_2-u_1)\,,
\end{align*}
\begin{align*}
  & 2(u_2U_2-u_1U_1)-2(u_1+u_2)\pd_xV-(U_1+U_2)(u_2-u_1)
\\ =&
(u_2+u_1)(U_2-U_1-2\pd_xV)=0\,,
\end{align*}
and
\begin{align*}
  &  e_{LKP}(u_2)U_2-e_{LKP}(u_1)U_1
  \\=&
\pd_x\{2(4\pd_t+\pd_x^3+3\pd_x^{-1}\pd_y^2)V+6(u_2U_2-u_1U_1)\}
  \\ =&
\pd_x^2\left\{3\pd_x(U_1-U_2)+3(\pd_x^{-1}\pd_y+2v_2)(U_1+U_2)
+ \nabla_v\widetilde{M}_-(u_1,v_2)V+\nabla_v\widetilde{M}_+(u_2,v_2)V\right\}\,.
\end{align*}
Combining the above with \eqref{eq:LMiura-t1} and
\eqref{eq:LMiura-t2}, we see that $e_{LKP}(u_1)U_1=0$ and that $U_1$
is a solution of \eqref{eq:0ef-u1}.  Thus we complete the proof.
\end{proof}

Let 
\begin{gather*}
A_1=-\pd_X-2v_2\,,\quad
A_2= -4\left\{\pd_X^2+(3v_2-b_2)\pd_X-b_1+3v_2^2+\frac{3}{2}u_2-2b_2v_2\right\}\,,\\
A_3=\pd_X-2v_2\,,\quad  
A_4= -4\left\{\pd_X^2-(3v_2-b_2)\pd_X-b_1+3v_2^2+\frac{3}{2}u_1-2b_2v_2\right\}\,,\\  
\\
 \mathcal{A}_1=\begin{pmatrix}
A_1 & I \\ A_2 & O
\end{pmatrix}\,,\quad
 \mathcal{A}_2=\begin{pmatrix}
A_3 & I \\ A_4 & O
\end{pmatrix}\,,
\\
\mathcal{B}_1=\begin{pmatrix}I & O \\   -3A_1-4b_2  & 3I\end{pmatrix}\,,
\quad
\mathcal{B}_2=\begin{pmatrix}I & O \\   -3A_3-4b_2  & 3I\end{pmatrix}\,.
\end{gather*}
Formally, $A_1^*=A_3$ and $A_2^*=A_4$.
By \eqref{eq:LMiura-1}--\eqref{eq:LMiura-t2},
\begin{gather}
  \label{eq:calA1-B1}
\mathcal{A}_1\begin{pmatrix}  V \\ \pd_X^{-1}\pd_YV\end{pmatrix}
=\mathcal{B}_1\begin{pmatrix}  U_1 \\ \pd_X^{-1}\pd_YU_1\end{pmatrix}\,,
\\
  \label{eq:calA2-B2}
\mathcal{A}_2\begin{pmatrix}  V \\ \pd_X^{-1}\pd_YV\end{pmatrix}
=\mathcal{B}_2
\begin{pmatrix}  U_2 \\ \pd_X^{-1}\pd_YU_2\end{pmatrix}\,.
\end{gather}

\begin{lemma}
\label{lem:kerA-solU}
Let $\beta_1=b_2-2a_{14}$, $U^*(\beta)=h(\bx)\Phi^{1,*}(\bx,-i\beta)$ and
\begin{equation*}
V(\beta)=\Phi^2(\bx,-i\beta)\Phi^{2,*}_1(\bx)\,,\quad  
\widehat{V}(\beta)=(\beta-\k_1)\Phi^2_1(\bx)\Phi^{2,*}_1(\bx,-i\beta)\,.
\end{equation*}
  \begin{enumerate}
  \item If $\beta=\k_1$, $\k_4$ or $\beta_1$, then
$A_2V(\beta)=0$ and $A_4U^*(\beta)=0$. Moreover, $\pd_xV(\beta_1)$              and $\pd_x\widehat{V}(\beta)$ are solutions of \eqref{eq:0ef-u2}.

\item Let
$\widetilde{V}=(b_2-2\k_1-\k_4)\pd_\beta V(\k_4)
-(b_2-\k_1-2\k_4)\pd_\beta\widehat{V}(\k_1)$.
Then $\pd_X\widetilde{V}$ is a solution of \eqref{eq:0ef-u2}.
\item If $b_2=\k_1+2\k_4$ or $b_2=2\k_1+\k_4$, then
$A_2\pd_\beta V(\beta_1)=0$ and $A_4\pd_\beta U^*(\beta_1)=0$. 
\end{enumerate}
\end{lemma}
\begin{proof}
  In view of \eqref{eq:Miura-Lax1} and \eqref{eq:Miura-Laxt2},
  \begin{gather*}
    \nabla M_-(v_2)h^{-1}\Phi^2(\bx,-i\beta)=0\,,\quad
    \nabla_v\widetilde{M}_-(u_1,v_2)h^{-1}\Phi^2(\bx,-i\beta)=0\,,
\\
    \nabla M_+(v_2)h\Phi^1(\bx,-i\beta)=0\,,\quad
    \nabla_v\widetilde{M}_+(u_2,v_2)h\Phi^1(\bx,-i\beta)=0\,.
  \end{gather*}
By a simple computation, we have
\begin{align*}
& \beta x+\beta^2y-\beta^3t-\frac{1}{2}(\theta_1+\theta_4)
  \\ = &
(\beta-a_{14})X+\left\{\beta^2-\frac{1}{2}(\k_1^2+\k_4^2)\right\}Y
-\lambda(\beta)t\,,
\end{align*}
where $\lambda(\beta)=(\beta-\k_1)(\beta-\k_4)(\beta-b_2+2a_{14})$.
Since
  \begin{equation*}
f_1=\tau_1=2\exp\frac{\theta_2+\theta_3}{2}\cosh(X+2a_{23}Y)\,,\quad
f_2=2\exp\frac{\theta_1+\theta_4}{2}\sinh(X+2a_{14}Y)\,,
\end{equation*}
and $\exp\frac{\theta_1+\theta_4}{2}h^{-1}$ depends only on $X$ and
$Y$, we see that $h^{-1}\Phi^2$ and $h\Phi^{1,*}$ are functions of $X$
and $Y$ multiplied by $e^{-\lambda(\beta)t}$ and
$e^{\lambda(\beta)t}$, respectively and that
  \begin{equation}
    \label{eq:Ub}
    A_2V(\beta)+4\lambda(\beta)\pd_x^{-1}V(\beta)=0\,,\quad
    A_4U^*(\beta)-4\lambda(\beta)\pd_x^{-1}U^*(\beta)=0\,.
  \end{equation}
  Thus we have $A_2V(\beta)=A_4U^*(\beta)=0$
  if $\beta=\k_1$, $\k_4$, $\beta_1$
  and $A_2\pd_\beta V(\beta_1)=A_4\pd_\beta U^*(\beta_1)=0$
  if $\beta=2\k_1+\k_4$ or $\beta=\k_1+2\k_4$.
  Note that $\beta=\beta_1$ is a double root of $\lambda(\beta)=0$ provided
 $\beta=2\k_1+\k_4$ or $\beta=\k_1+2\k_4$.
 \par
 Since $h^{-1}=\Phi^{2,*}_4(\bx)$, it follows from  Lemma~\ref{lem:prodJdJ} that
 $\pd_xV(\beta)$ is a solution of \eqref{eq:linequ2P}.
 If $\beta=\k_1$, $\k_4$, $\beta_1$, then $\pd_xV(\beta)$ is a stationary
 solution of \eqref{eq:linequ2P} in the moving coordinate $(X,Y)$ and
$\pd_xV(\beta)$ is a solution of \eqref{eq:0ef-u2}.
\par
Let $A=\pd_X^2-4b_1+6u_2$ and
$\overline{V}(\eps)=(\tilde{\beta}-\k_1)\Phi^2(\bx,-i\beta)\Phi^{2,*}(\bx,-i\tilde{\beta})$ with
$\beta=\k_4+\eps\lambda'(\k_1)$ and
$\tilde{\beta}=\k_1+\eps\lambda'(\k_4)$. Then
$\pd_\eps\overline{V}(0)=(\k_4-\k_1)\widetilde{V}$ and
  \begin{gather*}
\left(3\pd_Y^2-4b_2\pd_X\pd_Y+\pd_XA\pd_X\right)\overline{V}(\eps)
=4\left(\lambda(\beta)-\lambda(\tilde{\beta})\right)\pd_X\overline{V}(\eps)\,.
\end{gather*}
Since $\lambda(\k_1)=\lambda(\k_4)=0$, we have
\begin{align*}
  & \left(3\pd_Y^2-4b_2\pd_X\pd_Y+\pd_XA\pd_X\right)\pd_\eps\overline{V}(0)
  \\ =&
 4\left(\lambda(\k_4)-\lambda(\k_1)\right)\pd_X\pd_\eps\overline{V}(0)
+4\left(\pd_\eps\beta\lambda'(\k_4)-\pd_\eps\tilde{\beta}\lambda'(\k_1)\right)
\pd_X\overline{V}(0)
= 0\,.
\end{align*}
Thus we complete the proof.
\end{proof}

\begin{lemma}
  \label{lem:bfUU*}
Let $\a>0$ and $U$ and $\pd_XW$ be solutions of \eqref{eq:0ef-u2} such that
 \begin{gather*}
\mathbf{U}=\!(U,\pd_x^{-1}\pd_yU)^T
\in C^1(\R_Y;L^2(\R;e^{2\a X}dX))\,,
\\
\mathbf{U}^*=(-\frac{4}{3}b_2W+\pd_X^{-1}\pd_YW,W)^T
\in C^1(\R_Y;L^2(\R;e^{-2\a X}dX))\,.
\end{gather*}
Then
$$\la \mathbf{U},\mathbf{U^*}\ra:=\int_\R\left(-\frac{4}{3}b_2UW
  +\pd_X^{-1}\pd_YUW\right)\,dX$$
does not depend on $Y$.
\end{lemma}
\begin{proof}
Let $\mathbf{A}=
  \begin{pmatrix}
    0 & \pd_X \\ -\frac{1}{3}\pd_XA & \frac{4}{3}\pd_X
  \end{pmatrix}$.
Since $\pd_Y\mathbf{U}=\mathcal{A}\mathbf{U}$ and
$\pd_Y\mathbf{U^*}=-\mathbf{A}^*\mathbf{U^*}$,
\begin{align*}
  \frac{d}{dY}\la \mathbf{U},\mathbf{U^*}\ra
  =& \la \pd_Y\mathbf{U},\mathbf{U^*}\ra
     +\la \mathbf{U},\pd_Y\mathbf{U^*}\ra
     =0\,.
\end{align*}
\end{proof}

Let 
\begin{gather*}
\mathbf{W_1}=(\pd_XV(\k_1),\pd_YV(\k_1))^T\,,\quad
\mathbf{W_1^*}=(-\frac{4}{3}b_2V(\k_1)+\pd_X^{-1}\pd_YV(\k_1),V(\k_1))^T\,,
\\
\mathbf{W_2}=(\pd_X\widetilde{V},\pd_Y\widetilde{V})^T\,,
\quad
\mathbf{W_2^*}=(-\frac{4}{3}b_2\widetilde{V}
  +\pd_X^{-1}\pd_Y\widetilde{V},\widetilde{V})^T\,,
  \\
\mathbf{W_3}=(\pd_XV(\beta_1),\pd_YV(\beta_1))^T\,,\quad
\mathbf{W_3^*}=(-\frac{4}{3}b_2V(\beta_1)+\pd_X^{-1}\pd_YV(\beta_1),
V(\beta_1))^T\,,
\\  
\mathbf{W_4}=(\pd_X\widehat{V}(\beta_1),\pd_Y\widehat{V}(\beta_1))^T\,,
\quad
  \mathbf{W_4^*}=(-\frac{4}{3}b_2\widehat{V}(\beta_1)
  +\pd_X^{-1}\pd_Y\widehat{V}(\beta_1),\widehat{V}(\beta_1))^T\,.
\end{gather*}

\begin{lemma}
  \label{lem:orth-bfUU*}
Assume
 \begin{equation}
   \label{eq:ass-alpha-t}
0<\a<\left||b_2-3a_{14}|-\sqrt{c_{14}}\right|\,.
 \end{equation}
\begin{enumerate}
\item Suppose that $|b_2-3a_{14}|<\sqrt{c_{14}}$.
Then $\mathbf{W_i}\in C^1(\R;e^{2\a X}dX)$ and
$\mathbf{W_j^*}\in C^1(\R;e^{-2\a X}dX)$ for $1\le i$, $j\le 4$, and
\begin{gather*}
\la \mathbf{W_i},\mathbf{W_i^*}\ra=0\quad\text{and}\quad
\la \mathbf{W_i},\mathbf{W_j^*}\ra=0\text{ if $|i-j|\ge2$,}
\\
\la \mathbf{W_1},\mathbf{W_2^*}\ra=-\la \mathbf{W_2},\mathbf{W_1^*}\ra\ne0\,,
\quad
\la \mathbf{W_3},\mathbf{W_4^*}\ra=-\la \mathbf{W_4},\mathbf{W_3^*}\ra\ne0\,.
\end{gather*}
\item Suppose that $|b_2-3a_{14}|>\sqrt{c_{14}}$. Then
$\mathbf{W_1}$,   $\mathbf{W_2}\in C^1(\R;e^{2\a X}dX)$ and
$\mathbf{W_1^*}$, $\mathbf{W_2^*}\in C^1(\R;e^{-2\a X}dX)$, and
\begin{gather*}
\la \mathbf{W_i},\mathbf{W_i^*}\ra=0\text{  for $i=1$, $2$,}\quad
\la \mathbf{W_1},\mathbf{W_2^*}\ra=-\la \mathbf{W_2},\mathbf{W_1^*}\ra\ne0\,.
\end{gather*}
\end{enumerate}
\end{lemma}
\begin{proof}
  Let $\gamma=\beta-a_{14}$.
 As in the proof of Lemma~\ref{lem:orthogonality}, we have
 as $(a_{23}-a_{14})Y\to\pm\infty$,
 \begin{multline}
\label{eq:approxV}    
V(\beta)=\frac{e^{\mu_{2,\pm}}}{2}
\frac{a_{23}\pm\sqrt{c_{23}}-\beta}{\k_4-a_{23}\mp\sqrt{c_{23}}}
e^{(\gamma^2-c_{14})Y}
\\ \times
\left(\gamma-\sqrt{c_{14}}\tanh Z_{2,\pm}\right)e^{\gamma z_2}\sech Z_{2,\pm}
\left(1+O(e^{-2\sqrt{c_{23}}|z_1|})\right)\,,
\end{multline}
\begin{multline}
\label{eq:approxhatV}    
\widehat{V}(\beta)=\frac{\sqrt{c_{14}}e^{-\mu_{2,\pm}}}{\beta-\k_4}
\frac{\k_4-a_{23}\mp\sqrt{c_{23}}}{\beta-a_{23}\mp\sqrt{c_{23}}}
e^{-(\gamma^2-c_{14})Y}
\\  \times
\left(\gamma+\sqrt{c_{14}}\tanh Z_{2,\pm}\right)e^{-\gamma z_2}\sech Z_{2,\pm}
\left(1+O(e^{-2\sqrt{c_{23}}|z_1|})\right)\,.
\end{multline}
Especially, as $(a_{23}-a_{14})Y\to\infty$,
\begin{equation*}
  V(\k_1)=-\frac{\k_3-\k_1}{\k_4-\k_3}\frac{c_{14}}{2}
\sech^2Z_{2,+}
  \left(1+O\left(e^{-2\sqrt{c_{23}}|z_1|}\right)\right)\,,  
\end{equation*}
 \begin{multline*}
   \overline{V}(\eps)=
   \frac{\beta-\k_3}{(\tilde{\beta}-\k_3)(\tilde{\beta}-\k_4)}
\exp\left((\beta-\tilde{\beta})X+(\beta^2-\tilde{\beta}^2)Y
   -(\lambda(\beta)-\lambda(\tilde{\beta}))t\right)
   \\ \times
   \left(\beta-a_{14}-\sqrt{c_{14}}\tanh Z_{2,+}\right)
   \left(\beta'-a_{14}+\sqrt{c_{14}}\tanh Z_{2,+}\right)
   \left(1+O\left(e^{-2\sqrt{c_{23}}|z_1|}\right)\right)\,,
 \end{multline*}
where $\beta=\k_4+\lambda'(\k_1)\eps$ and
$\tilde{\beta}=\k_1+\lambda'(\k_4)\eps$. Since
$\lambda(\beta)-\lambda(\tilde{\beta})=O(\eps^2)$,
  \begin{align*}
&    (\beta-\tilde{\beta})X+(\beta^2-\tilde{\beta}^2)Y -(\theta_4-\theta_1)
    \\=& \eps\left\{\left(\lambda'(\k_1)-\lambda'(\k_4)\right)X
         +2\left(\k_4\lambda'(\k_1)-\k_1\lambda'(\k_4)\right)Y\right\}+O(\eps^2)
    \\=&
4\eps\left\{(b_2-3a_{14})(Z_{2,+}-\mu_{2,+})+2c_{14}^{3/2}Y\right\}+O(\eps^2)\,,
  \end{align*}
and
  \begin{align*}
    & \frac{(\beta-\k_3)e^{\theta_4-\theta_1}}
      {(\tilde{\beta}-\k_3)(\tilde{\beta}-\k_4)}
  \left(\beta-a_{14}-\sqrt{c_{14}}\tanh Z_{2,+}\right)
   \left(\beta'-a_{14}+\sqrt{c_{14}}\tanh Z_{2,+}\right)
    \\=&
         \frac{\k_4-\k_3}{\k_3-\k_1}\left(-\frac{\sqrt{c_{14}}}{2}+C_1\eps\right)
         \left\{\sech^2Z_{2,+}
 +4\eps(b_2-3a_{14})(1+\tanh Z_{2,+})\right\}+O(\eps^2)\,,
  \end{align*}
where  $C_1$ is a constant, we have as $(a_{23}-a_{14})Y\to\infty$,
  \begin{align*}
    \widetilde{V}=&
-\frac{\k_4-\k_3}{\k_3-\k_1}
\bigl\{ \left((b_2-3a_{14})Z_{2,+}+2c_{14}^{3/2}Y\right)\sech^2Z_{2,+}
\\ & +(b_2-3a_{14})(1+\tanh Z_{2,+})-\frac{2C_1}{c_{14}}V(\k_1)\bigr\}
\left(1+O(e^{-2\sqrt{c_{14}}|z_1|})\right)\,.
  \end{align*}
Combining the above with Lemma~\ref{lem:bfUU*}, we have
  \begin{align*}
\la \mathbf{W_1},\mathbf{W_2^*}\ra =&-\la \mathbf{W_2},\mathbf{W_1^*}\ra
 \\=&
-\frac{4}{3}b_2\int_\R \pd_XV(\k_1)\widetilde{V}\,dX
+\int_\R \left(\pd_YV(\k_1)\widetilde{V}-V(\k_1)\pd_Y\widetilde{V}\right)\,dX
    \\= &
\frac{4}{3}(b_2-3a_{14})\int_\R V(\k_1)\pd_X\widetilde{V}
-\int_\R V(\k_1)(\pd_Y-2a_{14}\pd_X)\widetilde{V}\,dX
\\=& \frac{4\sqrt{c_{14}}}{3}(b_2-2\k_1-\k_4)(b_2-\k_1-2\k_4)\,.
  \end{align*}
  \par
 If $|b_2-3a_{14}|<\sqrt{c_{14}}$, then
it follows from \eqref{eq:approxV} and \eqref{eq:approxhatV} that
$V(\beta_1)$ and $\widehat{V}(\beta_1)$ decay exponentially as $z_2\to\pm\infty$,
and $e^{\a X}\mathbf{W_3}$, $e^{\a X}\mathbf{W_4}\in L^2(\R)$.
On the other hand, if $|b_2-3a_{14}|>\sqrt{c_{14}}$, then
$V(\beta_1)$ and $\widehat{V}(\beta_1)$ grow exponentially as $z_2\to\infty$
or as $z_2\to-\infty$, and
$e^{\a X}\mathbf{W_3}$, $e^{\a X}\mathbf{W_4}\not\in L^2(\R)$.
\par
Suppose that $|b_2-3a_{14}|<\sqrt{c_{14}}$.
By integration by parts, we have $\la \mathbf{W_i}, \mathbf{W_i^*}\ra=0$.
Let $\gamma_1=b_2-3a_{14}$.
Combining \eqref{eq:approxV} and \eqref{eq:approxhatV} with
Lemma~\ref{lem:orth-bfUU*}, we have
\begin{align*}
\la \mathbf{W_3},\mathbf{W_4^*}\ra=& -\la \mathbf{W_4},\mathbf{W_3^*}\ra
  \\=& -\frac{4}{3}b_2\int_\R \pd_XV(\beta_1)\widehat{V}(\beta_1)\,dX
       +\int_\R \left(\pd_YV(\beta_1)\widehat{V}(\beta_1)
       -V(\beta_1)\pd_Y\widehat{V}(\beta_1)\right)\,dX
  \\=& -\frac{4\gamma_1}{3}\int_\R \pd_XV(\beta_1)\widehat{V}(\beta_1)\,dX
       +2(\gamma_1^2-c_{14})\int_\R V(\beta_1)\widehat{V}(\beta_1)\,dX\,.
  \\=& -\frac{2}{3}(b_2-2\k_1-\k_4)(b_2-\k_1-2\k_4)^2\,.
\end{align*}
\par

Since $\la \mathbf{W_i},\mathbf{W_j^*}\ra$ is independent of $Y$
and $(\beta_1-a_{14})^2-c_{14}\ne0$, 
we have $\la \mathbf{W_i},\mathbf{W_j^*}\ra=0$ if $|i-j|\ge2$.
\end{proof}

\begin{lemma}
 \label{lem:A4-est}
 Let $\beta_1=b_2-2a_{14}$. Let $A_2$ and $A_4$ be
 operators on $L^2(\R;e^{2\a X}dX)$.
  \begin{enumerate}    
  \item If $0<\a<\sqrt{c_{14}}-|b_2-3a_{14}|$, then for every $Y\in\R$,
  $\ker(A_4)=\ker(A_2^*)=\{0\}$ and
  $\ker(A_2)=\ker(A_4^*)=\spann\{V(\k_1),V(\beta_1)\}$.
  Moreover, if $F\in L^2(\R;e^{2\a X}dX)$ and
  \begin{equation}
    \label{eq:orthF}
    \int_\R FV(\k_1)\,dX=\int_\R FV(\beta_1)\,dX=0\,,
  \end{equation}
then     $A_4U=F$ has a solution satisfying
\begin{equation}
  \label{eq:U-Fest1}
      \sum_{0\le i\le2}\|\pd_X^iU(\cdot,Y)\|_{L^2(\R;e^{2\a X}dX)}\le C
      \|F(\cdot,Y)\|_{L^2(\R;e^{2\a X}dX)}\,,
    \end{equation}
    where $C$ is a constant that is independent of $Y$ and $F$.
  \item 
If $0<\a<|b_2-3a_{14}|-\sqrt{c_{14}}$, then for every $Y\in\R$,
    $\ker(A_2)=\ker(A_4^*)=\spann\{V(\k_1)\}$ and
    $\ker(A_2^*)=\ker(A_4)=\{0\}$. Moreover, if $F\in L^2(\R;e^{2\a X}dX)$ and
    \begin{equation}
      \label{eq:orthF'}
     \int_\R FV(\k_1)\,dX=0\,, 
    \end{equation}
then $A_4U=F$ has a solution
    satisfying
    \begin{equation}
      \label{eq:U-Fest}
      \sum_{0\le i\le2}\|\pd_X^iU(\cdot,Y)\|_{L^2(\R;e^{2\a X}dX)}\le C
      \|F(\cdot,Y)\|_{L^2(\R;e^{2\a X}dX)}\,,
    \end{equation}
    where $C$ is a constant that is independent of $Y$ and $F$.
  \end{enumerate}
\end{lemma}
\begin{proof}
Suppose that $|\gamma_1|=|b_2-3a_{14}|<\sqrt{c_{14}}$.
Let $\bx=(X,Y,0)$, $k_1^*(X,Y)=V(\k_1)$, $k_2^*(X,Y)=V(\beta_1)$ and 
\begin{gather*}
k_1(X,Y)=\frac{k_2^*(X,Y)}{4\operatorname{Wr}(k_1^*,k_2^*)(X,Y)}\,,\quad
k_2(X,Y)=-\frac{k_1^*(X,Y)}{4\operatorname{Wr}(k_1^*,k_2^*)(X,Y)}\,,
\end{gather*}
where $\operatorname{Wr}(k_1^*,k_2^*)(X,Y)=
k_1^*(X,Y)\pd_Xk_2^*(X,Y)-k_2^*(X,Y)\pd_Xk_1^*(X,Y)$.
Formally, $A_4^*=A_2$ and  $A_2k_1^*=A_2k_2^*=0$, $A_4k_1=A_4k_2=0$
by Lemma~\ref{lem:kerA-solU}.
\par
Since 
\begin{gather}
  \label{eq:orderk1k2}
k_1=O\left(\cosh^2\sqrt{c_{14}}z_2\right)\,,\quad
k_2=O\left(e^{-(\gamma_1^2-c_{14})Y}e^{-\gamma_1z_2}\cosh\sqrt{c_{14}}z_2\right)\,,
\\
\label{eq:orderk1k2*}
k_1^*=O\left(\sech^2\sqrt{c_{14}}z_2\right)\,,\quad
k_2^*=O\left(e^{(\gamma_1^2-c_{14})Y}e^{\gamma_1z_2}\sech\sqrt{c_{14}}z_2\right)\,,
\end{gather}
we have $\ker(A_4)=\{0\}$ and $\ker(A_4^*)=\spann\{k_1,k_2\}$.
Let
\begin{equation}
  \label{eq:defU}
  U(X,Y)=\sum_{i=1,2}k_i(X,Y)\int_{-\infty}^Xk_i^*(X',Y)F(X',Y)\,dX'\,.
\end{equation}
Formally, $U(X,Y)$ is a solution of $A_4U=F$.
Let $z_2=X+2a_{14}Y$ and $z_2'=X'+2a_{14}Y$. If $z_2z_2'\ge0$
and $|z_2'|\ge |z_2|$, then
\begin{gather*}
|k_1(X,Y)k_1^*(X',Y)|\le Ce^{-2\sqrt{c_{14}}|X-X'|}\,,
  \\
|k_2(X,Y)k_2^*(X',Y)|\le Ce^{-\gamma_1(X-X')}e^{-\sqrt{c_{14}}|X-X'|}\,,
\end{gather*}
where $C$ is a constant independent of $X$, $X'$ and $Y$.
Since $\a<\sqrt{c_{14}}-|\gamma_1|$, we have \eqref{eq:U-Fest1}
for $F$ satisfying \eqref{eq:orthF}.
\par
If $b_2<3a_{14}-\sqrt{c_{14}}$,
then $\ker(A_4)=\{0\}$ and $\ker(A_4^*)=\spann\{k_1\}$.
Moreover for $F$ satisfying \eqref{eq:orthF'},
\eqref{eq:defU} is a solution of $A_4U=F$ satisfying \eqref{eq:U-Fest}.
\par
If $b_2>3a_{14}+\sqrt{c_{14}}$, 
then $\ker(A_4)=\{0\}$ and $\ker(A_4^*)=\spann\{k_1\}$. Moreover,
$$
|k_2(X,Y)k_2^*(X',Y)|\le Ce^{-\gamma_1(X-X')}e^{\sqrt{c_{14}}|X-X'|}\,,
$$
and for $F$ satisfying \eqref{eq:orthF'},
\begin{equation*}
  U(X,Y)=k_1(X,Y)\int_{-\infty}^Xk_1^*(X',Y)F(X',Y)\,dX'
  -k_2(X,Y)\int^{\infty}_Xk_2^*(X',Y)F(X',Y)\,dX'
\end{equation*}
is a solution of $A_4U=F$ satisfying \eqref{eq:U-Fest}.
\end{proof}
\begin{lemma}
\label{lem:calA2-est}  
Assume \eqref{eq:ass-alpha-t}. Let $\mathcal{A}_1$ and $\mathcal{A}_2$ be operators
on $(H^1\times L^2)(\R;e^{2\a X}dX)$.
\begin{enumerate}
\item Suppose that $|b_2-3a_{14}|<\sqrt{c_{14}}$. Then
$\ker(\mathcal{A}_1)=\spann\left\{\pd_X^{-1}\mathbf{W_1},
  \pd_X^{-1}\mathbf{W_3}\right\}$ and $\ker(\mathcal{A}_2)=\{0\}$.  
Moreover, 
if $\mathcal{A}_2\mathbf{V}=\mathcal{B}_2\mathbf{U_2}$ and
$\la \mathbf{U_2},\mathbf{W_1^*}\ra=\la \mathbf{U_2},\mathbf{W_3^*}\ra=0$,
then
\begin{equation}
\label{eq:calA12-est}
  \|\mathbf{V}(\cdot,Y)\|_{(H^2\times H^1)(\R;e^{2\a X}dX)}
  \le C\|\mathbf{U_2}(\cdot,Y)\|_{(H^1\times L^2)(\R;e^{2\a X}dX)}\,,
\end{equation}
where $C$ is a constant that does not depend on $Y$ and $\mathbf{U_2}$.
\item
Suppose that $|b_2-3a_{14}|>\sqrt{c_{14}}$.
Then
$\ker(\mathcal{A}_1)=\spann\left\{\pd_X^{-1}\mathbf{W_1}\right\}$ and
$\ker(\mathcal{A}_2)=\{0\}$.
Moreover, 
if $\mathcal{A}_2\mathbf{V}=\mathcal{B}_2\mathbf{U_2}$ and
$\la \mathbf{U_2},\mathbf{W_1^*}\ra=0$, then
\begin{equation*}
  \|\mathbf{V}(\cdot,Y)\|_{(H^2\times H^1)(\R;e^{2\a X}dX)}
  \le C\|\mathbf{U_2}(\cdot,Y)\|_{(H^1\times L^2)(\R;e^{2\a X}dX)}\,,
\end{equation*}
where $C$ is a constant that does not depend on $Y$ and $\mathbf{U_2}$.
\end{enumerate}
\end{lemma}
\begin{proof}
We will prove Lemma~\ref{lem:calA2-est} assuming that  
$|b_2-3a_{14}|<\sqrt{c_{14}}$.
We can prove the other case in exactly the same way.
If  $\mathbf{V}=(V_1,V_2)^T\in\ker(\mathcal{A}_2)$,
then $V_2=-A_3V_1$ and $A_4V_1=0$, and it follows that $V_1=V_2=0$
since $\ker(A_4)=\{0\}$ by Lemma~\ref{lem:A4-est}.
If $\mathbf{V}=(V_1,V_2)^T\in\ker(\mathcal{A}_1)$,
then $V_2=-A_1V_1$ and $A_2V_1=0$.
Since $\ker(A_2)=\spann\{V(\k_1),V(\beta_1)\}$
by Lemma~\ref{lem:A4-est} 
and $\nabla M_-(v_2)V(\k_1)=\nabla M_-(v_2)V(\beta_1)=0$ by
Corollary~\ref{cl:kerM-P}, 
we have
$\ker(\mathcal{A}_1)=\spann\{\pd_X^{-1}\mathbf{W_1},\pd_X^{-1}\mathbf{W_3}\}$.
\par
Let $\mathbf{U_2}=(U_{21},U_{22})^T$ and $\mathbf{V}=(V_1,V_2)^T$.
Suppose that $\mathcal{A}_2\mathbf{V}=\mathcal{B}_2\mathbf{U_2}$
and that $\la\mathbf{U_2},\mathbf{W_1^*}\ra
=\la\mathbf{U_2},\mathbf{W_3^*}\ra=0$. Then
\begin{equation}
  \label{eq:pfcalA12-est2}
V_2=U_{21}-A_3V_1\,,\quad A_4V_1=-(3A_3+4b_2)U_{21}+3U_{22}\,.  
\end{equation}
Since  
$\mathcal{B}_2^*(0, V(\k_1))^T=3\mathbf{W_1^*}$ and
$\mathcal{B}_2^*( 0, V(\beta_1))^T=3\mathbf{W_3^*}$,
\begin{align}
    \label{eq:pfcalA12-est3}
\la -(3A_3+4b_2)U_{21}+3U_{22},V(\k_1)\ra=&
3\la \mathbf{U_2}, \mathbf{W_1^*}\ra=0\,,
\\
  \la  -(3A_3+4b_2)U_{21}+3U_{22},V(\beta_1)\ra
  =& 3\la \mathbf{U_2}, \mathbf{W_3^*}\ra=0\,.
\end{align}
Combining \eqref{eq:pfcalA12-est2} and \eqref{eq:pfcalA12-est3} with
Lemma~\ref{lem:A4-est}, we have \eqref{eq:calA12-est}.
\end{proof}
\begin{lemma}
\label{lem:kermapU1-U2}
Assume \eqref{eq:ass-alpha-t}.
Let $\mathcal{A}_1$ and $\mathcal{A}_2$ be operators
on $(H^1\times L^2)(\R;e^{2\a X}dX)$ and
 $\mathcal{A}_1\mathbf{V}=\mathcal{B}_1\mathbf{U_1}$,
$\mathcal{A}_2\mathbf{V}=\mathcal{B}_2\mathbf{U_2}$.
\begin{enumerate}
\item Suppose that $|b_2-3a_{14}|<\sqrt{c_{14}}$ and that
$\la \mathbf{U_2},\mathbf{W_i^*}\ra=0$ for $i=1,2,3,4$.
Then $\mathbf{U_1}=\mathbf{0}$ if and only if $\mathbf{U_2}=\mathbf{0}$.
\item
  Suppose that $|b_2-3a_{14}|>\sqrt{c_{14}}$ and that
$\la \mathbf{U_2},\mathbf{W_i^*}\ra=0$ for $i=1$, $2$.
  Then $\mathbf{U_1}=\mathbf{0}$ if and only if $\mathbf{U_2}=\mathbf{0}$.
\end{enumerate}
\end{lemma}
\begin{proof}
  Let $|b_2-3a_{14}|<\sqrt{c_{14}}$.
  Suppose that $\mathbf{U_2}=\mathbf{0}$. Then $\mathbf{V}=\mathbf{0}$
  and $\mathcal{B}_1\mathbf{U_1}=0$
  since $\ker(\mathcal{A}_2)=\{\mathbf{0}\}$ by Lemma~\ref{lem:calA2-est}.
  Since $\ker(\mathcal{B}_1)=\{\mathbf{0}\}$, we have
  $\mathbf{U_1}=\mathbf{0}$.
  \par
  Suppose that $\mathbf{U_1}=\mathbf{0}$. Then
  $\mathcal{A}_1\mathbf{V}=\mathbf{0}$ and it follows from
  Lemma~\ref{lem:calA2-est} that there exist $\delta_1$, $\delta_2\in\C$
such that
$\mathbf{V}=\delta_1\pd_X^{-1}\mathbf{W_1}+\delta_2\pd_X^{-1}\mathbf{W_3}$.
Since $\nabla M_-(v_2)V(\k_1)=\nabla M_-(v_2)V(\beta_1)=0$ and
$A_2V(\k_1)=A_2V(\beta_1)=0$,
\begin{align*}
  \mathcal{A}_2\mathbf{V}
  =&\delta_1   \begin{pmatrix}
      \nabla M_+(v_2)V(\k_1) \\ A_4V(\k_1) \end{pmatrix}
  +  \delta_2    \begin{pmatrix}
    \nabla M_+(v_2)V(\beta_1) \\ A_4V(\beta_1)\end{pmatrix}
  \\=&
    \delta_1 \begin{pmatrix}  2\pd_XV(\k_1) \\ (A_4-A_2)V(\k_1)
    \end{pmatrix}
  +  \delta_2 \begin{pmatrix}  2\pd_XV(\beta_1) \\ (A_4-A_2)V(\beta_1)
  \end{pmatrix}
  \\=& 2\mathcal{B}_2(\delta_1\mathbf{W_1}+\delta_2\mathbf{W_3})\,.
\end{align*}
Since $\ker(\mathcal{B}_2)=\{\mathbf{0}\}$,
$\mathbf{U_2}=2(\delta_1\mathbf{W_1}+\delta_2\mathbf{W_3})$,
and it follows from Lemma~\ref{lem:orth-bfUU*} that
$\delta_1=\delta_2=0$.
\par
Similarly, we can prove the case where $|b_2-3a_{14}|>\sqrt{c_{14}}$.
This completes the proof of Lemma~\ref{lem:kermapU1-U2}.
\end{proof}

Now we are in position to prove Proposition~\ref{prop:0notev}.
\begin{proof}[Proof of Proposition~\ref{prop:0notev}]
Let $\mathbf{U_2}=(U_2,\pd_X^{-1}\pd_YU_2)^T$.
By Proposition~\ref{lem:mLu=f},
$\pd_X^i\pd_Y^jU_2\in\calX_1$ for every $i\in\Z$ and $j\ge0$ and
\begin{gather*} 
e^{\a X}\mathbf{U_2}(X,0)\in H^1(\R)\times L^2(\R)\,,
  \quad
  \lim_{Y\to\infty}e^{2\a a_{23}Y}\|U_2(\cdot,Y)\|_{(H^1\times
L^2)(\R;e^{2\a X}dX)}=0\,.
\end{gather*}
Then it follows from Lemma~\ref{lem:bfUU*} that
\begin{equation}
  \label{eq:37}\left\{
    \begin{aligned}
& \la \mathbf{U_2},\mathbf{W_i^*}\ra=0\quad\text{if
    $1\le i\le 4$ and $|b_2-3a_{14}|<\sqrt{c_{14}}$,}
\\ &  \la \mathbf{U_2},\mathbf{W_i^*}\ra=0\quad\text{if
    $i=1$, $2$ and $|b_2-3a_{14}|>\sqrt{c_{14}}$,}    
    \end{aligned}\right.
\end{equation}
since 
$\|e^{-\a(\cdot+2a_{23}Y)}\mathbf{W_i^*}(\cdot,Y)\|_{L^2(\R)}\to0$ exponentially
as $Y\to\infty$ or as $Y\to-\infty$.
By \eqref{eq:37} and Lemma~\ref{lem:calA2-est}, there exist
$U_1$ and $V$ satisfying \eqref{eq:LMiura-1}--\eqref{eq:LMiura-t2}
and
\begin{equation}
 \label{eq:U1-U2bound}
\|\mathbf{U}_1(\cdot, Y)\|_{(H^1\times L^2)(\R;e^{2\a X}dX)}
\le C\|\mathbf{U}_2(\cdot, Y)\|_{(H^1\times L^2)(\R;e^{2\a X}dX)}\,,  
\end{equation}
where $\mathbf{U_1}=(U_1,\pd_X^{-1}\pd_YU_1)^T$ and
$C$ is a constant that does not depend on $Y$.
It follows from Lemma~\ref{lem:U1-U2} that
$U_1$ is a solution of \eqref{eq:0ef-u1}.
Moreover, $U_1\in\calX_1$ by \eqref{eq:U1-U2bound}.
Since $\lambda=0$ is not an eigenvalue
of $\mL_1$ on $\calX_1$ by Proposition~\ref{lem:wmL1-decay},  $U_1=0$
and it follows from Lemma~\ref{lem:kermapU1-U2} that $U_2=0$.
\end{proof}
\par

Combining Lemmas~\ref{lem:lambda2-etaI} and \ref{lem:anal-F},
Propositions~\ref{prop:specL2-X1} and \ref{prop:0notev} with Remark~\ref{rem:6.5},
we have the following.
\begin{proposition}
  \label{prop:condstability}
  Assume \eqref{ass-alpha'} and \eqref{eq:l2ne0}. Let $\eta_{0,1}$ be a sufficiently
 small positive number and let
$\tilde{\eta}_{0,2}\in (\tilde{\eta}_{*,+},\eta_+')$ if $a_{14}>a_{23}$ and
$\tilde{\eta}_{0,2}\in (\tilde{\eta}_{*,-},\eta_-')$ if $a_{14}<a_{23}$.
Then there exist positive constants $K$ and $b$ such that
$$\left\|\left(I-P_2'(\eta_{0,1},\tilde{\eta}_{0,2}\right)e^{t\mL_2}\right
\|_{B(\calX_1)}\le Ke^{-bt} \quad\text{for $t\ge0$.}$$
\end{proposition}
\bigskip

\subsection{Asymptotic linear stability of $2$-line solitons of P-type}
In this section, we investigate asymptotic behavior of solutions to
the linearized KP-II equation around $u_2=2\pd_x^2\log\tau$ with
\eqref{eq:P-type}. The $2$-line soliton $u_2$ can be expressed as
\begin{gather}
\label{eq:2-lsol-P}
u_2=\varphi_{\bc}(\bz)=2\pd_x^2\log\tilde{\tau}_2\,,
\\ \label{eq:ttau}
\tilde{\tau}_2=
\begin{vmatrix}
\cosh\sqrt{c_1}z_1 & \sinh\sqrt{c_2}z_2
\\ 
a_{23}\cosh\sqrt{c_1}z_1+\sqrt{c_1}\sinh\sqrt{c_1}z_1
&
a_{14}\sinh\sqrt{c_2}z_2+\sqrt{c_2}\cosh\sqrt{c_2}z_2
\end{vmatrix}\,,
\end{gather}
where $\bc=(c_1,c_2)=(c_{23},c_{14})$, $\bz=(z_1,z_2)$, 
$z_1=x+2a_{23}y-\omega_{23}t$, $z_2=x+2a_{14}y-\omega_{14}t$, and 
$a_{ij}$ and $\omega_{ij}$ are constants defined as \eqref{eqdef:a,c,omega}.
\par
Let $H_t(y)=(4\pi t)^{-1/2}\exp(-y^2/4t)$,
$W_t(y)=\frac12$ if $(2\k_1+\k_4)t\le y\le (\k_1+2\k_4)t$
and $W_t(y)=0$ otherwise.
\begin{theorem}
  \label{thm:AS}
Assume \eqref{eq:ass-alpha}. If
$u_0\in\calX_2\cap L^1_yL^2(\R;e^{2\a z_2}dz_2)$ and $u(t)$ is a solution of
\begin{equation}
  \label{eq:27}
4\pd_tu+\pd_x^3u+3\pd_x^{-1}\pd_y^2u+6\pd_x(u_2u)=0\,,\quad u(0)=u_0\,,
\end{equation}
in the class $ C([0,\infty);\calX_2)$. Then
$$\left\|e^{\a z_2}\left\{u(t,\cdot)-\left(H_{t/2\sqrt{c_2}}*W_t*f(y)\right)
\pd_{z_2}\varphi_\bc(\bz)\right\} \right\|_{L^2(\R^2)}=O(t^{-1/4})\quad
\text{as $t\to\infty$,}$$
where $$f(y)=-\frac{1}{4c_2}\int_\R u_0(x,y)
\left(\int_{-\infty}^x \pd_{z_2}\varphi_\bc  (x'+2a_{23}y,x'+2a_{14}y)\,dx'\right)\,dx\,.$$
\end{theorem}

Let $\widetilde{W}_t(y)=\frac12$
if $(2\k_2+\k_3)t\le y\le (\k_2+2\k_3)t$ and $\widetilde{W}_t(y)=0$ otherwise.
\begin{theorem}
  \label{thm:AS2}
  Assume \eqref{ass-alpha'}, \eqref{eq:l2ne0}, \eqref{eq:ass-alpha-t}
  and that $b_2\ne 2\k_1+\k_4$, $2\k_1+\k_4$.  
  Let $\eta_I=2\a(a_{23}-a_{14})$. There exists $\eta'>0$ such that if
$u(0)\in \calX_1\cap L^1(\R_y;L^2(\R;e^{2\a z_1}dz_1))$ and
  \begin{gather*}
    \int_{\R^2}u_0\overline{g^{2,*}_{2,k}(\cdot,\eta_R+i\eta_I)}\,dxdy=0
\quad\text{for $\eta_R\in[-\eta',\eta']$,}   
  \end{gather*}
then a solution  $u(t)$ of \eqref{eq:27} satisfies as $t\to\infty$,
$$\left\|e^{\a z_1}\left\{u(t,\cdot)
-\left(H_{t/2\sqrt{c_1}}*\widetilde{W}_t*f(y)\right) \pd_{z_1}\varphi_\bc(\bz)
  \right\}\right\|_{L^2(\R^2)}=O(t^{-1/4})\,,$$
where
$$f(y)=-\frac{1}{4c_1}\int_\R u_0(x,y)
\left(\int_{-\infty}^x \pd_{z_1}\varphi_\bc(x'+2a_{23}y,x'+2a_{14}y)\,dx'\right)\,dx\,.$$
\end{theorem}
\begin{remark}
If $b_2\ne 2\k_1+\k_4$, $\k_1+2\k_4$, then it follows from
\eqref{eq:pdlambdas-0} that $\pd_\eta\lambda^2_\pm(0)\ne0$,
and modulation of $[1,4]$-soliton gets away from the intersection
of $[1,4]$-soliton and $[2,3]$-soliton.
\end{remark}

To prove Theorems~\ref{thm:AS} and \ref{thm:AS2},
we need the first order asymptotics of $g^2_{2,k}(\eta)$
and $g^{2,*}_{2,k}(\eta)$ as $\eta\to0$. Let
\begin{gather*}
\gamma_1^\pm(c_2)
=\frac{1}{2\sqrt{c_1}}
\log\frac{a_{14}\pm\sqrt{c_2}-\k_3}{a_{14}\pm\sqrt{c_2}-\k_2}\,,
\quad
\gamma_2^\pm(c_1)
=\frac{1}{2\sqrt{c_2}}
\log\frac{\k_4-a_{23}\mp\sqrt{c_1}}{a_{23}\pm\sqrt{c_1}-\k_1}\,,
\\
\varphi^{1,\pm}(\bz)=\frac{d}{dc_1}
\varphi_\bc\left(z_1,z_2'-\gamma_2^\pm(c_1)\right)
\bigr|_{z_2'=z_2+\gamma_2^\pm(c_1)}\,,
\\
\varphi^{2,\pm}(\bz)=\frac{d}{dc_2}
\varphi_\bc\left(z_1'-\gamma_1^\pm(c_2),z_2\right)
\bigr|_{z_1'=z_1+\gamma_1^\pm(c_2)}\,.
\end{gather*}

\begin{lemma}
  \label{lem:gg*profile}
 For $k=1$ and $2$,
\begin{gather*}
g^2_{k,1}(x,y,\eta)= -\frac{1}{2\sqrt{c_k}}e^{iy\eta}
  \left(\pd_{z_k}\varphi_\bc(\bz)+O(\eta)\right)\,,
\\
g^2_{k,2}(x,y,\eta)=\frac12 e^{iy\eta}
\left(\varphi^{k,+}+\frac{1}{2c_k^{3/2}}\pd_{z_k}\varphi_\bc+O(\eta)\right)\,,
\\
g^{2,*}_{k,1}(x,y,\eta)=
\frac{1}{2}e^{iy\eta}\left(\int_{-\infty}^x\varphi^{k,-}+O(\eta)\right)\,,
\\
g^{2,*}_{k,2}(x,y,\eta)=
\frac{1}{2\sqrt{c_k}} e^{iy\eta}
\left(\int_{-\infty}^x\pd_{z_k}\varphi_\bc+O(\eta)\right)\,.
\end{gather*}
\end{lemma}
\begin{remark}
Lemma~\ref{cl:sgnpdlambda1} tells us that perturbations of
$[2,3]$-soliton propagate toward $y=-\infty$ along its crest.
The leading term of $g^2_{2,2}(\eta)$ suggests
that modulations of the speed parameter of $[1,4]$-soliton do not to
cause phase shifts of $[2,3]$-soliton for $y\gg1$.
\end{remark}
\begin{proof}[Proof of Lemma~\ref{lem:gg*profile}]
Since $\gamma_{14}(\eta)=\sqrt{c_2}+i(\k_4-\k_1)^{-1}\eta+O(\eta^2)$ and
\begin{gather*}
\frac{(\k_j-\k_1)(\k_4-\k_j)}{\k_j-\beta_{14}^-(-\eta)}
=\k_4-\k_j+ie_j\eta+O(\eta^2)\,,
\quad   e_j=\frac{\k_4-\k_j}{(\k_4-\k_1)(\k_j-\k_1)}\,,
\\
\frac{(\k_j-\k_1)(\k_4-\k_j)}{\beta_{14}^+(-\eta)-\k_j}
=\k_j-\k_1+ie_j'\eta +O(\eta^2)\,,
\quad   e_j'=\frac{\k_j-\k_1}{(\k_4-\k_1)(\k_4-\k_j)}\,,
\end{gather*}
it follows from \eqref{eq:g2pmpm*} and \eqref{eq:g2+j}--\eqref{eq:g2*-j} that
\begin{gather*}
g^2_{2,+}=
2e^{iy\eta}\pd_x^2\sum_{j=2,3} \frac{e^{\theta_1+\theta_j}}{\tau_2}
\Bigl\{\k_j-\k_1+i\eta\left(\frac{2(a_{14}-\k_j)}{(\k_4-\k_1)^2}
  -\frac{\k_j-\k_1}{\k_4-\k_1}z_2\right)+O(\eta)^2\Bigr\}\,,
\end{gather*}
\begin{multline*}
g^2_{2,-}=
2e^{iy\eta}\pd_x^2\sum_{j=2,3}\frac{e^{\theta_1+\theta_j}}{\tau_2}
\bigl\{\k_j-\k_1
\\ +i\eta
\left(e_j'+2\frac{\k_j-\k_1}{(\k_4-\k_1)^2}+\frac{\k_j-\k_1}{\k_4-\k_1}z_2\right)
+O(\eta^2) \bigr\}\,,  
\end{multline*}
\begin{gather*}
g^{2,*}_{2,+}=
e^{iy\eta}\pd_x\sum_{j=2,3} \frac{e^{\theta_4+\theta_j}}{\tau_2}
\Bigl\{\k_4-\k_j+i\bar{\eta}\left(e_j-\frac{\k_4-\k_j}{\k_4-\k_1}z_2\right)
  +O(\eta)^2\Bigr\}\,,\\
g^{2,*}_{2,-}=
e^{iy\eta}\pd_x\sum_{j=2,3}  \frac{e^{\theta_4+\theta_j}}{\tau_2}
\bigl\{\k_4-\k_j
+i\bar{\eta}\left(\frac{1}{\k_4-\k_1}+\frac{\k_4-\k_j}{\k_4-\k_1}z_2\right)
+O(\eta^2) \bigr\}\,,
\end{gather*}
and that
\begin{equation*}
g^2_{2,1}= 2e^{iy\eta}
\pd_x^2\sum_{j=2,3}\frac{e^{\theta_1+\theta_j}}{\tau_2}
\left(\k_j-\k_1+O(\eta)\right)\,,  
\end{equation*}
\begin{multline*}
g^2_{2,2}= \frac{-2e^{iy\eta}}{\k_4-\k_1}
\pd_x^2\sum_{j=2,3}\frac{e^{\theta_1+\theta_j}}{\tau_2}
\Bigl\{(\k_j-\k_1)z_2
  +\frac{\k_j-a_{14}}{\k_4-\k_j}+\frac{2(\k_j-\k_1)}{\k_4-\k_1}
+O(\eta)  \Bigr\}\,,  
\end{multline*}  
\begin{gather*}
  g^{2,*}_{2,1}=\frac{2e^{iy\eta}}{\k_4-\k_1}
\pd_x\sum_{j=2,3}\frac{e^{\theta_4+\theta_j}}{\tau_2}                    
\left\{\frac{\k_j-a_{14}}{\k_j-\k_1}+(\k_4-\k_j)z_2+O(\eta)\right\}\,,
  \\
g^{2,*}_{2,2}= 2e^{iy\eta}\pd_x\sum_{j=2,3}\frac{e^{\theta_4+\theta_j}}{\tau_2}  
\left(\k_4-\k_j+O(\eta)\right)\,.
\end{gather*}
Differentiating \eqref{eq:ttau} with respect to $z_1$, $z_2$, $c_1$ and $c_2$,
we have
\begin{equation}
  \label{eq:ttau-z}
\left\{  
\begin{aligned}
& 1+\frac{1}{\sqrt{c_1}}\frac{\pd_{z_1}\tilde{\tau}_2}{\tilde{\tau}_2}
 =\frac{2e^{\theta_3}}{\tau_2}\sum_{j=1,4}(-1)^j(\k_j-\k_3)e^{\theta_j}\,,
  \\ & 
1-\frac{1}{\sqrt{c_1}}\frac{\pd_{z_1}\tilde{\tau}_2}{\tilde{\tau}_2}
=\frac{2e^{\theta_2}}{\tau_2}\sum_{j=1,4}(-1)^j(\k_j-\k_2)e^{\theta_j}\,,
  \\   &
1+\frac{1}{\sqrt{c_2}}\frac{\pd_{z_2}\tilde{\tau}_2}{\tilde{\tau}_2}
 =\frac{2e^{\theta_4}}{\tau_2}\sum_{j=2,3}(\k_4-\k_j)e^{\theta_j}\,,
  \\ &
1-\frac{1}{\sqrt{c_2}}\frac{\pd_{z_2}\tilde{\tau}_2}{\tilde{\tau}_2}
=\frac{2e^{\theta_1}}{\tau_2}\sum_{j=2,3}(\k_j-\k_1)e^{\theta_j}\,,
\end{aligned}\right.
\end{equation}
\begin{equation}
  \label{eq:ttau-c}
\left\{\begin{aligned}
&  z_1+(\k_3-\k_2)\frac{\pd_{c_1}\tilde{\tau}_2}{\tilde{\tau}_2}
  = \frac{2e^{\theta_3}}{\tau_2}\sum_{j=1,4}(-1)^j(\k_j-\k_3)z_1e^{\theta_j}
-\frac{f_1^*f_2}{\tau_2}\,,
\\ &
  z_1-(\k_3-\k_2)\frac{\pd_{c_1}\tilde{\tau}_2}{\tilde{\tau}_2}
  = \frac{2e^{\theta_2}}{\tau_2}\sum_{j=1,4}(-1)^j(\k_j-\k_2)z_1e^{\theta_j}
+\frac{f_1^*f_2}{\tau_2}\,,
  \\ &
z_2+(\k_4-\k_1)\frac{\pd_{c_2}\tilde{\tau}_2}{\tilde{\tau}_2}
  = \frac{2e^{\theta_4}}{\tau_2}\sum_{j=2,3}(\k_4-\k_j)z_2e^{\theta_j}
+\frac{f_1f_2^*}{\tau_2}\,,
\\ &
  z_2-(\k_4-\k_1)\frac{\pd_{c_2}\tilde{\tau}_2}{\tilde{\tau}_2}
  = \frac{2e^{\theta_1}}{\tau_2}\sum_{j=2,3}(\k_j-\k_1)z_2e^{\theta_j}
-\frac{f_1f_2^*}{\tau_2}\,,
\end{aligned}\right.
\end{equation}
where $f_1^*=e^{\theta_3}-e^{\theta_2}$, $f_2^*=e^{\theta_4}+e^{\theta_1}$.
Combining the above with \eqref{eq:2-lsol-P}, we have
\begin{gather*}
 g^2_{2,1}= -\frac{1}{\sqrt{c_2}}e^{iy\eta}
 \left(\pd_{z_2}\varphi_\bc(\bz)+O(\eta)\right)\,,
\\
g^2_{2,2}= \frac12 e^{iy\eta}\left\{\pd_{c_2}\varphi_\bc
+\frac{1}{2c_2^{3/2}}\pd_{z_2}\varphi_\bc
-\frac{1}{2\sqrt{c_2}(\k_4-\k_3)(\k_4-\k_2)}
\pd_{z_1}\varphi_\bc+O(\eta)\right\}\,,
\end{gather*}
\begin{align*}
g^{2,*}_{2,1}=& \frac{e^{iy\eta}}{2\sqrt{c_2}}
 \left[1+\pd_x\left\{(\k_4-\k_1)\frac{\pd_{c_2}\tilde{\tau}_2}{\tilde{\tau}_2}
+\frac{1}{(\k_3-\k_1)(\k_2-\k_1)}\frac{\pd_{z_1}\tilde{\tau}_2}{\tilde{\tau}_2}
\right\}\right]+O(\eta)
  \\=&
       \frac{e^{iy\eta}}{2}\int_{-\infty}^x\left\{
\pd_{c_2}\varphi_\bc +\frac{1}{2\sqrt{c_2}(\k_3-\k_1)(\k_2-\k_1)}
\pd_{z_1}\varphi_\bc+O(\eta)\right\}\,,
\end{align*}
\begin{equation*}
g^{2,*}_{2,2}= \frac{e^{iy\eta}}{\sqrt{c_2}}\left(
  \pd_x\left\{\frac{\pd_{z_2}\tilde{\tau}_2}{\tilde{\tau}_2}\right)
  +O(\eta)\right\}
= \frac{e^{iy\eta}}{2\sqrt{c_2}}
\left(\int_{-\infty}^x\pd_{z_2}\varphi_\bc+O(\eta)\right)\,.
\end{equation*}
Since
\begin{gather*}
\varphi^{2,+}=
\pd_{c_2}\varphi_\bc-\frac{1}{2\sqrt{c_2}(\k_4-\k_3)(\k_4-\k_2)}
\pd_{z_1}\varphi_\bc\,,
\\
\varphi^{2,-}=
\pd_{c_2}\varphi_\bc +\frac{1}{2\sqrt{c_2}(\k_3-\k_1)(\k_2-\k_1)}
\pd_{z_1}\varphi_\bc\,,
\end{gather*}
we have Lemma~\ref{lem:gg*profile} for $k=2$.
\par

By \eqref{eq:g21pmpm*} and \eqref{eq:g21+j}--\eqref{eq:g21*-j},
\begin{gather*}
  g^2_{1,1}(X,Y,\eta)=
2e^{iy\eta}\pd_x^2\sum_{j=1,4}(-1)^j\frac{e^{\theta_2+\theta_j}}{\tau_2}
\left(\k_j-\k_2+O(\eta)\right)\,,
\\
g^2_{1,2}= \frac{2e^{iy\eta}}{\k_3-\k_2}
\pd_x^2\sum_{j=1,4}(-1)^j\frac{e^{\theta_2+\theta_j}}{\tau_2}
\Bigl\{-(\k_j-\k_2)z_1
+\frac{\k_j-a_{23}}{\k_j-\k_3}-2\frac{\k_j-\k_2}{\k_3-\k_2}
+O(\eta)  \Bigr\}\,,  
\\
g^{2,*}_{1,1}=\frac{2e^{iy\eta}}{\k_3-\k_2}
\pd_x\sum_{j=1,4}(-1)^j\frac{e^{\theta_3+\theta_j}}{\tau_2}                    
\left\{(\k_j-\k_3)z_1+  \frac{\a_{23}-\k_j}{\k_j-\k_2}+O(\eta)\right\}\,,
\\
g^{2,*}_{1,2}=
2e^{iy\eta}\pd_x\sum_{j=1,4}(-1)^j\frac{e^{\theta_3+\theta_j}}{\tau_2}
\left(\k_j-\k_3+O(\eta)\right)\,.
\end{gather*}
Combining the above with \eqref{eq:ttau-z} and \eqref{eq:ttau-c},
we have Lemma~\ref{lem:gg*profile} for $k=1$.
Thus we complete the proof.
\end{proof}

\begin{proof}[Proof of Theorems~\ref{thm:AS} and \ref{thm:AS2}]
 Let 
\begin{gather*}
a_k(t,\eta)=\frac{1}{\sqrt{2\pi}}
  \int_{\R^2}u(t,x,y)g^{2,*}_{2,k}(t,x,y)\,dxdy\quad
  \text{for $k=1$ and $2$,}\\  
D(\eta)=
\begin{pmatrix}
  d_1(\eta) & d_2(\eta) \\ -\eta^2d_2(\eta) & d_1(\eta)\end{pmatrix}\,,
\\
d_1(\eta)=\frac{\lambda_{2,+}(\eta)+\lambda_{2,-}(\eta)}{2}
=i(b_2-3a_{14})\eta-\eta\Im\gamma_{14}(\eta)\,,
\\
d_2(\eta)=\frac{\lambda_{2,+}(\eta)-\lambda_{2,-}(\eta)}{2i\eta}
= \Re\gamma_{14}(\eta)\,,
\end{gather*}
and $e(t,\eta)=\eta^2|a_1(t,\eta)|^2+|a_2(t,\eta)|^2$.
Since
\begin{equation*}
\pd_ta_k(t,\eta)=\int_{\R^2}u(t,x,y)\overline{\mL_2^*g^{2,*}_{2,k}(X,Y,\eta)}\,dxdy\,,  
\end{equation*}
it follows from Lemma~\ref{lem:gg*ev} that
\begin{gather*}
  \pd_t
  \begin{pmatrix}
    a_1(t,\eta) \\ a_2(t,\eta)
  \end{pmatrix}
  =  D(\eta)  \begin{pmatrix}
    a_1(t,\eta) \\ a_2(t,\eta)
  \end{pmatrix}\,,
\end{gather*}
and   $e(t,\eta)=e^{-2t\eta^2\Re\gamma_{14}(\eta)}e(0,\eta)$.
Thus we have for $j\ge0$,
\begin{align}
  \label{eq:e(t,eta)}
 \|\eta^je(t,\eta)\|_{L^2(-\eta_0,\eta_0)}\lesssim &
\|\eta^je^{-2t\eta^2\Re\gamma_{14}(\eta)}\|_{L^2(-\eta_0,\eta_0)}
\|e(0,\eta)\|_{L^\infty(-\eta_0,\eta_0)}
  \\ \lesssim & \notag
(1+t)^{-(2j+1)/4}\|u_0\|_{L^1_yL^2(\R;e^{2\a z_2}dz_2)}\,.
\end{align}
Moreover, we have
\begin{gather*}
D(\eta)=D_0(\eta)+
\begin{pmatrix}
  O(\eta^4) & O(\eta^2) \\ O(\eta^4) & O(\eta^4)
\end{pmatrix}\,,\\
D_0(\eta)=\{i(b_2-3a_{14})\eta -\frac{\eta^2}{2\sqrt{c_2}}\}I+
\sqrt{c_2} \begin{pmatrix}
0  & 1 \\ -\eta^2 & 0 \end{pmatrix}\,,
\\
e^{tD_0(\eta)}=e^{i(b_2-3a_{14})t\eta-t\eta^2/2\sqrt{c_2}}
\begin{pmatrix}
\cos(t\sqrt{c_2}\eta) & \eta^{-1}\sin(t\sqrt{c_2}\eta)
 \\
-\eta\sin(t\sqrt{c_2}\eta) & \cos(t\sqrt{c_2}\eta)
\end{pmatrix}\,.  
\end{gather*}
Using Lemma~\ref{lem:gg*profile}, \eqref{eq:e(t,eta)} and the variation of constants
formula, we have 
\begin{align}
\label{eq:a1-est}  
  &  \left\|a_1(t,\eta)+\sqrt{4 c_2}
    e^{i(b_2-3a_{14})t\eta-t\eta^2/2\sqrt{c_2}}
\frac{\sin(t\sqrt{c_2}\eta)}{\eta}\hat{f}(\eta)\right\|_{L^2(-\eta_0,\eta_0)}
  \\ \lesssim & \notag 
\sum_{j=1,2}\left(
\|e^{-\frac{t\eta^2}{2\sqrt{c_2}}}a_j(0,\eta)\|_{L^2(-\eta_0,\eta_0)}
  +\int_0^t \left\|\eta^{4-j}e^{-\frac{(t-s)\eta^2}{2\sqrt{c_2}}}
   a_j(s,\eta)\right\|_{L^2(-\eta_0,\eta_0)}\,ds\right)
  \\  & \notag
        +(1+t)^{-1/4}\|u_0\|_{L^1_yL^2(\R;e^{2\a z_2}dz_2)}
\lesssim (1+t)^{-1/4}\,.
\end{align}
Combining the Plancherel theorem with Lemma~\ref{lem:gg*2-asymp},
\eqref{eq:e(t,eta)} and \eqref{eq:a1-est}, we have
\begin{align*}
&  \left\|\int_{-\eta_0}^{\eta_0} \left\{a_1(t,\eta)+\sqrt{4 c_2}
e^{i(b_2-3a_{14})t\eta-t\eta^2/2\sqrt{c_2}}
 \frac{\sin(t\sqrt{c_2}\eta)}{\eta}\hat{f}(\eta)\right\}
g^2_{2,1}(x,y,\eta)\,d\eta\right\|_{\calX_2}    
  \\ & +
       \left\|\int_{-\eta_0}^{\eta_0}a_2(t,\eta)g^2_{2,2}(x,y,\eta)
       \,d\eta\right\|_{\calX_2}
\\ \lesssim & (1+t)^{-1/4}\,.       
\end{align*}
By Lemma~\ref{lem:gg*profile},
\begin{align*}
& e^{-iY\eta}g^2_{2,1}(X,Y,\eta)+(2\sqrt{c_2})^{-1}\pd_{z_2}\varphi_\bc(\bz)
=\eta\int_0^1\pd_{\eta}\left\{e^{-iY\eta}g^2_{2,1}(X,Y,\theta\eta)\right\}
\,d\theta\,.
\end{align*}
and it follows from Lemma~\ref{lem:gg*2-asymp} that
\begin{align*}
& \Biggl\|\int_{-\eta_0}^{\eta_0}e^{-t\eta^2/2\sqrt{c_2}}\sin(\sqrt{c_2}t\eta)
\int_0^1\pd_{\eta}\Bigl\{e^{-iy\eta}g^2_{2,1}(x,y,\theta\eta)
  \\ & \phantom{e^{-iy\eta}g^2_{2,1}(x,y,\theta\eta)}
  -\sum_{\pm}\tg_1(\frac{Z_{2,\pm}}{\sqrt{c_2}},\theta\eta)\chi_{1,R}^\pm(z_2)
\Bigr\}\,d\theta\hat{f}(\eta)e^{iy\eta}\,d\eta\Biggr\|_{\calX_2}
  \\ \lesssim &
\int_{-\eta_0}^{\eta_0}e^{-t\eta^2/2\sqrt{c_2}}|\hat{f}(\eta)|\,d\eta
 \lesssim  (1+t)^{-1/2}\|f\|_{L^1}\,.                
\end{align*}
By the Plancherel theorem,
\begin{align*}
& \left\|\int_{-\eta_0}^{\eta_0}e^{-t\eta^2/2\sqrt{c_2}}\sin(\sqrt{c_2}t\eta)
\int_0^1\pd_\eta\tg_1(\frac{Z_{2,\pm}}{\sqrt{c_2}},\theta\eta)\,d\theta
\chi_{1,R}^\pm(z_2)\hat{f}(\eta)e^{iy\eta}\,d\eta\right\|_{\calX_2}
  \\ \lesssim &
\left\|e^{-t\eta^2/2\sqrt{c_2}}\hat{f}(\eta)
\right\|_{L^2(-\eta_0,\eta_0)}
\left\|e^{\a z_2}\sup_{\eta\in[-\eta_0,\eta_0]}|\pd_\eta\tg_1(z_2,\eta)|\right\|_{L^2(\R)}
\\ \lesssim & (1+t)^{-1/4}\|f\|_{L^1}\,.
\end{align*}
By Proposition~\ref{prop:linearstability-1},
\begin{align*}
\left\|u(t,x+b_1t,y+b_2t)-\sum_{j=1,2}\int_{-\eta_0}^{\eta_0}
a_j(t,\eta)g^2_{2,j}(x,y,\eta)\,d\eta\right\|_{\calX_2} =O(e^{-bt}).
\end{align*}
Combining the above, we have Theorem~\ref{thm:AS}.
Using Proposition~\ref{prop:condstability}, we can prove
Theorem~\ref{thm:AS2} in exactly the same way as Theorem~\ref{thm:AS}.
Thus we compelte the proof.
\end{proof}

\bigskip

\section{$2$-line soliton of O-type}
\label{sec:LS-O}
In this section, we will prove linear stability of $2$-line soliton solutions
of $O$-type.

\subsection{The Darboux transformations and the
  linearized equation}
To begin with, we introduce Jost solutions for $[1,2]$-soliton and
$[3,4]$-soliton. Let
\begin{gather*}
\tau_{11}=f_1=e^{\theta_1}+e^{\theta_2}\,,\quad
\tau_{12}=f_2=e^{\theta_3}+e^{\theta_4}\,,\quad
\tau_2=\operatorname{Wr}(f_1,f_2)\,.
\end{gather*}
Let $u_2=2\pd_x^2\log\tau_2$, $L_2=-\pd_y+\pd_x^2+u_2$ and
for $j=1$ and $2$,
\begin{gather*}
u_{1j}=2\pd_x^2\log\tau_{1j}\,,\quad v_{2j}=\pd_x\log h_{2j}\,,\quad
h_{2j}=\frac{\tau_2}{\tau_{1j}}\,,
\\
L_{1j}=-\pd_y+\pd_x^2+u_{1j}\,,\quad
B_{1j}=4\pd_t+4\pd_x^3+6u_{1j}\pd_x+3\pd_xu_{1j}+3\pd_x^{-1}\pd_yu_{1j}\,,
\\
\mL_{1j}=\mL_0-\frac{3}{2}\pd_x(u_{1j}\cdot)\,.
\end{gather*}
Since $M_+(v_{2j})=u_2$ and $M_-(v_{2j})=u_{1j}$ for $j=1$, $2$,
it follows from \eqref{eq:Miura-Lax1} and \eqref{eq:Miura-Lax2} that
for $j=1$ and $2$,
\begin{gather}
  \label{eq:Mura-Lax1j}
\nabla M_+(v_{2j})=h_{2j}L_2^*h_{2j}^{-1}\pd_x^{-1}
=\pd_x^{-1}h_{2j}L_{1j}^*h_{2j}^{-1}\,,
\\  \label{eq:Miura-Lax2j}
\nabla M_-(v_{2j})=-h_{2j}^{-1}L_{1j}h_{2j}\pd_x^{-1}
=-\pd_x^{-1}h_{2j}^{-1}L_2h_{2j}\,.
\end{gather}
For $j=1$ and $2$, let
\begin{gather*}
\Phi^{1j}(\bx,k)=e^{ikx-k^2y+ik^3t}\frac{(ik-\pd_x)\tau_{1j}}{\tau_{1j}}\,,
 \quad \Phi^{1j,*}(\bx,k)=e^{-ikx+k^2y-ik^3t}\frac{(ik-\pd_x)\tau_{1j}}{\tau_{1j}}\,,
 \\
 \Phi^{11}_1(\bx)=-\Phi^{11}_2(\bx)
 =(\k_1-\k_2)\frac{e^{\theta_1+\theta_2}}{\tau_{11}}\,,
 \quad
  \Phi^{11,*}_1(\bx)=\Phi^{11,*}_2(\bx)=\frac{1}{\tau_{11}}\,,
  \\
\Phi^{12}_3(\bx)=-\Phi^{12}_4(\bx)
 =(\k_3-\k_4)\frac{e^{\theta_3+\theta_4}}{\tau_{12}}\,,
 \quad
 \Phi^{12,*}_3(\bx)=\Phi^{12,*}_4(\bx)=\frac{1}{\tau_{12}}\,.
\end{gather*}
Then $L_{1j}\Phi^{1j}=0$, $B_{1j}\Phi^{1j}=0$,
$L_{1j}^*\Phi^{1j,*}=0$, $B_{1j}^*\Phi^{1j,*}=0$, and it follows from
Lemma~\ref{lem:prodJdJ} that
$\pd_x(\Phi^{1j}\Phi^{1j,*})$ is a solution of
\begin{equation}
  \label{eq:mL1j}
  \pd_tu=\mL_{1j}u\,,
\end{equation}
and that $\Phi^{1j}\Phi^{1j,*}$ is a solution of the adjoint equation
of \eqref{eq:mL1j}.

Let $\Phi^2(\bx,k)$, $\Phi^{2,*}(\bx,k)$, $\Phi^2_n(\bx)$ and $\Phi^{2,*}_n(\bx)$ be as \eqref{def:Phi2}--\eqref{def:Phi2n}. That is,
\begin{gather*}
\Phi^2_1(\bx)=-\Phi^2_2(\bx)=\frac{e^{\theta_1+\theta_2}}{\tau_2}  
(\k_2-\k_1)(\pd_x-\k_1)(\pd_x-\k_2)f_2\,,
\\
\Phi^2_3(\bx)=-\Phi^2_4(\bx)=\frac{e^{\theta_3+\theta_4}}{\tau_2}  
(\k_3-\k_4)(\k_3-\pd_x)(\k_4-\pd_x)f_1\,,
\\
\Phi^{2,*}_1(\bx)=\Phi^{2,*}_2(\bx)=-\frac{f_2}{\tau_2}\,,
\quad
\Phi^{2,*}_3(\bx)=\Phi^{2,*}_4(\bx)=\frac{f_1}{\tau_2}\,.
\end{gather*}
Note that
$\pd_x(\Phi^2\Phi^{2,*})$ is a solution of $\pd_tu=\mL_2u$
and that $\Phi^2\Phi^{2,*}$ is a solution of $\pd_tu+\mL_2^*u=0$.

The Darboux transformations $\nabla M_\pm(v_{2j})$ connect
$\pd_x(\Phi^2\Phi^{2,*})$ with $\pd_x(\Phi^{1j}\Phi^{1j,*})$.

\begin{lemma}
 Let $\beta\in \C$ and $\beta'\in \C\setminus\{\k_1,\k_2,\k_3,\k_4\}$. Then
  \label{lem:M-pdPhiPhi*O}
  \begin{gather*}
\nabla M_+(v_{2j})\pd_x
\left(\Phi^{1j}(\bx,-i\beta)\Phi^{2,*}(\bx,-i\beta')\right)
=2\pd_x\left(\Phi^2(\bx,-i\beta)\Phi^{2,*}(\bx,-i\beta')\right)\,,
\\
\nabla M_-(v_{2j})\pd_x
\left(\Phi^{1j}(\bx,-i\beta)\Phi^{2,*}(\bx,-i\beta')\right)
=2\pd_x\left(\Phi^{1j}(\bx,-i\beta)\Phi^{1j,*}(\bx,-i\beta')\right)\,,
\\
\pd_x\nabla M_-(v_{2j})^*
\left(\Phi^{1j}(\bx,-i\beta)\Phi^{1j,*}(\bx,-i\beta')\right)
=2\pd_x\left(\Phi^2(\bx,-i\beta)\Phi^{1j,*}(\bx,-i\beta')\right)\,,
\\
\pd_x\nabla M_+(v_{2j})^*
\left(\Phi^2(\bx,-i\beta)\Phi^{2,*}(\bx,-i\beta')\right)
=2\pd_x\left(\Phi^2(\bx,-i\beta)\Phi^{1j,*}(\bx,-i\beta')\right)\,.
\end{gather*}
\end{lemma}
We can prove Lemma~\ref{lem:M-pdPhiPhi*O} in exactly the same way as
Lemmas~\ref{lem:M-pdPhiPhi*} and \ref{lem:M-PhiPhi*} by using
Claim~\ref{cl:L0*L0'} below.
\begin{claim}
  \label{cl:L0*L0'}
  For $j=1$, $2$ and $\beta\in\C$,
    \begin{gather*}
      \pd_x\left(h_{2j}^{-1}\Phi^{1j}(\bx,-i\beta)\right)
      =h_{2j}^{-1}\Phi^2(\bx,-i\beta)\,,\\
      \pd_x\left(h_{2j}\Phi^{2,*}(\bx,-i\beta)\right)
      =-h_{2j}\Phi^{1j,*}(\bx,-i\beta)\,,\\
      L_0^*\left(h_{2j}^{-1}\Phi^{1j}(\bx,-i\beta)\right)
      =2h_{2j}^{-1}\pd_x\Phi^2(\bx,-i\beta)\,,\\
      L_0\left(h_{2j}\Phi^{2,*}(\bx,-i\beta)\right)
      =-2h_{2j}\pd_x\Phi^{1j,*}(\bx,-i\beta)\,.
    \end{gather*}
  \end{claim}
We can prove Claim~\ref{cl:L0*L0'} in exactly the same way as
Claim~\ref{cl:L0*L0}.
\bigskip

\subsection{Resonant modes for $2$-soliton solutions of O-type}
\label{subsec:resonance-O}
For the $2$-line soliton solution $u=2\pd_x^2\log\tau$ with $\tau$
satisfying \eqref{eq:O-type}, we consider continuous resonant modes
which corresponds to the modulation of $[1,2]$-soliton and continuous
resonant modes which corresponds to the modulation of $[3,4]$-soliton.
Let $(b_1,b_2)$ and $(X,Y)$ be as in \eqref{eq:b1-b2O} and let $\mL_2$,
$Q^\pm_{ij}(\eta)$, $\widetilde{Q}^\pm_{ij}(\eta)$,
$\lambda_{ij}^\pm(\eta)$, $a_{ij}$ and $X_{ij}$ be as in
Section~\ref{subsec:resonance-P}. 
If $\a$ is a small positive number and $\eta\in\R$ is close
to $0$,  then $\pd_xQ_{12}^-(\eta)$, $\pd_x\widetilde{Q}_{12}^+(\eta)$
are continuous eigenfunctions of $\mL_2$  in $L^2(\R^2;e^{2\a X_{12}}dXdY)$ and
$\pd_xQ_{34}^-(\eta)$, $\pd_x\widetilde{Q}_{34}^+(\eta)$
are continuous eigenfunctions of $\mL_2$ in $L^2(\R^2;e^{2\a X_{34}}dXdY)$.
More precisely, we have the following.
\begin{equation}
  \label{eq:LQ12-34}
\left\{ \begin{gathered}
\mL_2 \pd_XQ_{12}^-(\eta)=\lambda_{12}^-(\eta)\pd_XQ_{12}^-(\eta)\,,
\quad \mL_2 \pd_X\widetilde{Q}_{12}^+(\eta)
=\lambda_{12}^+(\eta)\pd_X\widetilde{Q}_{12}^+(\eta)\,,
\\
\mL_2 \pd_XQ_{34}^-(\eta)=\lambda_{34}^-(\eta)\pd_XQ_{34}^-(\eta)\,,
\quad \mL_2 \pd_X\widetilde{Q}_{34}^+(\eta)
=\lambda_{34}^+(\eta)\pd_X\widetilde{Q}_{34}^+(\eta)\,,
\\
\mL_2^*Q_{12}^+(\eta)=\lambda_{12}^+(\eta)Q_{12}^+(\eta)\,,
\quad \mL_2^*\widetilde{Q}_{12}^-(\eta)
=\lambda_{12}^-(\eta)\widetilde{Q}_{12}^-(\eta)\,,
\\
\mL_2^*Q_{34}^+(\eta)=\lambda_{34}^+(\eta)Q_{34}^+(\eta)\,,
\quad \mL_2^*\widetilde{Q}_{34}^-(\eta)
=\lambda_{34}^-(\eta)\widetilde{Q}_{34}^-(\eta)\,.
 \end{gathered}
\right.  
\end{equation}
\par
Let $(X,Y)$ be as in \eqref{eq:b1-b2O},
$\lambda_{1,\pm}(\eta)=\lambda_{12}^\mp(\pm\eta)$,
$\lambda_{2,\pm}(\eta)=\lambda_{34}^\mp(\pm\eta)$ and 
\begin{align*}
& g^2_{1,+}(X,Y,\eta)=(\k_2-\k_1)d_{1,+}(\eta)\pd_x
\left(\Phi^2(\bx,-i\beta_{12}^-(\eta))\Phi^{2,*}_1(\bx)\right)\,,
\\ &
g^2_{1,-}(X,Y,\eta)=i\eta d_{1,-}(\eta)\pd_x\left(
\Phi^2_1(\bx)\Phi^{2,*}(\bx,-i\beta_{12}^+(-\eta))\right)\,,
\\ &
g^{2,*}_{1,+}(X,Y,\eta)=-\frac{i\bar{\eta}}{\k_2-\k_1}
\Phi^2_1(\bx)\Phi^{2,*}(\bx,-i\beta_{12}^-(-\bar{\eta}))\,,
\\ & 
g^{2,*}_{1,-}(X,Y,\eta)=\Phi^2(\bx,-i\beta_{12}^+(\bar{\eta}))\Phi^{2,*}_1(\bx)\,,
\\ &
g^2_{2,+}(X,Y,\eta)=(\k_4-\k_3) d_{2,+}(\eta)
\pd_x\left(\Phi^2(\bx,-i\beta_{34}^-(\eta)) \Phi^{2,*}_3(\bx)\right)\,,
\\ &
g^2_{2,-}(X,Y,\eta)=i\eta d_{2,-}(\eta)\pd_x\left(
\Phi^2_3(\bx)\Phi^{2,*}(\bx,-i\beta_{34}^+(-\eta)\right)\,,
\\ &
g^{2,*}_{2,+}(X,Y,\eta)=-\frac{i\bar{\eta}}{\k_4-\k_3}
\Phi^2_3(\bx)\Phi^{2,*}(\bx,-i\beta_{34}^-(-\bar{\eta}))\,,
\\ & 
g^{2,*}_{2,-}(X,Y,\eta)=\Phi^2(\bx,-i\beta_{34}^+(\bar{\eta}))\Phi^{2,*}_3(\bx)\,,
\end{align*}
where 
$d_{1,\pm}(\eta)=\pm 1/\gamma_{12}(\pm\eta)$ and
$d_{2,\pm}(\eta)=\pm 1/\gamma_{34}(\pm\eta)$.
We remark that $g^2_{1,\pm}$ and $g^{2,*}_{1,\pm}$ are continuous
eigenfunctions of $\mL_2$ and $\mL_2^*$ associated with
modulations of $[1,2]$-soliton and that 
$g^2_{2,\pm}$ and $g^{2,*}_{2,\pm}$ are continuous
eigenfunctions of $\mL_2$ and $\mL_2^*$ associated with
modulations of $[3,4]$-soliton, respectively.
\par
Let
\begin{align*}
& g^1_{1,+}(X,Y,\eta)=(\k_2-\k_1)d_{1,+}(\eta)\pd_x
\left(\Phi^{11}(\bx,-i\beta_{12}^-(\eta))\Phi^{11,*}_1(\bx)\right)\,,
\\ &
g^1_{1,-}(X,Y,\eta)=i\eta d_{1,-}(\eta)\pd_x\left(
\Phi^{11}_1(\bx)\Phi^{11,*}(\bx,-i\beta_{12}^+(-\eta))\right)\,,
\\ &
g^{1,*}_{1,+}(X,Y,\eta)=-\frac{i\bar{\eta}}{\k_2-\k_1}
\Phi^{11}_1(\bx)\Phi^{11,*}(\bx,-i\beta_{12}^-(-\bar{\eta}))\,,
\\ & 
g^{1,*}_{1,-}(X,Y,\eta)=\Phi^{11}(\bx,-i\beta_{12}^+(\bar{\eta}))
\Phi^{11,*}_1(\bx)\,,
\\ &
g^1_{2,+}(X,Y,\eta)=(\k_4-\k_3) d_{2,+}(\eta)
\pd_x\left(\Phi^{12}(\bx,-i\beta_{34}^-(\eta)) \Phi^{12,*}_3(\bx)\right)\,,
\\ &
g^1_{2,-}(X,Y,\eta)=i\eta d_{2,-}(\eta)\pd_x\left(
\Phi^{12}_3(\bx)\Phi^{12,*}(\bx,-i\beta_{34}^+(-\eta)\right)\,,
\\ &
g^{1,*}_{2,+}(X,Y,\eta)=-\frac{i\bar{\eta}}{\k_4-\k_3}
\Phi^{12}_3(\bx)\Phi^{12,*}(\bx,-i\beta_{34}^-(-\bar{\eta}))\,,
\\ & 
g^{1,*}_{2,-}(X,Y,\eta)=\Phi^{12}(\bx,-i\beta_{34}^+(\bar{\eta}))
\Phi^{12,*}_3(\bx)\,.
\end{align*}
\begin{remark}
We remark that $g^1_{1,\pm}$ and $g^1_{2,\pm}$ are continuous eigenfunctions of
$\mL_{11}$ and $\mL_{12}$, respectively and that
$g^{1,*}_{1,\pm}$ and $g^{1,*}_{2,\pm}$ are continuous eigenfunctions of $\mL_{11}^*$
and $\mL_{12}^*$, respectively.
Let $z_1=x+2a_{12}y$, $z_2=x+2a_{34}y$ and
  \begin{gather*}
g_1(z,\eta)=\sqrt{c_{12}}d_{1,+}(\eta)\pd_{z}^2\left\{
e^{-\gamma_{12}(\eta)z}\sech\sqrt{c_{12}}z\right\}\,\\
g_1^*(z,\eta)=\frac12\pd_{z}\left\{
e^{\gamma_{12}(-\bar{\eta})z}\sech\sqrt{c_{12}}z\right\}\,,
\\
g_2(z,\eta)=\sqrt{c_{34}}d_{2,+}(\eta)\pd_z^2\left\{
e^{-\gamma_{34}(\eta)z}\sech\sqrt{c_{34}}z\right\}\,\\
g_2^*(z,\eta)=\frac12\pd_z\left\{
e^{\gamma_{34}(-\bar{\eta})z}\sech\sqrt{c_{34}}z\right\}\,.
\end{gather*}
Then
\begin{gather*}
  g^1_{1,\pm}(x,y,\eta)=e^{iy\eta}g_1(z_1,\pm\eta)\,,\quad
  g^{1,*}_{1,\pm}(x,y,\eta)=e^{iy\bar{\eta}}g_1^*(z_1,\pm\eta)\,,
\\  
  g^1_{2,\pm}(x,y,\eta)=e^{iy\eta}g_2(z_2,\pm\eta)\,,\quad
  g^{1,*}_{2,\pm}(x,y,\eta)=e^{iy\bar{\eta}}g_2^*(z_2,\pm\eta)\,.
\end{gather*}
\end{remark}

By \eqref{eq:LQ12-34}, the definitions of $g^i_{j,\pm}$, $g^{i,*}_{j,\pm}$
and $\lambda_{j,\pm}$ and an argument in Section~\ref{subsec:resonance-1},
we have the following.
\begin{lemma}
  \label{lem:evmL2'}
For $j=1$ and $2$,
\begin{gather*}
\mL_2 g^2_{j,\pm}(x,y,\eta) =\lambda_{j,\pm}(\eta)g^2_{j,\pm}(x,y,\eta)\,,\quad
\mL_2^* g^{2,*}_{j,\pm}(x,y,\eta)  =
\overline{\lambda_{j,\pm}(\eta)}g^{2,*}_{j,\pm}(x,y,\eta)\,,
\\
\mL_{1j} g^1_{j,\pm}(x,y,\eta) =\lambda_{j,\pm}(\eta)g^1_{j,\pm}(x,y,\eta)\,,\quad
\mL_{1j}^* g^{1,*}_{j,\pm}(x,y,\eta)  =
\overline{\lambda_{j,\pm}(\eta)}g^{1,*}_{j,\pm}(x,y,\eta)\,,
\end{gather*}
where $\lambda_{1,\pm}(\eta)=i\eta\{b_2-3a_{12}\pm\gamma_{12}(\pm\eta)\}$
and $\lambda_{2,\pm}(\eta)=i\eta\{b_2-3a_{34}\pm\gamma_{34}(\pm\eta)\}$.
\end{lemma}
\par

Letting $\beta=\k_1$, $\k_3$ or calculating residues at $\beta'=\k_1$, $\k_3$
in Lemma~\ref{lem:M-pdPhiPhi*O} and combining the result with
Claim~\ref{cl:L0*L0'},
we see that $\nabla M_\pm(v_{2j})$
connect $g^2_{j,\pm}$ with $g^1_{j,\pm}$ and that
$\nabla M_\pm(v_{2j})^*$ connect $g^{2,*}_{j,\pm}$ with $g^{1,*}_{j,\pm}$.
Let $(X,Y)$ be a moving coordinate defined by \eqref{eq:b1-b2O} and
\begin{align*}
  & g^M_{1,+}(X,Y,\eta)=(\k_2-\k_1)d_{1,+}(\eta)\pd_x
\left(\Phi^{11}(\bx,-i\beta_{12}^-(\eta))\Phi^{2,*}_1(\bx)\right)\,,
\\ &
g^M_{1,-}(X,Y,\eta)=i\eta d_{1,-}(\eta)\pd_x\left(
\Phi^{11}_1(\bx)\Phi^{2,*}(\bx,-i\beta_{12}^+(-\eta))\right)\,,
\\ &
g^{M,*}_{1,+}(X,Y,\eta)=-\frac{i\bar{\eta}}{\k_2-\k_1}
\Phi^2_1(\bx)\Phi^{11,*}(\bx,-i\beta_{12}^-(-\bar{\eta}))\,,
\\ & 
g^{M,*}_{1,-}(X,Y,\eta)=\Phi^2(\bx,-i\beta_{12}^+(\bar{\eta}))\Phi^{11,*}_1(\bx)\,,
\\ &
g^M_{2,+}(X,Y,\eta)=(\k_4-\k_3) d_{2,+}(\eta)
\pd_x\left(\Phi^{12}(\bx,-i\beta_{34}^-(\eta)) \Phi^{2,*}_3(\bx)\right)\,,
\\ &
g^M_{2,-}(X,Y,\eta)=i\eta d_{2,-}(\eta)\pd_x\left(
\Phi^{12}_3(\bx)\Phi^{2,*}(\bx,-i\beta_{34}^+(-\eta)\right)\,,
\\ &
g^{M,*}_{2,+}(X,Y,\eta)=-\frac{i\bar{\eta}}{\k_4-\k_3}
\Phi^2_3(\bx)\Phi^{12,*}(\bx,-i\beta_{34}^-(-\bar{\eta}))\,,
\\ & 
g^{M,*}_{2,-}(X,Y,\eta)=
\Phi^2(\bx,-i\beta_{34}^+(\bar{\eta}))\Phi^{12,*}_3(\bx)\,.
\end{align*}

\begin{lemma}
  \label{lem:g-tg'}
For $j=1$ and $2$, 
\begin{gather*}
\nabla M_+(v_{2j})g^M_{j,\pm}=2g^2_{j,\pm}\,,\quad
\nabla M_-(v_{2j})g^M_{j,\pm}=2g^1_{j,\pm}\,,
\\
\nabla M_+(v_{2j})^*g^{2,*}_{j,\pm}=2g^{M,*}_{j,\pm}\,,\quad
\nabla M_-(v_{2j})^*g^{1,*}_{j,\pm}=2g^{M,*}_{j,\pm}\,,
\\
\nabla M_-(v_{2j})\pd_x^{-1}g^2_{j,\pm}
=\nabla M_+(v_{2j})\pd_x^{-1}g^1_{j,\pm}\,.
\end{gather*}
\end{lemma}

The Darboux transformations $\nabla M_\pm(v_{22})$ and $\nabla M_\pm(v_{21})$
connect $0$ with $g^2_{1,+}$ and $g^2_{2,+}$, respectively.
Let $(X,Y)$ be as \eqref{eq:b1-b2O} and 
\begin{align*}
& \tg^M_1(X,Y,\eta)=
(\k_2-\k_1)d_{1,+}(\eta)\Phi^2(\bx,-i\beta_{12}^-(\eta))\Phi^{2,*}_1(\bx)\,,
\\ &
\tg^{M,*}_1(X,Y,\eta)=\frac{1}{\k_2-\k_1}
     \Phi^2_1(\bx)\Phi^{12,*}(\bx,-i\beta_{12}^-(-\bar{\eta}))\,,
\\ & \tg^M_2(X,Y,\eta)=
(\k_4-\k_3)d_{2,+}(\eta)\Phi^2(\bx,-i\beta_{34}^-(\eta))\Phi^{2,*}_3(\bx)\,,
\\ &
\tg^{M,*}_2(X,Y,\eta)=\frac{1}{\k_4-\k_3}
     \Phi^2_3(\bx)\Phi^{11,*}(\bx,-i\beta_{34}^-(-\bar{\eta}))\,.
\end{align*}

\begin{lemma} \label{lem:eigenfunctions-mKP'}
\begin{gather}
  \label{eq:38}
 \nabla M_-(v_{22})\tg^M_1=0\,,\quad \nabla M_+(v_{22})^*g^{2,*}_{1,-}=0\,,
  \\ \label{eq:39}
  \nabla M_+(v_{22})\tilde{g}^M_1=2g^2_{1,+}\,,
  \quad \nabla M_+(v_{22})^*g^{2,*}_{1,1}=2\tilde{g}^{M,*}_1\,,
  \\   \label{eq:38'}
  \nabla M_-(v_{21})\tg^M_2=0\,,\quad \nabla  M_+(v_{21})^*g^{2,*}_{2,-}=0\,,
  \\ \label{eq:39'}
  \nabla M_+(v_{21})\tilde{g}^M_2=2g^2_{2,+}\,,
  \quad \nabla M_+(v_{21})^*g^{2,*}_{2,1}=2\tilde{g}^{M,*}_2\,.
\end{gather}
\end{lemma}
\begin{proof}
  Equations~\eqref{eq:38} and
  \eqref{eq:38'} follow from \eqref{eq:Miura-Lax2j} and the fact that
  $L_2\Phi^2=0$.  Letting $\beta=\beta_{12}^-(\eta)$ and taking a
  residue at $\beta'=\k_1$ or letting $\beta=\k_1$ and
  $\beta'=\beta_{12}^-(-\bar{\eta})$, we have \eqref{eq:39} from
  Lemma~\ref{lem:M-pdPhiPhi*O}.  We can prove 
  \eqref{eq:39'} in the same way.
\end{proof}

\bigskip

\subsection{Orthogonal relations}
Combining Lemma~\ref{lem:g-tg'} with Lemma~\ref{lem:orthgz} and
following the proof fo Lemma~\ref{lem:orthogonality[2,3]}, we have
the orthogonality relations for $g^i_{j,\pm}$ and $g^{i,*}_{j,\pm}$.
\begin{lemma}
  \label{lem:orthogonality'}
  Let $0<\eta_0<4\sqrt{2}\min\{\sqrt{c_{12}},\sqrt{c_{34}}\}$ and
  $\eta\in(-\eta_0,\eta_0)$.
Then for $i$, $j=1$, $2$ and $\varphi\in C^1(\R)$,
\begin{gather*}
\lim_{M\to\infty}\int_{-\eta_0}^{\eta_0}\varphi(\eta_1)
\left(\int_{-M}^M \int_{-\infty}^\infty
g^i_{j,\pm}(x,y,\eta)\overline{g^{i,*}_{j,\pm}(x,y,\eta_1)}\,dxdy\right)\,d\eta_1
=\pm2\pi i\eta\varphi(\eta)\,,
\\
\lim_{M\to\infty}\int_{-\eta_0}^{\eta_0}\varphi(\eta_1)\left(\int_{-M}^M
 \int_{-\infty}^\infty
g^k_{j,\pm}(x,y,\eta)\overline{g^{k,*}_{j,\mp}(x,y,\eta_1)}\,dxdy\right)\,d\eta_1=0\,.
\end{gather*}
\end{lemma}
Let $g^i_{j,k}(\eta)$ and $g^{i,*}_{j,k}(\eta)$ be as \eqref{def:gjk12}.
\begin{lemma}
  \label{lem:orthogonality2'}
 Let  $0<\eta_0<4\sqrt{2}\min\{c_{12},c_{34}\}$ and $\eta\in(-\eta_0,\eta_0)$.
Then for $i$, $j$, $k=1$, $2$ and and  $\varphi\in C^1(\R)$,
\begin{gather*}
\lim_{M\to\infty}\int_{-\eta_0}^{\eta_0}\varphi(\eta_1)\left(
  \int_{-M}^M \int_{-\infty}^\infty
g^i_{j,k}(x,y,\eta)\overline{g^{i,*}_{j,k}(x,y,\eta_1)}\,dxdy\right)\,d\eta_1
=2\pi\delta_{jk}\varphi(\eta)\,.
\end{gather*}
\end{lemma}
Lemma~\ref{lem:orthogonality2'} follows immediately from
Lemma~\ref{lem:orthogonality'}.
\bigskip

\subsection{Spectral projections}
Let $\a>0$ and
$$\calX_3=L^2(\R^2;e^{2\a z_1}dz_1dy)\,,\quad
\calX_4=L^2(\R^2;e^{2\a z_2}dz_2dy)\,.$$
For $i$, $j=1$ and $2$,
let $P_{ij}(\eta_0)$ be an operator defined by
\begin{gather}
  P_{ij}(\eta_0)f(x,y)=\frac{1}{\sqrt{2\pi}}\sum_{k=1,2}\int_{-\eta_0}^{\eta_0}
a_{ijk}(\eta)g^i_{j,k}(x,y,\eta)\,d\eta\,,
\\
a_{ijk}(\eta)=\frac{1}{\sqrt{2\pi}}\lim_{M\to\infty}\int_{-M}^M\int_\R f(x,y)
\overline{g^{i,*}_{j,k}(x,y,\eta)}\,dxdy\,.
\end{gather}
Let 
\begin{equation}
  \label{def:eta*'}
\eta_{*,3}=2(\sqrt{c_{12}}+\a)\sqrt{\a(\a+\k_2-\k_1)}\,,\quad
\eta_{*,4}=2(\sqrt{c_{34}}+\a)\sqrt{\a(\a+\k_4-\k_3)}\,.
\end{equation}
Then
$\Re\gamma_{12}(\eta_{*,3})=\sqrt{c_{12}}+\a$ and $\Re\gamma_{34}(\eta_{*,4})=\sqrt{c_{34}}+\a$.
We can prove that for $i$, $j=1$, $2$ and $\eta_0\in(0,\eta_{*,j})$,
the operator $P_{ij}(\eta_0)$ is a spectral projection on $\calX_j$
in the same way as Section~\ref{subsec:specproj},
\begin{lemma}
  \label{lem:spectral-proj'}
 Assume $0<\a<2\min\{\sqrt{c_{12}},\sqrt{c_{34}}\}$.
Let  $\eta'\in(0,\eta_{*,j})$ and
$\mL_2$ be a closed operator in $\calX_{j+2}$ for $j=1$ or $2$.
Suppose that $\eta_0\in(0,\eta']$. Then
\begin{enumerate}
\item
$\|P_{2j}(\eta_0)f\|_{\calX_j}\le C\|f\|_{\calX_j}$  for any $f\in \calX_j$,
where $C$ is a positive constant depending only on $\a$ and $\eta'$,
\item
$\mL_2 P_{2j}(\eta_0)f=P_{2j}(\eta_0)\mL_2 f$ for any $f\in D(\mL_2)$,
\item
$P_{2j}(\eta_0)^2=P_{2j}(\eta_0)$ on $\calX_{j+2}$,
\item
$e^{t\mL_2}P_{2j}(\eta_0)=P_{2j}(\eta_0)e^{t\mL_2}$ on $\calX_{j+2}$.
\end{enumerate}
\end{lemma}

To prove that $P_{2j}(\eta_0)$ is bounded on $\calX_{j+2}$,
we will investigate the asymptotic behavior of 
$g^2_{j,k}$ and $g^{2,*}_{j,k}$ as $y\to\pm\infty$.
Let $Z_{1,\pm}$ and $Z_{2,\pm}$ be as \eqref{eq:u2-asymptoticsO}.

\begin{lemma}  
  \label{lem:gg*-asymp-O}
As $Z_{2,+}\to\infty$,
  \begin{gather*}
 g^2_{1,\pm}(x,y,\eta)=c_{\pm,4}(\eta)
g_1(\frac{Z_{1,+}}{\sqrt{c_{12}}},\pm\eta)
e^{iy\eta}(1+O(e^{-2Z_{2,+}}))\,,
\\
g^{2,*}_{1,\pm}(x,y,\eta)=\frac{1}{c_{\pm,4}(-\bar{\eta})}
g^*_1(\frac{Z_{1,+}}{\sqrt{c_{12}}},\pm\eta)e^{iy\bar{\eta}}
(1+O(e^{-2Z_{2,+}}))\,,
\end{gather*}
and as $Z_{2,-}\to-\infty$,
\begin{gather*}
g^2_{1,\pm}(x,y,\eta)=c_{\pm,3}(\eta)g_1(\frac{Z_{1,-}}{\sqrt{c_{12}}},\pm\eta)
e^{iy\eta}(1+O(e^{2Z_{2,-}}))\,,
\\
g^{2,*}_{1,\pm}(x,y,\eta)
=\frac{1}{c_{\pm,3}(-\bar{\eta})}
g^*_1(\frac{Z_{1,-}}{\sqrt{c_{12}}},\pm\eta)e^{iy\bar{\eta}}
(1+O(e^{2Z_{2,-}}))\,.
\end{gather*}
where for $j=3$ and $4$,
  \begin{align*}
 & c_{+,j}(\eta)=
(\kappa_j-\beta_{12}^-(\eta))
(\kappa_j-\kappa_2)^{-\frac12(1-\gamma_{12}(\eta)/\sqrt{c_{12}})}
(\kappa_j-\kappa_1)^{-\frac12(1+\gamma_{12}(\eta)/\sqrt{c_{12}})}\,,
\\
& c_{-,j}(\eta)=
\frac{(\kappa_j-\kappa_1)^{\frac12(1-\gamma_{12}(-\eta)/\sqrt{c_{12}})}
(\kappa_j-\kappa_2)^{\frac12(1+\gamma_{12}(-\eta)/\sqrt{c_{12}})}}
{\k_j-\beta_{12}^+(-\eta)}\,.
  \end{align*}

As $Z_{1,+}\to\infty$,
  \begin{gather*}
g^2_{2,\pm}(x,y,\eta)=c_{\pm,3}(\eta)g_2(\frac{Z_{2,+}}{\sqrt{c_{34}}},\pm\eta)
e^{iy\eta}(1+O(e^{-2Z_{1,+}}))\,,
\\
g^{2,*}_{2,\pm}(x,y,\eta)=\frac{1}{c_{\pm,3}(-\bar{\eta})}
g^*_2(\frac{Z_{2,+}}{\sqrt{c_{34}}},\pm\eta)e^{iy\bar{\eta}}
(1+O(e^{-2Z_{1,+}}))\,,
\end{gather*}
and as $Z_{1,-}\to-\infty$,
\begin{gather*}
g^2_{2,\pm}(x,y,\eta)=c_{\pm,2}(\eta)g_2(\frac{Z_{2,-}}{\sqrt{c_{34}}},\pm\eta)
e^{iy\eta}(1+O(e^{2Z_{1,-}}))\,,
\\
 g^{2,*}_{2,\pm}(x,y,\eta)
 =\frac{1}{c_{\pm,2}(-\bar{\eta})}g_2^*(\frac{Z_{2,-}}{\sqrt{c_{34}}},\pm\eta)
 e^{iy\bar{\eta}}
(1+O(e^{2Z_{1,-}}))\,,
\end{gather*}
where for $j=1$ and $2$,
\begin{align*}
& c_{+,j}(\eta)=
(\beta_{34}^-(\eta)-\k_j)
(\kappa_3-\kappa_j)^{-\frac12(1+\gamma_{34}(\eta)/\sqrt{c_{34}})}
(\kappa_4-\kappa_j)^{-\frac12(1-\gamma_{34}(\eta)/\sqrt{c_{34}})}\,,
\\
& c_{-,j}(\eta)=
\frac{(\kappa_3-\kappa_j)^{\frac12(1-\gamma_{34}(-\eta)/\sqrt{c_{34}})}
(\kappa_4-\kappa_j)^{\frac12(1+\gamma_{34}(-\eta)/\sqrt{c_{34}})}}
{\beta_{34}^+(-\eta)-\k_j}\,.
\end{align*}
\end{lemma}

Using Lemma~\ref{lem:gg*-asymp-O}, we can prove Lemma~\ref{lem:spectral-proj'}
in the same manner as the argument in Section~\ref{subsec:specproj}.
To prove Lemma~\ref{lem:gg*-asymp-O}, we need the following.
\begin{claim}
  \label{cl:Phi22*s'}
  \begin{gather*}
\Phi^2(\bx,-i\beta)\Phi^{2,*}_1(\bx)=
\pd_x\left(\Phi^{12}(\bx,-i\beta)\Phi^{2,*}_1(\bx)\right)\,,
\\ 
\Phi^2(\bx,-i\beta)\Phi^{2,*}_3(\bx)=
\pd_x\left(\Phi^{11}(\bx,-i\beta)\Phi^{2,*}_3(\bx)\right)\,,
\\ 
\Phi^2_1(\bx)\Phi^{2,*}(\bx,-i\beta)=
\frac{\k_1-\k_2}{(\beta-\k_1)(\beta-\k_2)}
\pd_xJ_1\,,
\\
\Phi^2_3(\bx)\Phi^{2,*}(\bx,-i\beta)=
\frac{\k_4-\k_3}{(\beta-\k_3)(\beta-\k_4)}\pd_xJ_2\,,
\end{gather*}
where
\begin{align*}
& J_1=\Phi^{0,*}(\bx,-i\beta)\frac{e^{\theta_1+\theta_2}}{\tau_2}
\{(\k_3-\k_1)(\k_3-\k_2)\frac{e^{\theta_3}}{\beta-\k_3}
+(\k_4-\k_1)(\k_4-\k_2)\frac{e^{\theta_4}}{\beta-\k_4}\}\,,
\\ &
J_2=\Phi^{0,*}(\bx,-i\beta)\frac{e^{\theta_3+\theta_4}}{\tau_2}
\left\{(\k_3-\k_1)(\k_4-\k_1)\frac{e^{\theta_1}}{\beta-\k_1}
+(\k_3-\k_2)(\k_4-\k_2)\frac{e^{\theta_2}}{\beta-\k_2}\right\}\,.
\end{align*}
\end{claim}
Claim~\ref{cl:Phi22*s'} follows from  straightforward computations.

\begin{proof}[Proof of Lemma~\ref{lem:gg*-asymp-O}]
\begin{align*}
\tau_2e^{-(\theta_1+\theta_1)/2}=&
2\sqrt{(\k_3-\k_1)(\k_3-\k_2)}e^{\theta_3}\cosh Z_{1,-}
\\ & +2\sqrt{(\k_4-\k_1)(\k_4-\k_2)}e^{\theta_4}\cosh Z_{1,+} \,,        
\end{align*}
\begin{align*}
\tau_2e^{-(\theta_3+\theta_4)/2}=&
2\sqrt{(\k_3-\k_1)(\k_4-\k_1)}e^{\theta_1}\cosh Z_{2,-}
\\ & + 2\sqrt{(\k_3-\k_2)(\k_4-\k_2)}e^{\theta_2}\cosh Z_{2,+}\,,
\end{align*}
By \eqref{eq:btitj} and Claim~\ref{cl:Phi22*s'},
\begin{gather}
\label{eq:g21pm-O}  
g^2_{1,\pm}(x,y,\eta)=e^{iy\eta}\sum_{j=3,4}\pd_x^2g^2_{1,\pm,j}(x,y,\eta)\,,
\quad
g^{2,*}_{1,\pm}(x,y,\eta)=e^{iy\bar{\eta}}\sum_{j=3,4}
\pd_xg^{2,*}_{1,\pm,j}(x,y,\eta)\,,
\\ \label{eq:g22pm-O}
g^2_{2,\pm}(x,y,\eta)=e^{iy\eta}\sum_{j=1,2}\pd_x^2g^2_{2,\pm,j}(x,y,\eta)\,,
\quad
g^{2,*}_{2,\pm}(x,y,\eta)=e^{iy\bar{\eta}}\sum_{j=1,2}\pd_xg^{2,*}_{2,\pm,j}(x,y,\eta)\,,
\end{gather}
where
\begin{align*}
  & 
g^2_{1,+,j}=(\k_2-\k_1)d_{1,+}(\eta)(\k_j-\beta_{12}^-(\eta))
e^{-\gamma_{12}(\eta)z_1}
\frac{e^{\theta_j+(\theta_1+\theta_2)/2}}{\tau_2}\bigr|_{t=0}\,,
\\ &
 g^2_{1,-,j}=(\k_2-\k_1)d_{1,-}(\eta)
\frac{(\k_j-\k_1)(\k_j-\k_2)}{\beta_{12}^+(-\eta)-\k_j}
e^{-\gamma_{12}(-\eta)z_1}
\frac{e^{\theta_j+(\theta_1+\theta_2)/2}}{\tau_2}\bigr|_{t=0}\,
\\ &
  g^{2,*}_{1,+,j}=\frac{(\k_j-\k_1)(\k_j-\k_2)}{\k_j-\beta_{12}^-(-\bar{\eta})}
e^{\gamma_{12}(-\bar{\eta})z_1}
\frac{e^{\theta_j+(\theta_1+\theta_2)/2}}{\tau_2}\bigr|_{t=0}\,
\\ &
g^{2,*}_{1,-,j}=(\k_j-\beta_{12}^+(\bar{\eta}))e^{\gamma_{12}(\bar{\eta})z_1}
    \frac{e^{\theta_j+(\theta_1+\theta_2)/2}}{\tau_2}\bigr|_{t=0}\,,   
\end{align*}
\begin{align*}
  &
g^2_{2,+,j}=(\k_4-\k_3)d_{2,+}(\eta)(\beta_{34}^-(\eta)-\k_j)
e^{-\gamma_{34}(\eta)z_2}
\frac{e^{\theta_j+(\theta_3+\theta_4)/2}}{\tau_2}\bigr|_{t=0}\,,
\\ & 
 g^2_{2,-,j}=(\k_3-\k_4)d_{2,-}(\eta)
\frac{(\k_3-\k_j)(\k_4-\k_j)}{\beta_{34}^+(-\eta)-\k_j}
e^{-\gamma_{34}(-\eta)z_2}
\frac{e^{\theta_j+(\theta_3+\theta_4)/2}}{\tau_2}\bigr|_{t=0}\,
\\ & 
  g^{2,*}_{2,+,j}=\frac{(\k_3-\k_j)(\k_4-\k_j)}
  {\beta_{34}^-(-\bar{\eta})-}
e^{\gamma_{34}(-\bar{\eta})z_2}
    \frac{e^{\theta_j+(\theta_3+\theta_4)/2}}{\tau_2}\bigr|_{t=0}\,
\\ & 
g^{2,*}_{2,-,j}=(\beta_{34}^+(\bar{\eta})-\k_j)e^{\gamma_{34}(\bar{\eta})z_2}
    \frac{e^{\theta_j+(\theta_3+\theta_4)/2}}{\tau_2}\bigr|_{t=0}\,.
\end{align*}
Combining the above, we have Lemma~\ref{lem:gg*-asymp-O}.
\end{proof}

\bigskip

\subsection{Spectral  stability of $[1,2]$-soliton}
 We can prove spectral stability of $[1,2]$-soliton $u_{11}$ in $\calX_3$
in the same way as Proposition~\ref{lem:wmL1-decay}.
\begin{lemma}
  \label{lem:linear-stability[1,2]}
Assume $\a\in(0,\sqrt{c_{12}})$.
Let $\mL_0$ and $\mL_{11}$ be closed operators on $\calX_3$. Let
$\eta_{*,3}$ be as \eqref{def:eta*'} and $\eta_1\in(0,\eta_{*,3})$.
Then
$\sigma(\mL_{11})=\sigma(\mL_0)\cup\{\lambda_{1,\pm}(\eta)\mid
  \eta\in[-\eta_{*,3},\eta_{*,3}]\}$. Moreover,
 $\lambda_{1,\pm}(\eta)-\mL_{11}$ is invertible on
 $\left(I-P_{11}(\eta_0)\right)\calX_3$ if $\eta_0\in(0,\eta_{*,3})$ and
 $\eta\in(-\eta_0,\eta_0)$.
\end{lemma}

We will show spectral stability of $u_{11}$ in $\calX_4$.
\begin{lemma}
  \label{lem:mL-11inX2}
Assume
\begin{equation}
\label{eq:ass-alpha'}
0<\a<\min\left\{\sqrt{c_{12}}\,,\,
\frac{c_{12}}{2(a_{34}-a_{12})}\,,\,2(a_{34}-\k_2)\right\}\,.
\end{equation}
Suppose that $\mL_{11}$ is a closed operator in $\calX_4$.
Then there exists a positive number $b$ such that
\begin{equation*}
\sup_{\Re\lambda\ge -b}\|(\lambda-\mL_{11})^{-1}\|_{B(\calX_4)}<\infty\,.
\end{equation*}
\end{lemma}

Let $\eta=\eta_R+i\eta_I$, $\eta_R\in\R$ and $\eta_I=2(a_{34}-a_{12})\a$.
Then
 $\a z_2+iy\eta=\a z_1+iy\eta_R$ and
\begin{gather*}
\calX_4=e^{-y\eta_I}\calX_3=\bigoplus_{\eta=\eta_R+i\eta_I,\,\eta_R\in\R}
  e^{iy\eta}L^2(\R;e^{2\a z_1}dz_1)\,.
\end{gather*}
Let
\begin{align*}
\mL_{11}(\eta)v(z_1):=& e^{-iy\eta}\mL_{11}\left(v(z_1)e^{iy\eta}\right)
  \\ = &
-\frac{1}{4}\pd_{z_1}(\pd_{z_1}^2-4c_{12}+6u_{11})
+i(b_2-3a_{12})\eta+\frac{3}{4}\eta^2\pd_{z_1}^{-1}\,.
\end{align*}
In view of Lemma~\ref{lem:evmL2'},
\begin{equation*}
\mL_{11}(\eta)g_1(z_1,\eta)=\lambda_{1,+}(\eta)g(z_1,\eta)\,,
\quad
\mL_{11}(\eta)g_1(-z_1,\eta)=\tilde{\lambda}_{1,+}(\eta)g(z_1,\eta)\,,
\end{equation*}
where
$\tilde{\lambda}_{1,+}(\eta)=i\eta\{b_2-3a_{12}-\gamma_{12}(\eta)\}$.
Since
\begin{equation*}
  \|e^{\a z_2}g^1_{1,\pm}(x,y,\eta)\|_{L^\infty_yL^2_x}=
  \|e^{\a z_1}g_1(z_1,\pm\eta)\|_{L^2(\R_{z_1})}\,,
\end{equation*}
we see that 
$g^1_{1,\pm}(x,y,\eta)$ are continuous eigenfunctions of
$\mL_{11}$ if $\Re\gamma_{12}(\pm\eta)<\sqrt{c_{12}}+\a$.
On the other hand, we have
$\|e^{\a z_1}g_1(-z_1,\pm\eta)\|_{L^2(\R_{z_1})}<\infty$
provided $\Re\gamma_{12}(\pm\eta)<\sqrt{c_{12}}-\a$.

Suppose that $\a\in(0,\sqrt{c_{12}})$ and $\eta_I<c_{12}$.
Let
\begin{gather}
  \label{def:e*+'}
\eta_{*,+}=2(\sqrt{c_{12}}+\a)\sqrt{\a(\a+2a_{34}-2\k_1)}\,,
\\   \label{def:esh+'}
  \eta_{\#,+}=2(\sqrt{c_{12}}-\a)\sqrt{\a(\a+2a_{34}-2\k_2)}\,.
\end{gather}
Then
$\Re\gamma_{12}(\eta_{\#,+}+i\eta_I)=\sqrt{c_{12}}-\a$ and
$\Re\gamma_{12}(\eta_{*,+}+i\eta_I)=\sqrt{c_{12}}+\a$.
Since $\Re\gamma_{12}(\cdot+i\eta_I)$ is even and monotone increasing
on $(0,\infty)$, we have
$e^{\a z_1}g_1(z_1,\eta)\in L^2(\R)$
for $\eta_R\in(-\eta_{*,+},\eta_{*,+})$,
$e^{\a z_1}g_1(-z_1,\eta)\in L^2(\R)$ for $\eta_R\in(-\eta_{\#,+},\eta_{\#,+})$,
and
$$
\{\lambda_{1,+}(\eta)\mid \eta_R\in[-\eta_{*,+},\eta_{*,+}]\}
\cup \{\tilde{\lambda}_{1,+}(\eta)\mid \eta_R\in[-\eta_{\#,+},\eta_{\#,+}]\}
\subset \sigma_c(\mL_{11})\,.$$
On the other hand, if $0<\a<2(a_{34}-\k_2)$, then for every $\eta_R\in\R$,
$$\Re\gamma_{12}(-\eta)\ge \gamma_{12}(-i\eta_I)>\sqrt{c_{12}}+\a\,,$$
and $\{e^{iy\eta}g_1(\pm z_1,-\eta)\}$ are not continuous eigenfunctions
of $\mL_{11}$.

\par
Next, we will show that $\{\lambda_{1,+}(\eta)\}$ and
$\{\tilde{\lambda}_{1,+}(\eta)\}$ belong to the stable half plane.
\begin{lemma}
  \label{lem:lambda1+'}
Suppose that $\a$ satisfies \eqref{eq:ass-alpha'}.
Let $\eta=\eta_R+i\eta_I$ and $\eta_I=2\a(a_{34}-a_{12})$. Then 
$\lambda_{1,+}(i\eta_I)<\tilde{\lambda}_{1,+}(i\eta_I)$ and
\begin{gather}
  \label{eq:lem:lambda1+'}
\eta_R\pd_{\eta_R}\Re \lambda_{1,+}(\eta)=
-\eta_R\pd_{\eta_R}\Re \tilde{\lambda}_{1,+}(\eta)<0
\quad\text{for any $\eta_R\in\R\setminus\{0\}$,}
\\ \notag
\Re\lambda_{1,+}(\eta_R+i\eta_I)\le \Re \lambda_{1,+}(i\eta_I)\quad\text{for any $\eta_R\in\R$,}
\\ \notag
\Re\tilde{\lambda}_{1,+}(\eta_R+i\eta_I)\le
\Re\tilde{\lambda}_{1,+}(\eta_{\#,+}+i\eta_I)\le -C\a
\quad\text{for $\eta_R\in[-\eta_{\#,+},\eta_{\#,+}]$,}
\end{gather}
where $C$ is a positive constant that does not depend on $\a$.
\end{lemma}
To prove Lemma~\ref{lem:lambda1+'}, we need the following.
\begin{lemma}
  \label{lem:sgnpdlambda1'}
  Let $a_{23}$, $a_{14}$ and $b_2$ be as \eqref{eqdef:a,c,omega} and
  \eqref{eq:b1-b2O}. 
Then $b_2>\k_1+2\k_2$ and $b_2<2\k_3+\k_4$.
\end{lemma}
\begin{remark}
By Lemmas~\ref{lem:evmL2'} and \ref{lem:sgnpdlambda1'},
  \begin{equation*}
i^{-1}\pd_\eta\lambda_{1,\pm}(0)=b_2-3a_{12}\pm\sqrt{c_{12}}>0\,,\quad
i^{-1} \pd_\eta\lambda_{2,\pm}(0)=b_2-3a_{34}\pm\sqrt{c_{34}}<0\,.
  \end{equation*}
Thus perturbations of each line soliton of $u_2$ move
toward the rear of the other line soliton along its crest.
\end{remark}

\begin{proof}[Proof of Lemma~\ref{lem:sgnpdlambda1'}]
  By \eqref{eq:b1-b2O}and the fact that $\k_1<\k_2<\k_3<\k_4$,
  \begin{equation}
    \label{eq:b2-O}
b_2=2(a_{12}+a_{34})+\frac{\k_1\k_2-\k_3\k_4}{2(a_{34}-a_{12})}\,,    
  \end{equation}
\begin{align*}
  b_2-(\k_1+2\k_2)=& \k_3-\k_2+\frac{(\k_4-\k_1)(\k_4-\k_2)}{2(a_{34}-a_{12})}
                   >0\,,
\end{align*}
\begin{align*}
  b_2-(2\k_3+\k_4)=& \k_2-\k_3-\frac{(\k_3-\k_1)(\k_4-\k_1)}{2(a_{34}-a_{12})}
                   < 0\,.
\end{align*}
Thus we prove Lemma~\ref{lem:sgnpdlambda1'}.
\end{proof}
\begin{proof}[Proof of Lemma~\ref{lem:lambda1+'}]
Using Lemma~\ref{lem:sgnpdlambda1'},   
we can prove \eqref{eq:lem:lambda1+'} in exactly the same way as
Lemma~\ref{lem:lambda1-etaI}.
Hence it follows that $\Re\lambda_{1,+}(\eta_R+i\eta_I)\le  \lambda_{1,+}(i\eta_I)$
for every $\eta_R\in\R$. We have $\lambda_{1,+}(i\eta_I)<\tilde{\lambda}_{1,+}(i\eta_I)$
since $\eta_I>0$ and $\gamma_{12}(i\eta_I)>0$.
\par

By \eqref{eq:lem:lambda1+'}, we have
$\Re \tilde{\lambda}_{1,+}(\eta)\le
\Re \tilde{\lambda}_{1,+}(\eta_{\#,+}+i\eta_I)$
for $\eta_R\in[-\eta_{\#,+},\eta_{\#,+}]$.
Since
\begin{equation*}
\Re\gamma_{12}(\eta_{\#,+}+i\eta_I)=\sqrt{c_{12}}-\a\,,  \quad
\Im\gamma_{12}(\eta_{\#,+}+i\eta_I)=\sqrt{\a(\a+2a_{34}-2\k_2)}\,,
\end{equation*}
\begin{align*}
\Re\tilde{\lambda}_{1,+}(\eta_{\#,+}+i\eta_I)&=
-\eta_I(b_2-3a_{12}-\sqrt{c_{12}}+\a)+2(\sqrt{c_{12}}-\a)\a(\a+2a_{34}-2k_2)\,.
\end{align*}
Substituting $\eta_I=2\a(a_{34}-a_{12})$ and \eqref{eq:b2-O} into
the above, we have
\begin{align*}
  \Re\tilde{\lambda}_{1,+}(\eta_{\#,+}+i\eta_I)
=-\{3(\k_2-a_{34})^2+c_{34}\}\a+6(\k_2-a_{34})\a^2-2\a^3\,.  
\end{align*}
Thus we complete the proof.
\end{proof}
\begin{remark}
Let $\xi_{\#,+}=\sqrt{\a(\a+2a_{34}-2\k_2)}$.
If $\eta=\eta_{\#,+}+i\eta_I$, then $\gamma_{12}(\eta)=\sqrt{c_{12}}-\a+i\xi_{\#,+}$ and as $z_1\to\infty$,
  \begin{equation*}  
  g_1(-z_1,\eta)e^{iy\eta}\simeq e^{(\gamma_{12}(\eta)-\sqrt{c_{12}})z_1+iy\eta}
  =e^{ix(\xi_{\#,+}+i\a)+iy(\eta_{\#,+}+2a_{12}\xi_{\#,+}+2ia_{34}\a)}\,.
\end{equation*}
Since $u_{11}\to0$ as $z_1\to\infty$, we have
$\tilde{\lambda}_{1,+}(\eta_{\#,+}+i\eta_I)=ip(\xi_{\#,+}+i\a,\eta_{\#,+}+2a_{12}\xi_{\#,+}+2ia_{34}\a)
\in\sigma(\mL_0)$ in $\calX_4$.
\end{remark}

Let $v_{11}=\pd_x\tau_{11}/\tau_{11}|_{t=0}$.
Then $v_{11}=a_{12}+\psi_1$ with $\psi_1=\sqrt{c_{12}}\tanh\sqrt{c_{12}}z_1$ and
\begin{equation*}
\mM_{1,\pm}(\eta):=e^{-iy\eta}\nabla M_\pm(v_{11})e^{iy\eta}
=\pm\pd_{z_1}+i\eta\pd_{z_1}^{-1}-2\psi_1\,.
\end{equation*}
Let
\begin{gather*}
T^{1,-}_{12}(\eta)f(z_1)
=\int_\R k_3(z_1,z_1',\eta)f(z_1')\,dz_1'\,,
\\
T^{1,+}_{12}(\eta)f(z_1)
=\int_\R k_3(z_1',z_1,-\eta)f(z_1')\,dz_1'\,,
\end{gather*}
where
\begin{gather*}
k_3(z_1,z_1',\eta)=\frac{1}{2\gamma_{12}(\eta)}
\pd_{z_1}\left(e^{-\gamma_{12}(\eta)|z_1-z_1'|}
\sech\sqrt{c_{12}}z_1\cosh\sqrt{c_{12}}z_1'\right)\,.
\end{gather*}
We can prove bijectivity of $\mM_{1,\pm}(\eta)$ in exactly the same way
as Lemmas~\ref{lem:T1+}--\ref{lem:T1-}.
\begin{lemma}
  \label{lem:T1+''}
Assume \eqref{eq:ass-alpha'}.
Let  $\eta=\eta_R+i\eta_I$, $\eta_R\in\R$ and $\eta_I=2\a(a_{34}-a_{12})$.
Then
\begin{gather*}
  \mM_{1,+}(\eta)T^{1,+}_{12}(\eta)=T^{1,+}_{12}(\eta)\mM_{1,+}(\eta)=I
  \quad\text{on $L^2(\R;e^{2\a z_1}dz_1)$,}   
\\
\sum_{j=0,1,2}\la \eta\ra^{(2-j)/2}\|\pd_{z_1}^{j-1}T^{1,+}_{12}(\eta)f\|_{L^2(\R;e^{2\a z_1}dz_1)}
\le  C\|f\|_{L^2(\R;e^{2\a z_1}dz_1)}\,,  
\end{gather*}
where $C$ is a positive constant depending only on $\a$.
\end{lemma}

\begin{lemma}
  \label{lem:T1-''}
Assume \eqref{eq:ass-alpha'}.  
Let $\eta=\eta_R+i\eta_I$, $\eta_R\in\R$,
$\eta_I=2\a(a_{34}-a_{12})$ and $\eta_{*,+}$ be as \eqref{def:e*+'}.
Suppose that $\eta_0>\eta_{*,+}$.
Then for every  $\eta_R$ satisfying $|\eta_R|\ge \eta_0$,
\begin{gather*}
  \mM_{1,-}(\eta)T^{1,-}_{12}(\eta)=T^{1,-}_{12}(\eta)\mM_{1,-}(\eta)=I
  \quad\text{on $L^2(\R;e^{2\a z_1}dz_1)$,}   
  \\
  \sum_{j=0,1,2}\la \eta\ra^{(2-j)/2}\|\pd_{z_1}^{j-1}T^{1,-}_{12}(\eta)f\|_{L^2(\R;e^{2\a z_1}dz_1)}
\le C\|f\|_{L^2(\R;e^{2\a z_1}dz_1)}\,,
\end{gather*}
where $C$ is a positive constant depending only on $\a$ and $\eta_0$.
\end{lemma}
Since
$\mM_{1,-}(\eta)= \sech(\sqrt{c_{12}}z_1)(\gamma_{12}(\eta)^2-\pd_{z_1}^2)
\cosh(\sqrt{c_{12}}z_1)\pd_{z_1}^{-1}$,
\begin{gather*}
  \mM_{1,-}(\eta)\pd_x^{-1}g_1(\cdot,\eta)
  =\mM_{1,-}(\eta)\pd_x^{-1}g_1(-\cdot,\eta)=0\,.
\end{gather*}
As Corollary~\ref{cl:Mrelations-P}, we have
\begin{gather*}
\nabla M_+(v_{11})\pd_x
\left(\Phi^0(\bx,-i\beta)\Phi^{11,*}_1(\bx)\right)
=2\pd_x\left(\Phi^{11}(\bx,-i\beta)\Phi^{11,*}_1(\bx)\right)\,,
\\
\nabla M_+(v_{11})^*
\left(\Phi^{11}_1(\bx)\Phi^{11,*}(\bx,-i\beta')\right)
=2\Phi^{11}_1(\bx)\Phi^{0,*}(\bx,-i\beta')\,.
\end{gather*}
Letting $\beta=a_{12}\pm\gamma_{12}(\eta)$ and
$\beta'=a_{12}\pm\gamma_{12}(-\bar{\eta})$ in the above, we have
\begin{gather*}
\mM_{1,+}(\eta)\pd_x^{-1}g_1(\pm\cdot,\eta)=\pm 2g_1(\pm\cdot,\eta)\,,
\\
 \mM_{1,+}(\eta)^*g_1^*(\pm\cdot,\eta)=
\pm i\bar{\eta}e^{\pm\gamma_{12}(-\bar{\eta})z_1}\sech\sqrt{c_{12}}z_1\,.
\end{gather*}
Note that $\pd_{z_1}^{-1}g_1(z_1,\eta)\in L^2(\R;e^{2\a z_1}dz_1)$ 
if $|\eta_R|<\eta_{1,*}$ and that
$(\pd_{z_1}^{-1}g_1)(-z_1,\eta)\in L^2(\R;e^{2\a z_1}dz_1)$ 
if $|\eta_R|<\eta_{\#,*}$.
We can prove Lemma~\ref{lem:T1-'''} below in the same way as
Lemmas~\ref{lem:T1-'} and \ref{lem:T1-'k3small}.
\begin{lemma}
  \label{lem:T1-'''}
Assume \eqref{eq:ass-alpha'}.  
Let $\eta=\eta_R+i\eta_I$, $\eta_R\in\R$, $\eta_I=2\a(a_{34}-a_{12})$
and $\eta_{*,+}$ and $\eta_{\#,+}$ be as \eqref{def:e*+'} and
\eqref{def:esh+'}.
\begin{enumerate}
\item If $\eta_{\#,+}<|\eta_R|<\eta_{*,+}$, then
\begin{equation*}
\ker(\mM_{1,-}(\eta))=\spann\{\pd_x^{-1}g^1_{1,+}(\cdot,\eta)\}\,,\quad
\operatorname{Range}(\mM_{1,-}(\eta))=L^2(\R;e^{2\a z_1}dz_1)\,.
\end{equation*}
Suppose that $\eta_{\#,+}<\eta_{\#,0}\le|\eta_R|\le \eta_{*,0}<\eta_{*,+}$
and that $f\in L^2(\R;e^{2\a z_1}dz_1)$. Then there exists a unique solution of
$\mM_{1,-}(\eta)v=f$ satisfying
\begin{gather*}
\int_\R v(z_1)e^{\gamma_{12}(\eta)z_1}\sech\sqrt{c_{12}}z_1\,dz_1=0\,,
\\ 
\|v\|_{H^1(\R;e^{2\a z_1}dz_1)}+\|\pd_{z_1}^{-1}v\|_{L^2(\R;e^{2\a z_1}dz_1)}\le
 C\|f\|_{L^2(\R;e^{2\a z_1}dz_1)}\,,
\end{gather*}
where $C$ is a positive constant $C$ depending only on $\eta_{*,0}$,
$\eta_{\#,0}$ and $\a$.
\item
If $|\eta_R|<\eta_{\#,+}$, then
\begin{equation*}
\ker(\mM_{1,-}(\eta))=\spann\{\pd_x^{-1}g^1_{1,+}(\pm\cdot,\eta)\}\,,\quad
\operatorname{Range}(\mM_{1,-}(\eta))=L^2(\R;e^{2\a z_1}dz_1)\,.
\end{equation*}
Suppose that $|\eta_R|\le \eta_{\#,0}<\eta_{\#,+}$
and that $f\in L^2(\R;e^{2\a z_1}dz_1)$. Then there exists a unique solution of
$\mM_{1,-}(\eta)v=f$ satisfying
\begin{gather*}
\int_\R v(z_1)e^{\pm\gamma_{12}(\eta)z_1}\sech\sqrt{c_{12}}z_1\,dz_1=0\,,
\\ 
\|v\|_{H^1(\R;e^{2\a z_1}dz_1)}+\|\pd_{z_1}^{-1}v\|_{L^2(\R;e^{2\a z_1}dz_1)}\le
 C\|f\|_{L^2(\R;e^{2\a z_1}dz_1)}\,,
\end{gather*}
where $C$ is a positive constant $C$ depending only on $\eta_{\#,0}$ and $\a$.
\end{enumerate}
\end{lemma}

\begin{proof}[Proof of Lemma~\ref{lem:mL-11inX2}]
Using Lemmas~\ref{lem:lambda1+'} and \ref{lem:T1+''}--\ref{lem:T1-'''},
we can prove Lemma~\ref{lem:mL-11inX2} in exactly the same way as
Lemma~\ref{lem:mL-1inX2}.
\end{proof}
\bigskip

\subsection{Spectral stability of $[3,4]$-soliton}
We can prove spectral stability of $[3,4]$-soliton $u_{12}$ in $\calX_4$
in the same way as Proposition~\ref{lem:wmL1-decay}.

\begin{lemma}
  \label{lem:linear-stability[3,4]}
Assume $\a\in(0,\sqrt{c_{34}})$.
Let $\mL_0$ and $\mL_{12}$ be closed operators on $\calX_4$.
Let $\eta_{*,4}$ be as \eqref{def:eta*'}, $\eta_1\in(0,\eta_{*,4})$ and
Then
$\sigma(\mL_{12})=\sigma(\mL_0)\cup\{\lambda_{2,\pm}(\eta)\mid
  \eta\in[-\eta_{*,4},\eta_{*,4}]\}$. Moreover,
 $\lambda_{2,\pm}(\eta)-\mL_{12}$ is invertible on
 $\left(I-P_{12}(\eta_0)\right)\calX_4$ if $\eta_0\in(0,\eta_{*,4})$ and
 $\eta\in(-\eta_0,\eta_0)$.
\end{lemma}
Next, we will investigate spectral stability of $u_{12}$ in $\calX_3$.
\begin{lemma}
  \label{lem:mL-12inX1}
Assume
\begin{equation}
\label{eq:ass-alpha''}
0<\a<\min\left\{\sqrt{c_{34}}\,,\,
\frac{c_{34}}{2(a_{34}-a_{12})}\,,2(\k_3-a_{12})\right\}\,.
\end{equation}
Suppose that $\mL_{12}$ is a closed operator in $\calX_3$.
Then there exists a positive number $b$ such that
\begin{equation*}
\sup_{\Re\lambda\ge -b}\|(\lambda-\mL_{12})^{-1}\|_{B(\calX_3)}<\infty\,.
\end{equation*}
\end{lemma}

Let $\eta=\eta_R+i\eta_I$, $\eta_R\in\R$ and $\eta_I=2(a_{12}-a_{34})\a$.
Since  $\a z_1+iy\eta=\a z_2+iy\eta_R$,
\begin{gather*}
  \calX_3=e^{-y\eta_I}\calX_4=\bigoplus_{\eta=\eta_R+i\eta_I,\,\eta_R\in\R}
  e^{iy\eta}L^2(\R;e^{2\a z_2}dz_2)\,.
\end{gather*}
Let
\begin{align*}
  \mL_{12}(\eta)v(z_2):=& e^{-iy\eta}\mL_{12}\left(v(z_2)e^{iy\eta}\right)
  \\=&
-\frac{1}{4}\pd_{z_2}(\pd_{z_2}^2-4c_{34}+6u_{12})
+i(b_2-3a_{34})\eta+\frac{3}{4}\eta^2\pd_{z_2}^{-1}\,.  
\end{align*}
In view of Lemma~\ref{lem:evmL2'},
\begin{equation*}
\mL_{12}(\eta)g_2(z_2,-\eta)=\lambda_{2,-}(\eta)g(z_2,-\eta)\,,
\quad
\mL_{12}(\eta)g_2(-z_2,\eta)=\tilde{\lambda}_{2,-}(\eta)g(z_2,-\eta)\,,
\end{equation*}
where
$\tilde{\lambda}_{2,-}(\eta)=i\eta\{b_2-3a_{12}+\gamma_{34}(-\eta)\}$.
We will investigate whether $e^{iy\eta}g_2(\pm z_2,\pm\eta)$ are continuous
eigenfunctions of $\mL_{12}$ in $\calX_3$.
Since
\begin{gather*}
\|e^{\a z_1}e^{iy\eta}g_2(z_2,\pm\eta)\|_{L^\infty_yL^2_x}
=O\left(\|e^{\a z_2}e^{-\gamma_{34}(\pm\eta)z_2}\sech\sqrt{c_{34}}z_2\|_{L^2(\R)}
\right)\,,
\\  
\|e^{\a z_1}e^{iy\eta}g_2(-z_2,\pm\eta)\|_{L^\infty_yL^2_x}
=O\left(\|e^{\a z_2}e^{\gamma_{34}(\pm\eta)z_2}\sech\sqrt{c_{34}}z_2\|_{L^2(\R)}
\right)\,,
\end{gather*}
we see that $e^{iy\eta}g_2(z_2,\pm\eta)$ are continuous eigenfunctions of
$\mL_{12}$ if $\Re\gamma_{34}(\pm\eta)<\sqrt{c_{34}}+\a$ and that
$e^{iy\eta}g_2(-z_2,\pm\eta)$ are continuous eigenfunctions of
$\mL_{12}$ if $\Re\gamma_{34}(\pm\eta)<\sqrt{c_{34}}-\a$.
\par

If $0<\a<2(\k_3-a_{12})$, then 
$\Re\gamma_{34}(\eta)\ge \gamma_{34}(i\eta_I)>\sqrt{c_{34}}+\a$
and $e^{iy\eta}g_2(\pm z_2,\eta)$ are not continuous eigenfunctions
of $\mL_{12}$.
Let $\a\in(0,\sqrt{c_{34}})$, $c_{34}+\eta_I>0$ and
\begin{gather}
  \label{def:e*-'}
  \eta_{*,-}=2(\sqrt{c_{34}}+\a)\sqrt{\a(\a+2\k_4-2a_{12})}\,,
\\ \label{def:esh-'}
  \eta_{\#,-}=2(\sqrt{c_{34}}-\a)\sqrt{\a(\a+2\k_3-2a_{12})}  \,.
\end{gather}
If $|\eta_R|<\eta_{*,-}$, then $\Re\gamma_{34}(-\eta)<\sqrt{c_{34}}+\a$,
$g_2(z_1,-\eta)\in L^2(\R;e^{2\a z_2}dz_2)$ and $\lambda_{2,-}(\eta)$ is a
continuous eigenvalue of $\mL_{11}$.
If $|\eta_R|>\eta_{*,-}$, then $\Re\gamma_{34}(-\eta)>\sqrt{c_{12}}+\a$ and
$g^1_{2,-}(\eta)$ is not a continuous eigenfunction of $\mL_{12}$.
\par
If $|\eta_R|<\eta_{\#,-}$, then $\Re\gamma_{34}(-\eta)<\sqrt{c_{34}}-\a$,
$g_2(-z_1,-\eta)\in L^2(\R;e^{2\a z_2}dz_2)$ and
$\tilde{\lambda}_{2,-}(\eta)$ is a continuous eigenvalue of $\mL_{11}$.
If $|\eta_R|>\eta_{\#,-}$, then $\Re\gamma_{34}(-\eta)>\sqrt{c_{12}}-\a$ and
$e^{iy\eta}g_2(-z_2,-\eta)$ is not a continuous eigenfunction of $\mL_{12}$.
\par

Continuous eigenvalues $\{\lambda_{2,-}(\eta)\}$ and
$\{\tilde{\lambda}_{2,-}(\eta)\}$ belong to the stable half plane.
\begin{lemma}
  \label{lem:lambda2-'}
Assume \eqref{eq:ass-alpha''}.
Let $\eta=\eta_R+i\eta_I$ and $\eta_I=2\a(a_{12}-a_{34})$. Then
$\lambda_{2,-}(i\eta_I)< \tilde{\lambda}_{2,-}(i\eta_I)$,
\begin{gather*}
  \eta_R\pd_{\eta_R}\Re \lambda_{2,-}(\eta)
=-\eta_R\pd_{\eta_R}\Re \tilde{\lambda}_{2,-}(\eta) <0
\quad\text{for any $\eta_R\in\R\setminus\{0\}$,}
\\
\Re\lambda_{2,-}(\eta_R+i\eta_I)\le \Re \lambda_{2,-}(i\eta_I)
\quad\text{for any $\eta_R\in\R$,}
\\
\Re\lambda_{2,-}(\eta_R+i\eta_I)\le \Re \tilde{\lambda}_{2,-}(i\eta_I)
\le -C\a\quad\text{for any $\eta_R\in[-\eta_{\#,-},\eta_{\#,-}]$,}
\end{gather*}
where $C$ is a positive constant that does not depend on $\a$.
\end{lemma}
\begin{proof}
  We can prove Lemma~\ref{lem:lambda2-'} in exactly the same way as
  Lemma~\ref{lem:lambda1+'} by using Lemma~\ref{lem:sgnpdlambda1'} and
  the fact that $\eta_I=2(a_{12}-a_{34})\a<0$.
  Note that
  \begin{align*}
\Re\lambda_{2,-}(\eta_{\#,-}+i\eta_I)=&
-\eta_I(b_2-3a_{34}+\sqrt{c_{34}})+2(\sqrt{c_{34}}-\a)\a(\a+2\k_3-a_{12})
    \\ =&
-\a\{3(\k_3-a_{12})^2+c_{12}\}\a-6(\k_3-a_{12})\a^2-2\a^3\,.
  \end{align*}
\end{proof}
\begin{remark}
Let $\xi_{\#,-}=-\sqrt{\a(\a+2\k_3-2a_{12})}$.
If $\eta=\eta_{\#,-}+i\eta_I$, then $\gamma_{34}(-\eta)=\sqrt{c_{34}}-\a+i\xi_{\#,-}$ and
as $z_2\to\infty$,
  \begin{equation*}  
  g_2(-z_2,-\eta)e^{iy\eta}\simeq e^{(\gamma_{34}(-\eta)-\sqrt{c_{34}})z_2+iy\eta}
  =e^{ix(\xi_{\#,-}+i\a)+iy(\eta_{\#,-}+2a_{34}\xi_{\#,-}+2ia_{12}\a)}\,.
\end{equation*}
Since $u_{12}\to0$ as $z_2\to\infty$, we have
$\tilde{\lambda}_{2,-}(\eta_{\#,-}+i\eta_I)=ip(\xi_{\#,-}+i\a,\eta_{\#,-}+2a_{34}\xi_{\#,-}+2ia_{12}\a)
\in\sigma(\mL_0)$ in $\calX_3$.
\end{remark}

Let $v_{12}=\pd_x\tau_{12}/\tau_{12}|_{t=0}$.
Then $v_{12}=a_{34}+\psi_2$ with $\psi_2=\sqrt{c_{34}}\tanh\sqrt{c_{34}}z_2$ and
\begin{equation*}
\mM_{2,\pm}(\eta):=e^{-iy\eta}\nabla M_\pm(v_{12})e^{iy\eta}
=\pm\pd_{z_2}+i\eta\pd_{z_2}^{-1}-2\psi_2\,.
\end{equation*}
Let
\begin{gather*}
T^{1,-}_{34}(\eta)f(z_2)
=\int_\R k_4(z_2,z_2',\eta)f(z_2')\,dz_2'\,,
\\
T^{1,+}_{34}(\eta)f(z_2)
=\int_\R k_4(z_2',z_2,-\eta)f(z_2')\,dz_2'\,,
\end{gather*}
where
\begin{gather*}
k_4(z_2,z_2',\eta)=\frac{1}{2\gamma_{34}(\eta)}
\pd_{z_2}\left(e^{-\gamma_{34}(\eta)|z_2-z_2'|}
\sech\sqrt{c_{34}}z_2\cosh\sqrt{c_{34}}z_2'\right)\,.
\end{gather*}
We can prove bijectivity of $\mM_{2,\pm}(\eta)$ in exactly the same way
as Lemmas~\ref{lem:T1+}--\ref{lem:T1-'k3small}.

\begin{lemma}
  \label{lem:T1-34+}
Assume \eqref{eq:ass-alpha''}.  
Let $\eta=\eta_R+i\eta_I$, $\eta_R\in\R$ and
$\eta_I=2\a(a_{12}-a_{34})$. Let $\eta_{*,-}$ and $\eta_{\#,-}$ be as
\eqref{def:e*-'} and \eqref{def:esh-'}.
  \begin{enumerate}
  \item Suppose that $|\eta_R|\ge\eta_0 >\eta_{*,-}$. Then
\begin{gather*}
  \mM_{2,+}(\eta)T^{1,+}_{34}(\eta)=T^{1,+}_{34}(\eta)\mM_{2,+}(\eta)=I
  \quad\text{on $L^2(\R;e^{2\a z_2}dz_2)$,}   
\\
\sum_{j=0,1,2}\la \eta\ra^{(2-j)/2}
\|\pd_{z_2}^{j-1}T^{1,+}_{34}(\eta)u\|_{L^2(\R;e^{2\a z_2}dz_2)}
\le  C\|u\|_{L^2(\R;e^{2\a z_2}dz_2)}\,,  
\end{gather*}
where $C$ is a positive constant depending only on $\a$, $\eta_{*,0}$
and $\eta_{\#,0}$.
\item
  Suppose that $\eta_{\#,-}<\eta_{\#,0}\le |\eta_R|\le \eta_{*,0}<\eta_{*,-}$
  and that
\begin{equation*}
  \int_\R f(z_2)\overline{g^*_2(z_2,-\eta)}\,dz_2=0\,.
\end{equation*}
Then there exists a unique solution of $\mM_{2,+}(\eta)v=f$ satisfying
\begin{equation*}
\|v\|_{H^1(\R;e^{2\a z_2}dz_2)}+\|\pd_{z_2}^{-1}v\|_{L^2(\R;e^{2\a z_2}dz_2)}
\le  C\|f\|_{L^2(\R;e^{2\a z_2}dz_2)}\,, 
\end{equation*}
where $C$ is a positive constant depending only on $\eta_0$ and $\a$.
\item
Suppose that $|\eta_R|\le\eta_0 <\eta_{\#,-}$  and that
\begin{equation*}
  \int_\R f(z_2)\overline{g^*_2(\pm z_2,-\eta)}\,dz_2=0\,.
\end{equation*}
Then there exists a unique solution of $\mM_{2,+}(\eta)v=f$ satisfying
\begin{equation*}
\|v\|_{H^1(\R;e^{2\a z_2}dz_2)}+\|\pd_{z_2}^{-1}v\|_{L^2(\R;e^{2\a z_2}dz_2)}
\le  C\|f\|_{L^2(\R;e^{2\a z_2}dz_2)}\,, 
\end{equation*}
where $C$ is a positive constant depending only on $\eta_0$ and $\a$.
\end{enumerate}
\end{lemma}

\begin{lemma}
  \label{lem:T1-34-}
Assume \eqref{eq:ass-alpha''}.  
Let $\eta=\eta_R+i\eta_I$, $\eta_R\in\R$ and $\eta_I=2\a(a_{12}-a_{34})$.
  Then
\begin{gather*}
  \mM_{2,-}(\eta)T^{1,-}_{34}(\eta)=T^{1,-}_{34}(\eta)\mM_{2,-}(\eta)=I
  \quad\text{on $L^2(\R;e^{2\a z_2}dz_2)$,}   
  \\
  \sum_{j=0,1,2}\la \eta\ra^{(2-j)/2}
  \|\pd_{z_2}^{j-1}T^{1,-}_{34}(\eta)u\|_{L^2(\R;e^{2\a z_2}dz_2)}
\le C\|u\|_{L^2(\R;e^{2\a z_2}dz_2)}\,,
\end{gather*}
where $C$ is a positive constant depending only on $\a$.
\end{lemma}

\begin{proof}[Proof of Lemma~\ref{lem:mL-12inX1}]
Using Lemmas~\ref{lem:lambda2-'}--\ref{lem:T1-34-},
we can prove Lemma~\ref{lem:mL-12inX1} in exactly the same way as
Lemma~\ref{lem:mL-1inX2}.
\end{proof}
\bigskip

\subsection{The inverse of $\nabla M_\pm(v_{21})$}
To begin with, we introduce a Green function of $L_{11}$.
Let  $\tbx=(x,y)$, $\tbxd=(x',y')$, $\eta\in\R$, $z_2'=x'+2a_{34}y'$ and
\begin{gather*}
\mathfrak{g}^{11,-}_\pm(\tbx,\tbxd,\eta)=
-\frac{1}{2\gamma_{34}(\eta)}\Phi^{11}(x,y,0,-i\beta_{34}^\pm(\eta))
\Phi^{11,*}(x',y',0,-i\beta_{34}^\pm(\eta))\,,
\\
 \mathfrak{g}^{11,-}(\tbx,\tbxd,\eta)=\left\{
  \begin{aligned}
& \mathfrak{g}^{11,-}_-(\tbx,\tbxd,\eta)\quad\text{if $z_2>z_2'$,}
\\ &
\mathfrak{g}^{11,-}_+(\tbx,\tbxd,\eta)\quad\text{if $z_2<z_2'$,}
  \end{aligned}\right.
\\ 
\mathfrak{g}^{11,+}(\tbx,\tbxd,\eta)=\left\{
  \begin{aligned}
& \mathfrak{g}^{11,+}_+(\tbx,\tbxd,\eta)\quad\text{if $z_2>z_2'$,}
\\ & 
\mathfrak{g}^{11,+}_-(\tbx,\tbxd,\eta)\quad\text{if $z_2<z_2'$,}
  \end{aligned}\right.
\\  
\mathfrak{g}^{11,+}_\pm(\tbx,\tbxd,\eta)
=\mathfrak{g}^{11,-}_\pm(\tbxd,\tbx,-\eta)\,,
\end{gather*}
\begin{gather*}
\mathcal{G}^{11,-}(\tbx,\tbxd)=\frac{1}{2\pi}
\int_{-\infty}^\infty \mathfrak{g}^{11,-}(\tbx,\tbxd,\eta)\,d\eta\,,
\\
\mathcal{G}^{11,+}(\tbx,\tbxd)=
\frac{1}{2\pi}\int_{-\infty}^\infty\mathfrak{g}^{11,+}
(\tbx,\tbxd,\eta)\,d\eta\,.
\end{gather*}
The Green function $\mathfrak{g}^{11,-}$ can be rewritten as
\begin{align*}
& \mathfrak{g}^{11,-}(\tbx,\tbxd,\eta)= \mathfrak{g}^{11,-}_1(\tbx,\tbxd,\eta)
+\mathfrak{g}^{11,-}_2(\tbx,\tbxd,\eta)\,,
\\ &   
\mathfrak{g}^{11,-}_1(\tbx,\tbxd,\eta)
= -\frac{1}{2\gamma_{34}(\eta)}
e^{\frac{1}{2}\sum_{j=3,4}(\theta_{j,0}-\theta_{j,0}')
-\gamma_{34}(\eta)|z_2-z_2'|+i\eta(y-y')}\,,
\\ &
\mathfrak{g}^{11,-}_2(\tbx,\tbxd,\eta)= 
\mathfrak{g}^{11,-}_1(\tbx,\tbxd,\eta)
\sum_{j=1,2}
\left(e^{\theta_{j,0}'-\theta_{j,0}}
\frac{\Phi^{11}_j(\bx)\Phi^{11,*}_j(\bxd)}{\beta_{34}^\mp(\eta)-\k_j}\right)
\quad\text{if  $\pm(z_2-z_2')>0$.}
\end{align*}

We can prove the following in exactly the same way as
Lemma~\ref{lem:L-Green}.
\begin{lemma}
  \label{lem:L11-Green}
\begin{equation*}
L_{11}\mathcal{G}^{11,-}(\cdot,\tbxd)=\delta(\cdot-\tbxd)\,,\quad
L_{11}^*\mathcal{G}^{11,+}(\cdot,\tbxd)=\delta(\cdot-\tbxd)\,.
\end{equation*}
\end{lemma}

Next we introduce Green functions of $\nabla M_\pm(v_{21})$.
Let  $h_{21}=\tau_2/\tau_{11}|_{t=0}$ and
\begin{align*}
\tfg^{21,-}_\pm(\tbx,\tbxd,\eta)=-
\frac{h_{21}(\tbxd)}{h_{21}(\tbx)}\mathfrak{g}^{11,-}_\pm(\tbx,\tbxd,\eta)\,,
\\
\tfg^{21,+}_\pm(\tbx,\tbxd,\eta)
=\frac{h_{21}(\tbx)}{h_{21}(\tbxd)}\mathfrak{g}^{11,+}_\pm(\tbx,\tbxd,\eta)\,,
\\
\tfg^{21,-}=\tfg^{21,-}_\mp\,,\quad
\tfg^{21,+}=\tfg^{21,+}_\pm\quad\text{if $\pm(z_2-z_2')>0$.}
\end{align*}
By Lemma~\ref{lem:L11-Green} and \eqref{eq:Miura-Lax2j},
\begin{gather*}
\frac{1}{2\pi}\nabla M_-(v_{21})\pd_x
\int_\R \tfg^{21,-}(\tbx,\tbxd,\eta)\,d\eta
=\delta(\cdot-\tbxd)\,,
\\
\frac{1}{2\pi}\nabla M_+(v_{21})
\int_\R \tfg^{21,+}(\tbx,\tbxd,\eta)\pd_{x'}\,d\eta
=\delta(\cdot-\tbxd)\,.  
\end{gather*}
\par

For $\eta_1$ and $\eta_2$ satisfying $\eta_2\ge\eta_1\ge0$, let
\begin{align*}
& T^{21,-}(\eta_1,\eta_2)f=\frac{1}{2\pi}\pd_x
\int_{I_{\eta_1, \eta_2}}\int_{\R^2}
\tfg^{21,-}(\tbx,\tbxd,\eta)f(\tbxd)\,d\tbxd d\eta\,,
\\ & 
T^{21,+}(\eta_1,\eta_2)f=\frac{1}{2\pi}
\int_{I_{\eta_1,\eta_2}}\int_{\R^2}
\tfg^{21,+}(\tbx,\tbxd,\eta)\pd_{x'}f(\tbxd)\,d\tbxd d\eta\,.
\end{align*}

\begin{lemma}
  \label{lem:T21-pm}
Let $\a\in(0,2\sqrt{c_{34}})$ and $\eta_{*,4}$ be as \eqref{def:eta*'}.
 If $\eta_2\ge\eta_1\ge\eta_0>\eta_{*,4}$,
    \begin{gather*}
\|T^{21,\pm}(\eta_1,\eta_2)f\|_{\calX_4}\le C\|f\|_{\calX_4}\,,
    \end{gather*}
where $C$ is  a positive constant depending only on $\a$ and $\eta_0$.
\end{lemma}
\begin{proof}
Let
\begin{gather*}
\tfg^{21,-}_0(\tbx,\tbxd,\eta)=
-\frac{e^{\theta_{3,0}'}+e^{\theta_{4,0}'}}{e^{\theta_{3,0}}+e^{\theta_{4,0}}}
\mathfrak{g}^{11,-}_1(\tbx,\tbxd,\eta)\,,
\\  
\mathcal{T}_0^{21,-}(\eta_1,\eta_2)f=\frac{1}{2\pi}
\int_{I_{\eta_1,\eta_2}}\int_{\R^2}
\tfg^{21,-}_0(\tbx,\tbxd,\eta)f(\tbxd)\,d\tbxd d\eta\,,
\\
\mathcal{T}_j^{21,-}(\eta_1,\eta_2)f=\frac{1}{2\pi}
\int_{I_{\eta_1, \eta_2}}\int_{\R^2}
\frac{\tfg^{21,-}_0(\tbx,\tbxd,\eta)}{\sum_\pm\chi_\pm(z_2-z_2')
\beta_{34}^\mp(\eta)-\k_j}
f(\tbxd)\,d\tbxd d\eta\,,
\end{gather*}
and let
$$
\mathcal{W}_0(\tbx)=\frac{e^{\theta_{3,0}}+e^{\theta_{4,0}}}{h_{21}(\tbx)}\,,
\quad \mathcal{W}_1(\tbx)=(\k_1-\k_2)\frac{e^{\theta_{2,0}}}{\tau_{11}}\,,
\quad \mathcal{W}_2(\tbx)=\frac{e^{\theta_{1,0}}}{\tau_{11}}\,.
$$
Then
\begin{equation*}
  T^{21,-}=\pd_x\mathcal{W}_0\left\{\mathcal{T}_0^{21,-}+\mathcal{W}_1\mathcal{T}_1^{21,-}\mathcal{W}_2
    -\mathcal{W}_2\mathcal{T}_2^{21,-}\mathcal{W}_1\right\}\mathcal{W}_0^{-1}\,.
\end{equation*}
We remark that $\mathcal{W}_0$, $\mathcal{W}_1$, $\mathcal{W}_2$, $\mathcal{W}_0^{-1}$ and their derivatives
are bounded on $\calX_3$.  By \eqref{eq:btitj},
\begin{gather*}
\tfg^{21,-}_0(\tbx,\tbxd,\eta)=e^{i\eta(y-y')}k_{21}(z_2,z_2',\eta)\,,
\\
k_{21}(z,z',\eta)=\frac{1}{2\gamma_{34}(\eta)}e^{-\gamma_{34}(\eta)|z-z'|}
\sech\sqrt{c_{34}}z\cosh\sqrt{c_{34}}z'\,,
\\
|\pd_x^jk_{21}(z,z',\eta)| \lesssim
\la\eta\ra^{(j-1)/2}e^{(-\Re\gamma_{34}(\eta)+\sqrt{c_{34}})|z-z'|}
\quad\text{for $j\ge0$ and $z$, $z'\in\R$.}
\end{gather*}
Since $\Re\gamma_{34}(\eta)\ge \gamma_{34}(\eta_0)>\sqrt{c_{34}}+\a$ and
$$e^{\a(z-z')}k_{21}(z,z',\eta)\in L^\infty_zL^1_{z'}\cap L^\infty_{z'}L^1_z\,,$$
we have $\|T^{21,-}(\eta_1,\eta_2)f\|_{\calX_4}\le C\|f\|_{\calX_4}$.
\par

By integration by parts,
$T^{21,+}(\eta_1,\eta_2)=\sum_{1\le j\le 3}
T^{21,+}_j(\eta_1,\eta_2)$, where $f^\#(z_2,y)=f(x,y)$ and
\begin{align*}
& T^{21,+}_1(\eta_1,\eta_2)f(\tbx)  
= -\frac{1}{2\pi}\int_{I_{\eta_1,\eta_2}}\left(
\int_{-\infty}^\infty\int^{z_2}_{-\infty}\pd_{x'}\tfg^{21,+}_+(\tbx,\tbxd,\eta)
f(\tbxd)\,dz_2'dy'\right)\,d\eta\,,
  \\ &
T^{21,+}_2(\eta_1,\eta_2)f(\tbx)  
= -\frac{1}{2\pi}\int_{I_{\eta_1,\eta_2}}\left(
\int_{-\infty}^\infty\int_{z_2}^\infty\pd_{x'}\tfg^{21,+}_-(\tbx,\tbxd,\eta)
f(\tbxd)\,dz_2'dy'\right)\,d\eta\,,
\end{align*}
\begin{align*}
T^{21,+}_3(\eta_1,\eta_2)f(\tbx)
  = \frac{1}{2\pi}\sum_{j=1,2}
  & \int_{I_{\eta_1,\eta_2}}\int_\R
\frac{\mathcal{W}_0(\tbxd)\mathcal{W}_{3-j}(\tbx)\mathcal{W}_j(\tbxd)}{\mathcal{W}_0(\tbx)}\Bigr|_{z_2'=z_2}
\\ & \quad \times     \frac{(-1)^{j-1}e^{i(y-y')\eta}}
      {i\eta+(\k_4-\k_j)(\k_3-\k_j)}f^\#(z_2,y')\,dy'd\eta\,.
\end{align*}
Thus we can prove Lemma~\ref{lem:T21-pm} in exactly the same way as
Lemma~\ref{lem:T2pm}.  
\end{proof}

Next, we will modify $T^{21,+}$ and $T^{21,-}$ for small $\eta$ as in Section~\ref{subsec:invMpm2}.

Let $T^{21,+}_{low}(\eta_1,\eta_2)=T^{21,+}_{1,low}(\eta_1,\eta_2)
+\sum_{j=2,3}T^{21,+}_j(\eta_1,\eta_2)$, where
\begin{align*} 
 & T^{21,+}_{1,low}(\eta_1,\eta_2)f(\tbx)  
\\=& -\chi_+(z_2)\frac{1}{2\pi}\int_{I_{\eta_1,\eta_2}}\left(
\int_{-\infty}^\infty\int^\infty_{z_2}\pd_{x'}\tfg^{21,+}_+(\tbx,\tbxd,\eta)
f(\tbxd)\,dz_2'dy'\right)\,d\eta
\\ & +\chi_-(z_2)\frac{1}{2\pi}\int_{I_{\eta_1,\eta_2}}\left(
\int_{-\infty}^\infty\int_{-\infty}^{z_2}\pd_{x'}\tfg^{21,+}_+(\tbx,\tbxd,\eta)
f(\tbxd)\,dz_2'dy'\right)\,d\eta\,.
\end{align*}
By Claim~\ref{cl:L0*L0'},
\begin{align}
\label{eq:tfg21-+}
\pd_{x'} \tfg^{21,+}_+(\tbx,\tbxd,\eta)=
-\frac{h_{21}(\tbx)}{2\gamma_{34}(-\eta)}
\Phi^{11,*}(\tbx,0, -i\beta_{34}^+(-\eta))g^{2,*}_{2,-}(\tbxd,-\eta)\,,
\end{align}
and $T^{21,+}(\eta_1,\eta_2)f=T^{21,+}_{low}(\eta_1,\eta_2)f$ if
$f\in C_0^\infty(\R^2)$ and
\begin{equation}
\label{eq:T21+secular}
\int_{\R^2} f(\tbx)\overline{g^{2,*}_{2,-}(\tbx,\eta)}\,dxdy=0\quad\text{for $\eta\in I_{\eta_1,\eta_2}$.}
\end{equation}
Let
$T^{21,-}_{low}(\eta_1,\eta_2)=\sum_{j=1,3}T^{21,-}_j(\eta_1,\eta_2)
+T^{21,-}_{2,low}(\eta_1,\eta_2)$ and
\begin{gather*}
T^{2,-}_1(\eta_1,\eta_2)f(\tbx)=
\frac{1}{2\pi}\int_{I_{\eta_1,\eta_2}}\int_{-\infty}^\infty
\int_{z_2}^\infty \pd_x\tfg^{21,-}_+(\tbx,\tbxd,\eta)f(\tbxd)\,dz_2'dy'd\eta\,,
\\
T^{2,-}_2(\eta_1,\eta_2)f(\tbx)=
\frac{1}{2\pi}\int_{I_{\eta_1,\eta_2}}\int_{-\infty}^\infty
\int_{-\infty}^{z_2}\pd_x\tfg^{21,-}_-(\tbx,\tbxd,\eta)f(\tbxd)\,dz_2'dy'd\eta\,,
\\
T^{2,-}_{2,low}(\eta_1,\eta_2)f(\tbx)=
\frac{1}{2\pi}\int_{I_{\eta_1,\eta_2}}\int_{-\infty}^\infty
\int_0^{z_2}\pd_x\tfg^{21,-}_-(\tbx,\tbxd,\eta)f(\tbxd)\,dz_2'dy'd\eta\,,
\end{gather*}
\begin{align*}
T^{21,-}_3(\eta_1,\eta_2)f(\tbx)
=  \frac{1}{2\pi}\sum_{j=1,2}&\int_{I_{\eta_1,\eta_2}}\int_\R
\frac{\mathcal{W}_0(\tbx)\mathcal{W}_j(\tbx)\mathcal{W}_{3-j}(\tbxd)}{\mathcal{W}_0(\tbxd)}\Bigr|_{z_2'=z_2}
\\ & \quad\times      \frac{(-1)^je^{i(y-y')\eta}}
      {i\eta-(\k_3-\k_j)(\k_4-\k_j)}f^\#(z_2,y')\,dy'd\eta\,.  
\end{align*}
Note that $T^{21,-}(\eta_1,\eta_2)=\sum_{1\le j\le3}T^{21,-}_j(\eta_1,\eta_2)$. 
\par
By Claim~\ref{cl:L0*L0'},
\begin{equation*}
\pd_{x'}\tfg^{21,-}_-(\tbxd,\tbx,-\eta)=\frac{h_2(\tbx)}{2(\k_4-\k_3)}
\Phi^{11,*}(\bxd,-i\beta_{34}^-(-\eta))\tg^M_2(\tbxd,-\eta)\,.
\end{equation*}
Since
\begin{align*}
& T^{21,-}(\eta_1,\eta_2)^*f(\tbx)-T^{21,-}_{low}(\eta_1,\eta_2)^*f(\tbx)
  \\ = &
\frac{1}{2\pi}\chi_-(z_2)\int_{I_{\eta_1,\eta_2}}\left(
\int_{\R^2}\pd_{x'}\tfg^{21,-}_-(\tbxd,\tbx,-\eta)
   f(\tbxd)\,d\tbxd\right)\,d\eta\,,
\end{align*}
we have $T^{21,-}(\eta_1,\eta_2)^*f=T^{21,-}_{low}(\eta_1,\eta_2)^*f$  for $f\in C_0^\infty(\R^2)$ satisfying
\begin{equation}
  \label{eq:T21-*secular}
  \int_\R f(\tbx)\overline{\tg^M_2(\tbx,\eta)}\,dxdy=0
  \quad\text{for $\eta\in I_{\eta_1,\eta_2}$,}
\end{equation}
We can prove the following in the same way as Lemma~\ref{lem:T2lowpm}.
\begin{lemma}
  \label{lem:T21low+}
  Let $\a\in(0,2\sqrt{c_{34}})$ and $0\le\eta_1\le\eta_2\le\eta_0<\eta_{*,4}$.
  \begin{enumerate}
  \item There exists a positive constant $C$ depending only on
$\a$ and $\eta_0$ such that for $f\in \calX_4$,
\begin{gather*}
\|T^{21,\pm}_{low}(\eta_1,\eta_2)f\|_{\calX_4}
\le C\|f\|_{\calX_4}\,.
\end{gather*}
\item
There exists a positive constant $C$ depending only on
$\a$ and $\eta_0$ such that for $f\in\calX_4$ satisfying
\eqref{eq:T21+secular},
\begin{equation*}
\|T^{21,+}(\eta_1,\eta_2)f\|_{\calX_4}\le C\|f\|_{\calX_4}\,.
\end{equation*}
\item
There exists a positive constant $C$ depending only on
$\a$ and $\eta_0$ such that for $f\in\calX_4^*$ satisfying
\eqref{eq:T21-*secular},
\begin{equation*}
\|T^{21,-}(\eta_1,\eta_2)^*f\|_{\calX_4}\le C\|f\|_{\calX_4}\,.
\end{equation*}

\end{enumerate}
\end{lemma}

Finally, we will prove that $\ker \nabla M_+(v_{21})=\{0\}$.
\begin{lemma}
  \label{lem:functional-21}
  Assume $0<\a<\sqrt{c_{34}}$. If $\varphi\in\calX_4$ and
  $\nabla M_+(v_{21})\varphi=0$, then $\varphi=0$.
\end{lemma}
\begin{proof}
Let $a(x,y)=-\Phi^2_3(\tbx)\Phi^{2,*}_3(\tbx)$. Then $a(x,y)>0$ and
$a(x,y)=O\left(\sech^2\sqrt{c_{34}}z_2\right)$.
Moreover, it follows from \eqref{eq:Miura-Lax2j} that
$$(\pd_y-\pd_x^2-2\pd_xv_{21})a=\pd_x\nabla M_-(v_{21})a=0\,.$$
Thus we can prove Lemma~\ref{lem:functional-21} in exactly the same way as Lemma~\ref{lem:functional}.
\end{proof}
\bigskip

\subsection{The inverse of $\nabla M_\pm(v_{22})$}
Let $z_1'=x'+2a_{12}y'$ and
\begin{gather*}
\mathfrak{g}^{12,-}_\pm(\tbx,\tbxd,\eta)=
-\frac{1}{2\gamma_{12}(\eta)}\Phi^{12}(x,y,0,-i\beta_{12}^\pm(\eta))
\Phi^{12,*}(x',y',0,-i\beta_{12}^\pm(\eta))\,,
\\
 \mathfrak{g}^{12,-}(\tbx,\tbxd,\eta)=\left\{
  \begin{aligned}
& \mathfrak{g}^{12,-}_-(\tbx,\tbxd,\eta)\quad\text{if $z_1>z_1'$,}
\\ &
\mathfrak{g}^{12,-}_+(\tbx,\tbxd,\eta)\quad\text{if $z_1<z_1'$,}
  \end{aligned}\right.
\\ 
\mathfrak{g}^{12,+}(\tbx,\tbxd,\eta)=\left\{
  \begin{aligned}
& \mathfrak{g}^{12,+}_+(\tbx,\tbxd,\eta)\quad\text{if $z_1>z_1'$,}
\\ & 
\mathfrak{g}^{12,+}_-(\tbx,\tbxd,\eta)\quad\text{if $z_1<z_1'$,}
  \end{aligned}\right.
\\  
\mathfrak{g}^{12,+}_\pm(\tbx,\tbxd,\eta)
=\mathfrak{g}^{12,-}_\pm(\tbxd,\tbx,-\eta)\,,
\end{gather*}
\begin{gather*}
\mathcal{G}^{12,-}(\tbx,\tbxd)=\frac{1}{2\pi}
\int_{-\infty}^\infty \mathfrak{g}^{12,-}(\tbx,\tbxd,\eta)\,d\eta\,,
\\
\mathcal{G}^{12,+}(\tbx,\tbxd)=
\frac{1}{2\pi}\int_{-\infty}^\infty\mathfrak{g}^{12,+}
(\tbx,\tbxd,\eta)\,d\eta\,.
\end{gather*}
The Green function $\mathfrak{g}^{12,-}$ can be rewritten as
\begin{align*}
& \mathfrak{g}^{12,-}(\tbx,\tbxd,\eta)= \mathfrak{g}^{12,-}_1(\tbx,\tbxd,\eta)
+\mathfrak{g}^{12,-}_2(\tbx,\tbxd,\eta)\,,
\\ &   
\mathfrak{g}^{12,-}_1(\tbx,\tbxd,\eta)
= -\frac{1}{2\gamma_{12}(\eta)}
e^{\frac{1}{2}\sum_{j=1,2}(\theta_{j,0}-\theta_{j,0}')
-\gamma_{34}(\eta)|z_1-z_1'|+i\eta(y-y')}\,,
\\ &
\mathfrak{g}^{12,-}_2(\tbx,\tbxd,\eta)= 
\mathfrak{g}^{12,-}_1(\tbx,\tbxd,\eta)
\sum_{j=3,4}
\left(e^{\theta_{j,0}'-\theta_{j,0}}
\frac{\Phi^{12}_j(\bx)\Phi^{12,*}_j(\bxd)}{\beta_{12}^\mp(\eta)-\k_j}\right)
\quad\text{if  $\pm(z_1-z_1')>0$.}
\end{align*}

We can prove the following in exactly the same way as
Lemma~\ref{lem:L-Green}.
\begin{lemma}
  \label{lem:L12-Green}
\begin{equation*}
L_{12}\mathcal{G}^{12,-}(\cdot,\tbxd)=\delta(\cdot-\tbxd)\,,\quad
L_{12}^*\mathcal{G}^{12,+}(\cdot,\tbxd)=\delta(\cdot-\tbxd)\,.  
\end{equation*}
\end{lemma}

Next, we introduce Green functions of $\nabla M_\pm(v_{22})$.
Let  $h_{22}=\tau_2/\tau_{12}|_{t=0}$ and
\begin{align*}
\tfg^{22,-}_\pm(\tbx,\tbxd,\eta)=-
\frac{h_{22}(\tbxd)}{h_{22}(\tbx)}\mathfrak{g}^{12,-}_\pm(\tbx,\tbxd,\eta)\,,
\\
\tfg^{22,+}_\pm(\tbx,\tbxd,\eta)
=\frac{h_{22}(\tbx)}{h_{22}(\tbxd)}\mathfrak{g}^{12,+}_\pm(\tbx,\tbxd,\eta)\,,
\\
\tfg^{22,-}=\tfg^{22,-}_\mp\,,\quad
\tfg^{22,+}=\tfg^{22,+}_\pm\quad\text{if $\pm(z_1-z_1')>0$.}
\end{align*}
By Lemma~\ref{lem:L12-Green} and \eqref{eq:Miura-Lax2j},
\begin{gather*}
\frac{1}{2\pi}\nabla M_-(v_{22})\pd_x
\int_\R \tfg^{22,-}(\tbx,\tbxd,\eta)\,d\eta
=\delta(\cdot-\tbxd)\,,
\\
\frac{1}{2\pi}\nabla M_+(v_{22})
\int_\R \tfg^{22,+}(\tbx,\tbxd,\eta)\pd_{x'}\,d\eta
=\delta(\cdot-\tbxd)\,.  
\end{gather*}
\par

For $\eta_1$ and $\eta_2$ satisfying $\eta_2\ge\eta_1\ge0$, let
\begin{align*}
& T^{22,-}(\eta_1,\eta_2)f=\frac{1}{2\pi}\pd_x
\int_{I_{\eta_1,\eta_2}}\int_{\R^2}
\tfg^{22,-}(\tbx,\tbxd,\eta)f(\tbxd)\,d\tbxd d\eta\,,
\\ & 
T^{22,+}(\eta_1,\eta_2)f=\frac{1}{2\pi}
\int_{I_{\eta_1, \eta_2}}\int_{\R^2}
\tfg^{22,+}(\tbx,\tbxd,\eta)\pd_{x'}f(\tbxd)\,d\tbxd d\eta\,.
\end{align*}

\begin{lemma}
  \label{lem:T22-pm}
  Let $\a\in(0,2\sqrt{c_{12}})$ and $\eta_{*,3}$ be as \eqref{def:eta*'}.
 If $\eta_2\ge\eta_1\ge\eta_0>\eta_{*,3}$, then
\begin{gather*}
\|T^{22,\pm}(\eta_1,\eta_2)f\|_{\calX_3} \le C\|f\|_{\calX_3}\,,
\end{gather*}
where $C$ is  a positive constant depending only on $\a$ and $\eta_0$. 
\end{lemma}
\begin{proof}
Let
\begin{gather*}
\tfg^{22,-}_0(\tbx,\tbxd,\eta)=
-\frac{e^{\theta_{1,0}'}+e^{\theta_{2,0}'}}{e^{\theta_{1,0}}+e^{\theta_{2,0}}}
\mathfrak{g}^{12,-}_1(\tbx,\tbxd,\eta)\,,
\\  
\mathcal{T}_0^{22,-}(\eta_1,\eta_2)f=\frac{1}{2\pi}
\int_{I_{\eta_1, \eta_2}}\int_{\R^2}
\tfg^{22,-}_0(\tbx,\tbxd,\eta)f(\tbxd)\,d\tbxd d\eta\,,
\\
\mathcal{T}_j^{22,-}(\eta_1,\eta_2)f=\frac{1}{2\pi}
\int_{I_{\eta_1,\eta_2}}\int_{\R^2}
\frac{\tfg^{22,-}_0(\tbx,\tbxd,\eta)}{\sum_\pm\chi_\pm(z_1-z_1')
\beta_{12}^\mp(\eta)-\k_j}
f(\tbxd)\,d\tbxd d\eta\,,
\end{gather*}
and let
$$
\mathcal{W}_0(\tbx)=\frac{e^{\theta_{1,0}}+e^{\theta_{2,0}}}{h_{22}(\tbx)}\,,
\quad \mathcal{W}_1(\tbx)=(\k_3-\k_4)\frac{e^{\theta_{4,0}}}{\tau_{12}}\,,
\quad \mathcal{W}_2(\tbx)=\frac{e^{\theta_{3,0}}}{\tau_{12}}\,.
$$
Then
\begin{equation*}
  T^{22,-}=\pd_x\mathcal{W}_0\left\{\mathcal{T}_0^{22,-}+\mathcal{W}_1\mathcal{T}_1^{22,-}\mathcal{W}_2
    -\mathcal{W}_2\mathcal{T}_2^{22,-}\mathcal{W}_1\right\}\mathcal{W}_0^{-1}\,.
\end{equation*}
We remark that $\mathcal{W}_0$, $\mathcal{W}_1$, $\mathcal{W}_2$, $\mathcal{W}_0^{-1}$ and their derivatives
are bounded on $\calX_3$.  By \eqref{eq:btitj},
\begin{gather*}
\tfg^{22,-}_0(\tbx,\tbxd,\eta)=e^{i\eta(y-y')}k_{22}(z_1,z_1',\eta)\,,
\\
k_{22}(z,z',\eta)=\frac{1}{2\gamma_{12}(\eta)}e^{-\gamma_{12}(\eta)|z-z'|}
\sech\sqrt{c_{12}}z\cosh\sqrt{c_{12}}z'\,.
\end{gather*}
Thus we can prove Lemma~\ref{lem:T22-pm} in exactly the same way as
Lemma~\ref{lem:T2pm}.
\end{proof}

Let $T^{22,+}_{low}(\eta_1,\eta_2)=T^{22,+}_{1,low}(\eta_1,\eta_2)+\sum_{j=2,3}T^{22,+}_j(\eta_1,\eta_2)$,
where
\begin{align*}
& T^{22,+}_{1,low}(\eta_1,\eta_2)f(\tbx)  
\\=& \chi_+(z_1)\frac{1}{2\pi}\int_{I_{\eta_1,\eta_2}}\left(
\int_{-\infty}^\infty\int^\infty_{z_1}\pd_{x'}\tfg^{22,+}_+(\tbx,\tbxd,\eta)
f(\tbxd)\,dz_1'dy'\right)\,d\eta
\\ & -\chi_-(z_1)\frac{1}{2\pi}\int_{I_{\eta_1,\eta_2}}\left(
\int_{-\infty}^\infty\int_{-\infty}^{z_1}\pd_{x'}\tfg^{2,+}_+(\tbx,\tbxd,\eta)
f(\tbxd)\,dz_1'dy'\right)\,d\eta\,,
\end{align*}
\begin{align*}
T^{22,+}_2(\eta_1,\eta_2)f(\tbx)  
= -\frac{1}{2\pi}\int_{I_{\eta_1,\eta_2}}\left(
\int_{-\infty}^\infty\int_{z_1}^\infty\pd_{x'}\tfg^{22,+}_-(\tbx,\tbxd,\eta)
f(\tbxd)\,dz_1'dy'\right)\,d\eta\,,
\end{align*}
\begin{align*}
T^{22,+}_3(\eta_1,\eta_2)f(\tbx)
  = \frac{1}{2\pi}\sum_{j=1,2}
  & \int_{I_{\eta_1,\eta_2}}\int_\R
\frac{\mathcal{W}_0(\tbxd)\mathcal{W}_{3-j}(\tbx)\mathcal{W}_j(\tbxd)}{\mathcal{W}_0(\tbx)}\Bigr|_{z_1'=z_1}
\\ & \quad \times     \frac{(-1)^{j-1}e^{i(y-y')\eta}}
      {i\eta+(\k_{j+2}-\k_1)(\k_{j+2}-\k_2)}f^\#(z_1,y')\,dy'd\eta\,.
\end{align*}
Then
\begin{align*}
& T^{2,+}_{low}(\eta_1,\eta_2)f(\tbx) -T^{2,+}(\eta_1,\eta_2)f(\tbx)  
\\ =&
\frac{1}{2\pi}\chi_+(z_1)\int_{I_{\eta_1,\eta_2}}\int_{\R^2}
\pd_{x'}\tfg^{22,+}_+(\tbx,\tbxd,\eta)f(\tbxd)\,dz_1'dy'd\eta\,.
\end{align*}
Since
\begin{equation*}
  \pd_{x'}\tfg^{22,+}_+(\tbx,\tbxd,\eta)=\frac{h_{22}(\tbx)}{2\gamma_{12}(-\eta)}
  \Phi^{12,*}(\tbx,-i\beta_{12}^+(-\eta))\overline{g^{2,*}_{1,-}(\tbx,\eta)}\,, 
\end{equation*}
we have $T^{2,+}_{low}(\eta_1,\eta_2)f(\tbx)=T^{2,+}(\eta_1,\eta_2)f(\tbx)$ for $f\in C_0^\infty(\R^2)$ satisfying
\begin{equation}
  \label{eq:secular22+}
  \int_{\R^2} f(\tbx)\overline{g^{2,*}_{1,-}(\tbx,\eta)}\,dxdy=0
  \quad\text{for $\eta\in I_{\eta_1,\eta_2}$.}
\end{equation}
We remark that $\nabla M_+(v_{22})^*g^{2,*}_{1,-}(\cdot,\eta)=0$ and
$e^{-\a z_1}g^{2,*}_{1,-}(\cdot,\eta)\in L^\infty_yL^2_{z_1}$ if $\eta\in(-\eta_{*,3},\eta_{*,3})$.
We can prove the following in exactly the same way as Lemma~\ref{lem:T2lowpm}
\begin{lemma}
  \label{lem:T22+low}
  Let $\a\in(0,2\sqrt{c_{12}})$ and $\eta_{*,3}$ be as \eqref{def:eta*'}.
  Suppose that $0\le \eta_1\le\eta_2\le\eta_0<\eta_{*,3}$. If $f$ satisfies \eqref{eq:secular22+}, then
\begin{gather*}
\|T^{22,+}(\eta_1,\eta_2)f\|_{\calX_3} \le C\|f\|_{\calX_3}\,,
\end{gather*}
where $C$ is  a positive constant depending only on $\a$ and $\eta_0$. 
\end{lemma}

Next, we will investigate solvability of $\nabla M_-(v_{22})^*u=f$.
If $\eta_2\ge \eta_1>\eta_{*,3}$, then it follows from Lemma~\ref{lem:T22-pm} that
$T^{22,-}(\eta_1,\eta_2)^*$ is bounded on $\calX_3^*$.
On the other hand, we have $\nabla M_-(v_{22})\tg^M_1(\cdot,\eta)=0$ from Lemma~\ref{lem:eigenfunctions-mKP'}
and $e^{\a z_1}\tg^M_1(\eta)\in L^\infty_yL^2_{z_1}$ if $\eta\in(-\eta_{*,3},\eta_{*,3})$.
For $f\in C_0^\infty(\R^2)$, we have
\begin{equation*}
T^{22,-}(\eta_1,\eta_2)^*
=\sum_{1\le j\le 3}T^{22,-,*}_j(\eta_1,\eta_2)\,,
\end{equation*}
where
\begin{align*}
&  T^{22,-,*}_1(\eta_1,\eta_2)f(\tbx)
= \frac{1}{2\pi}\int_{I_{\eta_1,\eta_2}}\left(
\int_{-\infty}^\infty\int_{z_1}^\infty
\pd_{x'}\tfg^{22,-}_-(\tbxd,\tbx,-\eta)
f(\tbxd)\,dz_1'dy'\right)\,d\eta\,,
\\ &
T^{22,-,*}_2(\eta_1,\eta_2)f(\tbx)
= \frac{1}{2\pi}\int_{I_{\eta_1,\eta_2}}\left(
\int_{-\infty}^\infty\int^{z_1}_{-\infty}
\pd_{x'}\tfg^{22,-}_+(\tbxd, \tbx,-\eta)
f(\tbxd)\,dz_1'dy'\right)\,d\eta\,,
\end{align*}
and $T^{22,-,*}_3(\eta_1,\eta_2)=T^{22,+}_3(\eta_1,\eta_2)$.
Let
\begin{gather*}
 T^{22,-,*}_{low}(\eta_1,\eta_2)
=T^{22,-,*}_{2,low}(\eta_1,\eta_2)
+\sum_{j=1,3}T^{22,-,*}_j(\eta_1,\eta_2)f\,,
\\
\begin{split}
& T^{22,-,*}_{1,low}(\eta_1,\eta_2)f(\tbx)  
\\=& -\frac{1}{2\pi}\chi_-(z_1)\int_{I_{\eta_1,\eta_2}}\left(
\int_{-\infty}^\infty\int^{z_1}_{-\infty}
     \pd_{x'}\tfg^{22,-}_-(\tbxd,\tbx,-\eta)
     f(\tbxd)\,dz_1'dy'\right)\,d\eta
\\ & +
\frac{1}{2\pi}\chi_+(z_1)\int_{I_{\eta_1,\eta_2}}\left(
     \int_{-\infty}^\infty\int^\infty_{z_1}
     \pd_{x'}\tfg^{22,-}_-(\tbxd,\tbx,-\eta)
f(\tbxd)\,dz_1'dy'\right)\,d\eta\,.
\end{split}
\end{gather*}

By Claim~\ref{cl:Phi22*s'}, 
\begin{gather}
  \label{eq:tg22M}
  \pd_{x'}\tfg^{22,-}_-(\tbxd,\tbx,-\eta)
  =\frac{h_{22}(\tbx)}{2(\k_2-k_1)}
  \Phi^{12,*}(\tbx,-i\beta_{12}^-(-\eta))
  \tilde{g}^M_{1,+}(\tbxd,-\eta)\,,
\end{gather}
and we have
$T^{22,-}(\eta_1,\eta_2)^*f=T^{22,-,*}_{low}(\eta_1,\eta_2)f$
for  $f$ satisfying
\begin{equation}
  \label{eq:tgM2orth}
  \la f, \tilde{g}^M_{1,+}(\cdot,\eta) \ra=0
  \quad\text{for $\eta\in I_{\eta_1,\eta_2}$.}
\end{equation}
We can prove the following in exactly the same way as Lemmas~\ref{lem:T2-*} and \ref{lem:T2-*low}.

\begin{lemma}
  \label{lem:T22-*low}
  Let $\a\in(0,2\sqrt{c_{12}})$ and $0\le \eta_1\le \eta_2\le \eta_0<\eta_{*,3}$.
Suppose that $f$ satisfies \eqref{eq:tgM2orth}. Then
\begin{equation*}
\|T^{22,-}(\eta_1,\eta_2)^*f\|_{\calX_3^*}
\le C\|f\|_{\calX_3^*}\,,
\end{equation*}
where $C$ is a positive constant depending only on  $\a$ and $\eta_0$.
\end{lemma}

Finally, we will show that $\ker \nabla M_+(v_{22})=\{0\}$.
\begin{lemma}
  \label{lem:functional-22}
  Assume $0<\a<\sqrt{c_{12}}$. If $\varphi\in\calX_3$ and
  $\nabla M_+(v_{22})\varphi=0$, then $\varphi=0$.
\end{lemma}
\begin{proof}
Let $a(x,y)=-\Phi^2_1(\tbx)\Phi^{2,*}_1(\tbx)$. Then $a(x,y)>0$ and
$a(x,y)=O\left(\sech^2\sqrt{c_{12}}z_1\right)$.
Moreover, it follows from \eqref{eq:Miura-Lax2j} that
$$(\pd_y-\pd_x^2-2v_{22})a=\pd_x\nabla M_-(v_{212})a=0\,.$$
Thus we can prove Lemma~\ref{lem:functional-22} in exactly the same way as Lemma~\ref{lem:functional}.
\end{proof}

\bigskip

\subsection{Linear stability of $2$-line solitons of O-type}
The $2$-line soliton $u_2$ of O-type is linearly stable both in $\calX_3$ and
$\calX_4$.
\begin{proposition}
  \label{prop:linearstability-3}
  \begin{enumerate}
\item
  Assume \eqref{eq:ass-alpha''}, $\a<2\sqrt{c_{12}}$ and that $\eta'\in(0,\eta_{*,3})$.
Then there exist positive constants $K$ and $b$ such that
if $\eta_0\in(0,\eta']$,
$$\left\|e^{t\mL_2}\left(I-P_{21}(\eta_0)\right)\right\|_{B(\calX_3)}
\le Ke^{-bt}\quad\text{for $t\ge0$.}$$  
  \item 
  Assume \eqref{eq:ass-alpha'}, $\a<2\sqrt{c_{34}}$ and that $\eta'\in(0,\eta_{*,4})$.
Then there exist positive constants $K$ and $b$ such that
if $\eta_0\in(0,\eta']$,
$$\left\|e^{t\mL_2}\left(I-P_{22}(\eta_0)\right)\right\|_{B(\calX_4)}
\le Ke^{-bt}\quad\text{for $t\ge0$.}$$
\end{enumerate}
\end{proposition}

In view of Lemma~\ref{lem:evmL2'}, $\{\lambda_{1,\pm}(\eta)\}$ are
continuous eigenvalues of $\mL_{11}$ and $\mL_2$ and that
$\{\lambda_{2,\pm}(\eta)\}$ are continuous eigenvalues of $\mL_{12}$
and $\mL_2$.  We can prove the following in exactly the
same way as Lemmas~\ref{lem:lambda1-etaI} and \ref{lem:lambda2-est}.
\begin{lemma}
  \label{lem:lambdas-Otype}
  \begin{enumerate}
  \item For $i=1$, $2$ and $\eta\in\R\setminus\{0\}$,
$\eta\pd_\eta\Re\lambda_{i,\pm}(\eta)<0$
and $\Re\lambda_{i,\pm}(\eta)<\Re\lambda_{i,\pm}(0)=0$.

\item
Let $\mL_0$ and $\mL_2$ be operators on $\calX_3$, then
$\lambda_{1,\pm}(\eta_{*,3})\in\sigma(\mL_0)$ and 
\begin{equation*}
\{\lambda_{1,\pm}(\eta)\mid \eta\in[-\eta_{*,3},\eta_{*,3}]\}
\subset \sigma(\mL_2)\,.  
\end{equation*}
\item
Let $\mL_0$ and $\mL_2$ be operators on $\calX_4$, then
$\lambda_{2,\pm}(\eta_{*,4})\in\sigma(\mL_0)$ and 
\begin{equation*}
\{\lambda_{2,\pm}(\eta)\mid \eta\in[-\eta_{*,4},\eta_{*,4}]\}
\subset \sigma(\mL_2)\,.  
\end{equation*}
\end{enumerate}
\end{lemma}

If $\lambda\in\sigma(\mL_2)$ are resolvent point of linearized operators around
$[1,2]$-soliton or $[3,4]$-soliton, then $\lambda\in\sigma_p(\mL_2)$ and the corresponding
eigenfunctions are exponentially localized in $z_2$ or $z_1$, respectively.
\begin{lemma}
  \label{lem:specO-12}
Assume \eqref{eq:ass-alpha''} and $\a<\sqrt{c_{12}}$.
Let $\mL_{11}$, $\mL_{12}$ and $\mL_2$ be closed operator in $\calX_3$.
\begin{enumerate}
\item If $\lambda\in \sigma(\mL_2)\cap \rho(\mL_{11})\cap \rho(\mL_{12})$,
 then $\lambda$ is an isolated eigenvalue of $\mL_2$ with
finite multiplicity.
\item Suppose that $\mL_2u=\lambda u$ for $u\in\calX_3$ and
  $\lambda\in\rho(\mL_{12})$.  Then there exists an $\a'>0$ such that
  $e^{\a'|z_1|}u\in\calX_3$.
\end{enumerate}
\end{lemma}
\begin{lemma}
  \label{lem:specO-34}
  Assume \eqref{eq:ass-alpha'} and $\a<\sqrt{c_{34}}$.
Let $\mL_{11}$, $\mL_{12}$ and $\mL_2$ be closed operator in $\calX_4$.
\begin{enumerate}
\item If $\lambda\in \sigma(\mL_2)\cap \rho(\mL_{11})\cap \rho(\mL_{12})$,
then $\lambda$ is an isolated eigenvalue of $\mL_2$ with finite multiplicity.
\item
  Suppose that $\mL_2u=\lambda u$ for $u\in\calX_4$ and $\lambda\in\rho(\mL_{11})$.
Then there exists an $\a'>0$ such that  $e^{\a'|z_2|}u\in\calX_4$.
\end{enumerate}
\end{lemma}
We can prove Lemmas~\ref{lem:specO-12} and \ref{lem:specO-34} in exactly the
same way as Lemmas~\ref{lem:spec1} and \ref{lem:imev}.

Next, we consider the simplicity of continuous eigenvalues
$\lambda_{1,\pm}(\eta)$ and $\lambda_{2,\pm}(\eta)$.
\begin{lemma}
  \label{lem:anal-F-O}
  \begin{enumerate}
  \item   Assume \eqref{eq:ass-alpha''} and $\a<2\sqrt{c_{12}}$. Let
    $\mL_2$ be a closed operator in $\calX_3$, $\eta_0\in(0,\eta_{*,3})$
     and $\mathcal{C}_1=\{\lambda_{1,\pm}(\eta)\mid \eta\in(-\eta_0,\eta_0)\}$.
Then $\mathcal{C}_1\cap\sigma_p(\mL_2)$ is a discrete set and
$\lambda-\mL_2$ is invertible on $(I-P_{21}(\eta_0))\calX_3$
if $\lambda\in\mathcal{C}_1\setminus\sigma_p(\mL_2)$.
  
  \item   Assume \eqref{eq:ass-alpha'} and $\a<2\sqrt{c_{34}}$.
Let $\mL_2$ be a closed operator in $\calX_4$, $\eta_0\in(0,\eta_{*,4})$ and
$\mathcal{C}_2=\{\lambda_{2,\pm}(\eta)\mid \eta\in(-\eta_0,\eta_0)\}$.
Then $\mathcal{C}_2\cap\sigma_p(\mL_2)$ is a discrete set and
$\lambda-\mL_2$ is invertible on $(I-P_{22}(\eta_0))\calX_4$
if $\lambda\in\mathcal{C}_2\setminus\sigma_p(\mL_2)$.
  \end{enumerate}
\end{lemma}
We can prove Lemma~\ref{lem:anal-F-O} in exactly the same way as
Lemma~\ref{lem:anal-F1}.
Now we are in position to prove nonexistence of unstable eigenvalues.
\begin{lemma}
  \label{lem:Liouville21}
  Assume \eqref{eq:ass-alpha''} and $\a<2\sqrt{c_{12}}$. Let
  $\mL_{12}$ and $\mL_2$ be closed operators in $\calX_3$.
  If $\mL_2\varphi=\lambda\varphi$
  for $\varphi\in\calX_3$ and $\lambda\in\rho(\mL_{12})$,
  then $\varphi=0$.
\end{lemma}
\begin{lemma}
  \label{lem:Liouville22}
  Assume \eqref{eq:ass-alpha'} and $\a<2\sqrt{c_{34}}$.
  Let $\mL_{11}$ and $\mL_2$ be closed operators in $\calX_4$.
  If $\mL_2\varphi=\lambda\varphi$
  for $\varphi\in\calX_4$ and $\lambda\in\rho(\mL_{11})$,
  then $\varphi=0$.
\end{lemma}

\begin{proof}[Proof of Lemmas~\ref{lem:Liouville21} and ~\ref{lem:Liouville22}]
Proof follows the line of the proof of Lemma~\ref{lem:Liouville}.
  For $j=1$ and $2$, let
$$
\mL_{Mj}=\mL_0+\frac32\pd_x(v_{2j}^2\cdot)-
\frac{3}{2}(\pd_xv_{2j}\pd_x^{-1}\pd_y+\pd_x^{-1}\pd_yv_{2j}\pd_x)\,,$$
where $\pd_x^{-1}\pd_yv_{2j}=\pd_y\tau_2/\tau_2-\pd_y\tau_{1j}/\tau_{1j}$.
Then by \eqref{eq:bH+} and \eqref{eq:bH-},
we have intertwining properties
\begin{gather}
\label{eq:M2j-L2j}  
  \mL_2\nabla M_+(v_{2j})=\mL_{Mj}\nabla M_+(v_{2j})\,,\quad
  \mL_{1j}\nabla M_-(v_{2j})=\mL_{Mj}\nabla M_-(v_{2j})\,.
\end{gather}
\par
Suppose that $\lambda\in \sigma(\mL_2)\cap \rho(\mL_{12})$. 
By Lemmas~\ref{lem:specO-12}(2) and \ref{cl:m0-bound}, there exists an
$\a'>0$ such that $e^{\a'|z_1|}\pd_x^i\pd_y^j\varphi\in\calX_3$ for
any $i$, $j\ge0$.
It follows from Lemma~\ref{lem:spectral-proj'} that
$\mL_2P_{21}(\eta_0)=P_{21}(\eta_0)\mL_2$ for any
$\eta_0\in(0,\eta_{*,3})$. Moreover, $\lambda\ne\lambda_{1,\pm}(\eta_{*,3})$
since $\lambda\in\rho(\mL_{12})\subset \rho(\mL_0)$.
Therefore, we may assume that for an $\eps\in(0,\a')$,
\begin{equation}
  \label{eq:orth3*}
\la \varphi,g^{2,*}_{1,\pm}(\eta)\ra=0
\quad\text{for
$\eta\in[-\eta_{*,3}(\a+\eps),\eta_{*,3}(\a+\eps)]$ and $k=1$, $2$,}  
\end{equation}
where $\eta_{*,3}(\a)=2(\sqrt{c_{12}}+\a)\sqrt{\a(\a+\k_2-\k_1)}$.

Using Lemmas~\ref{lem:T22-pm} and \ref{lem:T22+low}, we can prove that
\begin{equation}
  \label{eq:phiM21}
\varphi_M=T^{22,+}(0,\infty)\varphi\in \calX_3\,,\quad
\nabla M_+(v_{22})\varphi_M=\varphi
\end{equation}
in the same way as Lemma~\ref{lem:Liouville}.
\par
By \eqref{eq:M2j-L2j} and \eqref{eq:phiM21},
\begin{gather*}
  \nabla M_+(v_{22})(\mL_{M2}-\lambda)\varphi_M=
  (\mL_2-\lambda)  \nabla M_+(v_{22})\varphi_M=0\,,
\end{gather*}
and it follows from Lemma~\ref{lem:functional-22}
that $\mL_{M2}\varphi_M=\lambda\varphi_M$.
Combining the above with \eqref{eq:M2j-L2j}, we have
\begin{equation}
  \label{eq:M-21=0}
(\mL_{12}-\lambda)\nabla M_-(v_{22})\varphi_M
=\nabla M_-(v_{22})(\mL_{M2}-\lambda)\varphi_M=0\,.  
\end{equation}
Since $\lambda\in\rho(\mL_{12})$, we have
\begin{equation}
  \label{eq:34}
\nabla M_-(v_{22})\varphi_M=0\,,
\end{equation}
and $\varphi_M=\frac12\pd_x^{-1}\varphi\in e^{-\eps|z_1|}\calX_3$.
Using Lemma~\ref{lem:eigenfunctions-mKP'},
\eqref{eq:orth3*} and \eqref{eq:phiM21}, we have
 for $\eta\in[-\eta_{*,3}(\a+\eps),\eta_{*,3}(\a+\eps)]$,
\begin{align*}
 2\la\varphi_M,\tg^{M,*}_1\ra
  =& \la \varphi_M,\nabla M_+(v_{22})^* g^{2,*}_{1,1}(\eta)\ra
= \la \varphi,g^{2,*}_{1,1}(\eta)\ra   =0\,.
\end{align*}
Now we will prove that $\varphi=0$.
Let $\eta_{*,3}(\a-\eps)<\eta_1<\eta_{*,3}<\eta_2<\eta_{*,3}(\a+\eps)$.
Let $\tg^{M,*}_{1r}(\tbx,\eta)$ be $0$
if $|\eta|\le \eta_1$ or $|\eta|\ge 2\eta_2-\eta_1$,
\begin{align*}
&  \tg^{M,*}_1(\tbx,\eta)e^{-iy\eta}
-\tg^{M,*}_1(\tbx,\pm\eta_1)e^{\mp iy\eta_1}
& \text{if $\pm\eta\in [\eta_1,\eta_2]$,}
\\ 
& \left(\tg^{M,*}_1(\tbx,\pm\eta_2)e^{\mp iy\eta_2}
-\tg^{M,*}_1(\tbx,\pm\eta_1)e^{\mp iy\eta_1}\right)
\frac{2\eta_2-\eta_1\mp\eta}{\eta_2-\eta_1}
& \text{if $\pm\eta\in[\eta_2,2\eta_2-\eta_1]$.}
\end{align*}
By Lemmas~\ref{lem:eigenfunctions-mKP'} and\ref{lem:orthogonality'},
\begin{equation*}
\label{eq:orthgMgM*'}
\la \tg^M_1(\cdot,\eta), \tg^{M,*}_1(\cdot,\eta_1)\ra= \delta(\eta-\eta_1)\,.
\end{equation*}
Following the proof of Claim~\ref{cl:Liouville-1} and
Lemma~\ref{lem:Liouville}, we can prove that there exists
an $a\in L^2(-\eta_2,\eta_2)$ such that
\begin{gather*}
\varphi_M^*=e^{2\a z_2}\varphi_M-\frac{1}{2\pi}\int_{-\eta_2}^{\eta_2}
a(\eta)\left(\tg^{M,*}_1(\tbx,\eta)-\tg^{M,*}_{1r}(\tbx,\eta)e^{iy\eta}\right)\,
d\eta
\end{gather*}
satisfies an orthogonal condition
\begin{equation}
  \label{eq:orthphi*'}
\int_{\R^2} \tg^M_1(\tbx,\eta)\overline{\varphi^*(\tbx)}\,d\tbx=0
\quad\text{for $\eta\in[-\eta_2,\eta_2]$.}
\end{equation}
Moreover, 
$$\la \varphi_M,\varphi_M^*\ra=
\|\varphi_M\|_{\calX_3}^2\left(1+O((\eta_2-\eta_1)^{1/2}\right)\,,$$
since $\la \varphi_M,\tg^{M,*}_1(\eta)\ra=0$ for
$\eta\in[-\eta_{*,3}(\a+\eps),\eta_{*,3}(\a+\eps)]$.
\par
On the other hand, it follows from Lemma~\ref{lem:T22-*low} and
\eqref{eq:orthphi*'} that there exists $w^*\in\calX_3^*$ satisfying
$\nabla M_-(v_{22})^*w^*=\varphi_M^*$ and
\begin{equation*}
\la \varphi_M,\varphi_M^*\ra=\la \nabla M_-(v_{22})\varphi_M,w^*\ra=0\,.
\end{equation*}
Thus we have $\varphi_M=0$ and $\varphi=M_-(v_{22})\varphi_M=0$.
\par
If $\lambda\in \sigma(\mL_2)\cap \rho(\mL_{11})$, we can prove
Lemma~\ref{lem:Liouville21} by using $\nabla M_\pm(v_{21})$.
This completes the proof of Lemma~\ref{lem:Liouville21}.
\par
Using Lemmas~\ref{lem:functional-21} and \ref{lem:specO-34}, we can prove 
Lemma~\ref{lem:Liouville22} in the same way.
\end{proof}

Now we are in position to prove Proposition~\ref{prop:linearstability-3}.
\begin{proof}[Proof of Proposition~\ref{prop:linearstability-3}]  
Let us prove the part (1) of Proposition~\ref{prop:linearstability-3}.
By Lemmas~\ref{lem:mL-12inX1}, \ref{lem:specO-12} and \ref{lem:Liouville21},
 $\sigma(\mL_2)\subset \sigma(\mL_{11})\cup\sigma(\mL_{12})$ and 
$\sigma(\mL_{12})\subset \{\lambda\in\C\mid \Re\lambda\le -b\}$ for
a $b>0$. Moreover, $\sup_{\eta_0\le|\eta|\le\eta_{*,3}}
\Re\lambda_{1,\pm}(\eta)=\Re\lambda_{1,\pm}(\eta_0)<0$ and it follows from
Lemmas~\ref{lem:anal-F-O} and \ref{lem:Liouville21} that
$\lambda_{1,\pm}(\eta)-\mL_2$ is invertible on $(I-P_{21}(\eta_0))\calX_3$.
Thus by the Gearhart-Pr\"{u}ss theorem (\cite{Gh,Pruss}),
we have the former part of Proposition~\ref{prop:linearstability-3}.

Using Lemmas~\ref{lem:mL-11inX2}, \ref{lem:specO-34}, \ref{lem:anal-F-O} and
\ref{lem:Liouville22}, we can prove the part (2) of
Proposition~\ref{prop:linearstability-3} in the same way.
This completes the proof of Proposition~\ref{prop:linearstability-3}.
\end{proof}

\bigskip

\subsection{Asymptotic linear stability of $2$-line solitons of O-type}
Let $\tau_2$ be as \eqref{eq:O-type}. Then $u_2=2\pd_x^2\log\tau$
can be rewritten as
\begin{gather}
\label{eq:2-lsol-O}
u_2=\varphi_{\bc}(\bz)=2\pd_x^2\log\tilde{\tau}_2\,,
\\ \label{eq:ttauO}
\tilde{\tau}_2=
\begin{vmatrix}
\cosh\sqrt{c_1}z_1 & \cosh\sqrt{c_2}z_2
\\ 
a_{12}\cosh\sqrt{c_1}z_1+\sqrt{c_1}\sinh\sqrt{c_1}z_1
&
a_{34}\cosh\sqrt{c_2}z_2+\sqrt{c_2}\sinh\sqrt{c_2}z_2
\end{vmatrix}\,,
\end{gather}
where $\bc=(c_1,c_2)=(c_{12},c_{34})$ and $\bz=(z_1,z_2)$ and
$z_1=x+2a_{12}y-\omega_{12}t$, $z_2=x+2a_{34}y-\omega_{34}t$, where
$a_{ij}$ and $\omega_{ij}$ are constants defined as \eqref{eqdef:a,c,omega}.
\par
We will investigate asymptotic behavior of solutions to
the linearized KP-II equation around $u_2$.
\par
Let $H_t(y)=(4\pi t)^{-1/2}\exp(-y^2/4t)$,
$W_t(y)=\frac12$ if $(2\k_1+\k_2)t\le y\le (\k_1+2\k_2)t$
and $W_t(y)=0$ otherwise.
\begin{theorem}
  \label{thm:AS3}
Assume \eqref{eq:ass-alpha''} and $\a<2\sqrt{c_1}$. Suppose
$u_0\in\calX_3\cap L^1_yL^2(\R;e^{2\a z_1}dz_1)$ and that
$u(t)$ is a solution of
\begin{equation}
  \label{eq:27aa}
4\pd_tu+\pd_x^3u+3\pd_x^{-1}\pd_y^2u+6\pd_x(u_2u)=0\,,\quad u(0)=u_0\,,
\end{equation}
in the class $ C([0,\infty);\calX_3)$. Then
$$\left\|e^{\a z_1}\left\{u(t,\cdot)-\left(H_{t/2\sqrt{c_1}}*
W_t*f(y)\right)
\pd_{z_1}\varphi_\bc(\bz)\right\} \right\|_{L^2(\R^2)}=O(t^{-1/4})\quad
\text{as $t\to\infty$,}$$
where $$f(y)=-\frac{1}{4c_1}\int_\R u_0(x,y)
\left(\int_{-\infty}^x \pd_{z_1}\varphi_\bc  (x'+2a_{12}y,x'+2a_{34}y)\,dx'\right)\,dx\,.$$
\end{theorem}
Let $\widetilde{W}_t(y)=\frac12$
if $(2\k_3+\k_4)t\le y\le (\k_3+2\k_4)t$ and $\widetilde{W}_t(y)=0$ otherwise.
\begin{theorem}
  \label{thm:AS4}
Assume \eqref{eq:ass-alpha'} and $\a<2\sqrt{c_2}$. Suppose
$u_0\in\calX_4\cap L^1_yL^2(\R;e^{2\a z_2}dz_2e)$ and that
$u(t)$ is a solution of \eqref{eq:27aa}
in the class $ C([0,\infty);\calX_4)$. Then
$$\left\|e^{\a z_2}\left\{u(t,\cdot)-\left(H_{t/2\sqrt{c_1}}*\widetilde{W}_t*f(y)\right)
\pd_{z_2}\varphi_\bc(\bz)\right\} \right\|_{L^2(\R^2)}=O(t^{-1/4})\quad
\text{as $t\to\infty$,}$$
where $$f(y)=-\frac{1}{4c_2}\int_\R u_0(x,y)
\left(\int_{-\infty}^x \pd_{z_2}\varphi_\bc  (x'+2a_{12}y,x'+2a_{34}y)\,dx'\right)\,dx\,.$$
\end{theorem}

We can prove Theorems\ref{thm:AS3} and \ref{thm:AS4}
in exactly the same way as Theorem~\ref{thm:AS} by
using Proposition~\ref{prop:linearstability-3} and
the asymptotic profiles of $g^2_{j,k}(\eta)$ and
$g^{2,*}_{j,k}(\eta)$ as $\eta\to0$ (Lemma~\ref{lem:gg*profile'} below).

Let 
\begin{gather*}
\gamma_1^\pm(c_2)
=\frac{1}{2\sqrt{c_1}}
\log\frac{a_{34}\pm\sqrt{c_2}-\k_2}{a_{34}\pm\sqrt{c_2}-\k_1}\,,
\quad
\gamma_2^\pm(c_1)
=\frac{1}{2\sqrt{c_2}}
\log\frac{\k_4-a_{12}\mp\sqrt{c_1}}{\k_3-a_{12}\mp\sqrt{c_1}}\,,
\\
\varphi^{1,\pm}(\bz)=\frac{d}{dc_1}
\varphi_\bc\left(z_1,z_2'-\gamma_2^\pm(c_1)\right)
\bigr|_{z_2'=z_2+\gamma_2^\pm(c_1)}\,,
\\
\varphi^{2,\pm}(\bz)=\frac{d}{dc_2}
\varphi_\bc\left(z_1'-\gamma_1^\pm(c_2),z_2\right)
\bigr|_{z_1'=z_1+\gamma_1^\pm(c_2)}\,.
\end{gather*}

\begin{lemma}
  \label{lem:gg*profile'}
For $k=1$ and $2$,
\begin{gather}
  \label{eq:35'}
  g^2_{k,1}(x,y,\eta)= -\frac{1}{2\sqrt{c_k}}e^{iy\eta}
  \left(\pd_{z_k}\varphi_\bc(\bz)+O(\eta)\right)\,,
  \\
\notag
g^2_{k,2}(x,y,\eta)=
\frac12 e^{iy\eta}\left(\varphi^{k,+}+\frac{1}{2c_k^{3/2}}\pd_{z_k}\varphi_\bc+O(\eta)\right)\,,       
\\ \notag
g^{2,*}_{k,1}(x,y,\eta)=
\frac{1}{2}e^{iy\eta}\left(\int_{-\infty}^x\varphi^{k,-}+O(\eta)\right)\,,
\\ \label{eq:28'}
g^{2,*}_{k,2}(x,y,\eta)=
\frac{1}{2\sqrt{c_k}}e^{iy\eta}
\left(\int_{-\infty}^x\pd_{z_k}\varphi_\bc+O(\eta)\right)\,.
\end{gather}
\end{lemma}
\begin{proof}[Proof of Lemma~\ref{lem:gg*profile'}]
Using \eqref{eq:g21pm-O} and \eqref{eq:g22pm-O},
we have
\begin{equation*}
g^2_{1,1}= 2e^{iy\eta}
\pd_x^2\sum_{j=3,4}\frac{e^{\theta_1+\theta_j}}{\tau_2}
\left(\k_j-\k_1+O(\eta)\right)\,,  
\end{equation*}
\begin{multline*}
g^2_{1,2}= \frac{-2e^{iy\eta}}{\k_2-\k_1}
\pd_x^2\sum_{j=3,4}\frac{e^{\theta_1+\theta_j}}{\tau_2}
\Bigl\{(\k_j-\k_1)z_1
  +\frac{a_{12}-\k_j}{\k_j-\k_2}+\frac{2(\k_j-\k_1)}{\k_2-\k_1}
+O(\eta) \Bigr\}\,,
\end{multline*}  
\begin{gather*}
  g^{2,*}_{1,1}=\frac{2e^{iy\eta}}{\k_2-\k_1}e^{iy\eta}
\pd_x\sum_{j=3,4}\frac{e^{\theta_2+\theta_j}}{\tau_2}                    
\left\{(\k_j-\k_2)z_2-\frac{\k_j-a_{12}}{\k_j-\k_1}+O(\eta)\right\}\,,
  \\
g^{2,*}_{1,2}= 2e^{iy\eta}\pd_x\sum_{j=3,4}\frac{e^{\theta_2+\theta_j}}{\tau_2}  
\left(\k_j-\k_2+O(\eta)\right)\,,
\end{gather*}
\begin{equation*}
g^2_{2,1}= 2e^{iy\eta}
\pd_x^2\sum_{j=1,2}\frac{e^{\theta_3+\theta_j}}{\tau_2}
\left(\k_3-\k_j+O(\eta)\right)\,,
\end{equation*}
\begin{multline*}
g^2_{2,2}= \frac{-2e^{iy\eta}}{\k_4-\k_3}
\pd_x^2\sum_{j=1,2}\frac{e^{\theta_3+\theta_j}}{\tau_2}
\Bigl\{(\k_3-\k_j)z_2
  +\frac{a_{34}-\k_j}{\k_4-\k_j}+\frac{2(\k_3-\k_j)}{\k_4-\k_3}
  \Bigr\}+O(\eta)\,,  
\end{multline*}  
\begin{gather*}
  g^{2,*}_{2,1}=\frac{2e^{iy\eta}}{\k_4-\k_3}
\pd_x\sum_{j=1,2}\frac{e^{\theta_4+\theta_j}}{\tau_2}                    
\left\{(\k_4-\k_j)z_2+\frac{a_{34}-\k_j}{\k_3-\k_j}+O(\eta)\right\}\,,
  \\
g^{2,*}_{2,2}= 2e^{iy\eta}\pd_x\sum_{j=1,2}\frac{e^{\theta_4+\theta_j}}{\tau_2}  
\left(\k_4-\k_j+O(\eta)\right)\,,
\end{gather*}
in the same way as Lemma~\ref{lem:gg*profile}.
\par
Differentiating \eqref{eq:ttauO} with respect to $z_1$, $z_2$, $c_1$ and $c_2$, we have
we have
\begin{align*}
& 1+\frac{1}{\sqrt{c_1}}\frac{\pd_{z_1}\tilde{\tau}_2}{\tilde{\tau}_2}
 =2\sum_{j=3,4}(\k_j-\k_2)\frac{e^{\theta_2+\theta_j}}{\tau_2}\,,
  \\ & 
1-\frac{1}{\sqrt{c_1}}\frac{\pd_{z_1}\tilde{\tau}_2}{\tilde{\tau}_2}
=2\sum_{j=3,4}(\k_j-\k_1)\frac{e^{\theta_1+\theta_j}}{\tau_2}\,,
  \\   &
1+\frac{1}{\sqrt{c_2}}\frac{\pd_{z_2}\tilde{\tau}_2}{\tilde{\tau}_2}
 =2\sum_{j=1,2}(\k_4-\k_j)\frac{e^{\theta_j+\theta_4}}{\tau_2}\,,
  \\ &
1-\frac{1}{\sqrt{c_2}}\frac{\pd_{z_2}\tilde{\tau}_2}{\tilde{\tau}_2}
=2\sum_{j=1,2}(\k_3-\k_j)\frac{e^{\theta_j+\theta_3}}{\tau_2}\,,
\end{align*}
\begin{align*}
&  z_1+(\k_2-\k_1)\frac{\pd_{c_1}\tilde{\tau}_2}{\tilde{\tau}_2}
  = 2\sum_{j=3,4}(\k_j-\k_2)z_1\frac{e^{\theta_2+\theta_j}}{\tau_2}
-\frac{f_1^*f_2}{\tau_2}\,,
\\ &
  z_1-(\k_2-\k_1)\frac{\pd_{c_1}\tilde{\tau}_2}{\tilde{\tau}_2}
  = 2\sum_{j=3,4}(\k_j-\k_1)z_1\frac{e^{\theta_1+\theta_j}}{\tau_2}
+\frac{f_1^*f_2}{\tau_2}\,,
  \\ &
z_2+(\k_4-\k_3)\frac{\pd_{c_2}\tilde{\tau}_2}{\tilde{\tau}_2}
  = 2\sum_{j=1,2}(\k_4-\k_j)z_2\frac{e^{\theta_j+\theta_4}}{\tau_2}
+\frac{f_1f_2^*}{\tau_2}\,,
\\ &
  z_2-(\k_4-\k_3)\frac{\pd_{c_2}\tilde{\tau}_2}{\tilde{\tau}_2}
 = 2\sum_{j=1,2}(\k_3-\k_j)z_2\frac{e^{\theta_j+\theta_3}}{\tau_2}
-\frac{f_1f_2^*}{\tau_2}\,,
\end{align*}
where $f_1^*=e^{\theta_2}-e^{\theta_1}$, $f_2^*=e^{\theta_4}-e^{\theta_3}$.
Combining the above, we can prove Lemma~\ref{lem:gg*profile'}
in exactly the same way as Lemma~\ref{lem:gg*profile}.
\end{proof}

\bigskip

\section*{Acknowledgments}
This work was supported by JSPS KAKENHI Grant Number JP21K03328.  A
part of this research was done in the Institute of Mathematics,
Academia Sinica in Taiwan. The author would like to express his
gratitude to Professor Shih-Hsien Yu and Professor Derchyi Wu for
their hospitality.  I would like to thank Professor Vladimir Varlamov
who suggested me to study KP equations during his lifetime.  Last but
not least, I would express my deep gratitude to Professor Robert
L.~Pego who kindly accepted the author as a visitor in 2006 and 2011
and taught many things to the author.  \bigskip

\appendix
\section{Proof of Lemmas~\ref{lem:M-pdPhiPhi*} and \ref{lem:M-PhiPhi*}}
\label{sec:lem:M-pdPhiPhi*}

To prove Lemmas~\ref{lem:M-pdPhiPhi*} and \ref{lem:M-PhiPhi*}, we need the following.
\begin{claim}
    \label{cl:L0*L0}
    \begin{gather}
\label{eq:Phi1-Phi2-1}
\pd_x\left(h^{-1}\Phi^1(\bx,-i\beta)\right)=h^{-1}\Phi^2(\bx,-i\beta)\,,\\
\pd_x\left(\tau_1^{-1}\Phi^0(\bx,-i\beta)\right)=\tau_1^{-1}\Phi^1(\bx,-i\beta)\,,\\
\label{eq:Phi1-Phi2-2}
\pd_x\left(h\Phi^{2,*}(\bx,-i\beta)\right)=-h\Phi^{1,*}(\bx,-i\beta)\,,\\
\pd_x\left(\tau_1\Phi^{1,*}(\bx,-i\beta)\right)=-\tau_1\Phi^{0,*}(\bx,-i\beta)\,,\\
\label{eq:Phi1*-Phi2*-1}
L_0^*\left(h^{-1}\Phi^1(\bx,-i\beta)\right)=2h^{-1}\pd_x\Phi^2(\bx,-i\beta)\,,\\
L_0^*\left(\tau_1^{-1}\Phi^0(\bx,-i\beta)\right)
=2\tau_1^{-1}\pd_x\Phi^1(\bx,-i\beta)\,,\\
\label{eq:Phi1*-Phi2*-2}
L_0\left(h\Phi^{2,*}(\bx,-i\beta)\right)=-2h\pd_x\Phi^{1,*}(\bx,-i\beta)\,,\\
L_0\left(\tau_1\Phi^{1,*}(\bx,-i\beta)\right)
=-2\tau_1\pd_x\Phi^{0,*}(\bx,-i\beta)\,.
    \end{gather}
\end{claim}
\begin{proof}
By a straightforward computation,
\begin{align*}
  \pd_x\left(h^{-1}\Phi^1(\bx,-i\beta)\right)=&
\Phi^0(\bx,-i\beta)\frac{\tau_1}{\tau_2^2}\left(
\beta^2\tau_2-\beta\pd_x\tau_2+\operatorname{Wr}(\pd_xf_1,\pd_xf_2)\right)\,,
\\=& h^{-1}\Phi^2(\bx,-i\beta)\,.                                                
\end{align*}
By \eqref{eq:th} and the fact that $L_1\Phi^1=0$,
\begin{align*}
 h L_0^*\left(h^{-1}\Phi^1(\bx,-i\beta)\right)
  =& (2\pd_x^2+u_1-u_2)\Phi^1(\bx,-i\beta)
     -2h^{-1}(\pd_xh)\pd_x\Phi^1(\bx,-i\beta)
  \\=& 2\pd_x\left\{h\pd_x\left(h^{-1}\Phi^1(\bx,-i\beta)\right)\right\}\,.
\end{align*}
Combining the above with \eqref{eq:Phi1-Phi2-1}, we have \eqref{eq:Phi1-Phi2-2}.
\par
Let $\tau_{2,\Phi^*}(\bx,-i\beta)
=\operatorname{Wr}((\beta-\pd_x)^{-1}f_1,(\beta-\pd_x)^{-1}f_2)$. Then
\begin{align*}
  \pd_x\left(h\Phi^{2,*}(\bx,-i\beta)\right)=&
  \Phi^{0,*}(\bx,-i\beta)\frac{-\tau_{2,\Phi^*}(\bx,-i\beta)
    (\beta+\pd_x)\tau_1 +\tau_1\pd_x\tau_{2,\Phi^*}(\bx,-i\beta)}{ \tau_1^2}
\\=&\frac{\Phi^{0,*}(\bx,-i\beta)}{\tau_1^2}\sum_{j=0,1}(-1)^{j+1}
\pd_x^jf_1 \begin{vmatrix} (\beta-\pd_x)^{-1}f_1 & (\beta-\pd_x)^{-1}f_2
\\ \pd_x^{1-j}f_1 & \pd_x^{1-j}f_2 \end{vmatrix}
  \\=& -h\Phi^{1,*}(\bx,-i\beta)\,.
\end{align*}
By \eqref{eq:h} and the fact that $L_2^*\Phi^{2,*}=0$,
\begin{align*}
  L_0\left(h\Phi^{2,*}(\bx,-i\beta)\right)=&
(L_0h)\Phi^{2,*}(\bx,-i\beta)+2(\pd_xh)\pd_x\Phi^{2,*}(\bx,-i\beta)+hL_0\Phi^{2,*}(\bx,-i\beta)
  \\=&
h(2\pd_x^2+u_2-u_1)\Phi^{2,*}(\bx,-i\beta)+2(\pd_xh)\pd_x\Phi^{2,*}(\bx,-i\beta)\,.
\end{align*}
Since $u_2-u_1=2\pd_x^2\log h$, it follows from \eqref{eq:Phi1-Phi2-2} that
\begin{align*}
  L_0\left(h\Phi^{2,*}(\bx,-i\beta)\right)=&
2\pd_x^2\left(h\Phi^{2,*}(\bx,-i\beta)\right)
 -2h^{-1}(\pd_xh)\pd_x\left(h\Phi^{2,*}(\bx,-i\beta)\right)
\\=& -2h\pd_x\Phi^{1,*}(\bx,-i\beta)\,.
\end{align*}
We can prove the rest in exactly the same way.
Thus we complete the proof.
\end{proof}

\begin{proof}[Proof of Lemmas~\ref{lem:M-pdPhiPhi*} and \ref{lem:M-PhiPhi*}]
  By \eqref{eq:Miura-Lax1} and the fact that $L_2^*\Phi^{2,*}(\bx,-i\beta')=0$,
  \begin{align*}
 &  \nabla M_+(v_2)\pd_x\left(\Phi^1(\bx,-i\beta)\Phi^{2,*}(\bx,-i\beta')\right)
\\ =& hL_2^*\left(h^{-1}\Phi^1(\bx,-i\beta)\Phi^{2,*}(\bx,-i\beta')\right)
    \\=&
2h\pd_x\left(h^{-1}\Phi^1(\bx,-i\beta)\right)\pd_x\Phi^{2,*}(\bx,-i\beta')
+\Phi^{2,*}(\bx,-i\beta')hL_0^*\left(h^{-1}\Phi^1(\bx,-i\beta)\right)\,.
  \end{align*}
Combining the above with \eqref{eq:Phi1-Phi2-1} and \eqref{eq:Phi1*-Phi2*-1},
we have \eqref{eq:M2+pdPhiPhi*}.
\par
By \eqref{eq:Miura-Lax1}, \eqref{eq:Phi1-Phi2-1}, \eqref{eq:Phi1*-Phi2*-1}
and the fact that $L_1^*\Phi^{1,*}=0$,
  \begin{align*}
    & h^{-1}\pd_x\nabla M_-(v_2)^*\left(\Phi^1(\bx,-i\beta)
      \Phi^{1,*}(\bx,-i\beta')\right)
\\=& L_1^*\left(h^{-1}\Phi^1(\bx,-i\beta)\Phi^{1,*}(\bx,-i\beta)'\right)
\\=& 2\pd_x\left(h^{-1}\Phi^1(\bx,-i\beta)\right)\pd_x\Phi^{1,*}(\bx,-i\beta')
     + L_0^*\left(h^{-1}\Phi^1(\bx,-i\beta)\right)\Phi^{1,*}(\bx,-i\beta')
    \\=& 2h^{-1}\pd_x\left(\Phi^2(\bx,-i\beta)\Phi^{1,*}(\bx,-i\beta')\right)\,.
  \end{align*}
  \par
By \eqref{eq:Miura-Lax2}, \eqref{eq:Phi1-Phi2-2}, \eqref{eq:Phi1*-Phi2*-2}
and the fact that $L_1\Phi^1=L_2\Phi^2=0$,
\begin{align*}
&  h\nabla M_-(v_2)\pd_x\left(\Phi^1(\bx,-i\beta)\Phi^{2,*}(\bx,-i\beta')\right)
\\=& -L_1\left(h\Phi^1(\bx,-i\beta)\Phi^{2,*}(\bx,-i\beta')\right)
  \\=& -2\pd_x\Phi^1(\bx,-i\beta)\pd_x\left(h\Phi^{2,*}(\bx,-i\beta')\right)
       -\Phi^1(\bx,-i\beta)L_0\left(h\Phi^{2,*}(\bx,-i\beta')\right)
  \\=& 2h\pd_x\left(       \Phi^1(\bx,-i\beta)\Phi^{1,*}(\bx,-i\beta')\right)\,,
\end{align*}
  \begin{align*}
& h\pd_x\nabla M_+(v_2)^*\left(\Phi^2(\bx,-i\beta)
 \Phi^{2,*}(\bx,-i\beta')\right)
\\=&
-L_2\left(h\Phi^2(\bx,-i\beta)\Phi^{2,*}(\bx,-i\beta')\right)
 \\=& 
-2(\pd_x\Phi^2(\bx,-i\beta))\pd_x\left(h\Phi^{2,*}(\bx,-i\beta')\right)
-\Phi^2(\bx,-i\beta))L_0\left(h\Phi^{2,*}(\bx,-i\beta')\right)
\\=& 2h\pd_x
     \left(\Phi^2(\bx,-i\beta)\Phi^{1,*}\bx,-i\beta')  \right)\,.
\end{align*}
\par
Using \eqref{eq:Miura-Lax3}, \eqref{eq:Miura-Lax4} and Claim~\ref{cl:L0*L0}, 
we can prove the rest in the same way.
\end{proof}

\section{Estimates for the Lax operators}
\label{sec:est-Lax}
For a standard bootstrap argument of $\nabla M_\pm(v_2)$ in Section~\ref{sec:LS-P},
we need the following.
\begin{lemma}
  \label{cl:m0-bound}
Let $\a\in\R\setminus\{0\}$, $a\in\R$ and $\calX=L^2(\R^2;e^{2\a(x+2ay)}dxdy)$.
Then $\pm\pd_x+\pd_x^{-1}\pd_y$ are closed operators on $\calX$. 
If $f$, $u\in \calX$ and $(\pm\pd_x+\pd_x^{-1}\pd_y)u=f$, then
\begin{align*}
 \|\pd_xu\|_\calX+\|\pd_x^{-1}\pd_yu\|_\calX \lesssim \|u\|_\calX+\|f\|_\calX\,.
\end{align*}
\end{lemma}
\begin{proof}
An operator $L_0^*=\pd_x^2+\pd_y$ is closed in $\calX$
since $\lambda-L_0^*$ has bounded inverse on $\calX$ for large $\lambda>0$.
Indeed,
  \begin{align*}
    \|(\lambda-L_0^*)^{-1}f\|_\calX
    =& \left\|\frac{\hat{f}(\xi+i\a,\eta+2ia\a)}
       {\lambda-(i\xi-\a)^2-i\eta+2a\a}\right\|_{L^2(\R^2)}
\\ \le &\frac{1}{\Re\lambda-\a^2+2a\a}\|f\|_\calX\,.  
  \end{align*}
Since $\pd_x+\pd_x^{-1}\pd_y=L_0^*\pd_x^{-1}$
and $\pd_x^{-1}$ is bounded on $\calX$,
$\pd_x+\pd_x^{-1}\pd_y$ is closed on $\calX$.
We can prove the closedness of $-\pd_x+\pd_x^{-1}\pd_y$ in the same way.
\par    
Let $m_{0,\pm}(\xi,\eta)=\pm i\xi+\eta/\xi$.
  Then
  $$m_{0,\pm}(\xi+i\a,\eta+2ia\a)=
  i\left\{\frac{\a(2a\xi-\eta)}{\xi^2+\a^2}\pm\xi\right\}
  \mp\a+\frac{\xi\eta+2\a^2a}{\xi^2+\a^2}\,,$$
  and there exists $R>0$ such that
  \begin{align*}
& |\Im m_{0,\pm}(\xi+i\a,\eta+2ia\a)|\gtrsim  \la \xi\ra
      \quad \text{if $|\xi|\ge R$ and $|\eta|\le \la\xi\ra^2$,}
\\ & |\Re m_{0,\pm}(\xi+i\a,\eta+2ia\a)|\gtrsim \la\xi\ra
      \quad \text{if $|\xi|\ge R$ and $|\eta|\ge \la\xi\ra^2$,}
\\ & |\Re m_{0,\pm}(\xi+i\a,\eta+2ia\a)|\gtrsim |\eta|\la\xi\ra^{-1}
      \quad \text{if $|\xi|\ge 1$ and $|\eta|\ge R\la\xi\ra$,}
\\ & |\Im m_{0,\pm}(\xi+i\a,\eta+2ia\a)|\gtrsim |\eta|
      \quad \text{if $|\xi|\le 1$ and $|\eta|\ge R\la\xi\ra$,}
  \end{align*}
\begin{equation}
      \label{eq:24}
\sup_{|\xi|+|\eta|\la \xi\ra^{-1}\ge 2R}
\frac{|\xi|+|\eta|\la\xi\ra^{-1}}{|m_{0,\pm}(\xi+i\a,\eta+2i\a a)|}
<\infty\,.    
\end{equation}
By the Plancherel theorem,
\begin{equation}
  \label{eq:28}
\|f\|_\calX=\|\hat{f}(\xi+i\a,\eta+2i\a a)\|_{L^2(\R^2)}\quad\text{for $f\in\calX$.}
\end{equation}
Combining \eqref{eq:24} and \eqref{eq:28}, we have Lemma~\ref{cl:m0-bound}.
Thus we complete the proof.
\end{proof}

\section{Resolvent estimates}
\par
To prove smoothness of eigenfunctions of $\mL_2$, we need the following.
\begin{proposition}
  \label{lem:mLu=f}
Let $\a\in\R\setminus\{0\}$ and $\calX$ be as in Lemma~\ref{cl:m0-bound}.
Suppose that $\mL_0$ is a closed operator on $\calX$ and that
$\mL_0u=f$. Then
$$\sum_{i=1,2}\|\pd_x^iu\|_{\calX}+\|\pd_yu\|_{\calX}
\le C(\|u\|_\calX+\|f\|_\calX)\,,$$
where $C$ is a constant that does not depend on $f$.
\end{proposition}
\begin{proof}
Since
\begin{align*}
 \Re p(\xi+i\a,\eta+2ia\a)=&
\frac14\xi(\xi^2-3\a^2+4b_1)
-\frac{3}{4}\frac{\xi}{\xi^2+\a^2}(\eta^2-4a^2\a^2)
\\ &  +\eta\left(b_2-\frac{3}{2}\frac{a\a^2}{\xi^2+\a^2}\right)\,,
\end{align*}
\begin{equation}
  \label{eq:Rep-est}
|\Re p(\xi+i\a,\eta+2ia\a)|
\gtrsim 
\begin{cases}
\la \xi\ra^3 \quad & \text{if $\la \xi\ra^2\gg \la \eta\ra$,}\\
|\xi|^{-1}\eta^2\quad & \text{if $|\eta|\gg \xi^2\ge1$.}
\end{cases}  
\end{equation}
Since
\begin{equation*}
\Im p(\xi+i\a,\eta+2ia\a)  =\frac{\a}{4}
\left(3\xi^2-\a^2+3\frac{(\eta-2a\xi)^2}{\xi^2+\a^2}\right)
+\a(b_1+2ab_2-3a^2)\,,
\end{equation*}
there exists a constant $C$ such that
\begin{equation}
\label{eq:Imp-est}
\Im p(\xi+i\a,\eta+2ia\a)+C\gtrsim \begin{cases}
\la\xi\ra^2\quad & \text{for every $(\xi,\eta)$,}
\\ 
\la\eta\ra^2\quad & \text{if $|\eta|\gg1\ge |\xi|$.} 
\end{cases}
\end{equation}
By \eqref{eq:Rep-est} and \eqref{eq:Imp-est},
\begin{equation}
  \label{eq:p-a14}
|p(\xi+i\a,\eta+2ia\a)+iC|\gtrsim \left\{
  \begin{aligned}
\la \xi\ra^3 \quad &\text{if $\la\xi\ra^2\gg \la\eta\ra$,}
\\
\la \xi\ra^2 \quad &\text{if $\la\xi\ra^2\simeq \la\eta\ra$,} 
\\
\la\eta\ra^2\la \xi\ra^{-1} \quad &\text{if $\la\xi\ra^2\ll \la\eta\ra$.}   
  \end{aligned}\right.
\end{equation}
Combining the above \eqref{eq:p-a14} with \eqref{eq:28},
we have Proposition~\ref{lem:mLu=f}.
Thus we complete the proof.
\end{proof}

\begin{lemma}
  \label{cl:com-L0}
  Let $\a>0$ and $\mL_0$ be a closed operator on $\calX_2$.
Suppose that  $\lambda\in\rho(\mL_0)$. Then for $i$, $j=1$, $2$,
$$\|[\mL_0,\chi_{i,R}(z_1)\chi_{j,R}(z_2)](\lambda-\mL_0)^{-1}\|_{B(\calX_2)}
\le CR^{-1}\,,$$
where $C$ is a constant independent of $R$ and locally uniform with
respect to $\a$.
\end{lemma}
\begin{lemma}
  \label{cl:com-L1}
Let $\a>0$ and $\mL_1$ be a closed operator on $\calX_2$.
Suppose that $\lambda\in \rho(\mL_0)\cap\rho(\mL_1)$. Then for $i$, $j=1$, $2$,
$$\|[\mL_1,\chi_{i,R}(z_1)\chi_{j,R}(z_2)](\lambda-\mL_1)^{-1}\|_{B(\calX_2)}
\le CR^{-1}\,,$$
where $C$ is a constant independent of $R$ and locally uniform with
respect to $\a$.
\end{lemma}
\begin{proof}[Proof of Lemmas~\ref{cl:com-L0} and \ref{cl:com-L1}]
  By Proposition~\ref{lem:mLu=f} and \eqref{eq:m(D)-bound},
  \begin{gather}
\label{eq:com-1}
\sum_{0\le i+2j\le2}\|\pd_x^i\pd_y^j(\lambda-\mL_0)^{-1}f\|_{\calX_2}
\lesssim \|f\|_{\calX_2}\,.
\end{gather}
Moreover, it follows from \eqref{eq:28} and \eqref{eq:p-a14} that
\begin{gather}
\label{eq:com-2}
 \|\pd_x^{-2}\pd_y^2(\lambda-\mL_0)^{-1}f\|_{\calX_2}
\lesssim \|f\|_{\calX_2}\,.
  \end{gather}
Let $\chi(x,y)=\chi_{i,R}(z_1)\chi_{j,R}(z_2)$. We have
  \begin{align}
\notag
 [\pd_x^{-1}\pd_y^2,\chi]=&\pd_x^{-1}[\pd_y^2,\chi]
-\pd_x^{-1}\chi_x\pd_x^{-1}\pd_y^2\,,
\\= & \label{eq:com-2''}
\pd_x^{-1}(\chi_{yy}+2\chi_y\pd_y)
- \left\{\chi_x-\pd_x^{-1}(\chi_{xx}\cdot)\right\}\pd_x^{-2}\pd_y^2\,.
\end{align}
Since $\|\pd_x^m\pd_y^n\chi\|_{L^\infty}=O(R^{-m-n})$ and
$\pd_x^{-1}$ is bounded on $\calX_2$, we see that
Lemma~\ref{cl:com-L0} follows from
\eqref{eq:com-1}--\eqref{eq:com-2''}.
\par
By \eqref{eq:com-1}, \eqref{eq:com-2} and the resolvent identity
$$(\lambda-\mL_1)^{-1}-(\lambda-\mL_0)^{-1}
=-\frac{3}{2}(\lambda-\mL_0)^{-1}\pd_x\{u_1(\lambda-\mL_1)^{-1}\}\,,$$
we have
  \begin{gather}
\label{eq:com-3}
\sum_{0\le i+2j\le2}\|\pd_x^i\pd_y^j(\lambda-\mL_1)^{-1}f\|_{\calX_2}
  \lesssim \|f\|_{\calX_2}\,,  
\quad
 \|\pd_x^{-2}\pd_y^2(\lambda-\mL_1)^{-1}f\|_{\calX_2}
\lesssim \|f\|_{\calX_2}\,.
  \end{gather}
Using \eqref{eq:com-2''}--\eqref{eq:com-3}, we can
prove Lemma~\ref{cl:com-L1} in the same way as Lemma~\ref{cl:com-L0}.
Thus we complete the proof.
\end{proof}

\begin{lemma}
\label{lem:imev-5}
Let $\a>0$, $a\in\R$ and $z=x+2ay$. Let $k=0$, $1$, $2$ and $\mL_k$ be
closed operators on $\calX=L^2(\R^2;e^{2\a z})$.  Suppose that
$\lambda\in\rho(\mL_k)$. There exists $\beta_0\in(0,\a)$ such that if
$(\lambda-\mL_k)u=f$, $u\in\calX$ and $(1+e^{\beta z})f\in \calX$
for $\beta\in[-\beta_0,\beta_0]$, then
$$\|e^{\beta z}u\|_\calX\le C\|e^{\beta z}f\|_{\calX}\,,$$
where $C$ is a positive constant independent of $u$ and $f$.
\end{lemma}
\begin{proof}
Suppose that $\beta>0$ and that $(\lambda-\mL_k)u=f$ for $u\in\calX$ and $(1+e^{\beta z})f\in \calX$.
Let  $$p_n(z)=e^{\beta n}(1+\tanh\frac{\beta(z-n)}{2})\,.$$
Then $\sup_{z\in\R}|\pd_z^jp_n/p_n|=O(\beta^j)$ for $j\ge1$, $
p_n(z)\uparrow 2e^{\beta z}$ as $n\to\infty$ and
\begin{gather*}
(\lambda-\mL_k)(p_n(z)u)+[\mL_k,p_n(z)]u= p_n(z)f\,,
  \\ \left\| [(\lambda-\mL_k)^{-1}[\mL_k,p_n]p_n^{-1}
  \right\|_{B(\calX)}=O(\beta)\,,
\end{gather*}
and $ \|p_n(z)u\|_{\calX}\lesssim \|e^{\beta z}f\|_{\calX}$
uniformly in $n$. Letting $n\to\infty$, we have from Beppo-Levi's theorem
$$\|e^{\beta z_2}u\|_{\calX_2}\lesssim \|e^{\beta z_2}f\|_{\calX_2}\,,$$
provided $\beta>0$ is sufficiently small.

If $\beta<0$ and $(\lambda-\mL_k)u=f$ for $u\in\calX$ and $(1+e^{\beta z})f\in \calX$,
we can prove $\|e^{\beta z}u\|_\calX\lesssim \|e^{\beta z}f\|_\calX$
by using $p_{-n}(z)=e^{-\beta n}(1+\tanh\frac{\beta(z+n)}{2})$.
\end{proof}

\end{document}